\documentclass[11pt]{article}
\usepackage[margin=0.75in]{geometry}

\usepackage{amsthm}
\usepackage{amssymb}
\usepackage{amsmath}
\usepackage{comment}
\usepackage{thm-restate}
\usepackage{url} 
\usepackage{hyperref}
\usepackage[noabbrev,capitalise]{cleveref}
\usepackage[utf8]{inputenc}

\usepackage{color}
\usepackage[normalem]{ulem}

\newcommand{\mc}[1]{\mathcal{#1}}

\newcounter{i}
\theoremstyle{plain}
\newtheorem{thm}{Theorem}[section]
\newtheorem{lem}[thm]{Lemma}
\newtheorem{claim}{Claim}[thm]
\newtheorem{proposition}[thm]{Proposition}

\newtheorem{cor}[thm]{Corollary}

\newtheorem{conj}[thm]{Conjecture}
\newtheorem{question}[thm]{Question}

\usepackage{color}

\usepackage[dvipsnames, table]{xcolor}

\newenvironment{proofclaim}[1][]%
{\noindent \emph{Proof.} {}{#1}{}}{\hfill
	$\Diamond$\vspace{1em}}

\usepackage{amsthm,thmtools}

\numberwithin{claim}{thm}

\usepackage{etoolbox}
 \AtEndEnvironment{proof}{\setcounter{claim}{0}}

\theoremstyle{plain} % just in case the style had changed
\newcommand{\thistheoremname}{}
\newtheorem{genericthm}{\thistheoremname}

\theoremstyle{definition}
\newtheorem{definition}[thm]{Definition}

\newtheorem{remark}[thm]{Remark}

\newcommand{\Prob}[1]{\ensuremath{%
\mathbb P\left[#1\right]
}}

\newcommand{\Expect}[1]{\ensuremath{%
\mathbb E\left[#1\right]
}}

\title{A Proof of Nash-Williams' Conjecture}

\author{
Michelle Delcourt
\thanks{Department of Mathematics, Toronto Metropolitan University (formerly named Ryerson University),
Toronto, Ontario M5B 2K3, Canada {\tt mdelcourt@torontomu.ca}. Research supported by NSERC under Discovery Grant No. 2019-04269 and a Sloan Research Fellowship.}
\and
Luke Postle
\thanks{Combinatorics and Optimization Department,
University of Waterloo, Waterloo, Ontario N2L 3G1, Canada {\tt lpostle@uwaterloo.ca}. Partially supported by NSERC
under Discovery Grant No. 2026-04610.}}
\date{\today}

\begin{document}

\maketitle

\begin{abstract} 
A central open question in extremal design theory is Nash-Williams' Conjecture from 1970 that every triangle-divisible graph on $n$ vertices (for $n$ large enough) with minimum degree at least $0.75 n$ has a triangle decomposition. In this paper, we prove this conjecture in full.

In 2016, Barber, K\"{u}hn, Lo, and Osthus~\cite{BKLO16} proved that if the fractional relaxation of Nash-Williams' Conjecture holds for minimum degree $cn$ for some constant $c\ge 0.75$, then Nash-Williams' Conjecture holds for any constant $c' > c$. The previously best-known bound on the fractional relaxation was due to Delcourt and Postle~\cite{DP2021progress} from 2021 with $c= \frac{7+\sqrt{21}}{14} \approx 0.82733$. This bound on the fractional relaxation has grown in importance over the years as it has been directly tied to bounds for a number of other problems in extremal design theory.
\vskip.05in

This paper consists of three parts.  {\bf In Part I, our first main result is a proof of the Fractional Nash-Williams' Conjecture}: if $G$ is a graph on $n$ vertices with minimum degree at least $\frac{3n}{4}$, then $G$ has a fractional triangle decomposition. 
\vskip.05in

{\bf In Part II, our second main result is a Fractional Stability Theorem for Nash-Williams' Conjecture}: if a graph $G$ on $n$ vertices has minimum degree close to $\frac{3n}{4}$ but no fractional $K_3$-decomposition, then $G$ is close (in edit distance) to the join of two $\frac{n}{4}$-regular graphs each on $\frac{n}{2}$ vertices.  We use this to prove that if a triangle-divisible graph $G$ on $n$ vertices has minimum degree close to $\frac{3n}{4}$ but no $K_3$-decomposition, then $G$ is close (in edit distance) to the join of two $\frac{n}{4}$-regular graphs each on $\frac{n}{2}$ vertices. 
\vskip.05in 
{\bf In Part III, our final main result is a proof of Nash-Williams' Conjecture in full.}
\vskip.05in
We further note that combined with other known work, even our proof of the Fractional Nash-Williams' Conjecture yields a number of results as immediate corollaries including:

\begin{itemize}
    \item {\bf Asymptotic ``Erd\H{o}s meets Nash-Williams" Conjecture:} more generally a triangle decomposition of girth at least $g$ exists for any  constants $c' > \frac{3}{4}$ and $g > 0$ provided $n$ is large enough.
    \item {\bf Decomposition Threshold for any $3$-chromatic $F$ is at most $\frac{3}{4}$:} an $F$-decomposition of $G$ exists for any constant $c' > \frac{3}{4}$ provided $F$ is $3$-chromatic and the host graph $G$ is $F$-divisible.
    \item {\bf Asymptotic Approximate Oberwolfach Nash-Williams' Conjecture:} namely, any collection of $2$-regular graphs with $(1-o(1))\cdot e(G)$ edges packs into any $d$-regular graph $G$ on $n$ vertices provided $d\ge (\frac{3}{4}+o(1))n$ and $n$ is large enough. 
    \item {\bf Maximum Triangle Packing in Minimum Degree Graphs:} if $G$ is a graph on $n$ vertices with minimum degree $dn$ for $d\in [\frac{1}{2},\frac{3}{4}]$, then $G$ admits a fractional triangle packing of value $\left(\frac{d}{2}-\frac{1}{4}\right)n^2$ and hence a set of $\left(\frac{d}{2}-\frac{1}{4} - o(1)\right)n^2$ edge-disjoint triangles as conjectured by Yuster.
\end{itemize}
\end{abstract}

\tableofcontents

\part*{Prelude}

\section{Introduction}

\subsection{Nash-Williams' Conjecture}

Combinatorial design theory has a rich history, dating back to the mid-1800s.  One of the most fundamental objects in this area is what is known in the literature as a \emph{Steiner triple system} on $n$ elements; formally, this is a set $S$ of triples of an $n$-element set $X$ such that every pair of elements of $X$ lives in exactly one triple of $S$. Viewing this through a graph theoretic lens, a Steiner triple system on $n$ elements is equivalent to a decomposition of the edge set of the complete graph $K_n$ into edge-disjoint triangles; in general such a decomposition of the edges of a graph $G$ is called a \emph{triangle decomposition} or \emph{$K_3$-decomposition} of $G$. Indeed, two obvious necessary conditions for any host graph $G$ to admit a $K_3$-decomposition consist of divisibility constraints, namely $3~|~e(G)$ and $2~|~d_G(v)$ for every $v\in V(G)$; such a graph is said to be \emph{$K_3$-divisible}. One of the oldest theorems in design theory is a result of Kirkman~\cite{K47} from 1847 that proves $K_n$ has a $K_3$-decomposition whenever $K_n$ is $K_3$-divisible (i.e.,~$n\equiv 1,3 \pmod 6$).

For host graphs other than $K_n$, $K_3$-divisibility alone may not be sufficient to guarantee the existence of a $K_3$-decomposition; consider for instance the cycle $C_6$ which is certainly $K_3$-divisible but does not contain $K_3$ as a subgraph, let alone an edge decomposition into copies of $K_3$. A natural question is if we impose a local density condition on $G$ such as bounding the minimum degree from below can we guarantee the existence of a $K_3$-decomposition? The natural extremal version of this question is what is the smallest minimum degree that guarantees a $K_3$-divisible graph has a $K_3$-decomposition. The most famous conjecture in this setting, and certainly the most central open conjecture in extremal design theory, is what is known in the literature as Nash-Williams' Conjecture from 1970~\cite{NW70}.

\begin{conj}[Nash-Williams~\cite{NW70}]\label{conj:NW}
Let $G$ be a $K_3$-divisible graph with $n$ vertices and
minimum degree $\delta (G) \geq \frac{3}{4}n$. If $n$ is sufficiently large, then $G$ admits a $K_3$-decomposition.
\end{conj}

As noted in an addendum to Nash-Williams' original article~\cite{NW70}, if true, then this conjecture would be best possible as Graham was able to produce a construction showing that the fraction $\frac{3}{4}$ is tight. We also note that Nash-Williams' Conjecture implies a corresponding weaker completion problem: if a \emph{triangle packing} (i.e.~edge-disjoint set of triangles) of $K_n$ has the property that every vertex is in at most $\frac{n}{8}$ of the triangles, then it can be extended to a Steiner triple system provided $n\equiv 1,3 \pmod 6$.  

Given the difficulty of Conjecture~\ref{conj:NW}, much of the partial progress in the last decade has come from viewing this problem from an optimization perspective; in this vein, it is natural to study the fractional relaxation as follows.  A \emph{fractional $K_3$-decomposition} of $G$ is an assignment of non-negative weights to each copy of $K_3$ in $G$ such that the sum of the weights along each edge is precisely one. The \emph{fractional $K_3$-decomposition threshold} $\delta^*_{K_3}$ is the infimum of all real numbers $c$ such that every graph $G$ with minimum degree at least $c\cdot v(G)$ has a fractional $K_3$-decomposition. Interestingly Graham's construction also admits no fractional $K_3$-decomposition, and hence, shows that $\delta_{K_3}^*\ge \frac{3}{4}$. Indeed, a more general form of that construction which admits no fractional $K_3$-decomposition is the join of any two $d$-regular graphs on $\frac{n}{2}$ vertices for any $d < \frac{n}{4}$ ~\cite{DHLP25b}, where the \emph{join} of graphs $H_1$ and $H_2$ is the graph obtained by taking the vertex-disjoint union of $H_1$ and $H_2$ and adding a complete bipartite graph between $V(H_1)$ and $V(H_2)$.

In 2005, Yuster~\cite{Yuster2005asymptotically} was the first to systematically study these fractional decomposition thresholds and to explicitly upper bound $\delta_{K_3}^*$ by a constant smaller than 1; although such a statement certainly would have been implicit in the thesis work of Gustavsson~\cite{gustavsson1991decompositions} from 1991. In a breakthrough result appearing in 2014 and published in 2016,  Barber, K\"{u}hn, Lo, and Osthus~\cite{BKLO16} show that the existence of a $K_3$-decomposition is fundamentally related to the existence of a fractional $K_3$-decomposition as follows.

\begin{thm}[Barber, K\"{u}hn, Lo, and Osthus~\cite{BKLO16}]\label{thm:TriangleDecomp}
Let $\varepsilon > 0$.  Any sufficiently large, $K_3$-divisible graph $G$ on $n$ vertices with minimum degree
$$\delta(G) \geq \left(\max\left\{\delta^*_{K_{3}}, \frac{3}{4}\right\}+\varepsilon\right)\cdot n$$
admits a $K_3$-decomposition.
\end{thm}

The result above naturally leads to the following fractional relaxation of Nash-Williams' Conjecture which would thus imply that Nash-Williams' Conjecture holds for any constant strictly greater than $3/4$.

\begin{conj}[Fractional Nash-Williams]\label{conj:FracNW}
If $G$ is a graph on $n$ vertices with
minimum degree $\delta (G) \geq \frac{3}{4}n$, then $G$ admits a fractional $K_3$-decomposition.
\end{conj}

Subsequent to the work of Barber, K\"{u}hn, Lo, and Osthus~\cite{BKLO16} there has been much interest in the fractional relaxation. In 2014, Garaschuk~\cite{garaschuk2014linear} showed $\delta_{K_3}^*\le 0.956$; in 2016, Dross~\cite{dross2016fractional} improved this to $0.9$; in 2020, Dukes and Horsley~\cite{dukes_minimum_2020} improved this further to $0.852$. Prior to this work, the previously best-known upper bound is due to the authors'~\cite{DP2021progress} previous work from 2021 showing that $\delta_{K_3}^*\leq \frac{7+\sqrt{21}}{14}\approx0.82733$. 

\subsection{Main Results}
 
In this paper, we introduce a novel approach for this problem which we will refer to as ``discharging in the dual'' which we will discuss in more detail in Section~\ref{sec:proof}. The main theorem of Part I of this paper is the resolution of the Fractional Nash-Williams' Conjecture as follows.

\begin{thm}\label{thm:fracNW} 
Fractional Nash-Williams' Conjecture is true.
\end{thm}

Combined with Theorem~\ref{thm:TriangleDecomp}, this resolves the asymptotic version of Nash-Williams' Conjecture as follows.

\begin{cor}\label{cor:main}
Let $\varepsilon > 0$. Let $G$ be a $K_3$-divisible graph with $n$ vertices and
minimum degree $\delta(G)\geq \left( \frac{3}{4} + \varepsilon\right)\cdot n$.  If $n$ is sufficiently large, then $G$ admits a $K_3$-decomposition.
\end{cor}

Theorem~\ref{thm:fracNW} has a number of other corollaries for extensions or variations of Nash-Williams' Conjecture which we detail later in Section~\ref{s:Applications} as the value of $\delta_{K_3}^*$ has been tied to many related problems in other recent works.

The main result of Part II is a fractional stability theorem. Recall from the literature that graphs $G$ and $H$ on the same number of vertices are \emph{$\alpha$-close} if one can add and/or delete at most $\alpha\cdot v(G)^2$ edges to transform $G$ into a graph isomorphic to $H$ (two graphs are \emph{$\alpha$-far} if they are not $\alpha$-close). In Part II, we further analyze our proof of the Fractional Nash-Williams' Conjecture to establish that the general form of Graham's construction mentioned above constitutes the extremal examples as follows.

\begin{thm}\label{thm:fracNWstability} 
For each sufficiently small $\sigma > 0$, there exists $\varepsilon > 0$ such that the following holds: If $G$ is a graph on $n$ vertices with $\delta(G) \ge \left(\frac{3}{4}-\varepsilon\right)n$ and $G$ has no fractional $K_3$-decomposition, then both of the following hold:
\begin{itemize}
    \item[(i)] at least $(1-\sigma)n$ vertices of $G$ have degree at most $(\frac{3}{4}+\sigma)n$, and
    \item[(ii)] the maximum cut of $G$ is at least $(1-\sigma)\cdot \frac{n^2}{4}$.
\end{itemize} 
\end{thm}

We note that although the above would have been natural to conjecture, to our knowledge, no one has previously done so in the literature.  This presumably is due both to the open state of Nash-Williams' Conjecture but also the difficulty in proving such a stability result (the two graphs $H_1$ and $H_2$ inside each set of the cut can be \emph{any} roughly $\frac{n}{4}$-regular graphs and so otherwise admit no structure). 

In Part III, we show that non-extremal graphs also admit the embedding of many $K_3$-absorbers and so using the method of refined absorption and Theorem~\ref{thm:fracNWstability}, we prove the following decomposition stability theorem. 

\begin{thm}\label{thm:NWstability} 
For each sufficiently small $\sigma > 0$, there exists $\varepsilon > 0$ such that the following holds for all sufficiently large $n$: If $G$ is a $K_3$-divisible graph on $n$ vertices with $\delta(G) \ge \left(\frac{3}{4}-\varepsilon\right)n$ and $G$ has no $K_3$-decomposition, then both of the following hold:
\begin{itemize}
    \item[(i)] at least $(1-\sigma)n$ vertices of $G$ have degree at most $(\frac{3}{4}+\sigma)n$, and
    \item[(ii)] the maximum cut of $G$ is at least $(1-\sigma)\cdot \frac{n^2}{4}$.
\end{itemize} 
\end{thm}

In the remainder of Part III, we then proceed to show that any such extremal graph with minimum degree at least $\frac{3}{4}n$ admits a $K_3$-decomposition. This will be done by first \emph{cleaning} the graph (that is decomposing atypical vertices in the extremal structure) and then proceeding to use only \emph{cross} triangles (triangles with exactly two edges across the bipartition); this is done using our new results on refined absorption for partitioned graphs. This works if there are many absorbers in this cross setting which we show exist unless both of $H_1$ and $H_2$ are each $\alpha$-close to either $K_{n/4,n/4}$ or $K_{n/4}\sqcup K_{n/4}$. Then the decomposition strategy divides into three more cases (up to symmetry) depending on which of $H_1$ and $H_2$ are close to which of these graphs. Again, showing the decompositions requires more cleaning and extensive use of our refined absorption for partitioned graphs (our graphs would now have $3$ or $4$ parts and the decompositions use triangles across the various parts in the appropriate numbers). 

The final main result of Part III is thus a proof of Nash-Williams' Conjecture as follows.

\begin{thm}\label{thm:NW} 
Nash-Williams' Conjecture is true.
\end{thm}

\subsection{Outline of Paper} 
In Section~\ref{sec:proof}, we discuss the main proof ideas. In Section~\ref{s:Applications}, we explicitly detail the further applications of the Fractional Nash-Williams' Conjecture and their relevant history. In Part I, we prove the Fractional Nash-Williams' Conjecture. In Part II, we prove our Fractional Stability Theorem. In Part III, we prove Nash-Williams' Conjecture. Each part contains its own proof overview with the statement of key lemmas followed by an outline of said part.

\section{On the Key Proof Ideas and their Novelty}\label{sec:proof}

First, let us discuss the techniques used in the various results to date on the Fractional Nash-Williams' Conjecture. In her thesis work from 2014~\cite{garaschuk2014linear}, Garaschuk used polyhedral cone methods which we discuss further in the next subsection. In a different vein, Dross~\cite{D15} (as well as Dukes and Horsley\cite{dukes_minimum_2020} who refined Dross' method) used the max-flow min-cut theorem to derive a duality contradiction. The authors in the previously best known work  on the Fractional Nash-Williams' Conjecture utilized `edge gadgets' and techniques from nonlinear optimization to directly construct a fractional $K_3$-decomposition (that is, they proved a certain explicit weighting of triangles that summed to $1$ on every edge had nonnegative triangle weights via optimization). 

For our proof of the Fractional Nash-Williams' Conjecture in Part I, we return to Garaschuk's approach but in a novel twist combine it with the method of discharging from graph coloring. We call this combined method \emph{discharging in the dual}. For our proof of a fractional stability theorem in Part II, we further extend this method via the use of \emph{layered discharging}. Finally in Part III, we use our new results on partitioned designs from~\cite{DP26} combined with a thorough analysis of all the extremal cases as well as cleaning procedures for each case to prove Nash-Williams' Conjecture in full. In the remainder of this section, we discuss the key proof ideas for each of the parts in more detail (technical proof overviews with definitions and lemma statements can be found at the start of each part), but first we crucially recall Farkas' Lemma~\cite{F02} and how it relates to fractional $K_3$-decompositions.

Now is also the time to address the complexity of our proof. In hypergraph decomposition problems (and similarly matching or embedding problems), two natural barriers arise and understanding them has proved crucial. The first is a \emph{space barrier} where the lack of a (robust) fractional decomposition is the issue; the second is a \emph{divisibility barrier} where an object lacks the required divisibility properties (or robust versions of them). These barriers reoccur often throughout the literature; indeed, one might \emph{assume} a robust fractional property and a robust divisibility property and show these properties guarantee the desired decomposition; Keevash and Mycroft~\cite{KM15} did just this in their theory of hypergraph matchings and Lang~\cite{lang2023tiling} did the same for his theory of tilings. The authors recently did the same for hypergraph decompositions of partitioned hypergraphs~\cite{DP26} which we will need here. That said, to prove Nash-Williams' Conjecture we do not have the luxury of such robust properties as we are proving the existence of a triangle decomposition exactly to the point it breaks. Moreover, Nash-Williams' Conjecture is especially difficult as $3n/4$ is \emph{both} a space barrier (being the threshold for a fractional $K_3$-decomposition) and a divisibility barrier (the threshold for the disappearance of $K_3$-absorbers). Even more, unlike $1$-uniform problems such as hypergraph matching, the barrier for designs prove to be extra cumbersome as modifying the decomposition plan to work around the space barrier needed us to resort to the theory of partitioned hypergraphs which involves understanding both \emph{integer divisibility} (the normal notion of residing in the integer lattice of the degree vectors) as well as \emph{fractional divisibility} (residing in the cone of the degree vectors). All of these factors contribute to the high complexity of the proof.

\subsection{Part I: Discharging in the Dual}

\subsubsection{Farkas' Lemma for Decompositions}

One of the foundational results in polyhedral combinatorics and linear programming is the classical Farkas' Lemma~\cite{F02} from 1902. The following is one standard version (it has many equivalent forms):

\begin{lem}[Farkas' Lemma]
The matrix equation $Ax=b$ has a nonnegative solution $x\ge 0$ if and only if there does not exist a vector $y$ such that $A^T y \ge 0$ and $b^T y < 0$.
\end{lem}

The connection to fractional $K_3$-decompositions is as follows. Let $G$ be a graph and let $A$ be the \emph{edge-triangle incidence matrix} of $G$, that is, the rows of $A$ correspond to the edges of $G$, the columns of $A$ to the triangles of $G$ and the cell $(e,T)$ is $1$ if $e\in T$ and $0$ otherwise. The key observation is that $G$ has a fractional $K_3$-decomposition if and only if $Ax = 1^{E(G)}$ has a non-negative solution $x$. Applying Farkas' Lemma to this matrix equation yields the following (this key fact formed the basis of Garaschuk's thesis). 

\begin{lem}[Farkas' Lemma applied to $K_3$-Decompositions]\label{lem:FarkasK3}
A graph $G$ does not admit a fractional $K_3$-decomposition if and only if there exists an assignment of weights to the edges of $G$ such that the weight of every triangle is non-negative but the global sum is negative.    
\end{lem}

We note that this holds more generally for (hypergraph) $F$-decompositions, and indeed even more generally for hypergraph matchings (decompositions can be conceived of as a hypergraph matching of the auxiliary \emph{design hypergraph}; see~\cite{DPI} for a definition). From this perspective, Farkas' Lemma has a rich history of use in hypergraph matchings and related problems as follows. 

\begin{remark}
Notably, Farkas' Lemma was used to approximately determine the minimum degree thresholds for hypergraph matchings by R\"odl, Ruci\'nski, and Szemer\'edi~\cite{RRS06b} in their seminal paper from 2006. Subsequently it has appeared in a number of related works.  For instance, Farkas' Lemma was used by Keevash and Mycroft~\cite{KM15} in their treatise on a geometric theory for hypergraph matching, by Han~\cite{H21} in 2021, and by Bowtell, Kathapurkar, Morrison, and Mycroft~\cite{B25} in 2025. Each of these latter papers directly seek a contradiction via counting arguments from which discharging can be seen as a most potent generalization. On the other hand, in her thesis Garaschuk endeavored to understand the cone generated by non-negative triangle weightings which while yielding partial results proved too complicated to master. 
\end{remark}

We thus instead return to finding a contradiction to the existence of such a weighting directly but with the use of much more sophisticated counting arguments, namely \emph{discharging}, a fundamental tool in graph coloring. Indeed, if we view the weights as charges the above looks quite similar to a usual discharging setup. 

\subsubsection{Discharging}

\emph{Discharging} has a long and storied history, most notably in the proof of the Four Color Theorem. The main idea behind discharging is that to show the sum of charges in a graph is in fact nonnegative (thus yielding a contradiction), one stipulates \emph{discharging rules} that preserve charge and then proves that each charged object (i.e.~vertex, edge or face) has final nonnegative charge. This latter step is usually shown to hold by noting that certain local configurations are forbidden in a minimal counterexample (called \emph{reducible configurations}).

Discharging was first introduced by Wernicke~\cite{W04} in 1904 in an attempt to prove the Four Color Theorem. Heesch in the 1960s and 1970s was the first to use computers with discharging to approach proving the Four Color Theorem. This vision was fulfilled when Appel and Haken~\cite{AH76} found a computer-based discharging proof of the Four Color Theorem in the 1970s. In 1996, Robertson, Seymour, Sanders and Thomas~\cite{RSST97} provided an alternate proof of the Four Color Theorem also using computers and discharging (notably their discharging rules are verified by hand but the millions of reducible configurations are still verified to be reducible via computer). 

Since the proof of the Four Color Theorem, discharging has been used extensively in the study of colorings of planar graphs and graphs on surfaces. For these settings, discharging usually assigns charges to vertices, edges and/or vertices and a contradiction is achieved to Euler's formula which asserts the sum should be negative (or at most $6(g-2)$ for a graph embeddable on a surface of genus $g$). Discharging has also been used in the study of properties of critical graphs via the so-called Potential Method. Various more complicated forms of discharging have been introduced over the years to solve ever more complicated problems such as \emph{layered discharging} (wherein discharging rules are administered in layered steps), \emph{triggered discharging} (wherein certain rules only fire when `triggered'), and \emph{global discharging} (wherein discharging rules can be global in nature sending to objects arbitrarily far away in the graph) to name a few.

As for our proof of the Fractional Nash-Williams' Conjecture we turn to discharging to derive a contradiction to the existence of the edge-weighting in Farkas' Lemma for Triangle Decompositions for graphs $G$ with $\delta(G)\ge \frac{3}{4}\cdot v(G)$. Our proof does not make use of reducible configurations (since induction is notoriously unhelpful for decomposition problems). Rather \emph{the only `local forbidden configuration' is a triangle whose sum of edge-weights is negative}. While this is seemingly not a lot of information to utilize, we will yet show it suffices to yield a contradiction precisely when $\delta(G)\ge \frac{3}{4} \cdot v(G)$. (Indeed there are many triangles in such graphs so that is at least something to work with). As far as we can tell, this might be the first use of discharging in design theory, and thus we have labeled this new combined method of Farkas' Lemma and discharging as \emph{`discharging in the dual'} since we believe the method shows promise for a number of other problems in design theory and the theory of graph and hypergraph decompositions.  

\subsubsection{On the Discharging Rules}

But what should the discharging rules be? How is a discharging perspective useful here over pure counting arguments? The freedom in discharging lies in the power of pure imagination - the possibility to conjure any set of discharging rules. For one cannot really disprove that there could exist the desired set of discharging rules. As for our problem of the Fractional Nash-Williams' Conjecture, we initially started with a very natural set of discharging rules. The key is to observe there are two types of triangles in our problem that contain negative edges: those with one negative edge and two positive edges (which we call \emph{acute}) and those with two negative edges and one positive edge (which we call \emph{obtuse}); this follows since three negative edges is a forbidden configuration by assumption. We note for the reader here that we assume there are no zero weight edges by slightly increasing the weight of any such edges (see Proposition~\ref{prop:NonZeroUntied} for a formal proof).

Our first natural discharging rule set consisted of one rule wherein every positive edge divides its weight in half between its end vertices and then uniformly distributes that weight to incident negative edges where we take care to re-aggregate the weight in obtuse triangles. This rule yielded two major cases each reducing to a 4-variable optimization program that we solved by computer; surprisingly this rule set proved that a minimum degree of $0.785n$ suffices to guarantee a fractional triangle decomposition (a marked improvement on our previous best-known result!). We continued to tinker with these degree-based rule sets until we had a proof of $0.77n$. However, all of these rule sets, which we would term as \emph{first-order} since they only involved the first neighborhood of edges, annoyingly came up short in proving the Fractional Nash-Williams' Conjecture as they failed to capture the full intricacies of the problem.

Not being satisfied with anything less than a full proof of the Fractional Nash-Williams' Conjecture, we turned to new rule sets, ones that we might term \emph{second-order} in that they took into account edges that were in the second neighborhood. Already our first-order discharging rule sets which had considered non-constant degree vertices seemed daunting (but analyzable via computer optimization); the new second-order rules seemed near impenetrable to analysis. To ease our study and conception of the different rule sets, we introduced a general framework to study many possible rule sets all at once. We would have every negative edge \emph{demand} a certain weight from positive edges it was in a triangle with (based only on the weights and type of said triangle) \emph{but then} each positive edge would distribute its weight to all demanding negative edges proportional to their demands. This \emph{demand setup} captured second-order behavior while crucially ensuring that each positive edge had nonnegative final charge. Now of course we had to find the right set of demands to ensure every negative edge had nonnegative final charge (not necessarily any easier than finding the right rule set). More problematically, how could we possibly analyze these second-order rule sets? How could we turn them into a finite (as in 4 or 15) variable optimization program? 

\subsubsection{Further Proof Ideas}

Now we move on to the next chapter of our story wherein we discuss the additional proof ideas needed to turn our dream of a second-order demand setup discharging rule set into a proof of the Fractional Nash-Williams' Conjecture. Again, here we assume we have an as yet unspecified demand function and it thus suffices to show that a negative edge $f=uv$ has a nonnegative final charge (equivalently weight) which we denote by $\textrm{ch}_F(f)$. 

\begin{itemize}
    \item {\bf Linear Approximation Lower Bound:} There are two main issues when trying to analyze the final charge of a negative edge arising from this demand setup, both concerning the nature of terms of charge sent from a positive edge $e$ that $f$ demands from: these terms are the weight of $e$ times the demand of $f$ to $e$ divided by the sum of all the demands of negative edges to $e$. The first issue is that such a term is naturally convex (whether it is actually convex would depend on the specific demand function). This is unfortunate as such functions are not likely to admit Zykov symmetrization (as for example concave functions do, see Balogh and Wigal~\cite{BW25}). The second and more problematic issue is that such a term is highly nonlinear, intertwining up to $\theta(n)$ demand functions in the denominator and so is extremely hard to analyze. 
    
    Our solution to these issues is to instead pursue a lower bound on the final charge which we denote by $\textrm{ch}_F'(f)$ (see Definition~\ref{def:chargeprime}), that is much more manageable to work with, specifically we will consider what might be conceived of as a `linear approximation' of said demand terms. One key, (since there are extremal examples for the conjecture) is to choose this linear approximation so as to not incur any loss for the extremal examples. This is accomplished with an appropriate rescaling and that we indeed obtain a lower bound is simply due to the fact that $1/(1-x) \ge 1+x$ for all real $x < 1$. This linear approximation then separates each of our previously intertwined denominators into a linear sum of terms. Of course we should not necessarily expect the minimum of this linear approximation lower bound $\textrm{ch}_F'(f)$ to be the same as for the original $\textrm{ch}_F(f)$ (namely the absolute value of the weight of $f$ as desired). That said, it is also not absurd to think the optimal value \emph{might} be the same. 
    \item {\bf Cleaning:} The next proof idea we use before proceeding with Zykov symmetrization of our linear approximation $\textrm{ch}_F'(f)$ is to worst case lower bound the linear terms for triangles outside the common neighborhood of $u$ and $v$ (we note this includes the terms for $u$ and $v$ themselves). We do this as there is not much control of the structure of these terms. We note the cleaned formula will involve a large positive term depending on the size of $N(u)\cap N(v)$ and the sum of many negative terms (second-order terms coming from the demands of negative $f'$ in $N(u)\cap N(v)$ to positive edges $e$ incident with $f$). In the remainder then it suffices to \emph{upper bound} the sum of these terms which we denote by $\textrm{ch}_F''(f)$ (after an appropriate rescaling), see Definition~\ref{def:chargedoubleprime}.   
    \item {\bf Zykov Symmetrization:} We could next proceed to use Zykov symmetrization on $\textrm{ch}_F''(f)$ at this stage by symmetrizing over positive edges in $G[N(u)\cap N(v)]$. We note before symmetrizing, we would actually first turn non-edges of $G[N(u)\cap N(v)]$ to be very positive edges. Then symmetrization would in theory would yield a structure of the common neighborhood having up to $4$ parts (cliques of positive edges): one part is $N^+(u)\cap N^-(v)$, that is the positive neighborhood of $u$ intersect the negative neighborhood of $v$ (which we note only contains positive edges within it as such an edge lies inside $N^-(v)$), similarly a second part is $N^-(u)\cap N^+(u)$, and there are up to two more parts in $N^+(u)\cap N^+(v)$ (but only two as there are no negative triangles). Now in theory this symmetrization has reduced our problem to a finite variable optimization program for a fixed demand function, specifically a 17-variable one (14 weights of edges inside a $K_6$ besides that of $f$ and 3 part sizes).
    \item {\bf A Value Formulation and Graph Lagrangians:} Now this 17-variable program would be daunting for a computer for even simple demand functions at least to provide a provable guarantee. Since we also require a fractional stability analysis for our proof of Nash-Williams' Conjecture, we pursue an alternate \emph{computer-free} approach. It is useful to conceive of the above problem as a \emph{Graph Lagrangian} problem wherein the \emph{value} of a pair of vertices $x,y\in N(u)\cap N(v)$ (denoted as $\textrm{Val}_f(xy)$, see Definition~\ref{def:Value}) is determined by the sum of the linear demand terms for triangles involving $x,y,u,v$. It is fairly straightforward to show this value is at most $2|\phi(f)|$ for any pair $x,y$ (where $\phi(f)$ denotes the weight of $f$), see Proposition~\ref{prop:EdgeWeightUpperBound}. Unfortunately, this simple bound is not sufficient to upper bound $\textrm{ch}''_F(f)$ to the desired number. It was at this point, we dreamt up an \emph{intermediary property} we call \emph{triangle-thinness} that we could prove yields the desired upper bound on $\textrm{ch}''_F(f)$. Namely triangle-thinness would mean that the sum of values over the pairs of three vertices $x,y,z\in N(u)\cap N(v)$ is at most $4|\phi(f)|$ (instead of $6|\phi(f)|$ as guaranteed by the simple bound). We prove in Section~\ref{s:ValueLemma} this key Value Lemma, Lemma~\ref{lem:Value}, that this triangle-thinness property yields the desired upper bound on $\textrm{ch}''_F(f)$. Our proof of the Value Lemma proceeds via `mass moment transfer' type arguments. 
    
    \item {\bf Proving Triangle-Thinness via Strict Sponsorship:} Now the real crux in our entire proof is finally determining what demand function would ensure the property of triangle-thinness in the common neighborhood. One natural choice we considered is the demand function that results in each positive edge dividing its weight equally among demanding triangles (and further dividing weight in obtuse triangles proportional to the triangle's negative edge weights). This demand function however admitted \emph{local counterexamples} (a local structure where one negative edge $f$ does not end up with nonnegative final charge); this made some intuitive sense as it does not judiciously balance demands. A next natural choice we considered is for negative edges to only demand what they need where for acute triangles this demand is divided proportionally between the two positive edges according to their weights. This too turned out to have local counterexamples at least for the linear approximation lower bound $\textrm{ch}'_F(f)$; the issue here was low weight edges sending charge across an acute triangle even though the other positive edge was much heavier than the negative edge. 
    
    This led to the idea of a \emph{sponsor} (an acute triangle is \emph{sponsored} if at least one of the positive edges has weight at least the absolute value of the negative edge weight). The easiest approach to analyze using sponsorship that yielded the triangle-thinness property was \emph{strict sponsorship} wherein the heaviest positive edge in a sponsored acute triangle pays all the negative weight (we can avoid ties by judiciously increasing positive weights ever so slightly). Now proving triangle-thinness, Lemma~\ref{lem:EdgeWeightTriangleThin}, for strict sponsorship is neither trivial nor easy (roughly 8 pages of structural graph theory analysis filling all of Section~\ref{s:FractionalCases}) but this concludes the proof ideas we used to prove the Fractional Nash-Williams' Conjecture. See Section~\ref{s:FracNW} for the start of the formal proof with the formal definitions and theorem statements. 
\end{itemize}

\subsection{Part II: Fractional Stability and Property Testing}

\subsubsection{On the Need for Fractional Stability}

One might naively think given the proof of the Fractional Nash-Williams' Conjecture that all that remains to consider is the existence of absorbers and a stability analysis for when absorbers do not exist. However, even in the case of a graph $G$ where absorbers exist in abundance, more is needed. In this case, one can attempt the refined absorber machinery, taking reserves $X$ and an omni-absorber $A$. One could then use the inheritance lemma and the almost fractional packings guaranteed from our fractional proof to generate a low-weight almost fractional packing of $G\setminus (X\cup A)$; this can even be turned into a low-weight almost fractional decomposition by redecomposing the neighborhoods of underweight vertices via Dirac's theorem. 

The issue is that to fix the almost here (so as to be able to use nibble), one needs both fractional absorbers (which would follow via absorbers) and for the fractional decomposition to be sufficiently \emph{seeded} to correct the almost via said fractional absorbers (meaning the cliques in the off decomposition of said fractional absorber would need to be able to accept negative weight). But there is no guarantee that the inherited fractional decompositions are seeded and indeed seem likely to use mostly `cross' triangles in the fractional extremal case. Trying to correct this issue by first seeding the fractional decomposition only creates a logical loop as then the almost would depend on the seeding which then could not fix the almost. This is all to say this naive approach suffers from a \emph{hierarchy issue}. 

From an alternate perspective, such a hierarchy issue seems inevitable since $\frac{3n}{4}$-regular graphs that have a complete bipartite balanced cut (the fractional extremal) can only be (fractionally) decomposed via `cross' triangles (triangles that cross the cut) due to counting arguments. So any proof that does not ensure that only cross triangles are used in the decomposition is doomed to fail. Thus we are led to proving a stability version of our fractional result (in Part II) to characterize these extremal cases. Whereas in Part III, we analyze these fractional extremal cases (and the absorber extremal cases) in much further depth to prove they have decompositions. 

\subsubsection{Key Proof Ideas}

We now turn to discussing the key proof ideas for Theorem~\ref{thm:fracNWstability}, our Stability Version of Fractional Nash-Williams' Conjecture. Recall that goal is to show a graph $G$ on $n$ vertices with minimum degree $\left(\frac{3}{4}-\varepsilon\right)n$ that has no fractional $K_3$-decomposition must be $\sigma$-close to being \emph{extremal}, that is the join of two $n/4$-regular graphs on $n/2$ vertices, or equivalently this means that {\bf almost all vertices of $G$ have degree close to $3n/4$ and the maximum cut of $G$ contains roughly $n^2/4$ edges}. 

\begin{remark}[On the Difficulty of our Extremal Problem]
Before proceeding with the key proof ideas, let us first provide some insight on why the problem is difficult. One might naively think we could just assume we are not $\alpha$-close to extremal and quickly derive a contradiction. It is true that not having the degree property is easy to characterize since that would mean having a decent set of vertices of much larger degree; on the other hand, not having such a max cut is a more amorphous property. More problematically, our proof of the Fractional Nash-Williams' Conjecture was to derive a contradiction via discharging in the dual formulation; this required considering triangle-covered edge-weightings. Thus a stability version of that proof would also necessitate understanding edge-weightings. However, we note for the reader that there are at least two very different extremal weightings: the natural one to consider is edges across the max cut being negative say $-1$ for normalization purposes and the edges inside each set being $+2$ -this is thus triangle covered and will sum to negative if the number of cut edges is strictly \emph{more} than twice the number of internal edges; the more surprising one arises when the subgraphs inside each side of the cut are complete bipartite graphs (so the whole graph is roughly a balanced blowup of $K_4$) - here a weighting of $-2$ to the internal and $+1$ to the cut edges yields a triangle-covered edge-weighting and the sum will be negative provided the number of cut edges is strictly \emph{less} than twice the number of internal edges. (Indeed in this later case, one can even take interpolations of it with its own permuted weightings resulting simply with values for the three matchings of $K_4$ that sum to $0$ with the requirement that the most positive valued matching having slightly fewer edges.)
\end{remark}

\noindent Now we discuss the key proof ideas as follows:

\begin{itemize}
    \item {\bf Value Stability:} The first key is to prove a stability version, Lemma~\ref{lem:ValueStability}, of our value stability lemma, Lemma~\ref{lem:Value}. Indeed we first prove a stability version, Lemma~\ref{lem:ValueStability2}, of our value lemma with 4 symmetrized parts, Lemma~\ref{lem:Value4Parts}. While the proof of Lemma~\ref{lem:Value4Parts} is roughly a page, that of Lemma~\ref{lem:ValueStability2} is roughly 5 pages since there actually four extremal families of step graphons even for this simplified problem. Similarly the proof of Lemma~\ref{lem:ValueStability} is a further four pages from Lemma~\ref{lem:ValueStability2} while Lemma~\ref{lem:Value} was just an extra one page from Lemma~\ref{lem:Value4Parts}. The proof idea for Lemma~\ref{lem:ValueStability2} is that of an algebraic contradiction somewhat typical of stability arguments but complicated by there being four different extremal families. The proof idea for Lemma~\ref{lem:ValueStability} is building a partition into four well-chosen parts (typical of Tur\'an stability arguments) but then coupling this with a random sampling argument of the induced step graphons and then arguing via the nature of the four families that such sampling forces a blowup structure. This latter argument is quite specific to our individual problem and does not in general seem to work for the extremal problem of maximizing the edge-weight of $K_r$-free graphs subject to prescribed upper bounds on sums of edges in cliques.
    \item {\bf Extremal Structure of the Common Neighborhood:} Armed with the statement of value stability which we will in fact need in all its technical glory, we then proceed to characterize the graph structure and edge-weights of the common neighborhoods of edges in a graph on $n$ vertices with minimum degree slightly less than $3n/4$ whose charge is still negative (or barely positive compared to its original weight) after our original discharging - what we will call an \emph{extremal edge}. This structure comes in two varieties, however one of these varieties has a saving grace if we turn to layered discharging to add extra discharging rules. 
    \item {\bf Layered Discharging to force Global Extremal Structure:} One variety has the property that there are linearly negative edges incident with our extremal edge whose absolute value of their weight is at least two times compared to that of our extremal edge.  Thus we add an extra discharging rule that non-extremal edges (those who after original discharging are decently positive) send $1/2$ their final charge uniformly to incident edges and $1/2$ uniformly to all other edges. Thus if we consider a maximum weight edge with small secondary charge it follows that is it of the second variety: it has a common neighborhood of size roughly $n/2$ with an almost-balanced almost-complete bipartite cut of negative edges almost all of which have similar weight to our chosen edge. 

    If the two sides of the cut are almost independent sets, then this is already close to our desired extremal structure (since we can also argue almost all vertices have degree roughly $3n/4$ or all edges are nonnegative after the first extra discharging due to our extremal structure results). So instead we assume there exists a decent number of edges in one of the set. From there, we proceed to argue these edges inside the set must be positive with roughly twice the weight of maximum still negative edges. By considering their remaining incident edges to the rest of the graph, these positive edges can be used to once more force the extremal structure. (Indeed this is quite satisfying as there were the two different extremal dual weightings mentioned before: the more restricted one was captured by our first use of extremal with the two quasi-independent sets while the more natural one was captured when we forced this more globally consistent weighting).
\end{itemize}

\subsubsection{Existence of Low-Weight Seeded Fractional Decompositions via  Property Testing}

The ideas described previously conclude those necessary for the proof of Theorem~\ref{thm:fracNWstability}. However in order to prove Theorem~\ref{thm:NWstability} it is not sufficient to just have a fractional decomposition for graphs with minimum degree close to $3n/4$ that are not extremal. Rather we require an ``almost" low-weight seeded fractional decomposition (meaning every every clique receives some but not too much weight and every edge sums to weight roughly one). This will follow by combining the Inheritance Lemma, wherein we find that almost every small random subgraph also has minimum degree at least $3n/4$, and results from \emph{property testing}. 

A fundamental result in property testing is that the proportional size of maximum cut in a graph $G$ is \emph{estimable} meaning that almost every subgraph of $G$ has a maximum cut of roughly the same proportion as $G$. It is obvious that almost every subgraph of $G$ should have a maximum cut at least as large as $G$; the much harder fact to prove is that they do not have larger cuts but this is what estimable tell us. Using this and inheritance, it will follow that if a graph $G$ is not close to extremal (extremal here means almost all degrees close to $3n/4$ and the maximum cut is close to $n^2/4$ edges), then in fact almost all small random subgraphs are also not close to extremal. This might seem counterintuitive but such is the power of Szemer\'edi's Regularity Lemma that if even some tiny proportion of small random subgraphs fail to have our extremal property then almost all fail to have it.

Here are the required definitions. We first recall that a \emph{fractional $F$-packing} of a graph $G$ is an assignment of non-negative weights to copies of $F$ in $G$ such that for every edge $e$ of $G$ the sum of the weights of copies of $F$ containing $e$ is at most $1$. 

\begin{definition}
Let $F$ be a graph and let $C \ge 1 \ge \sigma$. A fractional $F$-packing of a graph $G$ on $n$ vertices is \emph{$C$-low-weight} if every copy of $F$ in $G$ has weight at most $\frac{C}{n^{v(F)-2}}$. A fractional $F$-packing of a graph $G$ on $n$ vertices is \emph{$\sigma$-seeded} if every copy of $F$ in $G$ has weight at least $\frac{\sigma}{n^{v(F)-2}}$. A fractional $F$-packing of a graph $G$ on $n$ vertices is \emph{$C$-regular} if it is $C$-low-weight and $\frac{1}{C}$-seeded.
\end{definition}

Here is the required generalization of Theorem~\ref{thm:fracNWstability}. Its proof is given at the end of Part II. Starting with our fractional stability, the Inheritance Lemma and estimability of max-cut will provide the almost low-weight properties; however to prove we can also make this seeded requires coloring and averaging techniques similar to those we developed in~\cite{DHLP25} which in turn were inspired by those of Montgomery in~\cite{montgomery_fractional_2019}.

\begin{thm}\label{thm:RegularFracStability}
For each sufficiently small $\sigma > 0$, there exists $\varepsilon > 0$ and $C\ge 1$ such that the following holds for all sufficiently large $n$: If $G$ is a graph on $n$ vertices with $\delta(G) \ge \left(\frac{3}{4}-\varepsilon\right)n$ and $G$ does not admit a $C$-regular fractional $K_3$-decomposition, then both of the following hold:
\begin{itemize}
    \item[(i)] at least $(1-\sigma)n$ vertices of $G$ have degree at most $(\frac{3}{4}+\sigma)n$, and
    \item[(ii)] the maximum cut of $G$ is at least $(1-\sigma)\cdot \frac{n^2}{4}$.
\end{itemize}     
\end{thm}

\subsection{Part III: Decompositions of Extremal Cases}

Now we proceed to discuss how to prove Nash-Williams' Conjecture armed with fractional stability. The first key is to characterize \emph{absorber stability}, that is when do many absorbers exist. As alluded to before, the case where many absorbers do not exist is a subset of the extremal cases from fractional stability; namely, many absorbers exist unless each of $H_1$ and $H_2$ are each close to either $K_{n/4}\sqcup K_{n/4}$ or $K_{n/4,n/4}$. Thus the proof of Nash-Williams' Conjecture ultimately divides into 5 cases. When the graph is not close to extremal, the proof follows from our fractional stability and absorbers stability results via our refined absorption decomposition theorem. For the remaining cases, we need to invoke the theory of partitioned designs that we developed in~\cite{DP26}. For the second case where we still have absorbers, care will nonetheless be needed to ensure cross triangles (i.e.~across the cut) are used to the best of their ability to decompose. For the remaining cases, even more parts are required and much care must be given to ``clean" the graph so as to be integrally and fractionally divisible in the various partitioned settings; even characterizing when graphs are integrally and fractionally divisible proves to be highly non-trivial in this new partitioned setting.

\subsubsection{Absorber Stability}

The key idea to prove our absorber stability lemma, Lemma~\ref{lem:ManyHinges}, is to first prove a lemma (Lemma~\ref{lem:ManyPaths}) characterizing when there exist many paths of a prescribed parity (i.e.~odd or even) length in graphs on $n$ vertices of minimum degree roughly $n/2$. Namely there are two extremal families: first there are those graphs that are close in edit distance to $K_{n/2}\sqcup K_{n/2}$ and second there are those graphs that are close in edit distance to $K_{n/2,n/2}$. In the former, vertices that `map' to different cliques will not have many paths between them (no matter the length); in the latter, vertices that map to the same part of the complete bipartite graph have few odd length paths while those that map to different parts have few even length paths.

Our decomposition theorem actually does not require many absorbers but rather many hinges and anti-edges (see Definitions~\ref{def:AntiEdge} and~\ref{def:hinge}) which imply the existence of many absorbers. The existence of many anti-edges (which is equivalent here to paths of length two between two specified vertices) follows easily from our minimum degree condition. As for hinges, we need to show for any two triangles $T_1$ and $T_2$ intersecting in an edge, call it $uv$, that there is a size $i$ (for us this will be $2$, $4$ or $6$) such that number of subgraphs $B$ with $i$ new vertices which when joined to $T_1\cup T_2$ admit two disjoint triangle decompositions (what we call a \emph{hinge} for $T_1$ and $T_2$) is $\Theta(n^i)$.   

To show the existence of many hinges, we use relatively simple constructions for $B$. Namely for $i\in [2]$, we let $w_i := V(T_i)\setminus \{u,v\}$ and we then consider their common neighborhood, namely $H:= G[N(w_1)\cap N(w_2)]$. It would suffice to find many paths of odd length in $H$ between $u$ and $v$ since any such path $P$ plus the edges from the internal vertices of $P$ to both $w_1$ and $w_2$ yields a hinge. Thus we use our path lemma to determine that $H$ is one of the two extremal families. This will show that one half of the graph (the common neighborhood) is close to one of the two extremal families above. To show that the other half is also close to one of these extremal (and hence the whole graph is close to the join of two of them as desired), we `flip' the hinge construction. Specifically, we consider the many choices of a vertex $v'\in N(v)\cap N(w_1)\cap N(w_2)$. Then we consider $H':= G[N(u)\cap N(v')]$. Now it would suffice to find many even length paths $P'$ in $H'$ from $w_1$ to $w_2$ since any such path yields a hinge along with $v'$.

\subsubsection{On the Need for Partitioned Designs}

In light of both the space and divisibility barriers, and armed with both fractional and absorber stability results, we then utilize partitioned designs, that is decompositions of a \emph{partitioned} graph $G$ (graphs equipped with a vertex partition; here for us it is at most four parts) into copies of a partitioned graph $F$. This setting encompasses Latin squares, mutually orthogonal Latin squares, resolvable designs and much more. We will need to prove various decomposition results in this setting (namely one for each of our 5 cases). Unfortunately, Keevash's very general Designs II results are not amenable for use in this extremal setting. The method of iterative absorption has been used in these extremal settings for partitioned graphs but the method does not provide a black-box theorem (and so proofs would have to be redone at quite some length for each case) and has not been used for decomposing into families of partitioned graphs which we will need. So instead we will use our method of refined absorption, now extended by us in~\cite{DP26} to decomposing into families of partitioned graphs in a black-box manner as we discuss further in a bit.

First though, a word should be spared for why we need to consider these partitioned settings in what is seemingly a non-partitioned problem. The first and most serious issue arises from the space barrier since if the degree is permitted to fall below $3n/4$ then counterexamples arise (for example the join of two $d$-regular graphs $H_1$ and $H_2$ with $d < n/4$ which do not have enough interior edges to decompose the cross edges). To prove Nash-Williams' Conjecture down to $3n/4$ minimum degree \emph{exactly} instead of approximately, we must grapple with the case that $G$ is the join of two $n/4$-regular graphs each on $n/2$ vertices. In such a case the cross edges must be decomposed exactly via the interior edges, that is a triangle decomposition for such a graph must only use  \emph{cross triangles} (that is triangles with at least one vertex from each side). This is naturally a partitioned graph problem wherein $G$ has vertex partition $(V_1,V_2)$ and we need to decompose the edges of $G$ into copies from the family $\mc{F}_0:=\{T_{112},T_{122}\}$ where $T_{ijk}$ denotes the triangle with vertices in parts $i,j,k$. Even more problematically in the remaining three `special' cases that arise from absorber stability wherein each of $H_1$ or $H_2$ is close to either $K_{n/2}\sqcup K_{n/2}$ or $K_{n/2,n/2}$, we essentially have a partitioned graph decomposition problem where $G$ is partitioned into $4$ parts (let us say parts $1$ and $2$ partition $H_1$ and parts $3$ and $4$ partition $H_2$) and there are various allowed families of triangles, namely the families $\mc{F}_1=\{T_{113}, T_{114}, T_{223}, T_{224},T_{133},T_{233},T_{144},T_{244}\}$, $\mc{F}_2=\{T_{113}, T_{114}, T_{223}, T_{224},T_{134},T_{234}\}$, and $\mc{F}_3=\{T_{123}, T_{124}, T_{134}, T_{234}\}$.

\subsubsection{Key Proof Ideas}

We now discuss the key proof ideas remaining in our proof of Nash-Williams' Conjecture as follows.

\begin{itemize}
    \item {\bf Characterizing Fractional and Integral Divisibility:} Deriving a good characterization of when a partitioned graph $G$ satisfies the necessary divisibility conditions for an integral $\mc{F}$-decomposition (what we refer to as \emph{integrally $F$-divisible}) is significantly more complicated than the non-partitioned setting. This notion is related to being contained in the \emph{lattice} of certain generator vectors, a well-studied problem that can be characterized via the \emph{congruence relations} of the lattice which in turn can be solved using the \emph{Smith Normal Form}. Similarly, deriving a good characterization of when a partitioned graph $G$ satisfies the necessary divisibility conditions for a fractional $\mc{F}$-decomposition (what we refer to as \emph{fractionally $F$-divisible}) is significantly more complicated than the non-partitioned setting (where it is trivial). This notion is related to being contained in the \emph{cone} of certain generator vectors, well-studied problem that can be characterized algorithmically by finding the set of \emph{facet-defining inequalities} of the cone via the \emph{double description method}. That said, obtaining verifiable certificates can be trickier. Thus in order to provide a rigorous mathematical proof, and given the rather symmetric nature of partitioned families, we directly characterize the \emph{span} of our generator vectors and use this to deduce the congruence relations for the associated lattice and facet-defining inequalities for the associated cone.
    
    \item {\bf Using Refined Absorption for Partitioned Designs:} Crucial then to our proof is to have an absorption method that works in these partitioned settings. While iterative absorption could perhaps be used in an ad-hoc manner for decomposing into each different partitioned graph, we thankfully recently~\cite{DP26} extended our method of refined absorption both to the partitioned setting and for decomposing into \emph{families} of partitioned graphs. More than that, we created a \emph{black-box} decomposition theorem that only requires the existence of a robust (i.e.~low-weight and seeded) fractional decomposition and the existence of many absorbers (and indeed just the existence of two types of gadgets, \emph{anti-edges} and \emph{hinges}). The former is obtained by our fractional stability result. The latter by our absorber stability analysis. By further separating into cases based on when the standard hinges do not exist in abundance, and restricting to use only the decomposition theorem on the more restricted set of partitioned triangles, we guarantee that in our specialized cases that hinges always exist as desired. This comes at the cost of needing to understand divisibility for each of these cases and `cleaning' the graph to said divisibility conditions as follows.
    
    \item {\bf Cleaning to Divisibility:} All of our partitioned cases require a certain \emph{cross-divisibility}, specifically the number of \emph{cross} edges is exactly twice the number of \emph{interior} edges. Since our graphs might have more interior edges, this necessitates removing some number of interior triangles to achieve this condition. In the cases where a side contains a giant interior clique, this is not so problematic. However, when both sides are close to complete bipartite graphs, this is far more subtle. When we are far from that case, Tur\'an stability also provides the requisite number of interior triangles.
    
    \item {\bf Decomposing `Oddballs' via P\'osa's Theorem:} Another complication arises since in order to guarantee a robust fractional decomposition into our partitioned graphs, we must first decompose the neighborhoods of \emph{oddball} vertices - those vertices whose neighborhoods do not neatly fit the desired partitioned structure. This can be done by a clever use of P\'osa's Theorem about Hamilton cycles. However care must be taken to ensure that such oddball decompositions still achieve the cross divisibility condition. Thankfully in the three cases where we have access to many interior triangles, we can simply remove the appropriate number of interior triangles \emph{after} decomposing the oddball neighborhoods. The remaining case where both sides are near to complete bipartite graphs is much trickier and requires a subtle interplay of its more restrictive divisibility conditions (where the number of edges of each of the three matchings of the four parts must all be equal) and decomposition of oddball neighborhoods.
    
    \item {\bf Parity Fixers:} One technicality we glossed over above for the two other specialized cases (where at least one side is the disjoint union of two roughly equal cliques) is that these cases each require a separate parity condition to be met in addition to the cross divisibility condition. For each of the two cases, we then prove the existence of a \emph{parity fixer} for that case, namely a $4$-wheel whose spokes are from just the right parts (these just barely exist and only since $\delta(G)\ge \frac{3}{4}n$). For each case, we then first set aside the parity fixer before doing anything else. Then when we are prepared to invoke our black-box decomposition theorem, we check whether the parity condition is satisfied. If yes, then we proceed; if no, then we `flip' the parity-fixer to use the other triangle decomposition instead - thus correcting the parity.
\end{itemize}

\section{Further Applications of the Fractional Nash-Williams' Conjecture}\label{s:Applications}

In this section, we discuss the corollaries and applications that arise from us confirming the validity of the Fractional Nash-Williams' Conjecture. In the remainder of the introduction, we discuss further applications of the Fractional Nash-Williams' Conjecture. 

\subsection{High Girth Triangle Decompositions}

First let us return to Steiner triple systems and an old question of Erd\H{o}s on the existence of high girth Steiner triple systems as we now describe. 

Another fundamental problem in extremal design theory that has received much attention in recent years is the existence of high-girth Steiner triple systems. We define a \emph{$(j,i)$-configuration} in a $K_3$-decomposition as a set of $i$ triangles spanning at most $j$ vertices.  Erd\H{o}s observed that every Steiner triple system on $n$ elements contains an $(i+3,i)$-configuration for every $1 \le i \le n-3$.  The \textit{girth} of a $K_3$-decomposition as the smallest integer $g$ such that there are $g$ triangles in the decomposition spanning at most $g+2$ vertices.  In 1973 Erd\H{o}s conjectured in~\cite{E73} that for each integer $g$ and for all large enough $n\equiv 1,3 \pmod 6$, there exists a Steiner triple system on $n$ elements of girth at least $g \geq 3$. In 2019, Bohman and Warnke~\cite{BW19} and independently published in 2020, Glock, K\"uhn, Lo, and Osthus~\cite{GKLO20} proved the existence of approximate high girth Steiner triple systems. In 2024, Kwan, Sah, Sawhney, and Simkin \cite{KSSS2024STS} using iterative absorption confirmed Erd\H{o}s' High Girth Steiner Triple System Conjecture.

The common generalization of  Erd\H{o}s' High Girth Steiner Triple System Conjecture and Nash-Williams' Conjecture was dubbed the ``Erd\H{o}s meets Nash-Williams' Conjecture'' by Glock, K\"uhn, and Osthus~\cite{GKO20Survey} in 2021. 

\begin{conj}[``Erd\H{o}s meets Nash-Williams''~\cite{GKO20Survey}]\label{conj:E-NW}
    For every integer $g$, any sufficiently large $K_3$-divisible graph $G$ on $n$ vertices with minimum degree $\delta(G)\geq \frac{3n}{4}$ admits a $K_3$-decomposition with girth at least $g$.
\end{conj}

We note this is the common generalization of Erd\H{o}s' Conjecture (wherein $G=K_n$) and Nash-Williams' Conjecture (wherein we do not require a high girth triangle decomposition). Recently, Delcourt, Henderson, Lesgourgues, and Postle~\cite{DHLP25} were able to relate this to $\delta_{K_3}^*$ using the method of refined absorption (developed in~\cite{DPI}) as follows (thus generalizing Theorem~\ref{thm:TriangleDecomp} to high-girth triangle decompositions).

\begin{thm}[Delcourt, Henderson, Lesgourgues, and Postle~\cite{DHLP25}]\label{thm:Main_ErdosNashWilliams}
For every integer $g \geq 3$ and real $\varepsilon > 0$, any sufficiently large  $K_3$-divisible graph $G$ on $n$ vertices with minimum degree
\[\delta(G) \geq \Big(\max\Big\{\delta^*_{K_{3}}, \frac{3}{4}\Big\}+\varepsilon\Big)\cdot n\]
admits a $K_3$-decomposition with girth at least $g$.
\end{thm}
Thus, we immediately obtain the following corollary from our main result and Theorem~\ref{thm:Main_ErdosNashWilliams}.

\begin{cor}\label{cor:HighGirth}
    For every integer $g \geq 3$ and real $\varepsilon > 0$, every sufficiently large $K_3$-divisible graph $G$ on $n$ vertices with minimum degree $\delta(G) \ge \left(  \frac{3}{4} + \varepsilon \right) \cdot n$ admits a $K_3$-decomposition with girth at least $g$. 
\end{cor}

\subsection{Decomposition Thresholds for Graphs with Chromatic Number 3}

Another major area of design theory is the study of $F$-decompositions for a graph $F$. Indeed, one of the most studied objects in design theory is an $(n,q,2)$-Steiner system which is equivalent to a $K_q$-decomposition of $K_n$. In the 1970s, Wilson~\cite{W75} revolutionized by design theory by showing that $K_n$ admits a $K_q$-decomposition provided it satisfies the necessary divisibility conditions (this was the graph case of the infamous Existence Conjecture recently proved by Keevash~\cite{K14}); more generally, Wilson showed that for any graph~$F$, every sufficiently large $F$-divisible complete graph admits an $F$-decomposition, where a graph $G$ is said to be \emph{$F$-divisible} if $e(F)~|~e(G)$ and the greatest common divisor of the degrees of the vertices in~$F$, denoted by~$\gcd(F)$, divides $d_G(v)$ for every vertex $v$ of $G$.  

One might naturally wonder about generalizing Nash-Williams' Conjecture to other graphs $F$. To this end, Glock, K\"{u}hn, Lo, Montgomery, Osthus~\cite{GKLMO19} showed that the existence of $F$-decompositions is related to the existence of fractional $F$-decompositions as follows. A \emph{fractional $F$-decomposition} of $G$ is an assignment of non-negative weights to each copy of $F$ in $G$ such that the sum of the weights along each edge is precisely one. The \emph{fractional $F$-decomposition threshold} $\delta^*_{F}$ is defined as is the infimum of all real numbers $c$ such that every sufficiently large graph $G$ with minimum degree at least $c\cdot v(G)$ has a fractional $F$-decomposition. (Note the sufficiently large condition can be dropped when $F$ is a clique due to a blowup trick, see~\cite{DHLP25b}). 

A folklore generalization of Nash-Williams' Conjecture credited to Gustavsson~\cite{gustavsson1991decompositions}, 
asserts that for each integer $q\ge 3$, every sufficiently large $K_q$-divisible graph $G$ on $n$ vertices with $\delta(G)\geq \big(1-\frac{1}{q+1}\big)n$ admits a $K_q$-decomposition. In 1991, Gustavsson~\cite{gustavsson1991decompositions} in his thesis provided a construction to show that this minimum degree value would be tight (this construction was reiterated by Yuster~\cite{Yuster2005asymptotically} in 2005). Generalizing Theorem~\ref{thm:TriangleDecomp}, Glock, K\"{u}hn, Lo, Montgomery, Osthus~\cite{GKLMO19} tied the minimum degree threshold for $K_q$-decompositions to its fractional threshold as follows.

\begin{thm}[Glock, K\"{u}hn, Lo, Montgomery, and Osthus~\cite{GKLMO19}]\label{thm:CliqueDecomp}
Let $q\ge 3$ and $\varepsilon > 0$.  Any sufficiently large, $K_q$-divisible graph $G$ on $n$ vertices with minimum degree
$$\delta(G) \geq \left(\max\left\{\delta^*_{K_{q}}, \frac{q}{q+1}\right\}+\varepsilon\right)\cdot n$$
admits a $K_q$-decomposition.
\end{thm}

Surprisingly, Delcourt, Henderson, Lesgourgues, and Postle~\cite{DHLP25b} recently constructed infinitely many counterexamples to the folklore conjecture for each $q \geq 4$; indeed, they even showed that the fractional relaxation was false by showing that $\delta^*_{K_q} > 1-\frac{1}{q+1}$ for all $q\ge 4$ and even more that for every $\varepsilon>0$, we have $\delta^*_{K_q} > 1-\frac{1}{\left(\frac{1+\sqrt{2}}{2}-\varepsilon\right)\cdot q}$ provided $q$ is large enough (in terms of $\varepsilon$). This left $q=3$, corresponding to the original Nash-Williams' Conjecture, as the remaining open case - which we resolve in the affirmative in this paper thus showing that triangles are indeed special! 

Of course, determining the correct value of $\delta^*_{K_q}$ would still be of interest.  Currently the best-known upper bound on the fractional decomposition threshold of larger cliques is by Montgomery~\cite{montgomery_fractional_2019} who proved that for each $q \geq 4$, $\delta^*_{K_q} \leq 1 - \frac{1}{100 \cdot q}$.

Now we discuss what is known for general $F$. In 2012, Yuster~\cite{Y12} showed for any graph $F$ that $\delta_F^* \le \delta_{K_{\chi(F)}}$, where $\chi(F)$ denotes the chromatic number of $F$. In 2019 using Yuster's work and iterative absorption, Glock, K\"{u}hn, Lo, Montgomery, and Osthus~\cite{GKLMO19} generalized Wilson's theorem and their $K_q$-decomposition result mentioned above as follows.

\begin{thm}[Glock, K\"{u}hn, Lo, Montgomery, and Osthus~\cite{GKLMO19}]\label{thm:FDecomp}
Let $q\ge 3$, $\varepsilon > 0$, and $F$ be a graph with $q:= \chi(F)$.  Any sufficiently large, $F$-divisible graph $G$ on $n$ vertices with minimum degree
$$\delta(G) \geq \left(\max\left\{\delta^*_{K_{q}}, \frac{q}{q+1}\right\}+\varepsilon\right)\cdot n$$
admits an $F$-decomposition.
\end{thm}

Combining this with our main result, we obtain the following.

\begin{cor}
For each $\varepsilon > 0$ and graph $F$ with $\chi(F)=3$, the following holds for sufficiently large $n$: If $G$ is a $K_3$-divisible graph on $n$ vertices with $\delta(G)\geq \left(  \frac{3}{4} + \varepsilon\right)\cdot n$, then $G$ admits an $F$-decomposition.
\end{cor}

\subsection{Minimum Degree Thresholds for Cycle Packing}
Another very storied area of design theory is decomposing graphs into families of cycles.  In 1967 while attending a meeting at the Oberwolfach Research Institute for Mathematics, inspired by the rotating assigned seating arrangements during meals Ringel posed what is now known in the literature as the famous Oberwolfach Problem, see for instance~\cite{LR91}.

\begin{question}[Oberwolfach Problem] Let $n \in \mathbb{N}$ and $F$ be a 2-regular graph on $n$ vertices. For which (odd) $n$ and $F$ does $K_n$ decompose into edge-disjoint copies of $F$?
\end{question}

We note that the Oberwolfach Problem is a storied problem with a vast history and to that end we direct the reader to the survey of Burgess, Danziger, and Traetta~\cite{BDT22} and the paper of Glock, Joos, Kim, K\"uhn, and Osthus~\cite{GJKKO21} for more history.  
 
Summarizing the most recent relevant history, we note that Kim, K\"uhn, Osthus, and Tyomkyn~\cite{KKOT} and Ferber, Lee, and Mousset~\cite{FLM} were able to show that that such $K_n$ must contain $n/2 - o(n)$ edge-disjoint copies of any fixed 2-regular graph $F$ on $n$ vertices.  Subsequently in a celebrated result, appearing online in 2018 and published in 2021 Glock, Joos, Kim, K\"{u}hn, and Osthus~\cite{GJKKO21} resolved the Oberwolfach Problem by verifying the problem for $n$ sufficiently large.

Over 100 papers have been written on partial solutions to and variants of Ringel's original Oberwolfach Problem; one important variation is the so-called Generalized Oberwolfach Problem.

\begin{question}[Generalized Oberwolfach Problem] Let $n \in \mathbb{N}$ and $F_1, \ldots , F_{(n-1)/2}$ be (possibly distinct) 2-regular graphs each on $n$ vertices. For which (odd) $n$ and $F_1, \ldots , F_{(n-1)/2}$ does $K_n$ decompose into edge-disjoint $F_1, \ldots , F_{(n-1)/2}$?
\end{question}

In 2022 Keevash and Staden~\cite{KS} solved the Generalized Oberwolfach Problem (in particular that for large enough $n$, it is always possible to decompose); their results also extend to the settings of dense typical graphs and directed graphs. 

A natural variant of the Oberwolfach Problem is to ask about the minimum degree setting instead.  Glock, Joos, Kim, K\"{u}hn and Osthus~\cite{GJKKO21} conjectured the following generalization of Nash-Williams' Conjecture and the Oberwolfach Problem (though we note they conjectured only the approximate version in the sense of allowing an $\varepsilon$ in the minimum degree condition).

\begin{conj}[Approximate Oberwolfach Nash-Williams' Conjecture - Glock, Joos, Kim, K\"{u}hn and Osthus~\cite{GJKKO21}]\label{conj:ApproxOberwolfachNW}
For all $\varepsilon > 0$, the following holds for sufficiently large $n$: If $G$ be a $d$-regular graph on $n$ vertices with even $d \geq \left(\frac{3}{4}+\varepsilon \right)\cdot n$ and $F$ is a 2-regular graph on $n$ vertices, then $G$ decomposes into edge-disjoint copies of $F$.    
\end{conj} 

Note that if $F$ is a collection of $\frac{n}{3}$ triangles the fraction $\frac{3}{4}$ would be tight (due to the extremal examples for Nash-Williams' Conjecture), but if $F$ is a Hamilton cycle the fraction can be lowered to $\frac{1}{2}$~\cite{CKLOT16}.

In 2019, Condon, Kim, K\"uhn, and Osthus~\cite{CKKO19} showed the previously best known progress as follows.

\begin{thm}[Condon, Kim, K\"uhn, and Osthus~\cite{CKKO19}]
For all $\varepsilon > 0$, the following holds for sufficiently large $n$:  Suppose that $F$ is a collection of 2-regular $n$-vertex graphs and let $G$ be a $d$-regular $n$-vertex graph with $d \geq \Big(\max\Big\{\delta^*_{K_{3}}, \frac{3}{4}\Big\}+\varepsilon\Big)\cdot n$. If $e(F) \leq (1 -\varepsilon)\cdot e(G)$, then $G$ decomposes into edge-disjoint copies of $F$.
\end{thm}

As mentioned above an important special case of Conjecture~\ref{conj:ApproxOberwolfachNW} (and indeed the tight example) is when $F$ is a triangle factor. The full version would appropriately be called the Resolvable Nash-Williams' Conjecture as follows.

\begin{conj}[Resolvable Nash-Williams' Conjecture]\label{conj:ResolvableNW}
Let  $G$ be $d$-regular graph on $n$ vertices where $d\ge \frac{3}{4}n$ is even. If $n$ is sufficiently large and $3~|~n$, then $G$ admits a resolvable triangle-decomposition.
\end{conj}

In another forthcoming paper~\cite{DPresolve}, we combine the lemmas from this paper with our general black-box decomposition theorem from our forthcoming paper~\cite{DP26} on refined absorption in partitioned hypergraphs to prove Conjecture~\ref{conj:ResolvableNW} in full.

Returning to the results of this paper, we have that combined with Theorem~\ref{thm:TriangleDecomp}, our result gives the following progress on the generalization of the Oberwolfach problem and Nash-Williams' Conjecture.

\begin{cor}
For all $\varepsilon > 0$, the following holds for sufficiently large $n$:  Suppose that $F$ is a collection of 2-regular $n$-vertex graphs and let $G$ be a $d$-regular $n$-vertex graph with $d \geq \left(  \frac{3}{4} + \varepsilon\right)\cdot n$. If $e(F) \leq (1 -\varepsilon)\cdot e(G)$, then $G$ decomposes into edge-disjoint copies of $F$.
\end{cor}

\subsection{Maximum Triangle Packing in Minimum Degree Graphs}

First we recall some definitions.  Let $F$ be a graph; an \emph{$F$-packing} of a graph $G$ is a set of edge-disjoint subgraphs of $G$ each isomorphic to $F$. A \emph{fractional $F$-packing} of a graph $G$ is an assignment of nonnegative weights $\phi$ to each copy of $F$ in $G$ such that for each edge $e\in E(G)$, we have that $\sum_{F'\cong F: e\in E(F')} \phi(F) \le 1$. The \emph{$F$-packing number of $G$}, denoted $\nu_F(G)$ is the maximum number of copies of $F$ in an $F$-packing of $G$. The \emph{fractional $F$-packing number of $G$}, denoted $\nu_F^*(G)$ is the maximum of $\sum_{F\cong F'} \phi(F')$ is a fractional $F$-packing $\phi$ of $G$. 

A seminal result of Haxell and R\"odl~\cite{HR01} from 2001 showed that these two packing numbers are close as follows.

\begin{thm}[Haxell and R\"odl~\cite{HR01}]\label{thm:HR}
Let $F$ be a graph. If $G$ is a graph on $n$ vertices, then 
$$\nu_F(G) \ge \nu_F^*(G) - o(n^2).$$
\end{thm}

We note that Yuster~\cite{Yuster05} in 2005 provided a shorter proof of the above result and extended it to packings of families $\mc{F}$ of graphs. Another seminal result of Yuster~\cite{Y12} is that the packing number of a $q$-clique in a graph $G$ yields an asymptotically equivalent bound for the packing number of any $q$-chromatic graph as follows.

\begin{thm}[Yuster~\cite{Y12}]
Let $F$ be a graph with $\chi(F)=q$. If $G$ is a graph on $n$ vertices, then 
$$\nu_F^*(G) \ge \nu_{K_q}^*(G) - o(n^2),$$
and hence by Theorem~\ref{thm:HR}, we find that 
$$\nu_F(G) \ge \nu_{K_q}(G) - o(n^2).$$
\end{thm}

\begin{remark}
For the interested reader, we also note that Farkas' Lemma~\cite{F02}, a classical result from optimization, also relates to linear programming duality. In particular for fractional $K_3$-packings, it is natural to study the \emph{fractional $K_3$-packing program} of a graph $G$, denoted $\nu_{K_3}^*(G)$, which is the maximum number of triangles that can be fractionally packed into $G$. This is a linear program whose dual is the \emph{fractional $K_3$-covering program} of $G$, denoted $\tau_{K_3}^*(G)$, which is the minimum sum of an edge weighting that covers all triangles of $G$. This program can be shifted by $-1$ for each edge to yield (via linear programming duality) that the value of a maximum fractional $K_3$-packing of $G$ is equal to the $e(G)$ plus the minimum sum of an edge-weighting of $E(G)$ where every triangle is non-negative and every edge has value at least $-1$. Now the relation to Farkas' Lemma is apparent as $G$ has a fractional $K_3$-decomposition if and only if the maximum fractional $K_3$-packing equals $e(G)$. [The value at least $-1$ condition may then be dropped due to scaling.] See~\cite{KLMP19, BHKNW24, BW25} for recent papers utilizing the fractional $K_3$-packing program.
\end{remark}

It is natural to wonder about the (fractional) packing numbers of a graph $G$ on $n$ vertices with prescribed minimum degree. Nash-Williams' Conjecture asserts that $\delta(G)\ge 3n/4$ suffices to guarantee a triangle decomposition (i.e.~$\nu_{K_3}(G)\ge \frac{e(G)}{3}$) provided $G$ is $K_3$-divisible. Of course, for every $d \le \frac{n}{2}$, there exist $d$-regular bipartite graphs on $n$ vertices with no triangles. It remains then to consider $d\in [\frac{n}{2},\frac{3n}{4}]$. Modifying the extremal example for fractional Nash-Williams' by deleting $\frac{3n}{4}-d$ matchings using only edges inside each part yields a construction with at most $\frac{dn}{2} - \frac{n^2}{4}$ triangles (since that is the total number of edges inside each part and if the graphs $H_1$ and $H_2$ are triangle-free, this also is an upper bound for the maximum number of edge-disjoint triangles). In 2012, Yuster~\cite{Y12} conjectured that this upper bound is optimal. In 2014, Yuster~\cite{Y14}, showed that an approximately matching lower bound (namely up to a $-o(n^2)$ term) would be implied by the Fractional Nash-Williams' Conjecture. Hence as a corollary of our proof of the Fractional Nash-Williams' Conjecture, we have the following theorem. Indeed, it follows from Yuster's work that this more generally holds for all graphs $F$ of chromatic number at most $3$ as follows.

\begin{thm}[Corollary of Fractional Nash-Williams' Conjecture via Yuster's result~\cite{Y14}]
Let $F$ be a graph with $\chi(F)\le 3$. For every $c\in [\frac{1}{2},\frac{3}{4}]$, if $G$ is a graph on $n$ vertices with $\delta(G)\ge cn$, then 
$$\nu_{F}(G) \ge \frac{e(G)}{3}-\left(\frac{1}{4}-\frac{c}{3}\right)n^2 - o(n^2) \ge \left(\frac{c}{2}-\frac{1}{4}\right )n^2 - o(n^2).$$
\end{thm}

Yuster proved this implication by showing that the Fractional Nash-Williams' Conjecture implies that for every $c\in [\frac{1}{2},\frac{3}{4}]$, if $G$ is a graph on $n$ vertices with $\delta(G)\ge cn$, then $\nu_3^*(G) \ge \frac{e(G)}{3}-\left(\frac{1}{4}-\frac{c}{3}\right)n^2 - o(n^2)\ge (\frac{c}{2}-\frac{1}{4})n^2 - o(n^2)$. Combined with the theorem of Haxell and R\"odl this yields the result above. We in fact can derive the version of Yuster's fractional result for $F=K_3$ without the $o(n^2)$ term from our specific lemmas on Fractional Nash-Williams' Conjecture as follows.

\begin{thm}\label{thm:StrongerYuster}
For every $c\in [\frac{1}{2},\frac{3}{4}]$, if $G$ is a graph on $n$ vertices with $\delta(G)\ge cn$, then 
$$\nu_3^*(G) \ge \frac{e(G)}{3}-\left(\frac{1}{4}-\frac{c}{3}\right)n^2 \ge \left(\frac{c}{2}-\frac{1}{4}\right)n^2.$$
\end{thm}

\begin{proof}
Let $\phi$ be a triangle-covered edge-weighting of $G$ (see Section~\ref{s:FracNW} for these definitions) such that $\phi(f)\ge -1$ for all edges $f\in E(G)$. Using Proposition~\ref{prop:NonZeroUntied}, it suffices to consider $\phi$ that is nonzero and untied. Apply the discharging rules from Section~\ref{s:FracNW} to $G$.   

Let $\varepsilon:= \frac{3}{4}-c$. By Lemma~\ref{lem:MinDegFracLowerBound}, we have that for every negative edge $f$, 
$$\textrm{ch}_F(f) \ge \textrm{ch}_F'(f) \ge -4\varepsilon \cdot |\phi(f)| \ge -4\varepsilon,$$
where we used that $\phi(f)\ge -1$ and hence $|\phi(f)|\le 1$. Since $\phi$ is triangle-covered, there does not exist a triangle of $G$ with only negative edges, that is the subgraph consisting of the negative edges of $G$ is triangle-free. Hence by Mantel's Theorem, we have that there are at most $\frac{n^2}{4}$ negative edges of $(G,\phi)$. Since every positive edge $e$ of $G$ satisfies that $\textrm{ch}_F(e)\ge 0$, it follows that
\begin{align*}
\sum_{e\in E(G)} \phi(e) &= \sum_{e\in E(G)} \textrm{ch}_F(e) = \sum_{e\in E(G): \phi(e) > 0} \textrm{ch}_F(e) + \sum_{f\in E(G): \phi(f) <0} \textrm{ch}_F(f)  \\
&\ge |\{f\in E(G): \phi(f) < 0\}|\cdot(-4\varepsilon) \ge -4\varepsilon\cdot \frac{n^2}{4} = -\varepsilon\cdot n^2. 
\end{align*}
It follows by linear programming duality that 
$$3\nu_{K_3}^*(G) \ge e(G) - \varepsilon\cdot n^2.$$
Thus
\begin{align*}
\nu_{K_3}^*(G) &\ge \frac{e(G)}{3} - \frac{\frac{3}{4}-c}{3} \cdot n^2 =
\frac{e(G)}{3} - \left(\frac{1}{4}-\frac{c}{3}\right)\cdot n^2\\
&\ge \frac{cn^2/2}{3} - \left(\frac{1}{4}-\frac{c}{3}\right)\cdot n^2 = \left(\frac{c}{6} + \frac{c}{3} - \frac{1}{4}\right)n^2 = \left(\frac{c}{2}-\frac{1}{4}\right)n^2,
\end{align*}
as desired.
\end{proof}

\begin{remark}
For the interested reader, we explain now explain how one could strengthen Yuster's argument to remove the $o(n^2)$ term as follows. Yuster actually showed that it is always possible to add $n/2 + o(n)$ edges to increase the minimum degree of a graph whose minimum degree is below $3n/4$. Repeatedly applying this and assuming Fractional Nash-Williams' Conjecture then yields his implication. If we allow `half-edges' (edges with weight $1/2$), we can remove the error terms as follows. Namely, one considers $\overline{G}\times K_2$, embeds it inside a $(\frac{3}{4}-c)n$ regular bipartite graph $H$, and applies K\"onig's theorem to find $(\frac{3}{4}-c)n$ pairwise edge-disjoint perfect matchings of $H$; we then retain only the edges of the matchings in $\overline{G}\times K_2$, return them back to $\overline{G}$ and halve them (where if the same edge appears in two matchings we assign it weight one). This argument shows it is possible to add $(\frac{3}{4}-c)\cdot \frac{n}{2}$ edges and/or half-edges to raise the (weighted) minimum degree to $3n/4$. Since the Fractional Nash-Williams' Conjecture is equivalent to its more general weighted versions (by blowing up, see~\cite{DHLP25b,DHLP25}), the stronger implication then follows.  
\end{remark}

\part{Proof of Fractional Nash-Williams' Conjecture}

\section{Proof Modulo Key Lemmas}\label{s:FracNW}

\subsection{More on Farkas' Lemma for Triangle Decompositions}

We now proceed to formalize further the study of the implications of Farkas' Lemma for $K_3$-decompositions. Let us provide names to these notions. For the purposes of this paper, we make the following definitions.

\begin{definition}
An \emph{edge-weighting} $\phi$ of a graph $G$ is a function from $E(G)$ to $\mathbb{R}$. The \emph{value} of $\phi$, denoted $|\phi|$, is  $\sum_{e\in E(G)} \phi(e)$.    
\end{definition}

\begin{definition}
An \emph{edge-weighted graph} is a pair $(G,\phi)$ where $G$ is a graph and $\phi$ is an edge-weighting of $G$. We write $E_\phi(G)$ for the edges of $G$. We say an edge $e$ of $(G,\phi)$ is \emph{positive} if $\phi(e) > 0$ and \emph{negative} if $\phi(e) < 0$. We let $E^+_{\phi}(G)$ denote the set of positive edges of $(G,\phi)$ and $E^-_{\phi}(G)$ denote the set of negative edges of $(G,\phi)$. For all these, we drop the $\phi$ in the subscript if it is clear from context.
\end{definition}

\begin{definition}
We say an edge-weighting $\phi$ of a graph $G$ is \emph{triangle-covered} if $\sum_{e\in E(T)} \phi(e) \ge 0$ for each triangle $T$ of $G$.   
\end{definition} 

We note that we also naturally extend definitions for edge-weightings to edge-weighted graphs. Thus Farkas' Lemma for $K_3$-decompositions (Lemma~\ref{lem:FarkasK3}) may be restated as follows. 

\begin{lem}[Farkas Lemma for $K_3$-decompositions restated]\label{lem:FarkasK3Restated}
A graph $G$ has a fractional $K_3$-decomposition if and only if every triangle-covered edge-weighting of $G$ has nonnegative value.    
\end{lem}

Hence the Fractional Nash-Williams' Conjecture (Conjecture~\ref{conj:FracNW} in this paper, proved by our Theorem~\ref{thm:fracNW}) is equivalent to the following.

\begin{conj}[Fractional Nash-Williams' Conjecture equivalent form]\label{conj:fracNW2}
If $(G,\phi)$ is a triangle-covered edge-weighted graph on $n$ vertices with
$\delta(G) \geq \frac{3}{4}n$, then $|\phi|\ge 0$.
\end{conj}

Thus for our proof of the above, we will assume for a contradiction that there exists a triangle-covered edge-weighting with negative value and then use discharging rules to move charge (while preserving total charge) such that every edge has nonnegative weight, yielding a contradiction. 

Our first observation is that we will be able to assume there are no edges of zero weight and all edges have distinct weights as follows.

\begin{definition}
An edge-weighting $\phi$ of a graph $G$ is \emph{nonzero} if $\phi(e)\ne 0$ for all $e\in E(G)$ and is \emph{untied} if $\phi(e)\ne \phi(f)$ for all distinct $e, f \in E(G)$.   
\end{definition}

\begin{proposition}\label{prop:NonZeroUntied}
If $\phi$ is a triangle-covered edge-weighting of a graph $G$ with $e(G)\ge 1$ and $\varepsilon > 0$ is real, then there exists a nonzero untied triangle-covered edge-weighting $\phi'$ of $G$ such that $\phi(e)\le \phi'(e)$ for all $e\in E(G)$ and $|\phi|\le |\phi'|\le |\phi|+\varepsilon$.    
\end{proposition}
\begin{proof}
Let $\varepsilon':= \min \{|\phi(e)-\phi(f)|:e,f\in E(G): \phi(e)\ne \phi(f)\}$ (where we let this be infinite if there do not exist two edges of distinct weight). Let $\varepsilon'':=\min\{\varepsilon,\varepsilon'\}\cdot \frac{1}{3\cdot e(G)}$. Let $m:=e(G)$ and let $e_1,\ldots,e_m$ be an (arbitrary) ordering of the edges of $G$. For each $i\in [m]$, we define $\phi'(e_i)$ as $\phi(e_i) + 2i\cdot \varepsilon''$ unless $\phi(e_i) = -2i\cdot \varepsilon''$ in which case we define $\phi'(e_i)$ as $\phi(e_i) + (2i+1)\cdot \varepsilon''$. Thus $\phi'(e)\ge \phi(e)$ for all $e\in E(G)$ and hence $\phi'$ is triangle-covered. Now $\phi'$ is nonzero and $|\phi'| \le |\phi|+\varepsilon$ by construction. Finally $\phi'$ is untied since otherwise there exist $i\ne j$ with $\phi'(e_i)=\phi'(e_j)$ from which it follows that $|\phi(e_i)-\phi(e_j)| \ge \varepsilon''$ by construction but also that $|\phi(e_i)-\phi(e_j)|\le (2m+1)\varepsilon''$ contradicting the definitions of $\varepsilon'$ and $\varepsilon''$.
\end{proof}

In light of Proposition~\ref{prop:NonZeroUntied} above, the Fractional Nash-Williams' Conjecture is now equivalent to the following conjecture (since then if there exist a triangle-covered edge-weighting of negative value there would also exist a nonzero untied triangle-covered edge-weighting of negative value). 

\begin{conj}[Fractional Nash-Williams' Conjecture second equivalent form]\label{conj:fracNW3}
If $(G,\phi)$ is a nonzero untied triangle-covered edge-weighted graph on $n$ vertices with
$\delta(G) \geq \frac{3}{4}n$, then $|\phi|\ge 0$.
\end{conj}

\subsection{Discharging and Strict Sponsorship}

To state our discharging rules, we require some additional notions as follows.

\newcommand{\fat}{obtuse}
\newcommand{\skinny}{acute}

\begin{definition}
Let $(G,\phi)$ be a nonzero untied triangle-covered edge-weighted graph. A triangle $T$ of $G$ is \emph{\fat} under $\phi$ if $T$ has exactly one positive edge and is \emph{\skinny} if $T$ has exactly two positive edges. An \fat~triangle $T$ with positive edges $e_1, e_2$ and negative edge $f$ is \emph{sponsored} if $\max\{\phi(e_1),\phi(e_2)\} \ge -\phi(f)$ and in such case we say the \emph{sponsor} of $T$ is the edge $e_i$ such that $\phi(e_i)=\max\{\phi(e_1),\phi(e_2)\}$ (note such $i$ is unique since $\phi$ is untied); we say $T$ is \emph{unsponsored} if $T$ is not sponsored. 
\end{definition}

We note that such a $\phi$ may have a triangle with $3$ positive edges but we do not utilize such triangles for discharging and so do not require a name for them. Meanwhile since such a $\phi$ is triangle-covered, we crucially have that there is no triangle with $3$ negative edges.

We next require the following definition.

\begin{definition}
Let $G$ be a graph and $e$ an edge of $G$. We say an edge $f$ is \emph{triangle-incident} to $e$ if $f$ and $e$ are in a triangle together and we let $T_G(e)$ denote the set of edges triangle-incident to $e$. If $\phi$ is an edge-weighting of $G$, we further denote $\{f\in T_G(e): \phi(f) < 0\}$ as $T^{-}_{G,\phi}(e)$ and $\{f\in T_G(e): \phi(f) > 0\}$ as $T^+_{G,\phi}(e)$. [We drop the subscripts if $G$ and $\phi$ are clear from context.]    
\end{definition}
We are now prepared to define our demand function as follows.

\newcommand{\dem}{\textrm{dem}}

\begin{definition}[Strict Sponsorship]
Let $\phi$ be a nonzero untied triangle-covered edge-weighting of a graph $G$. The \emph{strict sponsorship demand function} for $\phi$ is defined as, for a positive edge $e=uv$ and edge $f=uw \in T^{-}_{G,\phi}(e)$ (and letting $e':=uw$ and $T:=uvw$): 
$$\dem(e,f) :=  \begin{cases}
      |\phi(f)| & \text{if $T$ is an \fat~triangle}\\
       |\phi(f)| \cdot \frac{ \phi(e) }{\phi(e) + \phi(e')}& \text{if $T$ is an unsponsored \skinny~triangle}\\
       |\phi(f)| & \text{if $T$ is a sponsored \skinny~triangle and e is the sponsor of $T$}\\
      0 & \text{if $T$ is a sponsored \skinny~triangle and e is not the sponsor of $T$}
    \end{cases}$$
\end{definition}

Define discharging rules. 

\noindent {\bf Discharging Rule (Demand Proportional)}: Let $G$ be a graph and $\phi$ a nonzero untied triangle-covered edge-weighting of $G$. If $e$ is a positive edge of $G$ and $f\in T^{-}_{G,\phi}(e)$, then $e$ sends $f$ the following amount of charge: $$\textrm{ch}(e,f):= \phi(e) \cdot \frac{\dem(e,f)}{\sum_{f' \in T^{-}_{G,\phi}(e)} \dem(e,f')}.$$
For a positive edge $e$, the \emph{final charge of $e$}, denoted $\textrm{ch}_{\textrm{F}}(e)$, is thus:
$$\textrm{ch}_{\textrm{F}}(e) := \phi(e) - \sum_{f\in T^{-}_{G,\phi}(e)} \textrm{ch}(e,f).$$
For a negative edge $f$, the \emph{final charge of $f$}, denoted $\textrm{ch}_{\textrm{F}}(f)$, is thus:
$$\textrm{ch}_{\textrm{F}}(f) := \phi(f) + \sum_{e\in T^{+}_{G,\phi}(f)} \textrm{ch}(e,f).$$

Note then that by definition the sum of charges is preserved, that is $\sum_{e\in E(G)} \textrm{ch}_F(e) = |\phi|$. Furthermore, it follows from the definition that $\textrm{ch}_F(e) \ge 0$ for every positive edge $e$ (in fact, it will equal $0$ if $T^-_{G,\phi}(e)\ne \emptyset$ but remains positive otherwise).  Noting these facts, we have that the Fractional Nash-Williams' Conjecture (specifically in the form of Conjecture~\ref{conj:fracNW3} above) is implied by the following stronger conjecture.

\begin{conj}[implies Fractional Nash-Williams' Conjecture]\label{conj:FracNW4}
If $(G,\phi)$ is a nonzero untied triangle-covered edge-weighted graph on $n$ vertices with
$\delta(G) \geq \frac{3}{4}n$ and $f$ is a negative edge of $G$, then $\textrm{ch}_F(f)\ge 0$.    
\end{conj}

\subsection{The Linear Approximation Lower Bound}

Now we proceed with the linear approximation lower bound mentioned in the proof ideas section. To that end, we make the following definitions.

\newcommand{\demh}{\widehat{\textrm{dem}}}

\begin{definition}
For a positive edge $e$ and $f' \in T^{-}_{G,\phi}(e)$, we define the \emph{normalized demand of $f'$ for $e$}, denoted $\demh(e,f')$, as $\frac{\dem(e,f')}{\phi(e)}$.    
\end{definition}

\begin{definition}[Linear Approximation Lower Bound]\label{def:chargeprime} Let $(G,\phi)$ be a nonzero untied triangle-covered edge-weighted graph on $n$ vertices. If $e$ is a positive edge of $G$ and $f\in T^{-}_{G,\phi}(e)$, then we define the following: $$\textrm{ch}'(e,f):= \frac{4}{n}\cdot \dem(e,f)\left(1- \frac{1}{n}\cdot \sum_{f' \in T^{-}_{G,\phi}(e)} \demh(e,f')\right),$$
and for a negative edge $f$, we define
$$\textrm{ch}'_F(f) := \phi(f) + \sum_{e\in T^+_{G,\phi}(f)} \textrm{ch}'(e,f).$$
\end{definition}

A simple calculation yields that this is a lower bound as follows.

\begin{proposition}Let $(G,\phi)$ be a nonzero untied triangle-covered edge-weighted graph on $n$ vertices. If $e$ is a positive edge of $G$ and $f\in T^{-}_{G,\phi}(e)$, then
$$\textrm{ch}(e,f) \ge \textrm{ch}'(e,f)$$ and hence $\textrm{ch}_F(f) \ge \textrm{ch}'_F(f)$.    
\end{proposition}
\begin{proof}
This follows as
\begin{align*}
\textrm{ch}(e,f) &:= \dem(e,f) \cdot \frac{\phi(e)}{\sum_{f' \in T^{-}_{G,\phi}(e)} \dem(e,f')} \\
&= \frac{2}{n}\cdot \dem(e,f) \cdot \frac{1}{ \frac{2}{n}\cdot \sum_{f' \in T^{-}_{G,\phi}(e)} \frac{\dem(e,f')}{\phi(e)}} \\
&= \frac{2}{n}\cdot \dem(e,f) \cdot \frac{1}{1 + \left(\frac{2}{n}\cdot \sum_{f' \in T^{-}_{G,\phi}(e)} \demh(e,f')-1\right)}\\
&\ge \frac{2}{n}\cdot \dem(e,f) \cdot \left(1-\left(\frac{2}{n}\cdot \sum_{f' \in T^{-}_{G,\phi}(e)} \demh(e,f')-1\right)\right) = \textrm{ch}'(e,f),\end{align*}
where for the inequality we used that $\frac{1}{1+x} \ge 1-x$ for all real numbers $x$. 
\end{proof}

In light of the proposition above, Conjecture~\ref{conj:FracNW4} is implied by the following stronger conjecture.

\begin{conj}[implies Fractional Nash-Williams' Conjecture]\label{conj:FracNW5}
If $(G,\phi)$ is a nonzero untied triangle-covered edge-weighted graph on $n$ vertices with $\delta(G) \geq \frac{3}{4}n$ and $f$ is a negative edge of $G$, then $\textrm{ch}_F'(f)\ge 0$.    
\end{conj}

\subsection{Cleaning}

To simplify the analysis that will come later in studying the value, we first introduce some further notation as follows.

\begin{definition}
Let $(G,\phi)$ be a nonzero untied triangle-covered edge-weighted graph. For a negative edge $f\in E^-_{\phi}(G)$, a positive edge $e\in T^{+}_{G,\phi}(f)$ and a negative edge $f' \in T^{-}_{G,\phi}(e)$, we define the \emph{normalized demand of $f'$ for $e$ with respect to $f$} as
$$\demh_f(e,f'):= \dem(e,f)\cdot \demh(e,f').$$    
\end{definition}

\begin{definition}
Let $(G,\phi)$ be a nonzero untied triangle-covered edge-weighted graph. For a negative edge $f=uv \in E^-_{\phi}(G)$, a positive edge $e=u'v'\in T^{+}_{G,\phi}(f)$ and a vertex $x\in N(u)\cap N(v)$ such that $T=u'v'x$ is a triangle of $G$, we define the \emph{normalized demand of $T$ for $e$ with respect to $f$} as
$$\demh_f(e,T):= \sum_{f'\in E(T): \phi(f') < 0} \demh_f(e,f').$$
\end{definition}

We note then that $\demh_f(e,T)$ thus equals $\demh_f(e,u'x)+\demh_f(e,v'x)$ if $T$ is obtuse, whereas if $T$ is acute then $\demh_f(e,T)$ equals $\demh_f(e,f')$ where $f'$ is the unique negative edge. [If $T$ consists only of positive edges, then $\demh_f(e,T)=0$.] We note the following useful proposition. 

\begin{proposition}\label{prop:demandBound}
Let $(G,\phi)$ be a nonzero untied triangle-covered edge-weighted graph and $f=uv \in E^-_{\phi}(G)$. If $x\in N(u)\cap N(v)$ such that $T=u'v'x$ is a triangle of $G$, then
$$\demh_f(e,T) \le \dem(e,f) \le |\phi(f)|.$$   
\end{proposition}
\begin{proof}
First suppose $T$ is acute. Then $\demh_f(e,T)=\demh_f(e,f')$ where $f'$ is the unique negative edge of $T$. By definition $\demh_f(e,f') = \dem(e,f)\demh(e,f')$. yet $\dem(e,f)\le |\phi(f)|$ by definition and $\demh(e,f') = \frac{\dem(e,f')}{|\phi(e)|} \le 1$ by construction. Hence $\demh_f(e,T) \le \dem(e,f)\demh(e,f') \le \dem(e,f)\le |\phi(f)|$ as desired. 

So we assume $T$ is obtuse. Then $\demh_f(e,T)=\demh_f(e,u'x)+\demh_f(e,v'x)$. By definition $\demh_f(e,u'x)= \dem(e,f)\cdot \demh(e,u'x)$ and similarly $\demh_f(e,v'x)= \dem(e,f)\cdot \demh(e,v'x)$. Yet by definition $\demh(e,u'x) = \frac{\dem(e,u'x)}{\phi(e)} \le \frac{|\phi(u'x)|}{\phi(e)}$ and similarly $\demh(e,v'x) = \frac{\dem(e,v'x)}{\phi(e)} \le \frac{|\phi(v'x)|}{\phi(e)}$. Since $T$ is obtuse and $\phi$ is triangle-covered, we find that $\phi(e) \ge |\phi(u'x)|+|\phi(v'x)|$ and hence $\demh_f(e,T) \le \dem(e,f)\cdot \frac{|\phi(e,u'x)|+|\phi(e,v'x)|}{\phi(e)} \le \dem(e,f)\le |\phi(f)|$ as desired. 
\end{proof}

As mentioned in the proof ideas section, we now proceed to worst-case bound the terms arising from triangles with vertices outside the common neighborhood of the ends of $f$. To that end, we define the following.

\begin{definition}\label{def:chargedoubleprime}
 Let $G$ be a graph on $n$ vertices and $\phi$ a nonzero untied triangle-covered edge-weighting of $G$. For a negative edge $f=uv$, we define
$$\textrm{ch}''_F(f) := \sum_{e\in T^+_{G,\phi}(f)}~~\sum_{f'\in T^{-}_{G,\phi}(e):~f'\setminus e \in N(u)\cap N(v)} \demh_f(e,f').$$
\end{definition}

We note for the reader that using the notation of $\demh_f(e,T)$, the above has an equivalent formulation, namely $\textrm{ch}''_F(f) := \sum_{e\in T^+_{G,\phi}(f)}~~\sum_{x\in N(u)\cap N(v): e\cup\{x\} \textrm{is a triangle of G}} \demh_f(e,T).$

Now we can state our worst-case bound as follows.

\begin{lem}\label{lem:ChargePrime}
If $G$ is a graph on $n$ vertices , $\phi$ is a nonzero untied triangle-covered edge-weighting of $G$, and $f=uv$ is a negative edge of $G$, then
\begin{align*}
\textrm{ch}'_F(f) &= \phi(f)\left(1-\frac{4}{n}\cdot |N(u)\cap N(v)|\right) - \frac{4}{n^2} \sum_{e\in T^+_{G,\phi}(f)}~~\sum_{f'\in T^{-}_{G,\phi}(e)}~~\demh_f(e,f')
\end{align*}
Furthermore if $\delta (G) \geq \frac{3}{4}n$, then
\begin{align*}
\textrm{ch}'_F(f) &\ge \phi(f) \left(1 - \frac{3}{n} \cdot |N(u)\cap N(v)|\right) - \frac{4}{n^2} \cdot \textrm{ch}''_F(f).
\end{align*}
\end{lem}
\begin{proof}
We calculate that
\begin{align*}
\textrm{ch}'_F(f) &= \phi(f) + \sum_{e\in T^+_{G,\phi}(f)} \textrm{ch}'(e,f) \\
&= \phi(f) + \sum_{e\in T^+_{G,\phi}(f)} \frac{4}{n}\cdot \dem(e,f)\left(1- \frac{1}{n}\cdot \sum_{f' \in T^{-}_{G,\phi}(e)} \demh(e,f')\right)\\
&= \phi(f) + \frac{4}{n} \sum_{e\in T^+_{G,\phi}(f)} \dem(e,f) - \frac{4}{n^2} \sum_{e\in T^+_{G,\phi}(f)}~~\sum_{f' \in T^{-}_{G,\phi}(e)} \demh_f(e,f').
\end{align*}
Now we note that $\sum_{e\in T^+_{G,\phi}(f)} \dem(e,f) = |N(u)\cap N(v)|\cdot (-\phi(f))$ by the definition of the demand function (namely $f$ demands $-\phi(f)$ in total from every triangle containing $f$ of which there are $|N(u)\cap N(v)|$ such triangles). Thus we find that
\begin{align*}
\textrm{ch}'_F(f) &= \phi(f) \left(1 - \frac{4}{n} \cdot |N(u)\cap N(v)|\right) - \frac{4}{n^2} \sum_{e\in T^+_{G,\phi}(f)}~~\sum_{f' \in T^{-}_{G,\phi}(e)} \demh_f(e,f')
\end{align*}
as desired. 

Now we prove the inequality as follows. First we turn to bounding the terms inside the rightmost sums with $f'\setminus e \not\in N(u)\cap N(v)$. Namely for $e\in T^+_{G,\phi}(f)$, we find that
\begin{align*}\sum_{f' \in T^{-}_{G,\phi}(e): f'\setminus e \not\in N(u)\cap N(v)} \demh_f(e,f') &= \sum_{x\not\in N(u)\cap N(v): e\cup\{x\} \textrm{ is a triangle of $G$}}  \demh_f(e,e\cup\{x\}).\\
&\le |\phi(f)|\cdot |\{x\not\in N(u)\cap N(v): e\cup\{x\} \textrm{ is a triangle of $G$}\}|\\
&\le |\phi(f)|\cdot \max\{|N(u)\setminus N(v)|, |N(v)\setminus N(u)|\} \\
&\le |\phi(f)|\cdot \frac{n}{4},
\end{align*}
where for the first inequality we used that $\demh_f(e,e\cup{x})\le |\phi(f)|$ 
by Proposition~\ref{prop:demandBound}, the second inequality follows as at exactly one of $u$ or $v$ is an end of $e$, and the third inequality follows since $\delta(G)\ge \frac{3n}{4}$. 

Hence we calculate that  
\begin{align*}
\sum_{e\in T^+_{G,\phi}(f)}~~\sum_{f' \in T^{-}_{G,\phi}(e):~f'\setminus e \not\in N(u)\cap N(v)} \demh_f(e,f') \le \frac{n}{4} \cdot \sum_{e\in T^+_{G,\phi}(f)} \dem(e,f) \le \frac{n}{4}\cdot |N(u)\cap N(v)| \cdot |\phi(f)|,
\end{align*}
where we used that $|\{e\in T^{+}_{G,\phi}(f)| \le |N(u)\cap N(v)|$. Substituting this bound yields that (recalling as $f$ is negative that $|\phi(f)|=-\phi(f)$)
\begin{align*}
\textrm{ch}'_F(f) &\geq \phi(f) \left(1 - \frac{3}{n} \cdot |N(u)\cap N(v)|\right) - \frac{4}{n^2} \cdot \textrm{ch}''_F(f),
\end{align*}
as desired.
\end{proof}

In light of the lemma above, Conjecture~\ref{conj:FracNW5} now follows from the following even stronger conjecture.

\begin{conj}[implies Fractional Nash-Williams' Conjecture]\label{conj:FracNW6}
If $(G,\phi)$ is a nonzero untied triangle-covered edge-weighted graph on $n$ vertices with $\delta(G) \geq \frac{3}{4}n$ and $f$ is a negative edge of $G$, then 
$$\textrm{ch}_F''(f)\le |\phi(f)|\cdot \frac{n^2}{4} \cdot \left(\frac{3}{n}\cdot |N(u)\cap N(v)| - 1 \right).$$    
\end{conj}

\subsection{Value Formulation}

Now as to how to prove Conjecture~\ref{conj:FracNW6} above, as mentioned in the proof ideas section, the next key is to formulate the `value' of each pair of vertices in the common neighborhood of the ends of $f$ as follows. 

\begin{definition}\label{def:Value}
Let $(G,\phi)$ be a nonzero untied triangle-covered edge-weighted graph. Let $f=u_1u_2$ be a negative edge of $G$. Then for an edge $g=w_1w_2\in E(G[N(u)\cap N(v)])$, we define the \emph{value of $g$ to $f$} as 
$$\textrm{Val}_{f}(g) := \sum_{i,j \in \{1,2\}: \phi(u_iw_j) > 0} \demh_f(u_iw_j,u_iw_1w_2).$$
\end{definition}

Note that it follows from the definitions that $\textrm{ch}''_F(f) = \sum_{g\in E(G[N(u)\cap N(v)])} \textrm{Val}_f(g)$. To understand how our intermediary theorem about values arises in this context, we define the following auxiliary graph.

\begin{definition}
Let $(G,\phi)$ be a nonzero untied triangle-covered edge-weighted graph. Let $f=uv$ be a negative edge of $G$. We define the \emph{filtered common neighborhood graph of $f$ in $(G,\phi)$}, denoted $\mathrm{Filter}_{G,\phi}(f)$, as the graph $H$ with $V(H):= N_G(u)\cap N_G(v)$ and $$E(H):= E(G[N(u)\cap N(v)]) \setminus \bigg( E(G[N^+(u)\cap N^-(v)]) \cup E(G[N^-(u)\cap N^+(v)]) \cup E^+(G[N^+(u)\cap N^+(v)])\bigg).$$     
\end{definition}

The key observation is that this auxiliary graph is $K_5$-free as follows.

\begin{proposition}\label{prop:K5freeFilter}
Let $(G,\phi)$ be a nonzero untied triangle-covered edge-weighted graph. If $f=uv$ is a negative edge of $G$, then $\mathrm{Filter}_{G,\phi}(f)$ is $K_5$-free.    
\end{proposition}
\begin{proof}
Suppose not and let $S$ be a copy of $K_5$ in $H:=\mathrm{Filter}_{G,\phi}(f)$. Then $|V(S)\cap N^+(u)\cap N^-(v)| \le 1$ by definition and similarly $|V(S)\cap N^-(u)\cap N^+(v)|\le 1$. Since $V(H) = N(u)\cap N(v)$ and $\phi$ is nonzero, it follows that $|V(S)\cap N^+(u)\cap N^+(v)]| \ge 3$. But then by definition of $\mathrm{Filter}_{G,\phi}(f)$, the edges of $S$ in $G$ are negative in $\phi$ and hence $G[V(S)]$ contains a triangle consisting only of negative edges contradicting that $\phi$ is triangle-covered.      
\end{proof}

The reason for defining such an auxiliary graph is that the edges deleted from the common neighborhood by the filter graph have value zero as follows.

\begin{proposition}
Let $(G,\phi)$ be a nonzero untied triangle-covered edge-weighted graph and let $f=uv$ be a negative edge of $G$. If $g\in  E(G[N^+(u)\cap N^-(v)]) \cup E(G[N^-(u)\cap N^+(v)]) \cup E^+(G[N^+(u)\cap N^+(v)])$, then  $\textrm{Val}_{f}(g) =0$.  
\end{proposition}
\begin{proof}
Let $g=wx$. We claim that for a positive edge $e$ with $|e\cap \{u,v\}|=1$ and $|e\cap \{w,x\}|=1$, $\demh_f(e,e\cup g)=0$ and hence the proposition follows as desired. To prove this claim, we assume without loss of generality (by symmetry) that $e=uw$. (Hence $g\in E(G[N^+(u)\cap N^-(v)])\cup E^+(G[N^+(u)\cap N^+(v)])$.) First suppose $w\in N^-(v)$ and hence $g\in E(G[N^+(u)\cap N^-(v)])$. Then the triangle $e\cup g$ consists of only positive edges  and hence $\demh_f(e,e\cup g)=0$ as desired. So we assume $w\in N^+(v)$ and hence $g\in E^+(G[N^+(u)\cap N^+(v)])$.  Once again the triangle $e\cup g$ consists of only positive edges and hence $\demh_f(e,e\cup g)=0$ as desired.
\end{proof}

It is also straightforward to upper bound the value of an edge as follows. 

\begin{proposition}\label{prop:EdgeWeightUpperBound}
Let $G$ be a graph, $\phi$ be a nonzero untied triangle-covered edge-weighting of $G$, and $f=uv$ be a negative edge of $G$. If $g\in \mathrm{Filter}_{G,\phi}(f)$, then $\mathrm{Val}_f(g) \le 2 \cdot |\phi(f)|$.
\end{proposition}
\begin{proof}
Let $g=wx$. The key observation is that for a positive edge $e$ with $|e\cap \{u,v\}|=1$ and $|e\cap \{w,x\}|=1$, the sum of the the normalized demands of negative edges $f'$ in the triangle $G[V(e)\cup V(g)]$ for $e$ is at most $1$ and hence the corresponding term in $\mathrm{Val}_f(g)$ is at most $\dem(e,f)$. But then the sum over $\dem(e,f)$ for positive edges in the triangle $G[V(f)\cup \{w\}]$ is at most $|\phi(f)|$ and similarly for $G[V(f)\cup \{x\}]$. Thus the total sum is at most $2\cdot |\phi(f)|$.
\end{proof}

The facts above and Tur\'an's theorem would be enough to upper-bound the total value of edges in the common neighborhood (i.e.~$\textrm{ch}''_F(f)$) by $\frac{4}{5}\cdot |N(u)\cap N(v)|^2$ (as the edges have weight at most $2$ and there would be at most $\frac{4}{5} \cdot \frac{|N(u)\cap N(v)|^2}{2} = \frac{2}{5}\cdot |N(u)\cap N(v)|^2$ edges in a $K_5$-free graph on $|N(u)\cap N(v)|$ vertices). However, in order to prove Conjecture~\ref{conj:FracNW6}, we need a better upper bound. This will follow by proving an additional property that there is not too much weight on any triangle motivating the following definition. 

\begin{definition}
A \emph{$[0,1]$-edge-weighting} $\varphi$ of a graph $H$ is a function $\varphi: E(H)\rightarrow [0,1]$. We say $\varphi$ is \emph{triangle-thin} if for every triangle $T$ of $H$, we have that $\sum_{e\in E(T)} \varphi(e) \le 2$.     
\end{definition}

We thus call a pair $(H,\varphi)$ a $[0,1]$-edge-weighted graph if $H$ is graph and $\varphi$ is $[0,1]$-edge weighting of $H$. We prove the following in Section~\ref{s:ValueLemma} using techniques from optimization (specifically this follows under the penumbra of graph Lagrangians as invented by Motzkin and Strauss~\cite{MS65}).

\begin{lem}\label{lem:Value}
If $(H,\varphi)$ is a $K_5$-free $[0,1]$-edge-weighted triangle-thin graph on $n$ vertices, then $|\varphi| \le \frac{n^2}{4}$.    
\end{lem}

Then in Section~\ref{s:FractionalCases}, we show the strict sponsorship yields the desired triangle-thinness (appropriately scaled) as follows. But first a definition.

\begin{definition}
Let $(G,\phi)$ be a nonzero untied triangle-covered edge-weighted graph and $f=uv$ be a negative edge of $G$. If $T=xyz$ is a triangle in $\mathrm{Filter}_{G,\phi}(f)$, we define the \emph{value of $T$ to $f$} as
$$\textrm{Val}_f(T):= \textrm{Val}_f(xy)+\textrm{Val}_f(xz)+\textrm{Val}_f(yz).$$
\end{definition}

\begin{lem}\label{lem:EdgeWeightTriangleThin}
Let $(G,\phi)$ be a nonzero untied triangle-covered edge-weighted graph and $f=uv$ be a negative edge of $G$. If $T=xyz$ is a triangle in $\mathrm{Filter}_{G,\phi}(f)$, then $\textrm{Val}_f(T) \le 4\cdot |\phi(f)|$. 
\end{lem}

Now we are prepared to prove Conjecture~\ref{conj:FracNW6} assuming the above two lemmas (Lemmas~\ref{lem:Value} and~\ref{lem:EdgeWeightTriangleThin}).

\begin{proof}[Proof of Conjecture~\ref{conj:FracNW6} and hence of the Fractional Nash-Williams' Conjecture]
We let $H:= \mathrm{Filter}_{G,\phi}(f)$ and define an edge-weighting $\varphi$ on $E(H)$ as $\varphi(g) := \frac{\textrm{Val}_f(g)}{2\cdot |\phi(f)|}$ for each $g\in E(H)$. Then it follows from Proposition~\ref{prop:EdgeWeightUpperBound} that $\varphi$ is a $[0,1]$-edge-weighting of $H$. Similarly, it follows from Lemma~\ref{lem:EdgeWeightTriangleThin} that $\varphi$ is triangle-thin. Hence by Lemma~\ref{lem:Value}, we have that $|\varphi| \le \frac{v(H)^2}{4}$. Since $v(H)=|N(u)\cap N(v)|$, it follows that
$$\textrm{ch}''_F(f) = \sum_{g\in E(G[N(u)\cap N(v)])} \textrm{Val}_f(g) = 2\cdot |\phi(f)| \cdot |\varphi| \le |\phi(f)| \cdot \frac{1}{2} \cdot |N(u)\cap N(v)|^2.$$
\begin{claim}\label{cl:QuadEquation}
$$\frac{1}{2} \cdot |N(u)\cap N(v)|^2 \le \frac{n^2}{4} \left(\frac{3}{n}\cdot |N(u)\cap N(v)| - 1 \right)$$    
\end{claim}
\begin{proofclaim}
Let $d:= \frac{|N(u)\cap N(v)|}{n}$ and note the claimed inequality is equivalent to the following (by dividing by $n^2$ and substituting $d$):
$$\frac{1}{2} \cdot d^2 \le \frac{3d-1}{4}.$$
Rewriting this is equivalent to
$$2d^2-3d+1\le 0.$$
Since $2d^2-3d+1 = (2d-1)(d-1)$, we find that this inequality holds precisely for all $d\in [\frac{1}{2},1]$. Since $|N(u)\cap N(v)|\le 1$, we have that $d\le 1$. Crucially since $\delta(G)\ge \frac{3n}{4}$, we also have that $|N(u)\cap N(v)| \ge \frac{n}{2}$ and hence $d\ge \frac{1}{2}$ as desired.
\end{proofclaim}

It now follows from Claim~\ref{cl:QuadEquation} and our previous bound on $\textrm{ch}''_F(f)$ that $\textrm{ch}''(f)\le |\phi(f)|\cdot \frac{n^2}{4} \cdot \left(\frac{3}{n}\cdot |N(u)\cap N(v)| - 1 \right)$ as desired. 
\end{proof}

\subsection{Outline of Remainder of Part I}

In the remainder of this part, we provide proofs of the key lemmas. In Section~\ref{s:ValueLemma}, we prove Lemma~\ref{lem:Value}. In Section~\ref{s:FractionalCases}, we prove Lemma~\ref{lem:EdgeWeightTriangleThin}.

\section{Proof of the Value Lemma}\label{s:ValueLemma}

In this section, we prove our value lemma, Lemma~\ref{lem:Value}. In order to prove Lemma~\ref{lem:Value}, we first prove the following lemma which essentially characterizes the extremal configurations value-wise for graphs with $4$ parts and fixed part sizes (where the vertices in each part are symmetric both in regard to neighborhood and weights to neighbors).

\begin{lem}\label{lem:ValueSizesFixed}
Let $x_1\ge x_2\ge x_3\ge x_4$ be nonnegative reals such that $\sum_{i\in [4]} x_i=1$. For $i< j \in [4]$, let $a_{ij} \in [0,1]$ be real such that for all $i<j<k \in [4]$, we have $a_{ij}+a_{ik}+a_{jk}\le 2$. Then $\sum_{i<j \in [4]} a_{ij} \cdot x_i \cdot x_j \le \max\{x_1(x_2+x_3+x_4), ~(x_1+x_4)(x_2+x_3)\}$.
\end{lem}

\begin{proof}
For $i< j \in [4]$, let $a_{ij}$ satisfy the assumptions and be chosen so as to maximize $\sum_{i<j\in [4]} a_{ij}\cdot x_i\cdot x_j$ and subject to this the vector $a:=(a_{12},a_{13},a_{23},a_{14},a_{24},a_{34})$ is lexicographically maximized. Let $\varepsilon := \min\{ a_{ij},~1 - a_{ij}:~a_{ij}\in (0,1),~i<j\in [4]\}$. For $i<j<k\in [4]$, let $T_{ijk}$ denote the triangle-thin inequality $a_{ij}+a_{ik}+a_{jk}\le 2$; we say $T_{ijk}$ is \emph{tight} if $a_{ij}+a_{ik}+a_{jk}= 2$. For the various cases throughout the remainder of this proof where we derive a contradiction to the choice of the vector $a$, we simply provide a contradicting vector $a'$. While we leave the formal verification that $a'$ is $[0,1]$-edge-weighting and triangle-thin to the reader, the former follows as we add or subtract at most $\varepsilon$ and only from $a_{ij} \ne 0$ or $1$ while the latter follows from the former for $T_{ijk}$ that are not tight and for those that are as we ensure their sum is preserved in those cases.  

First we claim that $a_{12}=1$. Suppose not. If neither $T_{123}$ nor $T_{124}$ is tight, then $a':=(a_{12}+\varepsilon,a_{13},a_{23},a_{14},a_{24},a_{34})$ yields a contradiction. If $T_{123}$ is tight and $T_{124}$ is not tight, then $a_{13},a_{23}>0$ and hence $a':=(a_{12}+\varepsilon,a_{13}-\varepsilon,a_{23},a_{14},a_{24},a_{34})$ yields a contradiction. Similarly if $T_{124}$ is tight and $T_{123}$ is not tight, then $a_{14},a_{24}>0$ and hence $(a_{12}+\varepsilon,a_{13},a_{23},a_{14}-\varepsilon,a_{24},a_{34})$ yields a contradiction. So we assume that both $T_{123}$ and $T_{124}$ are tight. Since $a_{12} < 1$ this implies that all of $a_{13},a_{14},a_{23},a_{24} > 0$. Next suppose $a_{34}<1$. Then we let $a':=(a_{12}+\varepsilon,a_{13}-\frac{\varepsilon}{2},a_{23}-\frac{\varepsilon}{2},a_{14}-\frac{\varepsilon}{2},a_{24}-\frac{\varepsilon}{2},a_{34}+\varepsilon)$ which yields a contradiction. So we assume $a_{34}=1$; but then $a_{13},a_{24} < 1$. Thus $a':=(a_{12},a_{13}+\varepsilon,a_{23}-\varepsilon,a_{14}-\varepsilon,a_{24}+\varepsilon,a_{34})$ yields a contradiction. This proves the claim.

Next we claim that $a_{13}=1$. Suppose not. If neither $T_{123}$ nor $T_{134}$ is tight, then $a':=(a_{12},a_{13}+\varepsilon,a_{23},a_{14},a_{24},a_{34})$ yields a contradiction. If $T_{123}$ is tight and $T_{134}$ is not tight, then $a_{23}>0$ and hence $a':=(a_{12},a_{13}+\varepsilon,a_{23}-\varepsilon,a_{14},a_{24},a_{34})$ yields a contradiction. Similarly if $T_{134}$ is tight and $T_{123}$ is not tight, then $a_{14},a_{34}>0$ and hence $(a_{12},a_{13}+\varepsilon,a_{23},a_{14}-\varepsilon,a_{24},a_{34})$ yields a contradiction. So we assume that both $T_{123}$ and $T_{134}$ are tight. Since $a_{13} < 1$ this implies that all of $a_{12},a_{14},a_{23},a_{34} > 0$. Note that since $a_{14}>0$ and $a_{12}=1$, we find that $a_{24}<1$. But then $a':=(a_{12},a_{13}+\varepsilon,a_{23}-\varepsilon,a_{14}-\varepsilon,a_{24}+\varepsilon,a_{34}+\varepsilon)$ yields a contradiction. This proves the claim.

Since $a_{12}=a_{13}=1$, we find that $a_{23}=0$. Now we see that $\sum_{i<j\in[4]} a_{ij}\cdot x_i\cdot x_j = x_1(x_2+x_3) + (a_{14}x_1 + a_{24}x_2+a_{34}x_3)x_4 \le x_1(x_2+x_3) + (a_{14}x_1 + (1-a_{14})x_2 + (1-a_{14})x_3) x_4 = x_1(x_2+x_3) + (x_2+x_3)x_4 + a_{14}(x_1-x_2-x_3)x_4$. 

Now next suppose that $x_1\ge x_2+x_3$. Then $(x_2+x_3)x_4 + a_{14}(x_1-x_2-x_3)x_4 =  (1-a_{14})(x_2+x_3)x_4 + a_{14} x_1x_4  \le x_1x_4$ as desired (since $a_{14}\le 1$). So we assume $x_1\le x_2+x_3$. But then $(x_2+x_3)x_4 + a_{14}(x_1-x_2-x_3)x_4\le (x_2+x_3)x_4$ as desired.
\end{proof}

Next we prove a lemma that upper bounds the value when the graph has $4$ parts (a bound that works for all part sizes) using the previous lemma.

\begin{lem}\label{lem:Value4Parts}
Let $x_i$ for $i\in [4]$ be nonnegative reals such that $\sum_{i\in [4]} x_i=1$. For $i< j \in [4]$, let $a_{ij} \in [0,1]$ be real such that for all $i<j<k \in [4]$, we have $a_{ij}+a_{ik}+a_{jk}\le 2$. Then $\sum_{i<j \in [4]} a_{ij} \cdot x_i \cdot x_j \le \frac{1}{4}$.
\end{lem}
\begin{proof}
By Lemma~\ref{lem:ValueSizesFixed}, we find that $\sum_{i<j \in [4]} a_{ij} \cdot x_i \cdot x_j  \le  \max\{x_1(x_2+x_3+x_4), ~(x_1+x_4)(x_2+x_3)\}$. Since $\sum_{i\in [4]}x_i=1$ and $x_i\ge 0$ for all $i\in [4]$, it follows that $x_1(x_2+x_3+x_4) = x_1(1-x_1) \le \frac{1}{4}$ and similarly that $(x_1+x_4)(x_2+x_3) = (x_1+x_4)( 1 - (x_1+x_4)) \le \frac{1}{4}$. Thus $\sum_{i<j \in [4]} a_{ij}\cdot x_i\cdot x_j \le \frac{1}{4}$ as desired.
\end{proof}

We now use symmetrization to prove Lemma~\ref{lem:Value} (reducing down to four parts and then using Lemma~\ref{lem:Value4Parts} above). But first a definition.

\begin{definition}
Let $(H,\varphi)$ be an edge-weighted graph. We say $v\in V(H)$ is a \emph{weight-clone} of a vertex $u\in V(H)$ if $N_H(u)=N_H(v)$ and for every $x\in N_H(u)$, we have that $\varphi(ux)=\varphi(vx)$.     
\end{definition}

Let $\sim$ denote the relation where $u\sim v$ for $u,v\in V(H)$ if and only if $u$ and $v$ are weight-clones. We note that $\sim$ is an equivalence relation on $V(H)$ and call the equivalence classes of $\sim$ the \emph{weight-clone equivalence classes} of $(H,\varphi)$. 

\begin{proof}[Proof of Lemma~\ref{lem:Value}]
Let $(H,\varphi)$ be a $K_5$-free $[0,1]$-edge-weighted triangle-thin graph on $n$ vertices such that $|\varphi|$ is maximized and subject to that the number $m$ of equivalence classes of weight-clones is minimized.

First suppose that $m\ge 5$. Since $H$ is $K_5$-free, there exist two equivalence classes of weight-twins $S_1,S_2$ with $e_H(S_1,S_2)=0$. For each $i\in [2]$, let $v_i\in S_i$. For each $i\in[2]$, we define a new edge-weighted graph $(H_i,\varphi_i)$ with $V(H_i):=V(H)$, 
$$E(H_i):= \{xy\in E(H): x,y\in V(H)\setminus (S_1\cup S_2)\}~\cup~\{xy: x\in S_1\cup S_2,~y\in V(H)\setminus (S_1\cup S_2),~v_iy\in E(H)\},$$ and 
for an edge $xy\in E(H_i)$, we define $\varphi_i(xy):= \varphi(xy)$ if $x,y\in V(H)\setminus (S_1\cup S_2)$ and $\varphi_i(xy):= \varphi(v_iy)$ if $x\in S_1\cup S_2$ and $y\in V(H)\setminus (S_1\cup S_2)$. (That is $(H_i,\varphi_i)$ is obtained from $(H,\varphi)$ by replacing the vertices of $S_1\cup S_2$ with weight-clones of $v_i$.)

Note that for each $i\in [2]$, $H_i$ is $K_5$-free and $\phi_i$ is a $[0,1]$-edge-weighting of $H_i$. For each $i\in [2]$, let $s_i:= \sum_{y:~v_iy\in E(H)} \varphi(v_iy)$. It follows that $|\varphi_1| = |\varphi| + s_1\cdot |S_2| - s_2\cdot |S_2|$ and similarly $|\varphi_2| = |\varphi| + s_2\cdot |S_1|-s_1\cdot |S_1|$. If $s_1 > s_2$, then $|\varphi_1| > |\varphi|$ and hence $(H_1,\varphi_1)$ contradicts the choice of $(H,\varphi)$. Similarly if $s_2 > s_1$, then $|\varphi_2| > |\varphi|$ and $(H_2,\varphi_2)$ contradicts the choice of $(H,\varphi)$. So we assume $s_1=s_2$ and hence $|\varphi|=|\varphi_1|=|\varphi_2|$. But then $(H_1,\varphi_1)$ (and also $(H_2,\varphi_2)$ as well) contradicts the choice of $(H,\varphi)$ as $H_i$ has strictly fewer weight-clone equivalence classes than $H$. 

So we assume $m\le 4$. Let $X_1,\ldots, X_m$ be the weight-clone equivalence classes of $(H,\varphi)$. For $i\in [m]$, let $x_i:= \frac{|X_i|}{n}$. Thus $\sum_{i\in [m]} x_i=1$. For $i<j\in [4]$, let $a_{ij}:=\varphi(v_iv_j)$ for some $v_i\in X_i$ and $v_j\in X_j$ (note the value is the same no matter the choice of $v_i$ and $v_j$ as $X_i$ and $X_j$ are weight-clone equivalence classes). Then 
$$|\varphi| = \sum_{i<j\in [4]} a_{ij}\cdot |X_i|\cdot |X_j| = n^2 \sum_{i<j\in [4]} a_{ij}\cdot x_i\cdot x_j \le \frac{n^2}{4}$$
as desired where the final inequality follows from Lemma~\ref{lem:Value4Parts}.
\end{proof}

\begin{remark}
The $K_5$-free condition is necessary since if $H$ is a graph, the edge-weighting that assigns $\frac{2}{3}$ to every edge is both a $[0,1]$-edge-weighting and triangle-thin. Hence if instead the hypothesis is that $H$ is $K_m$-free for some $m > 5$, then taking the Tur\'an graph $T(n,m-1)$ with this edge-weighting shows that one cannot achieve an upper bound better than $t(n,m-1)\cdot \frac{2}{3} \approx (1-\frac{1}{m-1})\frac{n^2}{2} \cdot \frac{2}{3}$ which for $m > 5$ is $\ge \frac{4}{15}n^2 > \frac{n^2}{4}$.  
\end{remark}

\section{Proof of the Triangle-Thin Lemma for Strict Sponsorship}\label{s:FractionalCases}

In this section, we prove Lemma~\ref{lem:EdgeWeightTriangleThin}. To that end, we provide a restatement of said lemma but with the various cases and subcases enumerated (up to symmetry).

\begin{lem}[Restatement of Triangle-Thin Lemma with Cases Enumerated]\label{lem:EdgeWeightTriangleThinCases}
Let $(G,\phi)$ be a nonzero untied triangle-covered edge-weighted graph and $f=uv$ be a negative edge of $G$. If $T=xyz$ is a triangle of $\mathrm{Filter}_{G,\phi}(f)$, then $\mathrm{Val}_f(T) \le 4\cdot |\phi(f)|$ if either of the following cases hold:
\begin{enumerate}
    \item[(1)] $x\in N^+(u)\cap N^-(v), y\in N^+(v)\cap N^-(u), z\in N^+(u)\cap N^+(v)$ and any of the following subcases hold: The negative edges of $T$ are precisely
    \begin{enumerate}
        \item[(a)] $xz,yz$
        \item[(b)] $xy,xz$
        \item[(c)] $xy$
        \item[(d)] $xz$
        \item[(e)] $\emptyset$
    \end{enumerate}
    
    \item[(2)] $x\in N^+(u)\cap N^-(v)$ and $y, z\in N^+(u)\cap N^+(v)$ and any of the following subcases hold: \\
    The negative edges of $T$ are precisely
    \begin{enumerate}
        \item[(a)] $xy,yz$
        \item[(b)] $yz$.
    \end{enumerate}
\end{enumerate}
\end{lem}

The hardest case by far is (1a). The second hardest case is (2a). The third somewhat involved case is (1b). The remaining cases are fairly straightforward.

\textbf{In all cases, we assume without loss of generality that $\phi(f)=-1$ (for normalization purposes). }

To remain consistent across all the cases, we adopt the following labeling scheme. We label the 10 edges of the $K_5$ with vertices $u,v,x,y,z$ as follows:

$$g_1 := ux,~~g_2:= uz,~~g_3:= vz,~~g_4:= vy,~~g_5:= xy$$
$$g_6:= xz,~~g_7:= f,~~g_8:= yz,~~g_9:= vx,~~g_{10}:= uy.$$ 

We now label all triangles as follows:
$$T_1:=g_1g_2g_6,~~T_2:=g_2g_3g_7,~~T_3:=g_3g_4g_8,~~T_4:=g_4g_5g_9,~~T_5:=g_5g_1g_{10}$$
$$T_6:=g_1g_7g_9,~~T_7:=g_2g_8g_{10},~~T_8:=g_3g_9g_6,~~T_9:=g_4g_{10}g_7,~~T_{10}:=g_5g_6g_8$$ 

For ease of comprehending each case, we further denote $g_i$ as $e_i$ in cases where $g_i$ is positive and $g_i$ as $f_i$ in cases where $g_i$ is negative. Similarly we will write $T_i$ in cases where $T_i$ is acute and $T_i'$ in cases where $T_i$ is obtuse.

We note that the triangles that could possibly be tallied in the value sum are those with one vertex in $\{u,v\}$ and two in $\{x,y,z\}$, which are $T_1$ ($=uxz$), $T_3$ ($=vyz$), $T_4$ ($=vxy$), $T_5$( $=uxy$), $T_7$ ($=uyz$), $T_8$ ($=vxz$) where as the other triangles $T_2$ ($=uvz$), $T_6$ ($=uvx$), $T_9$ ($=uvy$), $T_{10}$ ($=xyz$) will be not be tallied. 

We also note that the sign $g_2=f$ is always negative by assumption. In Case (1), we also have that $g_1$, $g_2$, $g_3$, $g_4$ are positive and $g_{9}$, $g_{10}$ are negative (with the signs of $g_5,g_6,g_8$ varying among the subcases). In Case (2) is the same except that $g_{10}$ is instead positive. 

\subsection{Proof of Lemma\ref{lem:EdgeWeightTriangleThinCases}(1)}

Recall that in this case the positive edges include $e_1$, $e_2$, $e_3$, $e_4$ and the negative edges include $f_7,f_{9}$, $f_{10}$ (with the signs of $g_5,g_6,g_8$ varying among the subcases). 

\subsubsection{Proof of Lemma\ref{lem:EdgeWeightTriangleThinCases}(1a)}

This subcase assumes the remaining negative edges are $f_6$ and $f_8$. We thus have positive edges $e_1,e_2,e_3,e_4,e_5$ that form a $5$-cycle $e_1e_2e_3e_4e_5$ and negative edges $f_6,f_7,f_8,f_9,f_{10}$ which form a $5$-cycle  $f_6f_8f_{10}f_7f_9$. Similarly, here is a breakdown of the acute/obtuse triangles in this case:

$$\textrm{Acute Triangles:}~~T_1=e_1e_2f_6,~~T_2=e_2e_3f_7,~~T_3=e_3e_4f_8,~~T_4=e_4e_5f_9,~~T_5=e_5e_1f_{10}$$
$$\textrm{Obtuse Triangles:}~~T_6'=e_1f_7f_9,~~T_7'=e_2f_8f_{10},~~T_8'=e_3f_9f_6,~~T_9'=e_4f_{10}f_7,~~T_{10}'=e_5f_6f_8$$ 

For this subcase, we calculate that

\begin{align*}
\textrm{Val}_f(xy) &= \demh_f(e_1,T_5) + \demh_f(e_4,T_4),\\
\textrm{Val}_f(xz) &= \demh_f(e_1,T_1) + \demh_f(e_2,T_1) + \demh_f(e_3,T_8'),\\
\textrm{Val}_f(yz) &= \demh_f(e_4,T_3) + \demh_f(e_3,T_3) + \demh_f(e_2,T_7').  
\end{align*}

Note there are two terms above containing $e_i$ for each $i\in [4]$. Each term is at most $1$ by Proposition~\ref{prop:demandBound}; moreover as $\dem(e_2,f)+\dem(e_3,f)\le 1$, it follows from Proposition~\ref{prop:demandBound} that any term containing $e_2$ and any term containing $e_3$ sum to at most $1$. The following claim then follows immediately by inspection from the facts listed above.

\begin{claim}\label{cl:CycleCaseContradictions}
If any of the following hold:
\begin{itemize}
    \item[(i)] $\min\{\demh_f(e_1,T_1),~\demh_f(e_1,T_5)\}=0$ and $\min\{\demh_f(e_4,T_3),~\demh_f(e_4,T_4)\}=0$,
    \item[(ii)] $\max\{\demh_f(e_1,T_1),~\demh_f(e_1,T_5)\}=0$ or $\max\{\demh_f(e_4,T_3),~\demh_f(e_4,T_4)\}=0$, or
    \item[(iii)] $\min\bigg\{\demh_f(e_1,T_1),~\demh_f(e_1,T_5),~\demh_f(e_4,T_3),~\demh_f(e)4,T_4)\bigg\}=0$, $\min\{\demh_f(e)2,T_1),\demh_f(e_2,T_7')\}=0$ and $\min\{\demh_f(e_3,T_3),~\demh_f(e_3,T_8')\}=0$,
\end{itemize}   
then $\textrm{Val}_f(T) \le 4$.
\end{claim} 
\begin{proofclaim}
If (i) holds, then the sum of the terms containing $e_1$ and $e_4$ is at most $2$, while as noted before the sum of the terms containing $e_2$ and $e_3$ is at most $2$ and hence $\textrm{Val}_f(T)\le 4$ as desired.

If (ii) holds, then once again the sum of the terms containing $e_1$ and $e_4$ is at most $2$, while as noted before the sum of the terms containing $e_2$ and $e_3$ is at most $2$ and hence $\textrm{Val}_f(T)\le 4$ as desired.

If (iii) holds, then the sum of the terms containing $e_1$ and $e_4$ is at most $3$, while using what was noted before the sum of the terms containing $e_2$ and $e_3$ is at most $1$ and hence $\textrm{Val}_f(T)\le 4$ as desired.
\end{proofclaim}

Let $\sigma(i)$ denote the index of the $i$th heaviest edge among $(e_j:j\in [5]))$ so that we have $\phi(e_{\sigma(1)}) > \phi(e_{\sigma(2)}) > \phi(e_{\sigma(3)}) > \phi(e_{\sigma(4)}) > \phi(e_{\sigma(5)})$ (this is well-defined as the $i$th heaviest is unique since $\phi$ is untied). We make the following observation.

\begin{claim}\label{cl:5cycleObservation}
For $i\in [5]$ if $\sigma(i) < \sigma(i+1)$(i.e.~$e_i$ is the heaviest positive edge in $T_i$) and $\sigma(i) < \max\{\sigma(i-1),\sigma(i+2)\}$, then $e_i$ sponsors $T_i$. Symmetrically for $i\in [5]$, if $\sigma(i) > \sigma(i-1)$(i.e.~$e_i$ is the heaviest positive edge in $T_{i-1}$) and $\sigma(i) > \min\{\sigma(i-2),\sigma(i+1)\}$, then $e_i$ sponsors $T_{i-1}$.  
\end{claim}
\begin{proofclaim}
In what follows we use modular arithmetic and omit the taking of moduli unless needed for clarity. By symmetry, it suffices to prove only the first statement. Note $T_{5+(i-1 \mod 5)}'$ and $T_{5+(i+2 \mod 5)}'$ are obtuse. As $\phi$ is triangle-covered, this implies that $\phi(e_{i-1}),\phi(e_{i-2}) \ge |\phi(f_{5+i})|$. Since $\sigma(i) < \max\{\sigma(i-1),\sigma(i+2)\}$, we have that $\phi(e_{i}) > \min\{\phi(e_{i-1}),~\phi(e_{i+2})\}$ and hence $\phi(e_{i}) \ge |\phi(f_{5+i})|$. As $\sigma(i) < \sigma(i+1)$, we also have that $\phi(e_i) > \phi(e_{i+1})$ and hence $e_i$ is the sponsor of $T_i$ as desired. 
\end{proofclaim}

Next we utilize the above observation to prove the following crucial claim.
\begin{claim}\label{cl:5cycleNotes}
All of the following hold:
\begin{itemize}
    \item[(a)] $e_{\sigma(1)}$ sponsors both $T_{\sigma(1)}$ and $T_{\sigma(1)-1 \mod 5}$, 
    \item[(b)] $e_{\sigma(2)}$ sponsors at least one of $T_{\sigma(1)}$ and $T_{\sigma(1)-1 \mod 5}$ and both unless $e_{\sigma(1)}\cap e_{\sigma(2)}\ne \emptyset$ [~more particularly it would obviously be the sponsor of the one $e_{\sigma(1)}$ is not~]
    \item[(c)] if $e_{\sigma(1)}\cap e_{\sigma(2)}\ne \emptyset$, then $e_{\sigma(3)}$ sponsors all incident triangles in which it is the heaviest positive edge and thus $e_{\sigma(3)}$ sponsors at least one of $T_{\sigma(1)}$ and $T_{\sigma(1)-1 \mod 5}$.
\end{itemize}
\end{claim}
\begin{proofclaim}
By Claim~\ref{cl:5cycleObservation}, it follows that for $i\in [5]$, if $\sigma(i)\in \{1,2\}$, then $e_j$ sponsors all incident triangles in which it is the heaviest positive edge (since the associated maximum in the claim is at least $3$ as $\sigma(i)\in \{1,2\}$). From this, (a) and (b) follow immediately. 

Finally we prove that (c) holds. Suppose without loss of generality that $\sigma(3)=i$. Since $e_{\sigma(1)}\cap e_{\sigma(2)}\ne \emptyset$, we find that $\max\{\sigma(i+1),\sigma(i-1)\} \ge 4 > 3 =\sigma(i)$. If $\sigma(i+1) \ge 4$ (and thus $e_i$ is the heaviest positive edge of $T_i$), then $\max\{\sigma(i-1),\sigma(i+2)\} \ge 4 > 3=\sigma(i)$ (since $\{\sigma(i+1),\sigma(i-2)\}\ne \{1,2\}$ as $e_{\sigma(1)}\cap e_{\sigma(2)} \ne \emptyset$ by assumption) and hence by Claim~\ref{cl:5cycleObservation}, $e_i$ is the sponsor of $T_i$. The same holds for $\sigma(i-1)$ and $T_{i-1}$ by symmetry and hence (c) follows as desired.
\end{proofclaim}

We note for the reader that Claim~\ref{cl:5cycleNotes}(a)-(c) imply that at least four of $(T_i: i\in [5])$ have sponsors (though we do not require this fact specifically).

\begin{claim}\label{cl:Not5}
If $\sigma(1)= 5$ or $\sigma(2)= 5$, then $\textrm{Val}_f(T) \le 4$
\end{claim}
\begin{proofclaim}
First suppose $\sigma(1)=5$ (i.e.~$e_5$ is heaviest). Then by Claim~\ref{cl:5cycleNotes}(a), $e_5$ sponsors both $T_5$ and $T_4$ and hence by definition of strict sponsorship we find that $\demh_f(e_1,T_5)=0$ and $\demh_f(e_4,T_4)=0$. Thus Claim~\ref{cl:CycleCaseContradictions}(i) holds as desired.

So we assume $\sigma(2)=5$ (i.e.~$e_5$ is second heaviest). Once more if $e_5$ sponsors both $T_5$ and $T_4$, then the result follows. Thus by Claim~\ref{cl:5cycleNotes}(b), we have that $\sigma(1)\in \{1,4\}$. Suppose without loss of generality that $\sigma(1)=1$. Hence by Claim~\ref{cl:5cycleNotes}(a), $e_1$ sponsors both $T_1$ and $T_5$ and hence $\demh_f(e_2,T_1)=0$. Similarly, we find by Claim~\ref{cl:5cycleNotes}(b) that $e_5$ sponsors $T_4$ and hence $\demh_f(e_4,T_4)=0$. Then $\textrm{Val}_f(T) \le 4$ as desired unless $\demh_f(e_3,T_3), \demh_f(e_3,T_8') > 0$. The latter implies that $e_2$ is not the sponsor of $T_3$. Since $e_{\sigma(1)}\cap e_{\sigma(2)}\ne \emptyset$, this then implies by Claim~\ref{cl:5cycleNotes}(c) that $\sigma(3)\ne 2$. Similarly we find that $e_4$ is not the sponsor of $T_3$ and hence $\sigma(3)\ne 4$. Thus $\sigma(3)=3$. But then by Claim~\ref{cl:5cycleNotes}(c), we have that $e_3$ sponsors both $T_3$ and $T_2$ and hence $\demh_f(e_4,T_3)=0$ and Claim~\ref{cl:CycleCaseContradictions}(ii) holds for $e_4$ as desired.
\end{proofclaim}

\begin{claim}\label{cl:Not14}
If $\sigma(1)\in \{1,4\}$, then $\textrm{Val}_f(T) \le 4$.    
\end{claim}
\begin{proofclaim}
Without loss of generality we assume that $\sigma(1)=1$. By Claim~\ref{cl:5cycleNotes}(a), we find that $e_1$ sponsors $T_1$ and $T_5$. Hence $\demh_f(e_2,T_1)=0$.

First suppose $\sigma(2)=2$. By Claim~\ref{cl:5cycleNotes}(b) $e_2$ sponsors $T_2$ and hence $\dem(e_3,f)=0$ and thus $\demh_f(e_3,T_3)=0$ and $\demh_f(e_3,T_8')=0$. Note here that we have $e_{\sigma(1)}\cap e_{\sigma(2)}\ne \emptyset$, that is the hypothesis of Claim~\ref{cl:5cycleNotes}(c) holds. Further suppose $\sigma(3)=3$. By Claim~\ref{cl:5cycleNotes}(c), we find that $e_3$ sponsors $T_3$ and hence $\demh_f(e_4,T_3)=0$ and hence Claim~\ref{cl:CycleCaseContradictions}(iii) holds as desired. So instead suppose $\sigma(3)=5$. By Claim~\ref{cl:5cycleNotes}(c), we find that $e_5$ sponsors $T_4$ and hence $\demh_f(e_4,T_4)=0$ and Claim~\ref{cl:CycleCaseContradictions}(iii) holds as desired. So we assume $\sigma(3)=4$. By Claim~\ref{cl:5cycleNotes}(c), we find that $e_4$ sponsors both $T_4$ and $T_3$. But then we find that $\phi(e_1),\phi(e_4) > \phi(e_5) \ge |\phi(f_6)|+|\phi(f_8)|$ (where we used that $\phi$ is triangle-covered for the triangle $T_{10}'$). Hence $\demh_f(e_1,T_1) + \demh_f(e_4,T_3) =\frac{|\phi(f_6)|}{\phi(e_1)} + \frac{|\phi(f_8)|}{\phi(e_4)} \le \frac{|\phi(f_6)|+|\phi(f_8)|}{\phi(e_5)} \le 1$ and hence $\textrm{Val}_f(T) \le 4$ as desired. Thus we assume $\sigma(2)\ne 2$.

Next suppose $\sigma(2)=3$. By Claim~\ref{cl:5cycleNotes}(b), $e_3$ sponsors $T_3$ and $T_2$. As $e_3$ sponsors $T_2$, we find that $\dem(e_2,f)=0$ and hence $\demh_f(e_2,T_1)=0$ and $\demh_f(e_2,T_7')=0$. As $e_3$ sponsors $T_3$, we find that $\demh_f(e_4,T_3)=0$. We also have that $\phi(e_1),\phi(e_3) > \phi(e_5) \ge |\phi(f_6)|+|\phi(f_8)|$ (where we used that $\phi$ is triangle-covered for the triangle $T_{10}'$). Hence $\demh_f(e_1,T_1) + \demh_f(e_3,T_3) =\frac{|\phi(f_6)|}{\phi(e_1)} + \frac{|\phi(f_8)|}{\phi(e_3)} \le \frac{|\phi(f_1)|+|\phi(f_3)|}{\phi(e_5)} \le 1$, and hence $\textrm{Val}_f(T) \le 4$ as desired. Thus we assume $\sigma(2)\ne 3$.

So by Claim~\ref{cl:Not5}, we assume $\sigma(2)=4$. So by Claim~\ref{cl:5cycleNotes}(b), $e_4$ sponsors $T_4$ and $T_3$. Hence $\demh_f(e_3,T_3)=0$. But then we use that $\demh_f(e_2,T_7')+ \demh_f(e_3,T_8') \le \dem(e_2,f)+\dem(e_3,f) \le 1$. Meanwhile, we also have that $\phi(e_1),\phi(e_4) > \phi(e_5) \ge |\phi(f_6)|+|\phi(f_8)|$ (where we used that $\phi$ is triangle-covered for the triangle $T_{10}'$). Hence $\demh_f(e_1,T_1) + \demh_f(e_4,T_3) =\frac{|\phi(f_6)|}{\phi(e_1)} + \frac{|\phi(f_8)|}{\phi(e_4)} \le \frac{|\phi(f_6)|+|\phi(f_8)|}{\phi(e_5)} \le 1$ and hence $\textrm{Val}_f(T) \le 4$ as desired.
\end{proofclaim}

\begin{claim}\label{cl:Not23}
If $\sigma(1)\in \{2,3\}$, then $\textrm{Val}_f(T) \le 4$.
\end{claim}
\begin{proofclaim}
Without loss of generality, we assume that $\sigma(1)=2$. By Claim~\ref{cl:5cycleNotes}(a), we find that $e_2$ sponsors $T_2$ and $T_1$. As $e_2$ sponsors $T_2$, we find that $\dem(e_3,f)=0$; hence $\demh_f(e_3,T_3)=0$ and $\demh_f(e_3,T_8')=0$. As $e_2$ sponsors $T_1$, we find that $\demh_f(e_1,T_1)=0$.

First suppose $\sigma(2)=3$. By Claim~\ref{cl:5cycleNotes}(b), we find that $e_3$ sponsors $T_3$ and hence $\demh_f(e_4,T_3)=0$ and hence Claim~\ref{cl:CycleCaseContradictions}(i) holds as desired. So we assume $\sigma(2)\ne 3$.

Next suppose $\sigma(2)=1$. Note here that $e_{\sigma(1)}\cap e_{\sigma(2)}\ne \emptyset$ and hence the hypothesis of Claim~\ref{cl:5cycleNotes}(c) holds. Further suppose $\sigma(3)=5$. By Claim~\ref{cl:5cycleNotes}(c), we find that $e_5$ sponsors $T_4$ and hence $\demh_f(e_4,T_4)=0$; thus Claim~\ref{cl:CycleCaseContradictions}(i) holds as desired. So instead suppose $\sigma(3)=3$. By Claim~\ref{cl:5cycleNotes}(c), we find that $e_3$ sponsors $T_3$ and hence $\demh_f(e_4,T_3)=0$; thus Claim~\ref{cl:CycleCaseContradictions}(i) holds as desired. So we assume $\sigma(3)=4$. But then we have that $\phi(e_2),\phi(e_4) > \phi(e_5) \ge |\phi(f_6)|+|\phi(f_8)|$ (where we used that $\phi$ is triangle-covered for the triangle $T_{10}'$). Hence $\demh_f(e_2,T_1) + \demh_f(e_4,T_3) =\frac{|\phi(f_6)|}{\phi(e_2)} + \frac{|\phi(f_8)|}{\phi(e_4)} \le \frac{|\phi(f_6)|+|\phi(f_8)|}{\phi(e_5)} \le 1$ and hence $\textrm{Val}_f(T) \le 4$ as desired.

So by Claim~\ref{cl:Not5}, we assume $\sigma(2)=4$. By Claim~\ref{cl:5cycleNotes}(b) (as $e_2\cap e_4=\emptyset$), we find that $e_4$ sponsors $T_4$ and $T_3$. But then we also have that $\phi(e_2),\phi(e_4) > \phi(e_5) \ge |\phi(f_6)|+|\phi(f_8)|$ (where we used that $\phi$ is triangle-covered for the triangle $T_{10}'$). Hence $\demh_f(e_2,T_1) + \demh_f(e_4,T_3) =\frac{|\phi(f_6)|}{\phi(e_2)} + \frac{|\phi(f_8)|}{\phi(e_4)} \le \frac{|\phi(f_6)|+|\phi(f_8)|}{\phi(e_5)} \le 1$ and hence $\textrm{Val}_f(T) \le 4$ as desired.
\end{proofclaim}

Finally, we conclude the proof of this subcase via the claims above. By definition $\sigma(1)\in [5]$. Thus by Claim~\ref{cl:Not5} (if $\sigma(1)= 5$), by Claim~\ref{cl:Not14} (if $\sigma(1)\in\{1,4\}$) and by Claim~\ref{cl:Not23} (if $\sigma(1)\in \{2,3\}$), we find that $\textrm{Val}_f(T)\le 4$ as desired. 

\subsubsection{Proof of Lemma\ref{lem:EdgeWeightTriangleThinCases}(1b)}

This subcase assumes the remaining negative edges are $f_5$ and $f_6$. Thus we have positive edges $e_1,e_2,e_3,e_4,e_8$ and negative edges $f_5,f_6,f_7,f_9,f_{10}$. Similarly, here is a breakdown of the acute/obtuse triangles in this case:

$$\textrm{All Positive}:~~T_3=e_3e_4e_8,$$
$$\textrm{Acute Triangles:}~~T_1=e_1e_2f_6,~~T_2=e_2e_3f_7,~~T_7=e_2e_8f_{10},$$
$$\textrm{Obtuse Triangles:}~~T_4'=e_4f_5f_9,~~ T_5'=f_5e_1f_{10},~~T_6'=e_1f_7f_9,~~T_8'=e_3f_9f_6,~~T_9'=e_4f_{10}f_7,~~T_{10}'=f_5f_6e_8$$ 

For this subcase, we calculate that

\begin{align*}
\textrm{Val}_f(xy) &= \demh_f(e_1,T_5') + \demh_f(e_4,T_4'),\\
\textrm{Val}_f(xz) &= \demh(e_1,T_1)+\demh(e_2,T_1)+\demh_f(e_3,T_8'),\\
\textrm{Val}_f(yz) &= \demh_f(e_2,T_7).  
\end{align*}

Suppose for a contradiction that $\textrm{Val}_f(T) = \textrm{Val}_f(xy)+\textrm{Val}_f(xz)+\textrm{Val}_f(yz) > 4$. Since $\dem(e_2,f)+\dem(e_3,f)\le 1$, we have that $\demh_f(e_2,T_i)+\dem(e_3,T_8')\le 1$ for any $i\in \{1,7\}$. Since the sum of values above is strictly greater than $4$, it follows that all of $\demh_f(e_1,T_5')$,~$\demh_f(e_4,T_4')$,~$\demh_f(e_1,T_1)$,~\\$\demh_f(e_2,T_1)$,~$\demh_f(e_2,T_7)$ are nonzero.

Since $\demh_f(e_1,T_1),\demh_f(e_2,T_1) > 0$, it follows that $T_1$ is unsponsored and hence $\phi(e_1),\phi(e_2) < |\phi(f_6)|$. Since $\phi$ is triangle-covered (in particular for the triangle $T_{10}'$), we find that $|\phi(e_8) \ge |\phi(f_5)|+|\phi(f_6)|$ and hence $\phi(e_8) > \phi(e_1)$ and $\phi(e_8) > \phi(e_2)$. Since $\phi$ is triangle-covered (in particular for the triangle $T_{5}'$), we find that $\phi(e_1)\ge |\phi(f_5)|+|\phi(f_{10})|$. Hence $\phi(e_8) \ge |\phi(f_{10})|$. It follows that $e_8$ is the sponsor of $T_7=e_2e_8f_{10}$. But then $\demh_f(e_2,T_7)=0$, a contradiction.

\subsubsection{Proof of Lemma\ref{lem:EdgeWeightTriangleThinCases}(1c)}

This subcase assumes the only remaining negative edge is $f_5$. Thus we have have positive edges $e_1,e_2,e_3,e_4,e_6,e_8$ and negative edges $f_5,f_7,f_9,f_{10}$. Similarly, here is a breakdown of the acute/obtuse triangles in this case:

$$\textrm{All Positive}:~~T_1=e_1e_2e_6,~~T_3=e_3e_4e_8,~~T_{10}=e_5e_6e_8$$
$$\textrm{Acute Triangles:}~~T_2=e_2e_3f_7,~~T_4=e_4e_5f_9,~~T_5=e_5e_1f_{10},~~T_7=e_2e_8f_{10},~~T_8=e_3f_9e_6,$$
$$\textrm{Obtuse Triangles:}~~T_6'=e_1f_7f_9,~~T_9'=e_4f_{10}f_7.$$ 

For this subcase, we calculate that

\begin{align*}
\textrm{Val}_f(xy) &= \demh_f(e_1,T_5) + \demh_f(e_4,T_4),\\
\textrm{Val}_f(xz) &= \demh_f(e_3,T_8),\\
\textrm{Val}_f(yz) &= \demh_f(e_2,T_7').  
\end{align*}

Since $\dem(e_2,f)+\dem(e_3,f)\le 1$, we actually find that the sum of the above is at most $3$. 

\subsubsection{Proof of Lemma\ref{lem:EdgeWeightTriangleThinCases}(1d)}

This subcase assumes the only remaining negative edge is $f_6$. In this subcase we thus have positive edges $e_1,e_2,e_3,e_4,e_5,e_8$ and negative edges $f_6,f_7,f_9,f_{10}$. Similarly, here is a breakdown of the acute/obtuse triangles in this case:

$$\textrm{All Positive:}~~T_3=e_3e_4e_8,$$
$$\textrm{Acute Triangles:}~~T_1=e_1e_2f_6,~~T_2=e_2e_3f_7,~~T_4=e_4e_5f_9,~~T_5=e_5e_1f_{10}~~T_7=e_2e_8f_{10},~~T_{10}=e_5f_6e_8$$
$$\textrm{Obtuse Triangles:}~~T_6'=e_1f_7f_9,~~T_8'=e_3f_9f_6,~~T_9'=e_4f_{10}f_7,$$ 

For this subcase, we calculate that

\begin{align*}
\textrm{Val}_f(xy) &= \demh_f(e_1,T_5) + \demh_f(e_4,T_4),\\
\textrm{Val}_f(xz) &= \demh_f(e_1,T_1) + \demh_f(e_2,T_1) + \demh_f(e_3,T_8'),\\
\textrm{Val}_f(yz) &= \demh_f(e_2,T_7).  
\end{align*}

Suppose for a contradiction that $\textrm{Val}_f(T) = \textrm{Val}_f(xy)+\textrm{Val}_f(xz)+\textrm{Val}_f(yz) > 4$. It then follows that all of $\demh_f(e_4,T_4)$, $\demh_f(e_2,T_1)$, $\demh_f(e_2,T_7)$, $\demh_f(e_1,T_5)$, $\demh_f(e_1,T_1)$ are nonzero (since $\demh_f(e_2,T_i)+\demh_f(e_3,T_8') \le 1$ for each $i\in \{1,7\}$). This implies that neither $e_1$ nor $e_2$ is the sponsor of $T_1$ and hence $|\phi(f_6)| > \phi(e_1),\phi(e_2)$. But then $\phi(e_3)\ge |\phi(f_6)|+\phi(f_9)|$. So $\phi(e_3)\ge \phi(e_1),\phi(e_2)$ but then $\phi(e_3)\ge \phi(e_1) \ge 1+|\phi(f_9)|\ge \phi(f)$. So $e_3$ sponsors $T_2$. Hence $\dem(e_2,f)=0$ and therefore $\demh_f(e_2,T_1)=0$, a contradiction.

\subsubsection{Proof of Lemma\ref{lem:EdgeWeightTriangleThinCases}(1e)}

This subcase assumes there are no remaining negative edges (that is $e_5,e_6,e_8$ are all positive). Thus we have positive edges $e_1,e_2,e_3,e_4,e_5,e_6,e_8$ and negative edges $f_7,f_9,f_{10}$. Similarly, here is a breakdown of the acute/obtuse triangles in this case:

$$\textrm{All Positive:}~~T_1=e_1e_2e_6,~~T_3=e_3e_4e_8,~~T_{10}=e_5e_6e_8$$
$$\textrm{Acute Triangles:}~~T_2=e_2e_3f_7,~~T_4=e_4e_5f_9,~~T_5=e_5e_1f_{10}~~T_7=e_2e_8f_{10},~~T_8=e_3f_9e_6,$$
$$\textrm{Obtuse Triangles:}~~T_6'=e_1f_7f_9,~~T_9'=e_4f_{10}f_7,$$ 

For this subcase, we calculate that

\begin{align*}
\textrm{Val}_f(xy) &= \demh_f(e_1,T_5) + \demh_f(e_4,T_4),\\
\textrm{Val}_f(xz) &= \demh_f(e_3,T_8),\\
\textrm{Val}_f(yz) &= \demh_f(e_2,T_7).  
\end{align*}

Since $\dem(e_2,f)+\dem(e_3,f)\le 1$, we actually find that the sum of the above is at most $3$. 

\subsection{Proof of Lemma\ref{lem:EdgeWeightTriangleThinCases}(2)}

Recall that in this case the positive edges include $e_1$, $e_2$, $e_3$, $e_4, e_{10}$ and the negative edges include $f_7,f_{9}$ (with the signs of $g_5,g_6,g_8$ varying among the subcases). 

\subsubsection{Proof of Lemma\ref{lem:EdgeWeightTriangleThinCases}(2a)}

This subcase assumes the remaining negative edges are $f_5$ and $f_8$. Thus we thus have positive edges $e_1,e_2,e_3,e_4,e_6,e_{10}$ and negative edges $f_5,f_7,f_8,f_9$. Similarly, here is a breakdown of the acute/obtuse triangles in this case:

$$\textrm{All Positive:}~~T_1=e_1e_2e_6,$$
$$\textrm{Acute Triangles:}~~T_2=e_2e_3f_7,~~T_3=e_3e_4f_8,~~T_5=f_5e_1e_{10},~~T_7=e_2f_8e_{10},~~T_8=e_3f_9e_6,~~T_9=e_4e_{10}f_7,$$
$$\textrm{Obtuse Triangles:}~~T_4'=e_4f_5f_9,~~T_6'=e_1f_7f_9,~~T_{10}'=f_5e_6f_8$$ 

For this subcase, we calculate that

\begin{align*}
\textrm{Val}_f(xy) &= \demh_f(e_1,T_5) + \demh(e_{10},T_5)+\demh_f(e_4,T_4),\\
\textrm{Val}_f(xz) &= \demh_f(e_3,T_8),\\
\textrm{Val}_f(yz) &= \demh_f(e_4,T_3) + \demh_f(e_3,T_3) + \demh_f(e_2,T_7)+\demh(e_{10},T_7).  
\end{align*}

Suppose for a contradiction that $\textrm{Val}_f(T) = \textrm{Val}_f(xy)+\textrm{Val}_f(xz)+\textrm{Val}_f(yz) > 4$. Recall that $\dem(e_2,f)+\dem(e_3,f)\le 1$ and hence $\demh_f(e_2,T_7)+\demh_f(e_3,T_i)\le 1$ for any $i\in \{3,8\}$. Similarly $\dem(e_4,f)+\dem(e_{10},f)\le 1$ and hence $\demh_f(e_4,T_i)+\demh_f(e_{10},T_j)\le 1$ for any $i\in \{3,4\}$, $j\in\{5,7\}$. (This already implies the sum of values above is at most $5$). Since the sum of values above is strictly greater than $4$, this implies that all of the following hold:
\begin{itemize}
    \item[(i)] there exists $i\in \{4,10\}$ such that both of the terms containing $e_i$ are nonzero,
    \item[(ii)] $\demh_f(e_3,T_3),~\demh_f(e_3,T_8) > 0$, and
    \item[(iii)] $\demh_f(e_1,T_5) > 0$.
\end{itemize}

Since $\demh_f(e_3,T_8) > 0$, it follows that $\dem(e_3,f_9)>0$ and hence $e_6$ is not the sponsor of $T_8$. Since $e_6$ is not the sponsor of $T_8$, it follows by definition of sponsor (and $\phi$ being untied) that $\phi(e_6) < \max\{\phi(e_3),|\phi(f_9)|\}$. Since $\phi$ is triangle-covered (in particular for the triangle $T_4'$), we find that $\phi(e_4)\ge |\phi(f_5)|+|\phi(f_9)|$ and hence $\phi(e_4) \ge \phi(e_6).$ Since $\phi$ is triangle-covered (in particular for the triangle $T_{10}'$), we find that $\phi(e_6)\ge |\phi(f_5)|+|\phi(f_8)|$ and hence $\phi(e_4) \ge \phi(e_6) \ge |\phi(f_8)|.$ It follow that $T_3=e_4e_4f_8$ has a sponsor. Since $\demh_f(e_3,T_3) > 0$ by (ii), we find that $e_3$ is the sponsor of $T_3$. This implies that $\phi(e_3) > \phi(e_4)$ by definition of sponsor. 

As $e_3$ is the sponsor of $T_3$, we also find that $\demh_f(e_4,T_3)=0$ and hence (i) does not hold for $i=4$. Thus (i) holds for $i=10$ and we have that $\demh_f(e_{10},T_5),\demh_f(e_{10},T_7) > 0$. Since $\demh_f(e_1,T_5)>0$ by (iii) and $\demh_f(e_{10},T_5)>0$, we find that $T_5=f_5e_1e_{10}$ is unsponsored. Hence by definition of sponsor, we find that $|\phi(f_5)| > \phi(e_1),\phi(e_{10})$. Thus $\phi(e_4) \ge |\phi(f_5)| > \phi(e_{10})$. Since $\phi$ is triangle-covered (in particular for the triangle $T_6'$), we find that $\phi(e_1)\ge |\phi(f_7)|+|\phi(f_9)|$. Hence $\phi(e_4)\ge |\phi(f_7)|$. Thus $\phi(e_4)\ge \phi(e_1)\ge |\phi(f_7)|$. It now follows that $e_4$ is the sponsor of $T_9=e_4e_{10}f_7$. Since $f=f_7$, this implies that $\dem(e_{10},f)=0$ and hence $\demh_f(e_{10},T_5)=0$, a contradiction.

\subsubsection{Proof of Lemma\ref{lem:EdgeWeightTriangleThinCases}(2b)}

This subcase assumes the only remaining negative edge is $f_8$. In this subcase we thus have positive edges $e_1,e_2,e_3,e_4,e_5,e_6,e_{10}$ and negative edges $f_7,f_8,f_9$. Similarly, here is a breakdown of the acute/obtuse triangles in this case:

$$\textrm{All Positive:}~~T_1=e_1e_2e_6,~~T_5=e_5e_1e_{10},$$
$$\textrm{Acute Triangles:}~~T_2=e_2e_3f_7,~~T_3=e_3e_4f_8,~~T_4=e_4e_5f_9,~~T_7=e_2f_8e_{10},~~T_8=e_3f_9e_6,~~T_9=e_4e_{10}f_7,~~T_{10}=e_5e_6f_8$$
$$\textrm{Obtuse Triangles:}~~T_6'=e_1f_7f_9$$ 

For this subcase, we calculate that

\begin{align*}
\textrm{Val}_f(xy) &= \demh_f(e_4,T_4),\\
\textrm{Val}_f(xz) &= \demh_f(e_3,T_8),\\
\textrm{Val}_f(yz) &= \demh_f(e_3,T_3) + \demh_f(e_4,T_3) + \demh_f(e_2,T_7)+\demh(e_{10},T_7).  
\end{align*}

Since $\dem(e_2,f)+\dem(e_3,f)\le 1$, we find that $\demh_f(e_3,T_8) + \demh_f(e_2,T_7)\le \dem(e_3,f)+\dem(e_2,f)\le 1$. Similarly since $\dem(e_4,f)+\dem(e_{10},f)\le 1$, we find that $\demh_f(e_4,T_3) + \demh(e_{10},T_7)\le \dem(e_4,f)+\dem(e_{10},f)\le 1$. Thus the total sum of values above is at most $4$ as desired. This concludes the proof of Lemma~\ref{lem:EdgeWeightTriangleThinCases}(2b) and hence of Lemma~\ref{lem:EdgeWeightTriangleThinCases} and hence of the Fractional Nash-Williams' Conjecture.

\part{Fractional Stability}

\section{Proof Overview}

Our proof of Theorem~\ref{thm:fracNWstability} consists of three main parts (in addition to reusing many of the lemmas from the proof of the Fractional Nash-Williams' Conjecture). First we prove stability versions of Lemmas~\ref{lem:Value4Parts} and~\ref{lem:Value}. Second we use these value stability lemmas and the bounds from the discharging lemmas in Section~\ref{s:FracNW} to prove lemmas about the local structure of edges whose charge after discharging is negative or not very positive; here local means of the common neighborhood of their endpoints and structure refers to both graph structure but also the edge-weights in a supposed negative value triangle-covered edge-weighting of the graph $G$. Third we introduce two additional discharging rules which will be triggered in sequence after the discharging of Section~\ref{s:FracNW} (so this is an example of layered discharging); the first additional rule sends charge from edges with very positive $\textrm{ch}_F$, half of the charge is spread uniformly to incident negative edges and the other half of the charge is spread uniformly to all edges; the second additional rule applies the same discharging scheme to those edges that are very positive after the first additional rule fires. We show using the local structure lemmas that all edges have nonnegative charge after these new additional discharging rules provided the graph is not close to extremal. 

The trickiness in needing two additional rules arises as we need not just the local structure of one common neighborhood (which would be roughly $n/2$ vertices, that is half the graph) but rather the extremal structure throughout the entire graph. This becomes quite complex then since while we can with a great deal of effort characterize the common neighborhood of negative edges that are extremal (i.e.~do not become very positive after our initial discharging), we need to show such neighborhoods not only admit a nice graph structure but also have negative weights roughly similar (or larger) to that of an extremal edge so as to propagate structure throughout the graph.

\subsection{Value Stability}

The first step in our proof of Theorem~\ref{thm:fracNWstability} is to prove a stability version of Lemma~\ref{lem:Value4Parts}, namely Lemma~\ref{lem:ValueStability2}. We delay the statement until the next section as it turns out there are four extremal families of $(x_i, a_{ij})$ which as we will see can be viewed as step graphons, see Definition~\ref{def:4extremal} for the definition of the four families. While the two base families are essentially complete balanced bipartite cuts, the latter two arise when the $x_i$ admit symmetry (either reflective or fully symmetric) since then it is possible to average two or even three bipartite cuts to create new extremal families. The stability lemma then shows that if the value is close to the extremal value of $\frac{1}{4}$, then the step graphon is close to a member of one of these four extremal families.

The next step then is to prove a stability version of Lemma~\ref{lem:Value}, namely Lemma~\ref{lem:ValueStability}. This stability lemma essentially says that if a $K_5$-free $[0,1]$-edge-weighted triangle-thin graph on $n$ vertices has value close to $\frac{n^2}{4}$, then it is close to being a blowup of a member of one of the four extremal families from Definition~\ref{def:4extremal}. In general, it would be possible to show for this type of problem (with the assumption of being $K_r$-free) that a nearly extremal example has a partition into $r-1$ `quasi-independent' sets such that almost all of the transversal $K_{r-1}$'s are nearly extremal for the step graphon version of the problem. However, whether that forces the edge weights of the $K_{r-1}$'s to have similar structure depends on the specific problem. For our problem, we show that it does and hence the nearly extremal graphs are close to being blowups but it is conceivable that in other similar problems different types of extremal $K_{r-1}$'s can simultaneously coexist and the nearly extremal examples would not necessarily need to be a blowup. Given this, we make no attempt to build a general theory of the maximum total edge-weight of $K_r$-free graphs satisfying certain clique weight upper bounds but rather instead proceed to prove just the stability version of our own specific problem ($K_5$-free triangle-thin). 

\subsection{Local Structure}

Next we proceed to characterize the local structure of extremal edges (those negative edges $f$ with $\textrm{ch}_F(f) \le \sqrt{\varepsilon}\cdot |\phi(f)|$). First we prove some lemmas about $\textrm{Val}_f(xy)$; the first such lemma, Lemma~\ref{lem:ObtuseAcutePair} characterizes the behavior of this function when one of $x$ and $y$ forms an obtuse triangle with $f$ and the other forms an acute triangle with $f$; the second value lemma, Lemma~\ref{lem:ExtremalPairValue}, characterizes the nature of pairs $xy$ where $\textrm{Val}_f(xy)$ is close to its maximum of $2\cdot |\phi(f)|$. We then use these value lemmas to prove our main structural lemma about the common neighborhood of extremal edges, Lemma~\ref{lem:NeighborhoodStructure}, which essentially says that the common neighborhood is small (close to $n/2$) and additionally falls into one of two categories: either $f=uv$ satisfies that both $u$ and $v$ have $\alpha\cdot n$ incident negative edges $f'$ with $|\phi(f')|\ge 2\cdot |\phi(f)|$ (what we call $\alpha$-locally-small) or $N(u)\cap N(v)$ admits a cut with at least $(1-\alpha)\cdot \frac{|N(u)\cap N(v)|^2}{2}$ negatives $f'$ with $|\phi(f')|\ge (1-\alpha)\cdot |\phi(f)|$ (what we call $\alpha$-bipartite).

\subsection{Forcing the Global Structure: Additional Discharging Rules}

Finally we proceed with forcing the global extremal structure. First to that end, we prove a ``cleaning lemma", Lemma~\ref{lem:FracCleaning}, which is useful in forcing structure out from a common neighborhood of an extremal edge $f$ for vertices that are in an obtuse triangle with $f$. 

Finally, we turn towards proving our main fractional stability result. To that end, we need the promised additional layer discharging rules. 
\vskip.1in

\noindent \textbf{Extra Discharging Rules:}
\vskip.1in
\begin{itemize}
    \item[\textbf{Rule 1:}] Every negative edge $f=uv$ with $\textrm{ch}_F(f)\ge \sqrt{\varepsilon}\cdot |\phi(f)|$ sends $\frac{1}{2}\cdot \textrm{ch}_F(f)\cdot \frac{1}{|N^-(u)\cup N^-(v)|}$ to each edge in $N^-(u)\cup N^-(v)$ and $\frac{1}{2}\cdot \textrm{ch}_F(f)\cdot \frac{1}{e(G)}$ to each edge of $G$. We denote the charge after this rule by $\textrm{ch}_1$. 
    \item[\textbf{Rule 2:}] Every negative edge $f=uv$ with $\textrm{ch}_1(f)\ge \sqrt{\varepsilon}\cdot |\phi(f)|$ sends $\frac{1}{2}\cdot \textrm{ch}_F(f)\cdot \frac{1}{|N^-(u)\cup N^-(v)|}$ to each edge in $N^-(u)\cup N^-(v)$ and $\frac{1}{2}\cdot \textrm{ch}_F(f)\cdot \frac{1}{e(G)}$ to each edge of $G$. We denote the charge after this rule by $\textrm{ch}_2$.
\end{itemize}

Recall our main fractional stability result is that if $G$ with $\delta(G)\ge (\frac{3}{4}-\varepsilon)n$ has no fractional $K_3$-decomposition, then we want to show $(1-\sigma)n$ of the vertices have degree at most $(\frac{3}{4}+\sigma)n$ and the max cut is at least $(1-\sigma)\cdot \frac{n^2}{4}$. Yet if $G$ has no fractional $K_3$-decomposition then by Lemma~\ref{lem:FarkasK3Restated}, $G$ has a triangle-covered edge-weighting $\phi$ with $|\phi| < 0$ and so applying our original and new discharging rules there exists $f^*$ with $\textrm{ch}_2(f^*) < 0$. We show that this leads to the extremal structure.

Roughly the proof idea is as follows. We consider a heaviest edge $f_0=u_0v_0$ that is still extremal after Rule 1 (i.e. $\textrm{ch}_F(f_0),\textrm{ch}_1(f_0) <\sqrt{\varepsilon}\cdot |\phi(f_0)|$). We show all negative edges close to $f_0$ in weight that do not discharge in Rules 1 or 2 (including $f_0$ itself) are actually $\alpha$-bipartite (for some small $\alpha$ depending on $\varepsilon$). 

Thus we consider $N(u_0)\cap N(v_0)$ and case accordingly. If the two sides $A_0,B_0$ of the promised bipartite cut (from $f_0$ being $\alpha$-bipartite) each have very few edges, then that means $G$ already has desired max cut size (namely between $A:= N(u_0)\cap N(v_0)$ and $B:=V(G)\setminus A$). Thus it remains to show most vertices have degree close to $3n/4$. This is already true for most vertices in $N(u_0)\cap N(v_0)$ (since most are incident with many negative edges with weight close to $f_0$, and most of those must not discharge in Rule 1 or 2 as otherwise the whole graph discharges, and so most are $\alpha$-bipartite). So we carefully choose a `typical' negative edge from the bipartite cut, call if $f_0'$, argue its common neighborhood is essentially $B$ and then apply the same argument we made before to show most vertices of $B$ have degree roughly $3n/4$. 

So we instead assume one of the sides, say $A_0$, has a decent number of edges (say $\ge \beta\cdot n^2$ for some well-chosen $\beta$). We then instead abandon $f_0$ and consider a `typical' edge $f_1=u_1v_1$ of the bipartite cut with the additional property that is incident with $\beta\cdot n$ edges inside $A_0$. Most of these neighbors behave `normally' (both with respect to the cut, but also in not being incident with many heavy edges that discharge in Rule 1 or 2) and thus we find many vertices $w$ that form an obtuse triangle with $f_1$ (say with $u_1w$ is positive and thus has weight roughly twice that of $f_0$). 

We then use one of these obtuse vertices $u_2$ (say it is in $A_1$ of the new cut $(A_1,B_1)$ and our cleaning lemma, Lemma~\ref{lem:FracCleaning}, to argue that most of $B:= (N(u_1)\setminus N(v_1)) \cup B_1$ lies inside $N^-(u_2)$ and with edges roughly weight of $f_0$. Since these edges are also $\alpha$-bipartite, this implies that such vertices have degree at most roughly $n/4$ inside $B$. Hence they have degree roughly $3n/4$ and setting $A:=V(G)\setminus B$, we find our desired max cut. Lastly it remains to show most vertices of $A$ have degree at most $3n/4$ which we again achieve by considering $f_2=u_2v_1$ (and finding a well-chosen $v_1$-obtuse vertex $v_2$ to apply the cleaning lemma to). 

\section{Value Stability}

\begin{definition}
A \emph{normalized vertex-weighting} $\mu$ of a graph $G$ is a function $\mu:V(G)\rightarrow [0,1]$ such that $\sum_{v\in V(G)}\mu(v)=1$. If $\varphi$ is an edge-weighting of $G$, then the \emph{value} of $\varphi$ under $\mu$,denoted $|\varphi| _{\mu}$, is $\sum_{e=uv\in E(G)}\varphi(e)\cdot \mu(u)\cdot\mu (v)$.      
\end{definition}

In these terms, Lemma~\ref{lem:Value4Parts} says that if $G\cong K_4$, $\mu$ is a normalized vertex weighting of $G$ and $\varphi$ is a triangle-thin $[0,1]$-edge-weighting of $G$, then the value of $\varphi$ under $\mu$ is at most $1/4$. We seek to characterize the pairs of $\mu,\varphi$ that are (close to) extremal for this setting. To that end, let us recall some terminology from the theory of graphons.

\begin{definition}
A \emph{graphon} is a symmetric measurable function $W: [0,1]^2\rightarrow [0,1]$ and a \emph{step graphon} is a graphon $W$ that admits a partition of $[0,1]$ into a finite number of intervals $A_1,\ldots, A_k$ such that $W$ is constant on $A_i\times A_j$ for all $i,j\in [k]$. We regard a step graphon equivalently as a triple $(G,\mu,\varphi)$ where $G$ is a (finite but possibly with loops) graph, $\mu$ is a normalized vertex-weighting of $G$, and $\varphi$ is a $[0,1]$-edge-weighting of $G$. The equivalence is we let $V(G):=\{v_1,\ldots, v_k\}$, $\mu(v_i):=|A_i|$, $\varphi(v_iv_j):= W(a_i,a_j)$ where $a_i\in A_i, a_j\in A_j$. We let $|W| := |\varphi|_{\mu}$ denote the \emph{value} of $W$. 
\end{definition}

\begin{definition}
We say two step graphons $(G,\mu,\varphi)$ and $(G',\mu',\varphi')$ are \emph{isomorphic} if there exists an isomorphism $\sigma$ from $G$ to $G'$ such that $\mu(v)=\mu'(\sigma(v))$ for all $v\in V(G)$ and $\varphi(e)=\varphi'(\sigma(e))$ for all $e\in E(G)$. Furthermore, we say two step graphons $(G,\mu,\varphi)$ and $(G',\mu',\varphi')$ are \emph{$\delta$-close} if there exists an isomorphism $\sigma$ from $G$ to $G'$ such that $|\mu(v)-\mu'(\sigma(v))|\le \delta$ for all $v\in V(G)$ and $|\varphi(e)-\varphi'(\sigma(e))|\le \delta$ for all $e\in E(G)$.    
\end{definition}

We note that the notion of $\delta$-close defined above is similar in spirit to what would be the notion of edit distance between step graphons (or weighted edit distance for weighted graphs) except those would more naturally require that for $e=uv$, we have $|\varphi(e)\cdot \mu(u)\cdot \mu(v) - \varphi'(\sigma(e))\cdot \mu'(\sigma(v))\cdot \mu'(\sigma(u))|\le \delta$ say (or the same summed over all edges but this is roughly equivalent if we view $v(G)$ as constant). Note that our notion is a stronger assumption. For our purposes, we do not need the more precise but weaker notion as we will show that step graphons that are close to the extremal value for our problem are close in this stronger notion to one of the extremal families.

Returning to our problem, we need some additional definitions. Recall that an edge-weighting $\varphi$ of graph $G$ is \emph{triangle-thin} if for every triangle $T$ of $G$, we have that $\sum_{e\in E(T)} \varphi(e) \le 2$.

\begin{definition}
We say a step graphon $(G,\mu,\varphi)$ is \emph{triangle-thin} if $\varphi$ is a triangle-thin edge-weighting of $G$. We say a step graphon $(K_4,\mu,\varphi)$ is \emph{matched} if $\varphi(e)=\varphi(e')$ for every all $e,e' \in E(K_4)$ with $e\cap e'=\emptyset$ (i.e. for every matching $e,e'$ of $K_4$).     
\end{definition}

Now we turn to defining the step graphons that are extremal for our problem as follows.

\begin{definition}\label{def:4extremal}
Let $W=(K_4,\mu,\varphi)$ be a step graphon and let $x_1,x_2,x_3,x_4$ be an enumeration of $V(K_4)$ such that $\mu(x_1)\ge \mu(x_2) \ge \mu(x_3)\ge \mu(x_4)$. We say $W$ is 
\begin{itemize}
    \item \emph{star-cut} if $\mu(x_1)=1/2$, and for all distinct $i,j\in [4]$ with $\mu(x_i),\mu(x_j) > 0$ we have that $\varphi(x_ix_j)=1$ if $i=1,~j\in\{2,3,4\}$ and $0$ otherwise,
    \item \emph{cycle-cut} if $\mu(x_1)+\mu(x_4)=1/2$, $W$ is matched and $\varphi(x_1x_2)=\varphi(x_1x_3)=1$ and $\varphi(x_1x_4)=0$,
    \item \emph{reflective-cut} if $\mu(x_1)+\mu(x_4)=1/2$, $\mu(x_1)=\mu(x_2)$, $W$ is matched, $\varphi(x_1x_2)=1$ and $\varphi(x_1x_3)+\varphi(x_1x_4)=1$,
    \item \emph{quad-cut} if $\mu(x_i)=1/4$ for all $i\in [4]$, $W$ is matched and $\varphi(x_1x_2)+\varphi(x_1x_3)+\varphi(x_2x_3)=2$.
\end{itemize}
We let $\mc{S}$ denote the set of star-cut step graphons, $\mc{C}$ the set of cycle-cut, $\mc{R}$ the set of reflective-cut, $\mc{Q}$ the set of quad-cut. We say $W$ is \emph{value-extremal} if $W\in \mc{S}\cup \mc{C}\cup \mc{R}\cup \mc{Q}$. 
\end{definition}

\begin{remark}
It turns out the above are precisely the extremal graphons for our problem (those exactly attaining a value of $1/4$) but since we do not require this fact nor rely on it for our stability proof, we do not prove this. The reader should also note that in the latter three families, the conditions for $\mu$ become increasingly more restrictive while those for $\varphi$ become increasingly more lax. In fact, the latter two actually arise from the cycle-cut family via the symmetry of $\mu$; namely if two or even three cycle-cut graphons can be formed on the same part sizes, then their average is also a triangle-thin step graphon with value $1/4$.  
\end{remark}

\subsection{Value Stability for Four Symmetrized Parts}

We are now prepared to state our stability version of Lemma~\ref{lem:Value4Parts} as follows.

\begin{lem}[Value Stability for 4 Symmetrized Parts]\label{lem:ValueStability2}
There exists $\varepsilon_0 > 0$ such that for all $\varepsilon \in (0,\varepsilon_0]$ the following holds: If $W=(K_4,\mu,\varphi)$ is a triangle-thin step graphon with value at least $\frac{1}{4} - \varepsilon$, then $W$ is $\varepsilon^{0.01}$-close to a value-extremal step graphon.
\end{lem}
\begin{proof}
We choose $\varepsilon_0>0$ small enough as needed so that various inequalities throughout the proof regarding $\varepsilon$ hold true. Let $x_1,x_2,x_3,x_4$ be an enumeration of $V(K_4)$ such that $\mu(x_1)\ge \mu(x_2) \ge \mu(x_3)\ge \mu(x_4)$; for brevity, we let $m_i:= \mu(x_i)$ for $i\in [4]$ and $a_{ij}:= \varphi(x_ix_j)$ for distinct $i<j\in [4]$. (Recall that this can fundamentally be viewed as graph Lagrangian problem wherein we are characterizing the extremal of the quadratic form $m^T\cdot A \cdot m$ where $m$ is a vector with entries $m_i$ and $A$ is the matrix with entries $a_{ij}$ and $a_{ii}=0$ for $i\in [4]$). 

Define $b_3:=\varepsilon^{0.02}$, $b_4:= \varepsilon^{0.06}$, $d_2:=\varepsilon^{0.08}$, $d_4:=\varepsilon^{0.02}$. We will consider different cases based on whether $m_3$ is larger or smaller than $b_3$, whether $m_4$ is larger or smaller than $b_4$, whether $|m_1-m_2|$ is larger or smaller than $d_2$, and whether $|m_1-m_4|$ is larger or smaller than $d_4$. For each case, we will present a value-extremal step graphon $W^*=(K_4,\mu^*,\varphi^*)$ that is close to $W$; again for brevity, we write $m_i^*:= \mu^*(x_i)$ for $i\in [4]$ and $a_{ij}^*:= \varphi^*(x_ix_j)$ for distinct $i<j\in [4]$ (note we use the identity isomorphism in each case). For most cases, we will also first present a graphon $W'=(K_4, \mu',\varphi')$ with $\varphi':=\varphi$ and with $|\mu'(x_i)-\mu(x_i)|\le \sigma$ for some small $\sigma > 0$ (indeed, $\mu'(x_i)$ will equal $\mu^*(x_i)$). Note then that $W$ is $\sigma$-close to $W$. An important observation then for all these cases is that $|W'| \ge |W| - 2\sigma \cdot \sum_{i<j\in [4]} a_{ij} \ge |W|-12\sigma$. We then will show that $W'$ is close to $W^*$ from whence it follows that $W$ is also close to $W^*$ as desired. 

\begin{claim}\label{cl:valuetight}
$|m_1-\frac{1}{2}|\le \sqrt{\varepsilon}$ or $|m_1+m_4-\frac{1}{2}|\le \sqrt{\varepsilon}$.    
\end{claim}
\begin{proofclaim}
Suppose not. As $|m_1-\frac{1}{2}| > \sqrt{\varepsilon}$, we find that $m_1( m_2+m_3+m_4) < \frac{1}{4}-\varepsilon$. Similarly as $|m_1+m_4-\frac{1}{2}| > \sqrt{\varepsilon}$, we find that $(m_1+m_4)( m_2+m_3) < \frac{1}{4}-\varepsilon$. But then it follows from Lemma~\ref{lem:Value4Parts}, that the value of $W$ is strictly less than $\frac{1}{4}-\varepsilon$, a contradiction.    
\end{proofclaim}

Note that by Claim~\ref{cl:valuetight}, we find that $m_1\le \frac{1}{2}+\sqrt{\varepsilon}$ and hence $m_2+m_3+m_4 \ge \frac{1}{2}-\sqrt{\varepsilon}$. Since $m_2\ge m_3\ge m_4$, we find that $m_2\ge \frac{1}{6} - \frac{\sqrt{\varepsilon}}{3} > 0$ (since $\varepsilon < \frac{1}{4}$). 

\begin{claim}
$m_3\ge b_3$.
\end{claim}
\begin{proofclaim}
Suppose not, that is $m_3 < b_3$. We let $W^*:=(K_4,\mu',\varphi')$ be defined as $m_1^*=m_2^*:=\frac{1}{2}$, $m_3^*=m_4^*:=0$; $a_{12}^*=1$, $a_{ij}^*=a_{ij}$ for distinct $i < j\in [4]$ where $\{i,j\}\ne \{1,2\}$. Note that $W^*$ is a value-extremal step graphon (specifically it is star-cut). Note that $|m_i-m_i^*|\le b_3$ for $i\in \{3,4\}$. By Claim~\ref{cl:valuetight}, we have that $m_1+m_4\ge \frac{1}{2}-\sqrt{\varepsilon}$ and $m_1\le \frac{1}{2}+\sqrt{\varepsilon}$; hence $|m_1-m_1^*|=|m_1-\frac{1}{2}| \le \sqrt{\varepsilon}+b_3 \le 2b_3$. Similarly $m_2\le \frac{1}{2}+\sqrt{\varepsilon}$ and $m_2+m_3+m_4 \ge \frac{1}{2}-\sqrt{\varepsilon}$ and hence $|m_2-m_2^*|=|m_2-\frac{1}{2}| \le \sqrt{\varepsilon}+2b_3 \le 3b_3$. 

By assumption $|W| \ge \frac{1}{4}-\varepsilon$. Yet $|W| \le a_{12}m_1m_2 + m_3 + m_4 \le a_{12}\cdot \frac{1}{4} + 2b_3$ where we used that $m_4\le m_3\le b_3$ and $m_1m_2\le \frac{1}{4}$. By rearranging, we find that $a_{12}\ge 1 - 4\varepsilon - 8b_3$. Thus $|a_{12}-a_{12}^*|\le 4\varepsilon+8b_3 \le 12b_3$. Combining all our observations, we find that $W$ is $12b_3$-close to $W^*$; since $12b_3=12\cdot \varepsilon^{0.02} \le \varepsilon^{0.01}$ as $\varepsilon$ is small enough, we find that $W$ is $\varepsilon^{0.01}$-close to $W^*$ as desired. 
\end{proofclaim}

\begin{claim}\label{cl:b4}
$m_4\ge b_4$.
\end{claim}
\begin{proofclaim}
Suppose not, that is $m_4 < b_4$. We let $W^*=(K_4,\mu^*,\varphi^*)$ be defined as  $m_1^*:=\frac{1}{2}$, $m_2^*:= \frac{m_2}{2(m_2+m_3)}$, $m_3^*:= \frac{m_3}{2(m_2+m_3)}$, $m_4^*:=0$; $a_{12}^*:= a_{13}^*:=1$, $a_{23}^*:=0$ $a_{i4}^*=a_{i4}$ for $i \in [3]$. (Recall this is well-defined as $m_2 \ge m_3 \ge b_3 > 0$.) Note that $W^*$ is a value-extremal step graphon (specifically it is star-cut). Note that $|m_4-m_4^*|\le b_4$. By Claim~\ref{cl:valuetight}, we have that $m_1+m_4\ge \frac{1}{2}-\sqrt{\varepsilon}$ and $m_1\le \frac{1}{2}+\sqrt{\varepsilon}$; hence $|m_1-m_1^*|=|x_1-\frac{1}{2}| \le \sqrt{\varepsilon}+b_4 \le 2b_4$. Similarly $m_2+m_3 = 1 -(m_1+m_4)\le \frac{1}{2}+\sqrt{\varepsilon}+b_4$ and $m_2+m_3+m_4 \ge \frac{1}{2}-\sqrt{\varepsilon}$ and hence $|m_2+m_3-\frac{1}{2}| \le \sqrt{\varepsilon}+b_4 \le 2b_4$. Thus $|m_2-m_2^*| = |m_2+m_3- \frac{1}{2}| \cdot \frac{m_2}{m_2+m_3} \le \sqrt{\varepsilon}+b_4 \le 2b_4$ and similarly $|m_3-m_3^*| = |m_2+m_3- \frac{1}{2}| \cdot \frac{m_3}{m_2+m_3} \le \sqrt{\varepsilon}+b_4 \le 2b_4$. 

By assumption $|W|\ge \frac{1}{4}-\varepsilon$. Yet $|W|\le a_{12}m_1m_2+a_{13}m_2m_3 + a_{23}m_2m_3 + m_4$. Substituting that $m_4\le b_4$ by assumption and $a_{23} \le 2-a_{12}-a_{13}$ by triangle-thinness of $W$, we find that 
\begin{align*}
\frac{1}{4}-\varepsilon -b_4 &\le a_{12}m_1m_2+a_{13}m_1m_3 + (2-a_{12}-a_{13})m_2m_3 \\
&= m_1m_2 + m_1m_3 - (1-a_{12})m_2(m_1-m_3) - (1-a_{13})m_3(m_1-m_2) \\
&\le \frac{1}{4} - (1-a_{12})\cdot b_3\cdot \frac{b_3}{2}-(1-a_{13})\cdot b_3\cdot \frac{b_3}{2},
\end{align*}
where we used that $m_2\ge m_3\ge b_3$, that $m_1\ge \frac{1}{2}-\sqrt{\varepsilon}-b_4$ and $m_2 \le \frac{1}{2}+\sqrt{\varepsilon}-b_3$ and hence $m_1-m_3\ge m_1-m_2\ge b_3-b_4-2\sqrt{\varepsilon} \ge b_3-3b_4 \ge \frac{b_3}{2}$ as $b_3=\varepsilon^{0.02}$, $b_4=\varepsilon^{0.06}$, and $\varepsilon$ is small enough. Rearranging this gives
$$(2-a_{12}-a_{13})\le 2\cdot \frac{b_4+\varepsilon}{b_3^2} \le 4\cdot \frac{\varepsilon^{0.06}}{(\varepsilon^{0.02})^2} \le 4\varepsilon^{0.02} \le \varepsilon^{0.01},$$
where the last inequality follows since $\varepsilon$ is small enough.
Thus we find that $a_{12} \ge 1 - \varepsilon^{0.01}$, $a_{13}\ge 1 -\varepsilon^{0.01}$ and $a_{23} \le 2-a_{12}-a_{23} \le \varepsilon^{0.01}$. Thus $W$ is $\varepsilon^{0.01}$-close to $W^*$ as desired.
\end{proofclaim} 

\begin{claim}
$m_1 < \frac{1}{2} - \frac{b_4}{2}$    
\end{claim}
\begin{proofclaim}
Suppose not, that is $m_1 \ge \frac{1}{2}-\frac{b_4}{2}$. Since $m_4\ge b_4$ by Claim~\ref{cl:b4}, we have that $m_1+m_4 \ge \frac{1}{2} + \frac{b_4}{2} > \frac{1}{2} + \sqrt{\varepsilon}$ since $b_4 > 2\sqrt{\varepsilon}$ as $\varepsilon$ is small enough. But then by Claim~\ref{cl:valuetight}, we find that $m_1 \ge \frac{1}{2} - \sqrt{\varepsilon}$. Let $W'$ be defined as $m_1':= \frac{1}{2}$ and for $j\in \{2,3,4\}$, $m_j':= \frac{m_j}{2(m_2+m_3+m_4)}$, and $a_{ij}'=a_{ij}$ for all $i<j\in [4]$; we will use later the fact that here $m_1'-m_2'-m_3' = m_4'$. Hence $|m_1-m_1'|\le \sqrt{\varepsilon}$. This also implies that $|m_2+m_3+m_4-\frac{1}{2}|\le \sqrt{\varepsilon}$. But then $\frac{1}{2(m_2+m_3+m_4)} \ge \frac{1}{1+2\sqrt{\varepsilon}} \ge 1 -2\sqrt{\varepsilon}$ and similarly $\frac{1}{2(m_2+m_3+m_4)} \le \frac{1}{1-2\sqrt{\varepsilon}} \le 1+3\sqrt{\varepsilon}$ (since $\sqrt{\varepsilon} \ge 6\varepsilon$ as $\varepsilon\le 1/36$). Thus we find that $|m_i-m_i'|\le 3\sqrt{\varepsilon}$ for all $i\in \{2,3,4\}$. Thus $W$ is $3\sqrt{\varepsilon}$-close to $W'$. Moreover, we find that $|W'|\ge |W| - 36\sqrt{\varepsilon}$. 

Let $W^*$ be defined as $m_i^*:=m_i'$ for all $i\in [4]$, $a_{1i}^*:=1$ for $i\in\{2,3,4\}$ and $a_{23}^*=a_{24}^*=a_{34}^*:=0$. Note that $W^*$ is a value-extremal step graphon (specifically it is star-cut). Hence $|W^*|=\frac{1}{4}$. Thus $|W^*|-|W'| \le \varepsilon+36\sqrt{\varepsilon}$. On the other hand, we calculate that
\begin{align*}
|W^*|-|W'| &= \bigg(\sum_{i\in \{2,3,4\}} (1-a_{1i})m_1'm_i'\bigg) - a_{23}m_2'm_3'-a_{24}m_2'm_4'-a_{34}m_3'm_4'  \\
&\ge \bigg(\sum_{i\in \{2,3,4\}} (1-a_{1i})m_1'm_i'\bigg) - (2-a_{12}-a_{13})m_2'm_3'-(2-a_{12}-a_{14})m_2'm_4'-(2-a_{13}-a_{14})m_3'm_4'\\
&= (1-a_{12})m_2'(m_1'-m_3'-m_4') + (1-a_{13})m_3'(m_1'-m_2'-m_4')+(1-a_{14})m_4'(m_1'-m_2'-m_3')\\
&\ge (3-a_{12}-a_{13}-a_{14})b_4\cdot b_4 \cdot (1-2\sqrt{\varepsilon}),
\end{align*}
where we used triangle--thinness for upper bounds on $a_{23},a_{24},a_{34}$, and that $m_2'\ge m_3'\ge m_4'\ge b_4$ and that $m_1'-m_3'-m_4'\ge m_1'-m_2'-m_4'\ge m_1'-m_2'-m_3' = m_4' \ge m_4(1-2\sqrt{\varepsilon})\ge b_4(1-2\sqrt{\varepsilon})$. Thus we find that $a_{12}+a_{13}+a_{14} \ge 3-2\cdot \frac{\varepsilon+36\sqrt{\varepsilon}}{b_4^2} \ge 3 - 74 \cdot \frac{\sqrt{\varepsilon}}{(\varepsilon^{0.06})^2} \ge 3 - \varepsilon^{0.30}$ as $\varepsilon$ is small enough. Since $a_{ij}\le 1$ for all $i<j\in [4]$, we find that $a_{1i} \ge 1 - \varepsilon^{0.30}$. Similarly, it follows that $a_{23},a_{24},a_{34} \le \varepsilon^{0.30}$. Thus $W'$ is $ \varepsilon^{0.30}$-close to $W^*$. Hence $W$ is $(3\sqrt{\varepsilon} + \varepsilon^{0.30})$-close to $W^*$; since $3\sqrt{\varepsilon} +  \varepsilon^{0.30} \le \varepsilon^{0.01}$ as $\varepsilon$ is small enough we find that $W$ is $\varepsilon^{0.01}$-close to $W^*$ as desired. 
\end{proofclaim}

By Claim~\ref{cl:valuetight} since $m_1 < \frac{1}{2} -\sqrt{\varepsilon}$, we find that $|m_1+m_4-\frac{1}{2}| \le \sqrt{\varepsilon}$.

\begin{claim}\label{cl:CycleCut}
$m_1-m_2 \le d_2$.    
\end{claim}
\begin{proofclaim}
Suppose not, that is $m_1-m_2 > d_2$. Let $W'$ be defined as $m_1':= \frac{m_1}{2(m_1+m_4)}$, $m_4':=\frac{m_4}{2(m_1+m_4)}$, $m_2':= \frac{m_2}{2(m_2+m_3)}$, $m_3':=\frac{m_3}{2(m_2+m_3)}$ and $a_{ij}'=a_{ij}$ for all $i<j\in [4]$. Since $|m_1+m_4-\frac{1}{2}| \le \sqrt{\varepsilon}$, we find that $\frac{1}{2(m_1+m_4)} \ge \frac{1}{1+2\sqrt{\varepsilon}} \ge 1 -2\sqrt{\varepsilon}$ and similarly $\frac{1}{2(m_1+m_4)} \le \frac{1}{1-2\sqrt{\varepsilon}} \le 1+3\sqrt{\varepsilon}$ (since $\sqrt{\varepsilon} \ge 6\varepsilon$ as $\varepsilon\le 1/36$). It follows then $|m_4-m_4'|\le |m_1-m_1'| \le 3\sqrt{\varepsilon}$. Similarly $|m_2+m_3-\frac{1}{2}|\le \sqrt{\varepsilon}$ and via the same argument that $|m_3-m_3'|\le |m_2-m_2'|\le 3\sqrt{\varepsilon}$.  Thus $W$ is $3\sqrt{\varepsilon}$-close to $W'$ and we also find that $|W'|\ge |W|-36\sqrt{\varepsilon}$.

Let $W^*$ be defined as $m_i^*:=m_i'$ for $i\in [4]$ and $a_{12}^*=a_{13}^*=a_{24}^*=a_{34}^*:= 1$ and $a_{14}^*=a_{23}^*=0$. Note that $W^*$ is a value-extremal step graphon (specifically it is cycle-cut). Hence $|W^*|=\frac{1}{4}$. Thus $|W^*|-|W'| \le \varepsilon+36\sqrt{\varepsilon}$. On the other hand, we calculate that
\begin{align*}
|W^*|-|W'| &= \bigg(\sum_{i\in \{1,4\}} \sum_{j\in\{2,3\}} (1-a_{ij})m_i'm_j'\bigg) - a_{23}m_2'm_3'-a_{14}m_1'm_4'.
\end{align*}
Let $\delta:= \varepsilon^{0.32}$ and $\delta':=\varepsilon^{0.15}$. Now we case according to how $a_{23},a_{14}$ compare to $\delta$ and $\delta'$. First suppose that $a_{23},a_{14} < \delta'$. Then $a_{23}m_2'm_3'+a_{14}m_1'm_4' < 2\delta'$. It follows that for all $i\in \{1,4\}$ and $j\in \{2,3\}$ with $1-a_{ij}\le \frac{36\sqrt{\varepsilon} + 2\delta'}{b_4^2}$ as otherwise we obtain a contradiction. But then $W'$ is $\frac{36\sqrt{\varepsilon} + 2\delta'}{b_4^2}$-close to $W^*$ and hence $W$ is $\frac{72\sqrt{\varepsilon} + 2\delta'}{b_4^2}$-close to $W^*$; since $\frac{72\sqrt{\varepsilon} + 2\delta'}{b_4^2}\le 74\cdot \frac{\delta'}{b_4^2} = 74\cdot \frac{\varepsilon^{0.21}}{(\varepsilon^{0.099})^2} =74\cdot \varepsilon^{0.012} \le \varepsilon^{0.01}$ as $\varepsilon$ is small enough, we find that $W$ is $\varepsilon^{0.01}$-close to $W^*$ as desired. So we assume at least one of $a_{23}\ge \delta'$ or $a_{14}\ge \delta'$ holds.

Next suppose that $a_{23}\ge \delta'$ and $a_{14}\le \delta$. But then we calculate using that $m_i'\ge m_j'$ for $i<j\in [4]$ and using both that $a_{23}\le 2-a_{12}-a_{13}$ as well that $a_{23}\le 2-a_{24}-a_{34}$ by triangle-thinness that
\begin{align*}
|W^*|-|W'| &\ge (1-a_{12})m_1'm_3' + (1-a_{13})m_1'm_3' + (1-a_{24})m_3'm_4' + (1-a_{34})m_3'm_4' - a_{23}m_2'm_3'-a_{14}m_1'm_4' \\
&\ge (2-a_{12}-a_{13})m_1'm_3' + (2-a_{24}-a_{34})m_3'm_4' - a_{23}m_2'm_3' - \delta \\
&\ge a_{23}(m_1'+m_4'-m_2')m_3' - \delta \\
&\ge \delta'\cdot d_2\cdot \frac{b_3}{2} - \delta = \frac{1}{2} \cdot \varepsilon^{0.15}\cdot \varepsilon^{0.08}\cdot \varepsilon^{0.02} - \varepsilon^{0.32}> 36\sqrt{\varepsilon},
\end{align*}
a contradiction where we used that $m_3'\ge \frac{m_3}{2} \ge \frac{b_3}{2}$ and the last inequality follows since $\varepsilon$ is small enough.

Similarly we suppose that $a_{14}\ge \delta'$ and $a_{23}\le \delta$. But then we calculate using that $m_i'\ge m_j'$ for $i<j\in [4]$ and using both that $a_{14}\le 2-a_{12}-a_{24}$ and $a_{14}\le 2-a_{13}-a_{34}$ by triangle-thinness that
\begin{align*}
|W^*|-|W'| &\ge (1-a_{12})m_2'm_4' + (1-a_{13})m_3'm_4' + (1-a_{24})m_2'm_4' + (1-a_{34})m_3'm_4' - a_{23}m_2'm_3'-a_{14}m_1'm_4' \\
&\ge (2-a_{12}-a_{24})m_2'm_4' + (2-a_{13}-a_{34})m_3'm_4' - a_{14}m_1'm_4' - \delta \\
&\ge a_{14}(m_2'+m_3'-m_1')m_4' - \delta \\
&\ge \delta'\cdot \frac{b_4}{2}\cdot \frac{b_4}{2} - \delta = \frac{1}{4} \varepsilon^{0.15}\cdot (\varepsilon^{0.06})^2 - \varepsilon^{0.32}> 36\sqrt{\varepsilon},
\end{align*}
a contradiction where we used that $m_4'\ge \frac{b_4}{2}$ and hence $m_2'+m_3'-m_1' = m_4' \ge \frac{b_4}{2}$ and the last inequality follows since $\varepsilon$ is small enough.

Finally we assume that $a_{14},a_{23}\ge \delta$. Let $M_1 := \{14,23\}, M_2:= \{24,13\}, M_3:=\{34,12\}$. Let $S:= \{ij \in M_2\cup M_3 : 1-a_{ij} \ge \frac{\delta}{2} \}$. Since $2-a_{12}-a_{13}\ge a_{23}$ by triangle-thinness and $a_{23}\ge \delta$, it follows that $|S\cap \{12,13\}|\ge 1$. Similarly since $2-a_{24}-a_{34}\ge a_{23} \ge \delta$, it follows that $|S\cap \{24,34\}|\ge 1$. As $2-a_{12}-a_{14}\ge a_{14} \ge \delta$, it follows that $|S\cap \{12,24\}|\ge 1$. As $2-a_{13}-a_{34}\ge a_{14} \ge \delta$, it follows that $|S\cap \{13,34\}|\ge 1$. Combining these observations, we find that $M_2\subseteq S$ or $M_3\subseteq S$. 

Choose $i\in \{2,3\}$ with $M_i\subseteq S$. Let $W^+$ be a step graphon defined as $m^+_i=m'_i$ for $i\in [4]$, $a_{jk}^+:= a_{jk}-\frac{\delta}{2}$ for $jk \in M_1$, $a_{jk}^+:= a_{jk}+\frac{\delta}{2}$ for $jk \in M_i$, $a_{jk}^+:= a_{jk}$ for $jk \in M_\ell$ where $\ell \in \{2,3\}\setminus \{i\}$. Since $a_{14},a_{23}\ge \delta$, we have that $a_{14}^+,a_{23}^+\ge 0$. Since $M_i\subseteq S$, we also have that $a_{jk}^+\le 1$ for all $jk\in M_i$. Thus $W^+$ is indeed a step graphon. Furthermore since every triangle contains at least one element of each $M_1$, it follows that $W^+$ is triangle-thin. By Lemma~\ref{lem:Value4Parts}, we find that $|W^+|\le \frac{1}{4}$. 

Further suppose that $i=2$. Then
\begin{align*}
|W^+|-|W'| &= \frac{\delta}{2} \cdot ( m_1'm_3' + m_2'm_4' - m_1'm_4'-m_2'm_3') \\
&= \frac{\delta}{2} (m_1'-m_2')(m_3'-m_4') \\
&\ge \frac{\delta}{2}\cdot \frac{d_2}{2} \cdot \frac{d_2}{2} = \frac{1}{8}\cdot \varepsilon^{0.32} \cdot (\varepsilon^{0.08})^2 > 36\sqrt{\varepsilon},   
\end{align*}
a contradiction, where we used that $m_3'-m_4'=m_1'-m_2' \ge m_1-m_2 -6\sqrt{\varepsilon}\ge d_2-6\sqrt{\varepsilon}\ge \frac{d_2}{2}$ and the last inequality follows since $\varepsilon$ is small enough.

So we assume that $i=3$. Then
\begin{align*}
|W^+|-|W'| &= \frac{\delta}{2} \cdot ( m_1'm_2' + m_3'm_4' - m_1'm_4'-m_2'm_3') \\
&= \frac{\delta}{2} (m_1'-m_3')(m_2'-m_4') \\
&\ge \frac{\delta}{2}\cdot \frac{d_2}{2} \cdot \frac{d_2}{2} > 36\sqrt{\varepsilon},   
\end{align*}
a contradiction, where we used that $m_2'-m_4'=m_1'-m_3' \ge m_1'-m_2'\ge \frac{d_2}{2}$ as above and the last inequality follows as above.
\end{proofclaim}

\begin{claim}\label{cl:ReflectiveCut}
$m_1-m_4 < d_4$.
\end{claim}
\begin{proofclaim}
Suppose not, that is $m_1-m_4 \ge d_4$. By Claim~\ref{cl:CycleCut}, we have that $m_1-m_2 \le d_2$.  Let $W'$ be the step-graphon with $m'_1=m'_2 = \frac{m_1+m_2}{2}=:\overline{m}_{12}$ and $m'_3=m'_4=\frac{m_3+m_4}{2}=:\overline{m}_{34}$ and $a_{ij}':=a_{ij}$ for all $i<j\in[4]$. Since $m_1-m_2<d_2$, it follows that $|m_{\ell}-m_{\ell}'|\le \frac{d_2}{2}$ for $\ell\in [2]$. On the other hand, noting  that $\overline{m}_{12}+\overline{m}_{34} = \frac{1}{2}$ and therefore that $m_1' = \frac{1}{2}-m_4'$, then

$$|m_4-m_4'| = \bigg|m_4+\frac{1}{2}-m_4' - \frac{1}{2}\bigg| = \bigg|m_4+m_1'-\frac{1}{2}\bigg|\le \bigg|m_4+m_1-\frac{1}{2}\bigg| + |m_1-m_1'| \le \sqrt{\varepsilon}+\frac{d_2}{2}.$$
Hence $W$ is $\sqrt{\varepsilon}+\frac{d_2}{2}$-close to $W'$. In addition, we note that $|W'|\ge |W|-12(\sqrt{\varepsilon}+\frac{d_2}{2}) \ge \frac{1}{4}-12\sqrt{\varepsilon}-6d_2 \ge \frac{1}{4} - 18d_2$.

Let $M_1 := \{14,23\}, M_2:= \{24,13\}, M_3:=\{34,12\}$. For $i\in[3]$, let $\overline{a}_i:= \frac{1}{2} \sum_{jk\in M_i} a_{jk}$. Note that $\overline{a}_i\in [0,1]$. Let $W^*$ be defined as $m_i^*:=\frac{1}{4}$ for all $i\in [4]$; for $\ell\in [2]$ and $jk\in M_\ell$, define $a_{jk}^*:=\frac{\overline{a}_\ell}{\overline{a}_1+\overline{a}_2}$; for $jk\in M_3$, define $a_{jk}^* := 1$. Note that $W^*$ is a value-extremal step graphon (specifically it is reflective-cut) and hence $|W^*|=\frac{1}{4}$. Thus $|W^*|-|W'| \le 18\cdot d_2$. On the other hand, again noting that $\overline{m}_{12}+\overline{m}_{34} = \frac{1}{2}$, we calculate that
\begin{align*}
|W^*|-|W'| = \frac{1}{4}-|W'| &= (\overline{m}_{12}+\overline{m}_{34})^2 - |W'|\\
&= (1-a_{12})\overline{m}_{12}^2 + (2-2(\overline{a}_1+\overline{a}_2))\overline{m}_{12}\overline{m}_{34}+(1-a_{34})\overline{m}_{34}^2 \\
&\ge (\overline{a}_1+\overline{a}_2-1)\overline{m}_{12}^2  - 2 (\overline{a}_1+\overline{a}_2-1)\overline{m}_{12}\overline{m}_{34}+(\overline{a}_1+\overline{a}_2-1)\overline{m}_{34}^2 \\
&= (\overline{a}_1+\overline{a}_2-1)\cdot (\overline{m}_{12}-\overline{m}_{34})^2\\
&\ge (\overline{a}_1+\overline{a}_2-1) \cdot \left(\frac{d_4}{2}\right)^2
\end{align*}
where for the first inequality we used that for $jk\in M_3$, we have that $1-a_{jk}\ge \overline{a}_1+\overline{a}_2-1$ by the triangle-thinness of $W$ (applied to the average of two triangles), and for the last inequality we used that $\overline{m}_{12}-\overline{m}_{34} \ge \frac{m_1-m_4}{2} \ge \frac{1}{2}\cdot d_4$ (since $m_2\ge m_3$). It follows that $$\overline{a}_1+\overline{a}_2 \le 1 + \frac{4}{d_4^2} \cdot (18d_2).$$
Substituting this into the second line of the above set of equations, we find that
$$|W^*|-|W'| \ge (1-a_{12})\overline{m}_{12}^2 -144\cdot \frac{d_2}{d_4^2} \cdot \overline{m}_{12}\overline{m}_{34}+(1-a_{34})\overline{m}_{34}^2.$$
Since $\overline{m}_{12}\ge \frac{1}{4}$ and $\overline{m}_{12}\ge \overline{m}_{34}$, it follows that $1-a_{12}\le 18d_2+ \frac{144d_2}{d_4^2} \le 162\cdot \frac{d_2}{d_4^2} = 162\cdot \varepsilon^{0.08}\cdot (\varepsilon^{-0.02})^2 \le \varepsilon^{0.01}$ since $\varepsilon$ is small enough. Since $\overline{m}_{34}\ge \frac{b_3+b_4}{2} \ge \frac{b_3}{2}$, it follows that 
$$1-a_{34} \le \frac{1}{\overline{m}_{34}} \cdot \left(18d_2 + 144\cdot \frac{d_2}{d_4^2} \right)\le 162\cdot \frac{d_2}{d_4^2\cdot b_3} \le 162\cdot \frac{\varepsilon^{0.08}}{(\varepsilon^{0.02})^2\cdot \varepsilon^{0.02}} = 162 \cdot \varepsilon^{0.02} \le \varepsilon^{0.01}$$
since $\varepsilon$ is small enough. 
But since $a_{12},a_{34}\le 1$, we also have that $|W^*|-|W'| \ge 2(1-\overline{a}_1-\overline{a}_2)\overline{m}_{12}\overline{m}_{34} \ge \frac{1}{2}\cdot (1-\overline{a}_1-\overline{a}_2)\overline{m}_{34} $ (because $\overline{m}_{12}\ge \frac{1}{4}$) and hence 
$$\overline{a}_1+\overline{a}_2 \ge 1 - \frac{36d_2}{b_3} \ge 1-36\cdot \frac{\varepsilon^{0.08}}{\varepsilon^{0.02}} = 1 -36\cdot \varepsilon^{0.06}.$$ 
It then follows that $W'$ is $162 \varepsilon^{0.02}$-close to $W^*$ and hence $W$ is $\varepsilon^{0.01}$-close to $W^*$ as desired. 
\end{proofclaim}

By Claim~\ref{cl:ReflectiveCut}, we have that $m_1-m_4 < d_4$. Let $W'$ be the step-graphon with $m'_i:=\frac{1}{4}$ for all $i\in [4]$ and $a_{ij}':=a_{ij}$ for all $i<j\in[4]$. Since $m_1-m_4<d_4$, it follows that $|m_i-m_i'|\le d_4$ for all $i\in [4]$ and hence that $W$ is $d_4$-close to $W'$. In addition, we note that $|W'|\ge |W|-12d_4 \ge \frac{1}{4}-\varepsilon-12d_4$. 

Let $M_1 := \{14,23\}, M_2:= \{24,13\}, M_3:=\{34,12\}$. For $i\in[3]$, let $\overline{a}_i:= \frac{1}{2} \sum_{jk\in M_i} a_{jk}$. Note that $\overline{a}_i\in [0,1]$. Note that $|W'| = \frac{1}{8}\cdot \sum_{i\in [3]} \overline{a}_i$. Hence $\sum_{i\in [3]} \overline{a}_i \ge 2-8\varepsilon-96d_4$. Notice that triangle-thinness implies that for $jk\in M_i$ (applied to the average of two triangles), we find that $a_{jk} \le 2 - \sum_{\ell \in [3]\setminus\{i\}} \overline{a}_{\ell}$. Combining the two inequalities above yields that $a_{jk}\le \overline{a}_i + 8\varepsilon + 96d_4$ for $jk\in M_i$. Conversely that also implies (by applying to the other edge of the matching) that $a_{jk} \ge \overline{a}_i + 8\varepsilon + 96d_4$ and hence $|a_{jk}-\overline{a}_i|\le 8\varepsilon+96d_4$. 

Choose $i\in [3]$ such that $\overline{a}_i$ is minimized. Let $W^*$ be defined as $m_i^*:=\frac{1}{4}$ for all $i\in [4]$; for $\ell\in[3]\setminus \{i\}$ and $jk\in M_\ell$, define $a_{jk}^*:=\overline{a}_\ell$; for $jk\in M_i$, define $a_{jk}^* := 2 - \sum_{\ell\in [3]\setminus \{i\}} \overline{a}_{\ell}$. Since $\sum_{i\in [3]} \overline{a}_i \ge 2-8\varepsilon-96d_4$, we find by the choice of $i$ that $\sum_{\ell\in [3]\setminus \{i\}} \overline{a}_{\ell} \ge \frac{2}{3} \cdot (2-8\varepsilon-96d_4) \ge 1$ since $8\varepsilon+96d_4\le \frac{1}{2}$ as $\varepsilon$ is small enough. Hence $W^*$ is a step graphon; moreover is triangle-thin by construction. Indeed $W^*$ is a value-extremal graphon (specifically it is quad-cut). Finally note that $|\overline{a}_i - (2-\sum_{\ell\in [3]\setminus \{i\}} \overline{a}_{\ell})| \le 8\varepsilon+96d_4$ and hence for $jk\in M_i$, we have that $|a_{jk}-a_{jk}^*| \le 16\varepsilon+192d_4$. Hence we find that $W'$ is $(8\varepsilon+192d_4)$-close to $W^*$. But then $W$ is $(8\varepsilon+193d_4)$-close to $W^*$; since $8\varepsilon+193d_4 = 8\varepsilon+193\cdot \varepsilon^{0.02} \le \varepsilon^{0.01}$ as $\varepsilon$ is small enough, we find that $W$ is $\varepsilon^{0.01}$-close to $W^*$ as desired.
\end{proof}

\subsection{Full Value Stability}

Finally, we formulate our stability version of Lemma~\ref{lem:Value}. To that end, we make the following definition.

\begin{definition}
Let $(H,\varphi)$ be a $[0,1]$-edge-weighted graph on $n$ vertices. We say $(H,\varphi)$ is \emph{$\delta$-close} to a step graphon $(G,\mu,\varphi')$ on $n'$ vertices if, letting $x_1,\ldots,x_{n'}$ be an enumeration of $V(G)$, there exists a (not necessarily proper) $n'$-coloring $c$ of $V(H)$ with color classes $(X_i: i\in n')$ such that $\bigg|\mu(x_i)-\frac{|X_i|}{n}\bigg| \le \delta$ and at most $\delta n^2$ of the pairs $(u,v)$ with $u,v\in V(H)$ satisfy that $|\varphi(uv)-\varphi'(x_{c(u)}x_{c(v)})| > \delta$ where we also set $\varphi(uv):=0$ if $uv\not\in E(H)$ and $\varphi'(x_ix_i):=0$ for all $i\in [n']$.
\end{definition}

For the proof, we also need the following somewhat elementary proposition about random variables.

\begin{proposition}\label{prop:Expect}
Let $X$ be a random variable with values in $[0,1]$. Let $0\le a < b \le 1$. If \  $\Prob{ X \not\in [a,b]} \le p$, then $\Expect{X} \in [a-p,b+p]$ and hence $\Prob{~|X-\Expect{X}| > p+(b-a)~}\le p$.
\end{proposition}
\begin{proof}
We assume without loss of generality that $\Prob{X\not\in [a,b]}=p$. Let $X_1$ be the subset of outcomes of $X$ that lie in $[a,b]$ and $X_2$ be the set of outcomes of $X$ that do not, namely those in $[0,a) \cup (b, 1]$. Let $c:= \Expect{X_1}$; note that $c\in [a,b]$. Thus $\Expect{X} = \Prob{X\in [a,b]~}\cdot\Expect{X_1} + \Prob{X\not\in [a,b]~}\cdot \Expect{X_2} = (1-p)c + p\cdot \Expect{X_2} = c + p(\Expect{X_2}-c)$. Since $X$ lies in $[0,1]$, we have that $\Expect{X_2}\in [0,1]$ and hence $\Expect{X_2}-c \in [-1,+1]$. Thus $\Expect{X}\in [a-p,b+p]$ as desired.
\end{proof}

We are now prepared to state and prove our stability version of Lemma~\ref{lem:Value}.

\begin{lem}[Value Stability]\label{lem:ValueStability}
There exists $\varepsilon_0 > 0$ such that for all $\varepsilon \in (0,\varepsilon_0]$ the following holds: If $H$ is a $K_5$-free graph on $n$ vertices, $\varphi$ is a triangle-thin $[0,1]$-edge-weighting of $H$, and $|\varphi| \ge \frac{n^2}{4}-\varepsilon n^2$, then $(H,\varphi)$ is $\varepsilon^{0.001}$-close to some value-extremal step graphon. 
\end{lem}
\begin{proof}
We choose $\varepsilon_0$ small enough as needed so that various inequalities throughout the proof regarding $\varepsilon$ hold true. For $v\in V(H)$, let $d_{\varphi}(v) := \sum_{e\in E(H): v\in e} \varphi(e)$. Note by hand-shaking that $|\varphi| = \frac{1}{2} \sum_{v\in V(H)} d_{\varphi}(v)$. 

Let $H_1:= H$. For each $i\in [4]$, if $V(H_i)\ne\emptyset$, choose $v_i\in V(H_i)$ such that $d_{\varphi}(v_i)$ is maximized and set $X_i:= V(H_i)\setminus N_{H}(v_i)$ and $H_{i+1}:= H_i\setminus X_i$; if $V(H_i):= \emptyset$, set $X_i:=\emptyset$ and $H_{i+1}:=\emptyset$.

We claim that $V(H_5)=\emptyset$. Suppose not, that is there exists $v_5\in V(H_5)$. But then $v_1,v_2,v_3,v_4$ all exist and hence are pairwise adjacent by construction. But then as $v_5 \not\in X_1\cup X_2\cup X_3\cup X_4$, we find that $v_5$ is adjacent to all of $v_1, v_2, v_3, v_4$, contradicting that $H$ is $K_5$-free. This proves the claim. Now let $c$ be the (not necessarily proper) $4$-coloring of $V(H)$ with color classes $X_1,X_2,X_3,X_4$.

\begin{claim}\label{cl:InsideEdgeCount}
For each $i\in [4]$, $\sum_{e\in E(H[X_i])} \varphi(e) \le \varepsilon n^2$.    
\end{claim}
\begin{proofclaim}
Suppose towards a contradiction that the statement does not hold for $i\in [4]$. Let $s:= \sum_{e\in E(H[x_i])} \varphi(e)$. Hence $s > \varepsilon n^2$.

Define a new edge-weighted graph $(H_0,\varphi_0)$ with $V(H_0):=V(H)$, 
$$E(H_0):= \{xy\in E(H): x,y\in V(H)\setminus X_i\}~\cup~\{xy: x\in X_i,~y\in V(H)\setminus X_i,~v_iy\in E(H)\},$$ and 
for an edge $xy\in E(H_0)$, we define $\varphi_0(xy):= \varphi(xy)$ if $x,y\in V(H)\setminus X_i$ and $\varphi_0(xy):= \varphi(v_iy)$ if $x\in X_i$ and $y\in V(H)\setminus X_i$ (i.e.~by replacing vertices of $X_i$ by weight-clones of $v_i$). It follows that $H_0$ is $K_5$-free (since $H$ is) and that $\varphi_0$ is a triangle-thin $[0,1]$-edge-weighting of $H_0$ (since $\varphi$ is of $H$). Hence by Lemma~\ref{lem:Value}, we find that $|\varphi_0|\le \frac{n^2}{4}$. Yet 
$$|\varphi_0| = |\varphi| - \bigg(\sum_{v\in X_i} d_{\varphi}(v)\bigg) + |X_i|\cdot d_{\varphi}(v_i) + s \ge  |\varphi| + s > \frac{n^2}{4},$$
a contradiction - where for the first inequality we used that $d_{\varphi}(v_i) \ge d_{\varphi}(v)$ for all $v\in X_i$ by the choice of $v_i$.
\end{proofclaim}

Let $H^*$ be the complete $4$-partite graph with partition classes $X_1,X_2,X_3,X_4$ and let $\varphi^*$ be the $[0,1]$-edge-weighting of $H^*$ where for $e\in E(H^*)$, we define $\varphi^*(e):=\varphi(e)$ if $e\in E(H)$ and $\varphi^*(e):=0$ if $e\not\in E(H)$. We assume without loss of generality going forward that $|X_1|\ge |X_2|\ge |X_3|\ge |X_4|$ (as otherwise we may rename the classes accordingly for the remainder for the proof). By Claim~\ref{cl:InsideEdgeCount}, we find that $|\varphi^*|\ge \frac{n^2}{4}-5\varepsilon n^2$.

Now we define a step graphon $\overline{W} := (K_4,\overline{\mu},\overline{\varphi})$ where $V(K_4):=\{x_1,x_2,x_3,x_4\}$, $\overline{\mu}(x_i) := \frac{|X_i|}{n}$ for all $i\in [4]$ and $\overline{\varphi}(x_ix_j):= \frac{1}{|X_i|\cdot |X_j|}\cdot \sum_{u\in X_i, v\in X_j} \varphi(uv)$ and for all distinct $i,j\in [4]$. Note that the ranges of $\overline{\mu}$ and $\overline{\varphi}$ lie in $[0,1]$ by construction. Moreover, $|\overline{W}| = n^2\cdot |\varphi_{*}|$ by construction and hence $|\overline{W}| \ge \frac{1}{4}-5\varepsilon$.  Note that $\overline{W}$ is simply the average of all transversal step graphons of $H^*$. Since each of these is triangle-thin by assumption, so is $\overline{W}$. By Lemma~\ref{lem:ValueStability2}, we thus find that $\overline{W}$ is $(5\varepsilon)^{0.01}$-close to a value-extremal step graphon $W$. 

Let $\delta':= \frac{\varepsilon^{0.001}}{2}$. We now proceed to show in the remainder of the proof that $H$ is $\delta'$-close to $\overline{W}$ and hence $H$ is $\varepsilon^{0.001}$-close to $W$ as desired (since $(5\varepsilon)^{0.01}+\frac{\varepsilon^{0.001}}{2}\le \varepsilon^{0.001}$ as $\varepsilon$ is small enough). 

To that end, let $B:= \{(u,v)\in V(H): |\varphi(uv)-\overline{\varphi}(x_{i}x_{j})| > \delta', u\in X_i, v\in X_j \text{ for some $i,j\in [4]$ }\}$ where we also set $\varphi(uv):=0$ if $uv\not\in E(H)$ and $\overline{\varphi}(x_ix_i):=0$ for all $i\in [4]$. Since $\overline{\mu}(x_i) = \frac{|X_i|}{n}$ by construction, to show that $H$ is $\delta'$-close to $\overline{W}$, it suffices to show that $|B|\le \delta'n^2$. To that end, we define for $i\in [4]$, $B_i:= \{(u,v)\in B: u,v\in X_i\}$, and for $i<j\in [4]$, $B_{ij}:= \{(u,v)\in B: u\in X_i, v\in X_j\}$. Hence
$|B|=\sum_{i\in [4]} |B_i| + 2\cdot \sum_{i<j\in[4]} |B_{ij}|$.

From Claim~\ref{cl:InsideEdgeCount}, we can derive the following upper bound on $\sum_{i\in [4]} |B_i|$.

\begin{claim}\label{cl:InsideBadCount}
For each $i\in [4]$, $|B_i|\le \sqrt{\varepsilon}n^2$ and hence $\sum_{i\in [4]} |B_i| \le 4\sqrt{\varepsilon}n^2$.
\end{claim}
\begin{proofclaim}
Suppose for a contradiction that $|B_i| > \sqrt{\varepsilon}n^2$ for some $i\in [4]$. But then $\sum_{e\in E(H[X_i])} \varphi(e) \ge |B_i|\cdot \delta' \ge \sqrt{\varepsilon}n^2\cdot \frac{\varepsilon^{0.001}}{2} > \varepsilon n^2$ contradicting Claim~\ref{cl:InsideEdgeCount}, where the last inequality follows since $\varepsilon$ is small enough.
\end{proofclaim}

Now we proceed with upper bounding $|B_{ij}|$ for all $i<j\in [4]$ to complete the proof. To that end, consider the probability space wherein we choose a vertex $u_i \in X_i$ uniformly at random (if $X_i=\emptyset$, choose a dummy vertex whose incident edge-weights are zero). Consider the resulting random step graphon $Z:=(K_4,\overline{\phi},\varphi_Z)$ where $\varphi_Z(x_ix_j)=\varphi(u_iu_j)$ for $i<j\in [4]$. Note that $Z$ is triangle-thin as $\overline{W}$ is. Here is a useful claim.

\begin{claim}\label{cl:sameExp}
$\Expect{~|Z|~}=|\overline{W}|$.
\end{claim}
\begin{proofclaim}
We calculate that
\begin{align*}
\Expect{~|Z|~} &= \frac{1}{\prod_{k\in [4]} |X_k|} \cdot \sum_{u_1\in X_1,u_2\in X_2,u_3\in X_3,u_4\in X_4}~~ \sum_{i<j\in [4]} \varphi(u_iu_j) \cdot \overline{\mu}(x_i) \cdot \overline{\mu}(x_j)\\
&= \frac{1}{\prod_{k\in [4]} |X_k|} \cdot \sum_{i<j\in [4]} \sum_{u_i\in X_i, u_j\in X_j}~~  \varphi(u_iu_j) \cdot \overline{\mu}(x_i) \cdot \overline{\mu}(x_j) \prod_{\ell \in [4]\setminus \{i,j\}} |X_{\ell}| \\
&= \sum_{i<j\in [4]} \frac{1}{|X_i|\cdot |X_j|} \cdot \sum_{u_i\in X_i, u_j\in X_j}~~  \varphi(u_iu_j) \cdot \overline{\mu}(x_i) \cdot \overline{\mu}(x_j)\\
&= \sum_{i<j\in [4]} \overline{\varphi}(x_ix_j) \cdot \overline{\mu}(x_i) \cdot \overline{\mu}(x_j)\\
&= |\overline{W}|,
\end{align*} 
as desired.
\end{proofclaim}

Now consider the random variable $d:=\frac{1}{4}-|Z|$. It follows from Claim~\ref{cl:sameExp} that 

$$\Expect{d} = \frac{1}{4} - |\overline{W}| \le 5\varepsilon.$$

Let $\delta'':= (\sqrt{5\varepsilon})^{0.01}$. Let $\mc{Z}$ denote the family of triangle-thin step graphons on $4$ vertices which are $\delta''$-close to some value-extremal step graphon.  

\begin{claim}
$\Prob{ Z\not\in \mc{Z} } \le \sqrt{5\varepsilon}.$    
\end{claim}
\begin{proofclaim}
Since $d$ is nonnegative via Lemma~\ref{def:Value}, we have by Markov's inequality that $\Prob{d \ge \sqrt{5\varepsilon}} \le \sqrt{5\varepsilon}$. By Lemma~\ref{lem:ValueStability2}, if $z\in Z$ is not $\delta''$-close to any value-extremal step graphon, then $|z| < \frac{1}{4}-\sqrt{5\varepsilon}$ (where we used that $\sqrt{5\varepsilon}$ satisfies Lemma~\ref{lem:ValueStability2} as $\varepsilon$ is small enough). But then we have $\Prob{ Z\not\in \mc{Z} }\le \Prob{d\ge \sqrt{5\varepsilon}} \le \sqrt{5\varepsilon}$ as desired. 
\end{proofclaim}

For $i<j\in [4]$, let $p_{ij}:= \Prob{|\varphi_Z(x_ix_j)-\overline{\varphi}(x_ix_j)|> \delta'}$. It is useful to note that for all $i>j\in [4]$, we have that $\Expect{\varphi_Z(x_ix_j)} = \overline{\varphi}(x_ix_j)$. Before proceeding with the main proof, we record the following useful observation.

\begin{claim}\label{cl:EdgeProb}
For $i<j\in [4]$, $|B_{ij}|=p_{ij}\cdot |X_i|\cdot |X_j| \le p_{ij}\cdot n^2$.
\end{claim}
\begin{proofclaim}
Let $k,\ell$ be defined such that $\{i,j,k,\ell\}=[4]$. The claim then follows since
$$p_{ij} = \frac{|B_{ij}|\cdot |X_{k}|\cdot |X_{\ell}|}{|X_i|\cdot |X_j|\cdot |X_k|\cdot |X_{\ell}|} = \frac{|B_{ij}|}{|X_i|\cdot |X_j|},$$
and since $|X_i|,|X_j|\le n$.
\end{proofclaim}

Now we proceed somewhat along the lines of the proof of Lemma~\ref{lem:ValueStability2}. First we show that either $\mu(x_1)$ or $\mu(x_1)+\mu(x_4)$ is close to $\frac{1}{2}$ as follows.

\begin{claim}
$|\overline{\mu}(x_1)-\frac{1}{2}|\le \sqrt{5\varepsilon}$ or $|\overline{\mu}(x_1)+\overline{\mu}(x_4)-\frac{1}{2}|\le \sqrt{5\varepsilon}$.    
\end{claim}
\begin{proofclaim}
Suppose not. As $|\overline{\mu}(x_1)-\frac{1}{2}| > \sqrt{5\varepsilon}$, we find that $\overline{\mu}(x_1)( \overline{\mu}(x_2)+\overline{\mu}(x_3)+\overline{\mu}(x_4)) < \frac{1}{4}-5\varepsilon$. Similarly as $|\overline{\mu}(x_1)+\overline{\mu}(x_4)-\frac{1}{2}| > \sqrt{5\varepsilon}$, we find that $(\overline{\mu}(x_1)+\overline{\mu}(x_4))(\overline{\mu}(x_2)+\overline{\mu}(x_3)) < \frac{1}{4}-5\varepsilon$. But then it follows from Lemma~\ref{lem:Value4Parts}, that the value of $Z$ is strictly less than $\frac{1}{4}-5\varepsilon$ no matter the choice of $u_i$; hence $|\overline{W}|=\Expect{~|Z|~} < \frac{1}{4}-5\varepsilon$, a contradiction.    
\end{proofclaim}

We now proceed to case according to how large $\mu(x_4)$, $\mu(x_3)$, $\mu(x_1)$ are as follows. 

\begin{claim}\label{cl:FullValueX3}
$\overline{\mu}(x_3)> 4\delta''$.    
\end{claim}
\begin{proofclaim}
Suppose not, that is $\overline{\mu}(x_3) \le 4\delta''$. Recall that by Markov's inequality we have that $\Prob{d\ge \sqrt{5\varepsilon}} \le \sqrt{5\varepsilon}$. Let $z\in Z$ such that $|z| > \frac{1}{4} - \sqrt{5\varepsilon}$. Yet we also calculate that 
$$|z| \le \varphi_z(x_1x_2) \cdot \overline{\mu}(x_1)\cdot \overline{\mu}(x_2) + \overline{\mu}(x_3) + \overline{\mu}(x_4) \le \varphi_z(x_1x_2)\cdot \frac{1}{4} + 6\delta'',$$
and hence $\varphi_z(x_1x_2) \ge 1 - 4\sqrt{5\varepsilon}-24\delta''.$ Thus it follows that $\Prob{\varphi_Z(x_1x_2) < 1 - 4\sqrt{5\varepsilon}-24\delta''} \le \sqrt{5\varepsilon}.$ By Proposition~\ref{prop:Expect} applied to $\varphi_Z(x_1x_2)$ this implies that
$$\Prob{\bigg|\varphi_Z(x_1x_2) - \Expect{\varphi_Z(x_1x_2)}\bigg| >\sqrt{5\varepsilon} + 4\sqrt{5\varepsilon} + 24\delta''} \le \sqrt{5\varepsilon}.$$
Since $5\sqrt{5\varepsilon}+24\delta'' \le \delta'$ as $\varepsilon$ is small enough, it follows that $p_{12}\le \sqrt{5\varepsilon}$ and hence by Claim~\ref{cl:EdgeProb} that $|B_{12}|\le \sqrt{5\varepsilon}\cdot n^2$. On the other hand by Claim~\ref{cl:EdgeProb}, we find that $\sum_{i\in [3]} |B_{i4}| \le |X_4|\cdot n \le 4\delta''\cdot n^2$ and similarly that $|B_{13}|+|B_{23}|\le |X_3|\cdot n \le 4\delta''\cdot n^2$. Combining, we find that
$\sum_{i<j\in [4]} |B_{ij}| \le (\sqrt{5\varepsilon} + 8\delta'')\cdot n^2$.
Since $\sum_{i\in [4]} |B_i| \le 4\sqrt{\varepsilon}\cdot n^2$ by Claim~\ref{cl:InsideBadCount}, we find that $|B| \le 2\cdot (\sqrt{5\varepsilon} + 8\delta'')\cdot n^2 + 4\sqrt{\varepsilon}\cdot n^2 \le \delta'\cdot n^2$
as desired, where the inequality follows since $2\cdot(\sqrt{5\varepsilon} + 8(\sqrt{5\varepsilon})^{0.01}) + 4\sqrt{\varepsilon} \le \frac{\varepsilon^{0.001}}{2}$ as $\varepsilon$ is small enough.
\end{proofclaim}

\begin{claim}\label{cl:FullValueX4}
$\overline{\mu}(x_4) > 2\delta''$.    
\end{claim}
\begin{proofclaim}
Suppose not, that is $\overline{\mu}(x_4) \le 2\delta''$.  By Claim~\ref{cl:FullValueX3}, $\overline{\mu}(x_3) > 4\delta''$. Thus for all $z\in Z$, it follows that $z$ is not $\delta''$-close to a reflective-cut step graphon or a quad-cut graphon. Let $z\in Z\cap \mc{Z}$. Thus we have that $z$ is $\delta''$-close to a star-cut or cycle-cut step graphon. By the definition of star-cut and cycle-cut, we thus find that $\varphi_z(x_1x_2),\varphi_z(x_1x_3) \ge 1 -\delta''$ and $\varphi_z(x_2x_3) \le \delta''$. Hence $\Prob{\varphi_Z(x_1x_2) < 1-\delta''} \le \sqrt{5\varepsilon}$. By Proposition~\ref{prop:Expect} applied to $\varphi_Z(x_1x_2)$ this implies that
$$\Prob{\bigg|\varphi_Z(x_1x_2) - \Expect{\varphi_Z(x_1x_2)}\bigg| >\sqrt{5\varepsilon} + \delta''} \le \sqrt{5\varepsilon}.$$
Since $\sqrt{5\varepsilon}+\delta'' \le \delta'$ as $\varepsilon$ is small enough, it follows that $p_{12}\le \sqrt{5\varepsilon}$ and hence by Claim~\ref{cl:EdgeProb} that $|B_{12}|\le \sqrt{5\varepsilon}\cdot n^2$. Similar arguments applied to $\varphi_Z(x_1x_3)$ and $\varphi(x_2x_3)$ imply that $|B_{13}|,|B_{23}|\le \sqrt{5\varepsilon}\cdot n^2$. On the other hand by Claim~\ref{cl:EdgeProb}, we find that $\sum_{i\in [3]} |B_{i4}| \le |X_4|\cdot n \le 2\delta''\cdot n^2$. Combining, we find that
$\sum_{i<j\in [4]} |B_{ij}| \le 3\cdot (\sqrt{5\varepsilon} + 2\delta'')\cdot n^2$.
Since $\sum_{i\in [4]} |B_i| \le 4\sqrt{\varepsilon}\cdot n^2$ by Claim~\ref{cl:InsideBadCount}, we find that $|B| \le 6\cdot (\sqrt{5\varepsilon} + 2\delta'')\cdot n^2 + 4\sqrt{\varepsilon}\cdot n^2 \le \delta'\cdot n^2$
as desired, where the inequality follows since $6\cdot (\sqrt{5\varepsilon} + 2(\sqrt{5\varepsilon})^{0.01}) + 4\sqrt{\varepsilon} \le \frac{\varepsilon^{0.001}}{2}$ as $\varepsilon$ is small enough.
\end{proofclaim}

\begin{claim}\label{cl:FullValueX1}
$\overline{\mu}(x_1) < \frac{1}{2} - \delta''$.    
\end{claim}
\begin{proofclaim}
Suppose not, that is $\overline{\mu}(x_1) \ge \frac{1}{2}-\delta''$. Yet $\overline{\mu}(x_4) > 2\delta''$ by Claim~\ref{cl:FullValueX4}. Hence $\overline{\mu}(x_1)+\overline{\mu}(x_4) > \frac{1}{2} + \delta''$. Thus for all $z\in Z$, it follows that $z$ is not $\delta''$-close to any cycle-cut, reflective-cut or quad-cut graphon.

Let $S:= \{ i\in \{2,3,4\}: \overline{\mu}(x_i) > 0\}$. Let $z\in Z\cap \mc{Z}$. Thus we have that $z$ is $\delta''$-close to some star-cut step graphon. By the definition of star-cut, we find that $\varphi_z(x_1x_i) \ge 1-\delta''$ for all $i\in S$ and $\varphi_z(x_ix_j)$ for all $i<j\in S$. Hence $\Prob{\varphi_Z(x_1x_2) < 1-\delta''} \le \sqrt{5\varepsilon}$. By Proposition~\ref{prop:Expect} applied to $\varphi_Z(x_1x_2)$ this implies that
$$\Prob{\bigg|\varphi_Z(x_1x_2) - \Expect{\varphi_Z(x_1x_2)}\bigg| >\sqrt{5\varepsilon} + \delta''} \le \sqrt{5\varepsilon}.$$
Since $\sqrt{5\varepsilon}+\delta'' \le \delta'$ as $\varepsilon$ is small enough, it follows that $p_{12}\le \sqrt{5\varepsilon}$ and hence by Claim~\ref{cl:EdgeProb} that $|B_{12}|\le \sqrt{5\varepsilon}\cdot n^2$.
Similar arguments applied to $\varphi_Z(x_ix_j)$ for $i<j\in S\cup \{1\}$ imply that $|B_{ij}|\le \sqrt{5\varepsilon}\cdot n^2$. On the other hand for $i\in [4]$ and $j\in [4]\setminus S$, we have that $|B_{ij}|=0$.

Combining, we find that
$\sum_{i<j\in [4]} |B_{ij}| \le 6 \sqrt{5\varepsilon}\cdot n^2$.
Since $\sum_{i\in [4]} |B_i| \le 4\sqrt{\varepsilon}\cdot n^2$ by Claim~\ref{cl:InsideBadCount}, we find that $|B| \le 12\sqrt{5\varepsilon} \cdot n^2 + 4\sqrt{\varepsilon}\cdot n^2 \le \delta'\cdot n^2$
as desired, where the inequality follows since $12\sqrt{5\varepsilon} + 4\sqrt{\varepsilon} \le \frac{\varepsilon^{0.001}}{2}$ as $\varepsilon$ is small enough.
\end{proofclaim}

By Claim~\ref{cl:FullValueX1}, we have that $\overline{\mu}(x_1) < \frac{1}{2} - \delta''$. Thus it follows if $z\in Z$, then $z$ is not $\delta''$-close to a star-cut step graphon. Hence if $z\in Z\cap \mc{Z}$, we find that $z$ is $\delta''$-close to some cycle-cut, reflective-cut, or quad-cut step graphon. Since each of these are matched, this implies that for all $z\in Z\cap \mc{Z}$, we have that 
$$|\varphi_z(x_1x_2)-\varphi_z(x_3x_4)|, |\varphi_z(x_1x_3)-\varphi_z(x_2x_4)|, |\varphi_z(x_1x_4) -\varphi_z(x_2x_3)| \le 2\delta''.$$

For $i< j\in [4]$, let $b_{ij}$ denote a median of $\varphi_Z(x_ix_j)$. Since $\Prob{ |\varphi_Z(x_1x_2)-\varphi_Z(x_3x_4)| > 2\delta''} \le \sqrt{5\varepsilon}$ from above, we find that 
$$\Prob{\varphi_Z(x_3x_4) \ge b_{34}}\cdot \Prob{\varphi_Z(x_1x_2) \le b_{34}-2\delta''} = \Prob{\varphi_Z(x_3x_4) \ge b_{34} \text{ and } \varphi_Z(x_1x_2) \le b_{34}-2\delta''} \le \sqrt{5\varepsilon}.$$
Since $\Prob{\varphi_Z(x_3x_4) \ge b_{34}}\ge \frac{1}{2}$ by definition of median, we find that $\Prob{\varphi_Z(x_1x_2) \le b_{34}-2\delta''} \le 2\sqrt{5\varepsilon}.$ A symmetric argument implies that $\Prob{\varphi_Z(x_1x_2) \ge b_{34}+2\delta''} \le 2\sqrt{5\varepsilon}.$ Hence $\Prob{\varphi_Z(x_1x_2) \not\in [b_{34}-2\delta'',b_{34}+2\delta'']} \le 4\sqrt{5\varepsilon}$. By Proposition~\ref{prop:Expect} applied to $\varphi_Z(x_1x_2)$ this implies that
$$\Prob{\bigg|\varphi_Z(x_1x_2) - \Expect{\varphi_Z(x_1x_2)}\bigg| >4\sqrt{5\varepsilon} + 4\delta''} \le 4\sqrt{5\varepsilon}.$$
Since $4\sqrt{5\varepsilon}+4\delta'' \le \delta'$ as $\varepsilon$ is small enough, it follows that $p_{12}\le 4\sqrt{5\varepsilon}$ and hence by Claim~\ref{cl:EdgeProb} that $|B_{12}|\le 4\sqrt{5\varepsilon}\cdot n^2$.
Symmetric arguments applied to $\varphi_Z(x_ix_j)$ for $i<j\in [4]$ imply that $|B_{ij}|\le 4\sqrt{5\varepsilon}\cdot n^2$. Combining, we find that
$\sum_{i<j\in [4]} |B_{ij}| \le 24 \sqrt{5\varepsilon}\cdot n^2$.
Since $\sum_{i\in [4]} |B_i| \le 4\sqrt{\varepsilon}\cdot n^2$ by Claim~\ref{cl:InsideBadCount}, we find that $|B| \le 48\sqrt{5\varepsilon} \cdot n^2 + 4\sqrt{\varepsilon}\cdot n^2 \le \delta'\cdot n^2$
as desired, where the inequality follows since $48\sqrt{5\varepsilon} + 4\sqrt{\varepsilon} \le \frac{\varepsilon^{0.001}}{2}$ as $\varepsilon$ is small enough.
\end{proof}

\section{The Local Structure}

\subsection{Value Lemmas}

\begin{definition}
Let $(G,\phi)$ be a nonzero untied edge-weighted graph and let $f=uv$ be a negative edge of $G$. A vertex $x\in N(u)\cap N(v)$ is said to be \emph{obtuse} with respect to $f$ if $uvx$ is obtuse; furthermore we say $x$ is \emph{$u$-obtuse} with respect to $f$ if $x$ is obtuse and $ux$ is a positive edge of $G$; similarly we say $x$ is \emph{$v$-obtuse} with respect to $f$ if $x$ is obtuse and $vx$ is a positive edge of $G$. We say $x$ is \emph{acute} if $uvx$ is acute. For brevity we omit the with respect to $f$ if $f$ is clear from context.
\end{definition}

\begin{lem}[Obtuse-Acute Pair Value]\label{lem:ObtuseAcutePair}
Let $(G,\phi)$ be a nonzero untied triangle-covered edge-weighted graph and let $f=uv$ be a negative edge of $G$. Let $xy\in E(G[N(u)\cap N(v)])$. If
$x$ is $u$-obtuse and $y$ is acute and
    \begin{itemize} 
        \item[(i)] $xy$ is a positive edge, then $\textrm{Val}_f(xy) \le \dem(vy,f)$
        \item[(ii)] $xy$ is a negative edge, then $\textrm{Val}_f(xy)\le \min\left\{\frac{|\phi(xy)|}{\phi(ux)},~1\right\}\cdot |\phi(f)| + \dem(vy,f)$ 
    \end{itemize}
\end{lem}
\begin{proof}
First suppose $xy$ is a positive edge. Then $$\textrm{Val}_f(xy) = \demh_f(vy,vxy) = \demh_f(vy,vx) = \dem(vy,f) \cdot \demh(vy,vx) \le \dem(vy,f).$$ 
This proves (i).

So we assume $xy$ is a negative edge. Then 
$$\textrm{Val}_f(xy) = \demh_f(ux,uxy)+\demh_f(uy,uxy)+\demh_f(vy,vxy).$$
Since $xy$ is negative, we have that $vxy$ is obtuse. Hence
$$\demh(vy,xy)+ \demh(vy,vx) = \frac{\dem(vy,xy)}{\phi(vy)} + \frac{\dem(vy,vx)}{\phi(vy)} = \frac{|\phi(xy)|}{\phi(vy)} + \frac{|\phi(vx)|}{\phi(vy)} \le 1.$$
Thus
$$\demh_f(vy,vxy) = \demh_f(vy,xy) + \demh_f(vy,vx) = \dem(vy,f) \cdot (\demh(vy,xy) + \demh(vy,vx) ) \le \dem(vy,f).$$ 

Since $uxy$ is acute, we find that 
$$\demh_f(ux,uxy)+\demh_f(uy,uxy) = \demh_f(ux,xy)+\demh_f(uy,xy) = \dem(ux,f)\cdot \demh(ux,xy) + \dem(uy,f)\cdot \demh(uy,xy).$$
We claim the above is at most $\min\left\{\frac{|\phi(xy)|}{\phi(ux)},~1\right\}\cdot |\phi(f)|$ which will prove (ii) and hence conclude the proof. 

To see this, first suppose $uxy$ is sponsored. Further suppose $ux$ is the sponsor of $uxy$. Then $\dem(uy,xy)=0$, $\dem(ux,xy) = \frac{|\phi(xy)|}{\phi(ux)} \le 1$, $\dem(ux,f) = |\phi(f)|$ and the claim follows. So we suppose $uy$ is the sponsor of $uxy$ (and hence $|\phi(uy)|\ge |\phi(ux)|$ and $|\phi(uy)| \ge |\phi(xy)|$. Then $\dem(ux,xy)= 0$, $\dem(uy,xy) = \frac{|\phi(xy)|}{\phi(uy)} \le \min\left\{\frac{|\phi(xy)|}{\phi(ux)},~1\right\}$, $\dem(uy,f) \le |\phi(f)|$ and the claim follows. 

So we assume $uxy$ is not sponsored. Hence $|\phi(xy)| > \phi(ux),\phi(uy)$. Since $vxy$ is obtuse, we have that $\phi(vy)\ge |\phi(vx)|+|\phi(xy)|$. This implies that $\phi(vy) > \phi(uy)$ and also that $\phi(vy) > \phi(ux) \ge |\phi(f)|$. Thus we find that $vy$ is the sponsor of $uvy$. Hence $\dem(uy,f)=0$. Thus the above is at most $\dem(ux,f)\cdot \demh(ux,xy) \le |\phi(f)| \cdot \frac{|\phi(xy)| \cdot \frac{\phi(ux)}{\phi(ux)+\phi(uy)}}{\phi(ux)} = |\phi(f)|\cdot \frac{|\phi(xy)|}{\phi(ux)+\phi(uy)} \le \min\left\{\frac{|\phi(xy)|}{\phi(ux)},~1\right\}\cdot |\phi(f)|$ as claimed (since $|\phi(xy)| \le \phi(ux)+\phi(uy)$ as $\phi$ is triangle-covered).

\end{proof}

\begin{lem}[Extremal Pair Value]\label{lem:ExtremalPairValue}
Let $\varepsilon \in [0,\frac{1}{3}]$. Let $(G,\phi)$ be a nonzero untied  triangle-covered edge-weighted graph and let $f=uv$ be a negative edge of $G$. Let $xy\in E(G[N(u)\cap N(v)])$. If $\textrm{Val}_f(xy) \ge (2-\varepsilon)\cdot |\phi(f)|$, then either
\begin{enumerate}
\item[(a)] $x$ is $u$-obtuse and $y$ is $v$-obtuse and $|\phi(vx)|, |\phi(uy)| \ge 2\cdot  |\phi(f)|$, or
\item[(b)] $x$ is $v$-obtuse and $y$ is $u$-obtuse and $|\phi(ux)|, |\phi(vy)| \ge 2\cdot |\phi(f)|$, or
\item[(c)] $xy$ is a negative edge and $|\phi(xy)| \ge (1-3\varepsilon)\cdot |\phi(f)|$
\end{enumerate}
\end{lem}
\begin{proof}
Let $s:= \frac{\textrm{Val}_f(xy)}{|\phi(f)|}$. By assumption $s\ge 2-\varepsilon$. We case according to whether $x$ or $y$ are acute or obtuse as follows.
\vskip.1in

\noindent \textbf{Case 1: One of $x,y$ is obtuse and one of $x,y$ is acute.}
\vskip.05in
Without loss of generality, we assume $x$ is obtuse and $y$ is acute; similarly without loss of generality we assume $x$ is $u$-obtuse. First suppose that $xy$ is a positive edge. Then by Lemma~\ref{lem:ObtuseAcutePair}(i), we have that $s\le \frac{\dem(vy,f)}{|\phi(f)|} \le 1$, a contradiction. So we assume that $xy$ is a negative edge. Then by Lemma~\ref{lem:ObtuseAcutePair}(ii), we find that 
$$s\le \min\left\{\frac{|\phi(xy)|}{\phi(ux)},~1\right\} + \frac{\dem(vy,f)}{|\phi(f)|} \le \min\left\{\frac{|\phi(xy)|}{\phi(ux)},~1\right\} + 1.$$
Since $s\ge 2-\varepsilon$, we find that 
$$|\phi(xy)| \ge (1-\varepsilon)\cdot \phi(ux) \ge (1-\varepsilon)\cdot |\phi(f)|,$$
and hence (c) holds as desired where the final inequality follows since $uvx$ is obtuse and $\phi$ is triangle-covered.
\vskip.1in

\noindent \textbf{Case 2: Both of $x$ and $y$ are acute.}
\vskip.05in
If $xy$ is a positive edge, then $\textrm{Val}_f(xy)=0$, a contradiction. So we assume $xy$ is a negative edge. Note that 

$$\textrm{Val}_f(xy) = \demh_f(ux,xy) + \demh_f(uy,xy)+\demh(vx,vy)+\demh_f(vy,xy).$$

\begin{claim}\label{cl:AllNonzero}
$\demh_f(ux,xy),\demh_f(uy,xy),\demh(vx,vy),\demh_f(vy,xy) > 0$.    
\end{claim}
\begin{proofclaim}
Suppose not. We assume without loss of generality that $\demh_f(vx,xy)=0$. Now we calculate that
\begin{align*}
\textrm{Val}_f(xy) &= \demh_f(ux,xy) + \demh_f(uy,xy)+\demh_f(vy,xy) \\
&= \dem(ux,f)\cdot \demh(ux,xy) + \dem(uy,f)\cdot \demh(uy,xy) + \dem(vy,f)\cdot \demh(vy,xy)  \\
&\le \phi(ux)\cdot \demh(ux,xy) + \dem(uy,f)+ \dem(vy,f)\\
&\le \phi(ux)\cdot \frac{|\phi(xy)|}{\phi(ux)} + |\phi(f)|\\
&= |\phi(xy)| + |\phi(f)|
\end{align*}
where we used that $\dem(ux,f)\le \phi(ux)$ and that $\demh(uy,xy)\le 1$ and that $\demh(vy,xy)\le 1$ and that $\dem(uv,f)+\dem(vy,f)\le |\phi(f)|$ and that $\demh(ux,xy)\le \frac{|\phi(xy)|}{\phi(ux)}$ (since $\dem(ux,xy) \le |\phi(xy)|$). Since $s\ge 2-\varepsilon$, it follows
that
$$|\phi(xy)|\ge (1-\varepsilon)\cdot |\phi(f)|,$$
and (c) holds as desired. 
\end{proofclaim}

\begin{claim}
None of $uxy, vxy, uvx, uvy$ are sponsored.    
\end{claim}
\begin{proofclaim}
Suppose not. First suppose that at least one of $uvx$ or $uvy$ is sponsored. We assume without loss of generality that $uvx$ is sponsored and that $ux$ is its sponsor. Thus we find that $\phi(ux)\ge |\phi(f)|$. But we also have that $\dem(vx,f)=0$ and hence $\demh_f(vx,xy) = \dem(vx,f)\cdot \demh(vx,vy) = 0$, contradicting Claim~\ref{cl:AllNonzero}.

So we assume that at least one of $uxy$ or $vxy$ is sponsored. We assume without loss of generality that $vxy$ is sponsored and that $vy$ is its sponsor. But then we have that $\dem(vx,xy)=0$ and hence $\demh_f(vx,xy) = \dem(vx,f)\cdot \demh(vx,vy) = 0$, contradicting Claim~\ref{cl:AllNonzero}. 
\end{proofclaim}

Now we calculate that
\begin{align*}
\textrm{Val}_f(xy) &= \demh_f(ux,xy) + \demh_f(vx,xy)+ \demh_f(uy,xy)+\demh_f(vy,xy) \\
&= \dem(ux,f)\cdot \demh(ux,xy) + \dem(vx,f)\cdot \demh(vx,xy) +\dem(uy,f)\cdot \demh(uy,xy) \\
&~~~~+ \dem(vy,f)\cdot \demh(vy,xy)  \\
&\le \phi(ux)\cdot \demh(ux,xy) + \phi(vx)\cdot \demh(vx,xy) +\phi(uy)\cdot \demh(uy,xy) \\
&~~~~+ \phi(vy)\cdot \demh(vy,xy)  \\
&= \phi(ux)\cdot \frac{|\phi(xy)|}{\phi(ux)+\phi(uy)} + \phi(vx)\cdot \frac{|\phi(xy)|}{\phi(xx)+\phi(vy)} +\phi(uy)\cdot \frac{|\phi(xy)|}{\phi(ux)+\phi(uy)} \\
&~~~~+ \phi(vy)\cdot \frac{|\phi(xy)|}{\phi(vx)+\phi(vy)}  \\
&= 2\cdot |\phi(xy)|.
\end{align*}
Since $s\ge 2-\varepsilon$, it follows that $|\phi(xy)| \ge \frac{1}{2}\cdot (2-\varepsilon) |\phi(f)|$ and (c) holds as desired.

\vskip.1in

\noindent \textbf{Case 3: Both of $x$ and $y$ are obtuse.}
\vskip.05in
We assume without loss of generality that $x$ is $u$-obtuse (since outcomes (a) and (b) are symmetric). If $y$ is $u$-obtuse, then $\textrm{Val}_f(xy)=0$, a contradiction. So we assume that $y$ is $v$-obtuse. First suppose that $xy$ is a negative edge. Then we find that 
\begin{align*}
\textrm{Val}_f(xy) &= \demh_f(ux,xy) + \demh_f(ux,uy)+ \demh_f(vy,xy)+\demh_f(vy,vx) \\
&= \dem(ux,f)\cdot \demh(ux,xy) + \dem(ux,f)\cdot \demh(ux,uy) +\dem(vy,f)\cdot \demh(vy,xy) \\
&~~~~+ \dem(vy,f)\cdot \demh(vy,vx)  \\
&= |\phi(f)|\cdot \left(\frac{|\phi(xy)|+|\phi(uy)|}{\phi(ux)}\right)+ |\phi(f)|\cdot \left(\frac{|\phi(xy)|+|\phi(vx)|}{\phi(vy)}\right).\\
\end{align*}
Since $\phi$ is triangle-covered, we have that $|\phi(xy)|+|\phi(uy)| \le \phi(ux)$ and similarly $|\phi(xy)|+|\phi(vx)|\le \phi(vy)$. Using this and that $s\ge 2-\varepsilon$, it follows that $|\phi(xy)|+|\phi(uy)| \ge (1-\varepsilon)\cdot \phi(ux)$ and similarly $|\phi(xy)|+|\phi(vx)|\ge (1-\varepsilon) \cdot \phi(vy)$. Note that $\phi(ux)\ge |\phi(f)|+|\phi(vx)|$ and $\phi(vy)\ge |\phi(f)|+|\phi(uy)|$ since $\phi$ is triangle-covered. Now if (c) holds, the lemma follows. So we assume (c) does not hold and hence $|\phi(xy)| < (1-3\varepsilon)\cdot |\phi(f)|$. But then we find that
$$(1-3\varepsilon)\cdot |\phi(f)| + |\phi(uy)| \ge (1-\varepsilon)\cdot \phi(ux) \ge (1-\varepsilon) \cdot (|\phi(f)|+|\phi(vx)|).$$
Rearranging this yields that 
$$|\phi(uy)|\ge 2\varepsilon\cdot |\phi(f)| + (1-\varepsilon)\cdot |\phi(vx)|$$
and symmetrically by a similar logic we have that
$$|\phi(vx)|\ge 2\varepsilon\cdot |\phi(f)| + (1-\varepsilon)\cdot |\phi(uy)|.$$
Combining the two inequalities yields that 
$$|\phi(uy)|\ge 2\varepsilon\cdot |\phi(f)| + (1-\varepsilon)\cdot \bigg(2\varepsilon\cdot |\phi(f)| + (1-\varepsilon)\cdot |\phi(uy)|\bigg) = (4\varepsilon-2\varepsilon^2)\cdot |\phi(f)| +(1-2\varepsilon+\varepsilon^2)\cdot |\phi(uy)|.$$
Rearranging once more yields that
$$|\phi(uy)| \ge \frac{4\varepsilon-2\varepsilon^2}{2\varepsilon-\varepsilon^2} \cdot |\phi(f)| = 2\cdot |\phi(f)|$$
and symmetrically we find that $|\phi(vx)|\ge 2\cdot |\phi(f)|$ and hence (a) holds as desired. 

So we finally assume that $xy$ is a positive edge. Then we find that 
\begin{align*}
\textrm{Val}_f(xy) &= \demh_f(ux,uy)+\demh_f(vy,vx) \\
&= \dem(ux,f)\cdot \demh(ux,uy) + \dem(vy,f)\cdot \demh(vy,vx)\\
&= |\phi(f)|\cdot \left(\frac{|\phi(uy)|}{\phi(ux)}\right)+ |\phi(f)|\cdot \left(\frac{|\phi(vx)|}{\phi(vy)}\right).
\end{align*}
Since $s\ge 2-\varepsilon$, it follows that $|\phi(uy)| \ge (1-\varepsilon)\cdot \phi(ux)$ and similarly $|\phi(vx)|\ge (1-\varepsilon) \cdot \phi(vy)$.
Note that $\phi(ux)\ge |\phi(f)|+|\phi(vx)|$ and $\phi(vy)\ge |\phi(f)|+|\phi(uy)|$ since $\phi$ is triangle-covered. Using this, we find that
$$|\phi(uy)| \ge (1-\varepsilon)\cdot \phi(ux) \ge (1-\varepsilon) \cdot (|\phi(f)|+|\phi(vx)|).$$
and symmetrically by a similar logic we have that
$$|\phi(vx)|\ge (1-\varepsilon)\cdot (|\phi(f)| + |\phi(uy)|).$$
Combining the two inequalities yields that 
$$|\phi(uy)|\ge (1-\varepsilon)\cdot |\phi(f)| + (1-\varepsilon)^2\cdot \bigg(|\phi(f)| + |\phi(uy)|)\bigg) = (2-3\varepsilon+\varepsilon^2)\cdot |\phi(f)| +(1-2\varepsilon+\varepsilon^2)\cdot |\phi(uy)|.$$
Rearranging once more yields that
$$|\phi(uy)| \ge \frac{2-3\varepsilon+\varepsilon^2}{2\varepsilon-\varepsilon^2} \cdot |\phi(f)| \ge 2\cdot |\phi(f)|$$
since $2-3\varepsilon+\varepsilon^2 \ge 2(2\varepsilon-\varepsilon^2)$ as $\varepsilon \le \frac{1}{3}$. Symmetrically we find that $|\phi(vx)|\ge 2\cdot \phi(f)$ and hence (a) holds as desired. 
\end{proof}

\subsection{Extremal Neighborhood Structure Lemma}

\begin{definition}
Let $(G,\phi)$ be a nonzero untied triangle-covered edge-weighted graph on $n$ vertices. A negative edge $f=u_1u_2$ of $G$ is 
\begin{itemize}
    \item \emph{$\alpha$-locally-small} if for both $i\in \{1,2\}$, we have that $\bigg|\bigg\{f'=u_iw\in E_{\phi}^-(G): |\phi(f')| \ge 2\cdot| \phi(f)|~\bigg\}\bigg| \ge \alpha\cdot n$,
    \item \emph{$\alpha$-bipartite} if there exists a bipartition $\mc{P}=(A,B)$ of $N(u_1)\cap N(u_2)$ such that at least $(1-\alpha)\cdot \frac{|N(u_1)\cap N(u_2)|^2}{4}$ of the pairs $(a,b)$ with $a\in A$ and $b\in B$ satisfy that $ab\in E_{\phi}^-(G)$ and $|\phi(ab)| \ge (1-\alpha)\cdot |\phi(u_1u_2)|$. 
 \end{itemize}    
\end{definition}

\begin{lem}[Extremal Neighborhood Structure]\label{lem:NeighborhoodStructure}
There exists $\varepsilon_0$ such that the following holds for each $\varepsilon\in (0,\varepsilon_0]$: Let $\alpha \in (0,1)$. Let $(G,\phi)$ be a nonzero untied triangle-covered edge-weighted graph on $n$ vertices with $\delta(G)\ge (\frac{3}{4}-\varepsilon)n$ and let $f=uv$ be a negative edge of $G$. If $\sum_{g\in E(G[N(u)\cap N(v)])} \textrm{Val}_f(g) \ge (1-\varepsilon) \frac{|N(u)\cap N(v)|^2}{2}\cdot |\phi(f)|$, then $f$ is
\begin{enumerate}
    \item[(i)] $\alpha$-locally-small, or
    \item[(ii)] $(36\alpha+4900\cdot \varepsilon^{0.0005})$-bipartite.
\end{enumerate}
\end{lem}
\begin{proof}
We choose $\varepsilon_0$ small enough as needed to satisfy various inequalities throughout the proof. Suppose not. We assume without loss of generality that $|\phi(f)|=1$. Note then that $H:=(\textrm{Filter}_{G,\phi}(f),~\frac{1}{2}\cdot\textrm{Val}_f)$ is a $[0,1]$-edge-weighted graph; by Lemma~\ref{lem:EdgeWeightTriangleThin}, $H$ is triangle-thin and by Proposition~\ref{prop:K5freeFilter}, $H$ is $K_5$-free. Let $\beta:= \varepsilon^{0.001}$. Since $\sum_{g\in E(G[N(u)\cap N(v)])} \textrm{Val}_f(g) \ge (1-\varepsilon) \frac{|N(u)\cap N(v)|^2}{2}\cdot |\phi(f)|$, we find by Lemma~\ref{lem:ValueStability} that $(\textrm{Filter}_{G,\phi}(f),~\frac{1}{2}\cdot \textrm{Val}_f)$ is $\beta$-close to some value-extremal step graphon $W=(K_4,\mu,\varphi)$.

Let $S_u$ be the set of $u$-obtuse vertices $x$ in $V(H)$ with $|\phi(vx)|\ge 2\cdot |\phi(f)|$; similarly let $S_v$ be the set of $v$-obtuse vertices $v$ in $V(H)$ with $|\phi(uy)|\ge 2\cdot |\phi(f)|$. Let $w\in \{u,v\}$ such that $|S_w|=\min\{|S_u|,~|S_v|\}$ and let $S:= S_w$. 

\begin{claim}\label{cl:LocallySmall}
$|S|\le \alpha\cdot n$.     
\end{claim}
\begin{proofclaim}
Suppose not. But then $f$ is $\alpha$-locally-small and (i) holds, a contradiction.
\end{proofclaim}

\begin{claim}\label{cl:NotStarCycle}
Let $\delta:=61\sqrt{\beta}$. Then $H$ is not $\delta$-close to a star-cut or cycle-cut step graphon.    
\end{claim}
\begin{proofclaim}
Suppose not. That is, $H$ is $\delta$-close to a cycle-cut or star-cut step graphon $W'$. Then $H$ has a bipartition $(A,B)$  where $(\frac{1}{2}-3\delta)\cdot v(H) \le |A|,|B| \le (\frac{1}{2}+3\delta)\cdot v(H)$ and all but at most $\delta n^2$ of the pairs $a\in A, b\in B$ satisfy $\frac{1}{2} \cdot \textrm{Val}_f(ab) \ge 1-\delta$. By Lemma~\ref{lem:ExtremalPairValue} (with $\varepsilon$ of that lemma set to $2\delta$), each such pair satisfies either Lemma~\ref{lem:ExtremalPairValue}(a), (b), or (c).  Thus the pairs that satisfy (a) or (b) necessarily contain exactly one vertex from $S_u$ and one vertex from $S_v$. Since $|S| \le \alpha\cdot n$ by Claim~\ref{cl:LocallySmall}, we find that there are at most $\alpha\cdot n^2$ pairs satisfying (a) or (b). Hence the number of pairs satisfying (c), that is $ab\in E(G)$ and $|\phi(ab)|\ge (1-6\delta)\cdot |\phi(f)|$, is at least 
$$\left((\frac{1}{2}-3\delta\right)\cdot \left(\frac{1}{2}+3\delta\right)\cdot v(H)^2 - \delta\cdot n^2 - \alpha\cdot n^2 \ge \frac{v(H)^2}{4} (1-36\delta -36(\delta+\alpha)) ,$$  where we used that $|N(u)\cap N(v)|\ge n/3$ since $\varepsilon\le \frac{1}{12}$. It follows that $f$ is $(36\alpha+72\delta)$-bipartite and (ii) holds (since $4900\sqrt{\beta} \ge 4392\sqrt{\beta} =72\delta$), a contradiction.
\end{proofclaim}

Thus by Claim~\ref{cl:NotStarCycle} we assume $H$ is not $\delta$-close to a cycle-cut or star-cut step graphon $W'$. It follows that $W$ is reflective-cut or quad-cut. 

Let $x_1,x_2,x_3,x_4$ be an enumeration of $V(K_4)$ such that $\mu(x_1)\ge\mu(x_2)\ge \mu(x_3)\ge \mu(x_4)$. 

\begin{claim}
Let $b_4:= \frac{\delta-\beta}{2}$ and $w_0:= \delta-\beta$. For all $i\in [4]$, $\mu(x_i)\ge b_4$. For all $i<j\in [4]$, $\varphi(x_ix_j)\ge w_0$.
\end{claim}
\begin{proofclaim}
Suppose not. First suppose there exists $i\in [4]$ with $\mu(x_i)\le b_4$. But then in this case we find that $W$ is reflective-cut and $W$ is $(\delta-\beta)$-close to a a cycle-cut step graphon $W'$ and so $H$ is $\delta$-close to $W'$, contradicting Claim~\ref{cl:NotStarCycle}. 

So we assume there exist $i<j\in [4]$ such that $\varphi(x_ix_j) \le \delta-\beta$. But then $W$ is $(\delta-\beta)$-close to a cycle-cut step graphon $W'$ and hence $H$ is $\delta$-close to $W'$, contradicting Claim~\ref{cl:NotStarCycle}. 
\end{proofclaim}

Let $c$ be the $4$-coloring of $V(H)$ as guaranteed in the definition of $\beta$-close and let $V_1,V_2,V_3,V_4$ be the color classes (corresponding to the above enumeration $x_1,\ldots,x_4$). 

\begin{claim}\label{cl:NotTooManyParts}
All of the following hold:
\begin{itemize}
    \item[(i)] there exists at most one $i\in [4]$ such that $V_i$ has more than $\sqrt{\beta}\cdot n$ $u$-obtuse vertices,
    \item[(ii)] there exists at most one $j\in [4]$ such that $V_j$ has more than $\sqrt{\beta}\cdot n$ $v$-obtuse vertices,
    \item[(iii)] there exist at most two $k\in [4]$ such that $V_k$ has more than $2\sqrt{\beta}\cdot n$ acute vertices.
\end{itemize}
\end{claim}
\begin{proofclaim}
First suppose there exist $i\ne j\in [4]$ such that $V_i,V_j$ each have more than $\sqrt{\beta}\cdot n$ $u$-obtuse vertices. Let $O_i$ be the set of $u$-obtuse vertices in $V_i$ and let $O_j$ be the set of $u$-obtuse vertices in $V_j$. Since $\phi$ is triangle covered we have for each $x\in O_i,~y\in O_j$ that either $xy\not\in E(G)$ or $xy$ is a positive edge (since $xv,yv$ are negative edges as $x$ and $y$ are both $u$-obtuse) and hence in either case we find that $\textrm{Val}_f(xy)=0$. But then each pair satisfies that $|\frac{1}{2}\textrm{Val}_f(xy) - \varphi(x_ix_j)| > \beta$ (since $\varphi(x_ix_j) \ge w_0 > \beta$); yet there are at least $|O_i|\cdot |O_j| > \beta\cdot  n^2$, contradicting that $H$ is $\beta$-close to $W$. This proves (i). A symmetric argument proves (ii).

Finally suppose there exist three distinct indices $i_1,i_2,i_3\in [4]$ such that $V_{i_1},V_{i_2},V_{i_3}$ each have more than $2\sqrt{\beta}\cdot n$ acute vertices. For each $j\in [3]$, let $A_{i_j}$ be a subset of $V_{i_j}$ of size $2\sqrt{\beta}\cdot n$ consisting of acute vertices. Let $J$ be the graph with $V(J):= \bigcup_{j\in [3]} V_{i_j}$ and $E(J):= \bigcup_{j<j'\in [3]} \{ xy: x\in V_{i_j}, y\in V_{i_j'}, xy\in E^-_{\phi}(G)\}$. Since $\phi$ is triangle-covered, we find that $J$ is triangle-free. Hence by Mantel's Theorem, $e(J) \le \frac{1}{4}\cdot v(J)^2 = \frac{1}{4}\cdot (6\sqrt{\beta}\cdot n)^2 = 9\beta\cdot n^2$. but then there exist at least $3\cdot (2\sqrt{\beta}\cdot n^2) - 9\beta\cdot n^2 = 3\beta\cdot n^2$ $x,y\in V(J)$ with $xy\not \in E(J)$ (that is $xy\not\in E(G)$ or $xy\in E^+_{\phi}(G)$) and $x\in V_{i_j}, y\in V_{i_j'}$ for some $j<j'\in [3]$. But for each such pair we find that that $\textrm{Val}_f(xy)=0$ since $x$ and $y$ are acute and $xy\not\in E(G)$ or $xy$ is positive. But then each pair satisfies that $|\frac{1}{2}\textrm{Val}_f(xy) - \varphi(x_ix_j)| > \beta$ (since $\varphi(x_ix_j) \ge w_0 > \beta$); yet there are at least $|O_i|\cdot |O_j| > \beta\cdot  n^2$, contradicting that $H$ is $\beta$-close to $W$. This proves (iii). 
\end{proofclaim}

\begin{claim}\label{cl:Perm}
There exists a permutation $\sigma:[4]\rightarrow[4]$ such that all of the following hold:
\begin{itemize}
    \item[(i)] all but at most $3\sqrt{\beta}\cdot n$ vertices of $V_{\sigma(1)}$ are $u$-obtuse,
    \item[(ii)] all but at most $3\sqrt{\beta}\cdot n$ vertices of $V_{\sigma(2)}$ are $v$-obtuse,
    \item[(iii)] for $i\in \{3,4\}$, all but at most $2\sqrt{\beta}\cdot n$ of each of $V_{\sigma(i)}$ are acute.
\end{itemize}
\end{claim}
\begin{proofclaim}
For each $i\in [4]$, we have that  
$$|V_i| \ge \mu(x_i)\cdot v(H) - \beta \cdot n \ge \frac{\delta-\beta}{2}\cdot \left(\frac{1}{2}-2\varepsilon\right)n  -\beta \cdot n \ge \frac{\delta-\beta}{2}\cdot \frac{n}{3} -\beta\cdot n> 4\sqrt{\beta}\cdot n,$$
where we used that $\mu(x_i) \ge \frac{\delta-\beta}{2}$, that $v(H)=|N(u)\cap N(v)|\ge (\frac{1}{2}-2\varepsilon)n \ge \frac{n}{3}$ that $H$ is $\beta$-close to $W$, and that $\delta \ge 7\beta+24\sqrt{\beta}$. By the pigeonhole principle it follows that for each $i\in [4]$, $V_i$ satisfies at least one of the following: (a) $V_i$ has more than $\sqrt{\beta}\cdot n$ $u$-obtuse vertices; (b) $V_i$ has more than $\sqrt{\beta}\cdot n$ $v$-obtuse vertices; (c) $V_i$ has more than $2\sqrt{\beta}\cdot n$ acute vertices. 

However by Claim~\ref{cl:NotTooManyParts}, at most one $i\in [4]$ satisfies (a), at most one $j\in [4]$ satisfies (b), and at most two indices $k,\ell\in [4]$ satisfy (c). Combining these facts, we find that all of these indices exist and are distinct as desired.
\end{proofclaim}

We say a pair $y\ne z\in V(H)$ is \emph{abnormal} if $|\frac{1}{2}\textrm{Val}_f(yz) -\varphi'(x_{c(y)}x_{c(z)})| > \beta$; we say the pair is \emph{normal} if it is not abnormal.  Since $H$ is $\beta$-close to $W$, we have by definition that there are at most $\beta\cdot v(H)^2$ ($\le \beta\cdot n^2$) pairs that are abnormal. We say a vertex $y\in V(H)$ is \emph{abnormal} if there exist at least $\sqrt{\beta}\cdot n$ vertices $z\in V(H)$ such that $yz$ is abnormal; we say $y$ is \emph{normal} if it is not abnormal. It follows that there are at most $\sqrt{\beta}\cdot n$ abnormal vertices of $V(H)$. Let $\mc{N}$ denote the set of normal vertices of $H$.

Let $\sigma$ be as in Claim~\ref{cl:Perm}. Let $V_{\sigma(1)}'$ be the set of vertices in $V_{\sigma(1)}\cap \mc{N}$ that are $u$-obtuse. Let $V_{\sigma(2)}'\cap \mc{N}$ be the set of vertices in $V_{\sigma(2)}$ that are $v$-obtuse. For $i\in \{3,4\}$, let $V_{\sigma(i)}'$ be the set of vertices in $V_{\sigma(i)}\cap \mc{N}$ that are acute. By Claims~\ref{cl:Perm} and~\ref{cl:LocallySmall} and that fact that $|V(H)\setminus \mc{N}| \le \sqrt{\beta}\cdot n$, we have for each $i\in [4]$ that 
$$|V_{\sigma(i)}'| \ge |V_{\sigma(i)}|-4\sqrt{\beta}\cdot n \ge b_4\cdot \frac{n}{3} - \beta\cdot n - 4\sqrt{\beta}\cdot n \ge 4\sqrt{\beta}\cdot n,$$
where the final inequality follows as $\delta\ge 13\beta + 48\sqrt{\beta}$. 

Let $J$ be the graph with $V(J):= \bigcup_{i\in [4]} V_{\sigma(i)}'$ and $E(J):= \bigcup_{i<j\in [4]} \{yy': y\in V_{\sigma(i)},~y'\in V_{\sigma(j)}',~yy' \textrm{ is normal}\}$. Thus $v(J)\ge v(H) - 20\sqrt{\beta}\cdot n$. Moreover, by the definition of normal, we have that for each $i\ne j\in [4]$ and $y\in V_{\sigma(i)}'$ that $|N_J(y)\cap V_{\sigma(j)}'|\ge |V_{\sigma(j)'}|-\sqrt{\beta}\cdot n$.

For $i<j\in [4]$, define $m_{ij}:= \varphi(x_{\sigma(i)}x_{\sigma(j)})$. Recall that $m_{ij} \in [w_0,1]$. Since $W$ is matched as $W$ is reflective-cut or quad-cut, we have that $m_{12}=m_{34}$, $m_{13}=m_{24}$, $m_{14}=m_{23}$; similarly we have that $m_{12}+m_{13}+m_{23}=2$.

\begin{claim}\label{cl:MatchingSumSmall}
$m_{13} +m_{23}\le 1+2\beta$ and hence $m_{12}\ge 1-2\beta$.
\end{claim}
\begin{proofclaim}
Since $|V_{\sigma(3)}'| > 0$, there exists a vertex $y\in V_{\sigma(3)}'$ (that is $y$ is in $V_{\sigma(3)}$, is normal and acute and not in $S$). Since $|V_{\sigma(4)}'| > \sqrt{\beta}\cdot n$ from above, there exists a vertex $y'\in V_{\sigma(4)}'$ such that $yy'$ is normal. Since $|V_{\sigma(1)}'| > 2\sqrt{\beta}\cdot n$, there exists a vertex $z\in V_{\sigma(1)}'$ such that $yz,y'z$ are both normal.  Similarly, there exists a vertex $z'\in V_{\sigma(2)}'$ such that $yz',y'z'$ are both normal. 

Since $yy'$ is normal, we find that $\frac{1}{2}\textrm{Val}_f(yy') \ge \varphi(x_{\sigma(3)}x_{\sigma(4)}) -\beta \ge w_0-\beta > 0$. Since $y$ and $y'$ are both acute, it then follows that $yy'$ is a negative edge. Since $\phi$ is triangle-covered, it follows that at most one of $yz,y'z$ are negative. We thus assume without loss of generality that $yz$ is a positive edge. Since $y$ is acute and $z$ is $u$-obtuse, we have by Lemma~\ref{lem:ObtuseAcutePair}(i) that $\textrm{Val}_f(yz) \le \dem(vy,f)$. Since $y$ is acute and $z'$ is $v$-obtuse by Lemma~\ref{lem:ObtuseAcutePair}(i) and (ii), we also find that $\textrm{Val}_f(yz')\le \min\left\{\frac{|\phi(yz')|}{\phi(vz')},~1\right\}\cdot |\phi(f)| + \dem(uy,f)$. Since $\dem(vy,f)+\dem(uv,y)\le |\phi(f)|\le 1$, we find that 
$$\textrm{Val}_f(yz) + \textrm{Val}_f(yz') \le 1\cdot |\phi(f)| + \dem(uy,f)+\dem(vy,f) \le 2\cdot |\phi(f)| = 2.$$
Yet since $yz$ is normal we have that $|\frac{1}{2}\textrm{Val}_f(yz) - m_{13}| \le \beta$; similarly since $yz'$ is normal, we have that $|\frac{1}{2}\textrm{Val}_f(yz') - m_{23}| \le \beta$. Substituting these yields that 
$$m_{13} +m_{23} \le \frac{1}{2}\textrm{Val}_f(yz) + \beta + \frac{1}{2}\textrm{Val}_f(yz') + \beta \le 1 + 2\beta,$$
as desired. 
\end{proofclaim}

Let $\mc{E}:= \{yz\in E(J)\cap E^-_{\phi}(G): |\phi(yz)| \ge (1-(36\alpha+4900\sqrt{\beta}))\cdot |\phi(f)|\}$.

\begin{claim}\label{cl:NormalK4}
If $w_1,w_2,w_3,w_4\in V(J)$ such that $w_i\in V_{\sigma(i)}'$ for each $i\in [4]$ and $J[\{w_1,w_2,w_3,w_4\}]\cong K_4$, then $w_3w_4\in \mc{E}$; if $w_1,w_2\not\in S$, then $w_1w_2\in \mc{E}$; and exactly one of the following holds:
\begin{itemize}
    \item[(i)] $m_{13} > \frac{1}{2}+\beta$ and $w_1w_3,w_2w_4 \in \mc{E}$,
    \item[(ii)] $m_{23} > \frac{1}{2}+\beta$ and $w_2w_3,w_1w_4 \in \mc{E}$,
    \item[(iii)] $\frac{1}{2}-\beta \le m_{13},m_{23}\le \frac{1}{2}+\beta$, $w_1w_3,w_2w_4 \in \mc{E}$ and $\dem(uw_3,f),\dem(vw_4,f)\ge 1-6\beta$, or
    \item[(iv)] $\frac{1}{2}-\beta \le m_{13},m_{23}\le \frac{1}{2}+\beta$, $w_2w_3,w_1w_4 \in \mc{E}$ and $\dem(uw_4,f),\dem(vw_3,f)\ge 1-6\beta$.   
\end{itemize}
\end{claim}
\begin{proofclaim}
Note that for each $i<j\in [4]$, we have that $w_iw_j\in E(G)$ where this follows since $w_iw_j$ is normal and hence $\textrm{Val}_f(w_iw_j) \ge 2(m_{ij}-\beta)\ge 2(w_0-\beta) > 0$.

Since $m_{34}=m_{12}\ge 1-2\beta$ by Claim~\ref{cl:MatchingSumSmall}, it follows that $\textrm{Val}_f(w_1w_2), \textrm{Val}_f(w_3w_4) \ge 2-6\beta$. Suppose that $w_1,w_2\not\in S$. Then since $w_1$ is $u$-obtuse and $w_2$ is $v$-obtuse, it follows that $\min\{\phi(uw_1),\phi(vw_2)\} < 2\cdot |\phi(f)|$. Hence Lemma~\ref{lem:ExtremalPairValue}(a) or (b) do not hold for $w_1w_2$ and thus Lemma~\ref{lem:ExtremalPairValue}(c) holds for $w_1w_2$, that is $w_1w_2$ is a negative edge and $|\phi(w_1w_2)| \ge (1-18\beta)\cdot |\phi(f)|$. Thus $w_1w_2\in \mc{E}$.

Since $w_3$ and $w_4$ are acute, we similarly find that Lemma~\ref{lem:ExtremalPairValue}(c) holds for $w_3w_4$, that is $w_3w_4$ is a negative edge and $|\phi(w_3w_4)|\ge (1-18\beta)\cdot |\phi(f)|$. Thus $w_3w_4\in \mc{E}$.

Now we first suppose that $m_{13} > \frac{1}{2}+\beta$ thus we find that $\textrm{Val}_f(w_1w_3) \ge 2(m_{13}-\beta) > 1 = |\phi(f)| \ge \dem(vw_3,f)$. So Lemma~\ref{lem:ObtuseAcutePair}(i) does not hold for $w_1w_3$. Hence Lemma~\ref{lem:ObtuseAcutePair}(ii) holds for $w_1w_3$, and thus $w_1w_3$ is a negative edge. By symmetry, we also find that $w_2w_4$ is a negative edge. Furthermore, we also have that $w_3w_4$ is a negative edge since $\textrm{Val}_f(w_3w_4) \ge 2(w_0-\beta) > 0$. Since $\phi$ is triangle-covered, it follows that $w_1w_4$ and $w_2w_3$ are positive edges. By Lemma~\ref{lem:ObtuseAcutePair}(i) since $w_1$ is $u$-obtuse and $w_4$ is acute, we find that $\dem(vw_4,f) \ge \textrm{Val}_f(w_1w_4) \ge 2(m_{14}-\beta) = 2(m_{23}-\beta)$. But then $\dem(uw_4,f) \le 1 - \dem(vw_4,f) = 1-2(m_{23}-\beta)$. By Lemma~\ref{lem:ObtuseAcutePair}(ii) since $w_2$ is $v$-obtuse and $w_4$ is acute and $\phi(vw_2)\ge|\phi(f)|=1$, we find that
$$|\phi(w_2w_4)| \ge \textrm{Val}_f(w_2w_4)-\dem(uw_4,f) \ge 2(m_{24}-\beta) - (1 - 2(m_{23}-\beta)) \ge 2(m_{23}+m_{24})-4\beta -1\ge 1-4\beta,$$
where we used that $m_{23}+m_{24}\ge 2-m_{12} \ge 1$ as $m_{12}\le 1$. It follows that $w_2w_4\in \mc{E}$. 

Similarly by Lemma~\ref{lem:ObtuseAcutePair}(i) since $w_2$ is $v$-obtuse and $w_3$ is acute, we find that $\dem(uw_3,f) \ge \textrm{Val}_f(w_2w_3) \ge 2(m_{23}-\beta)$ and hence $\dem(vw_3,f) \le 1-\dem(uw_3,f) \le 1-2(m_{23}-\beta)$. By Lemma~\ref{lem:ObtuseAcutePair}(ii) since $w_1$ is $u$-obtuse and $w_3$ is acute and $\phi(uw_1)\ge|\phi(f)|=1$, we similarly find that $|\phi(w_1w_3)|\ge 1-4\beta$ and $w_1w_3\in \mc{E}$. Thus (i) holds as desired.

Next suppose that $m_{23} > \frac{1}{2}+\beta$. A symmetric argument to the one given above implies that (ii) holds as desired.

So we assume that $m_{13},m_{23} \le \frac{1}{2}+\beta$. Since $m_{13}+m_{23}\ge 2-m_{12}\ge 1$, we also find that $m_{13},m_{23} \ge \frac{1}{2}-\beta$.  Thus we find for each $i\in\{1,2\}$ and $j\in \{3,4\}$ that $\textrm{Val}_f(w_iw_j) \in [2m_{ij}-2\beta,2m_{ij}+2\beta] \subseteq [1-3\beta,1+3\beta]$. Since $\phi$ is triangle-covered, we have that at least one of $w_1w_3$ or $w_1w_4$ is positive (since $w_1w_2$ is negative as already demonstrated). 

We assume without loss of generality that $w_1w_3$ is positive (given the symmetry of outcomes (iii) and (iv)). Since $\textrm{Val}_f(w_1w_3)\ge 1-6\beta$ and $w_1w_3$ is positive, we find by Lemma~\ref{lem:ObtuseAcutePair}(i) that $\dem(vw_3,f) \ge \textrm{Val}_f(w_1w_3)\ge 1-6\beta$. But then $\dem(uw_3,f)\le 1-\dem(vw_3,f) \le 6\beta$. Since $1-6\beta > 6\beta$ as $\beta < \frac{1}{12}$ (since $\varepsilon$ is small enough), we thus find that Lemma~\ref{lem:ObtuseAcutePair}(i) does not hold for $w_2w_3$. Thus Lemma~\ref{lem:ObtuseAcutePair}(ii) holds for $w_2w_3$; it follows that $w_2w_3$ is negative and that $|\phi(w_2w_3)|\ge \textrm{Val}_f(w_2w_3)-\dem(uw_3,f) \ge 1 - 12\beta$. Since $\alpha\ge 12\beta$, it follows that $w_2w_3\in \mc{E}$. Moreover since $w_2w_3$ is negative and $\phi$ is triangle-covered, we also find that $w_2w_4$ is positive. Then a symmetric argument to the one given above yields that $\dem(uw_4,f)\ge 1-6\beta$, $w_1w_4\in \mc{E}$ and (iv) holds as desired (where we note for the reader that outcome (iii) holds if we had assumed $w_1w_4$ is positive).
\end{proofclaim}

Since $|V_{\sigma(i)}'| > 3\sqrt{\beta}\cdot n$ and given the degree condition of $J$ mentioned before, it follows that every $xy\in E(J)$ is an edge of some copy of $K_4$ of $J$ (and such $K_4$ necessarily contains exactly one vertex from each $V_{\sigma(i)}'$) and hence Claim~\ref{cl:NormalK4} applies to every edge of $J$.

Note that since $W$ is reflective-cut or quad-cut, we have in either case that $\mu(x_1)=\mu(x_2)$ and $\mu(x_3)=\mu(x_4)$ and thus $\mu(x_1)+\mu(x_3)=\mu(x_1)+\mu(x_4)=\mu(x_2)+\mu(x_3)=\mu(x_2)+\mu(x_4)=\frac{1}{2}$. Recall that since $H$ is $\beta$-close to $W$, we have for each $i\in [4]$ that $\bigg| |V_i|-\mu(x_i)\cdot v(H)\bigg| \le \beta\cdot n$ and hence that $\bigg||V_{\sigma(i)}'|-\mu(x_{\sigma(i)})\cdot v(H)\bigg| \le (\beta+5\sqrt{\beta})\cdot n$.

Now first suppose that $m_{13} > \frac{1}{2}+\beta$. Then it follows from Claim~\ref{cl:NormalK4}(i) and the remark above that every edge of $J$ is in a $K_4$ in $J$ that $E(J[V_{\sigma(i)},V_{\sigma(j)}]) \subseteq \mc{E}$ for all $\{i,j\} \in \bigg\{ \{1,3\},\{2,4\},\{3,4\}\bigg\}$ and that all but at most $\alpha \cdot n^2$ edges in $E(J[V_{\sigma(1)},V_{\sigma(2)}])$ are in $\mc{E}$. Moreover, we also find that $m_{23} \ne 1$ and hence whether $W$ is reflective-cut or quad-cut, we find that $\mu(x_{\sigma(2)})+\mu(x_{\sigma(3)})$. Hence $\bigg| |V_{\sigma(2)}'|+|V_{\sigma(3)}'| - \frac{v(H)}{2}\bigg| \le 12\sqrt{\beta}\cdot n$. Similarly $\bigg| |V_{\sigma(1)}'|+|V_{\sigma(4)}'| - \frac{v(H)}{2}\bigg| \le 12\sqrt{\beta}\cdot n$. Thus we find that 

\begin{align*}
|\mc{E}| &\ge (|V_{\sigma(2)}'|+|V_{\sigma(3)}'|)(|V_{\sigma(1)}'|+|V_{\sigma(4)}'|-\sqrt{\beta}\cdot n) -\alpha\cdot n^2\\
&\ge \left(\frac{v(H)}{2}-12\sqrt{\beta}\cdot n\right)\left(\frac{v(H)}{2}-13\sqrt{\beta}\cdot n\right) -\alpha n^2\\
&\ge \frac{|N(u)\cap N(v)|^2}{4} - 25\sqrt{\beta}\cdot n\cdot |N(u)\cap N(v)|-\alpha\cdot n^2\\
&\ge (1-300\sqrt{\beta}-36\alpha)\cdot \frac{|N(u)\cap N(v)|^2}{4} 
\end{align*}
where used that all vertices of $J$ are normal, that $v(H)=|N(u)\cap N(v)| \ge \frac{n}{3}$, and $\alpha \ge 300\sqrt{\beta}$. But then it follows that $f$ is $(36\alpha+4900\sqrt{\beta})$-bipartite as desired by fixing any bipartition $(A,B)$ of $N(u)\cap N(v)$ with $V_{\sigma(2)}'\cup V_{\sigma(3)}'\subseteq A$ and $V_{\sigma(1)}'\cup V_{\sigma(4)}'\subseteq B$.

Next suppose that $m_{23} > \frac{1}{2}+\beta$. Then the symmetric argument (using Claim~\ref{cl:NormalK4}(ii) instead) to the one given above yields that $f$ is $(36\alpha+4900\sqrt{\beta})$-bipartite as desired.

Finally we assume that $m_{12},m_{23}\le \frac{1}{2}+\beta$ (and hence also that $m_{12},m_{23}\ge \frac{1}{2}-\beta$). Thus for each $K_4$ in $J$, one of Claim~\ref{cl:NormalK4}(iii) or (iv) holds for it. We claim that the same outcome (either (iii) or (iv)) holds for all $K_4$ in $J$. Suppose not. Then there exists $w_3\in V_{\sigma(3)}'$ with $\dem(uw_3,f)\ge 1-6\beta$ (since outcome (iii) holds for some $K_4$) and similarly there exists $w_3'\in V_{\sigma(3)}'$ with $\dem(vw_3',f)\ge 1-6\beta$. Since $|V_{\sigma(4)}'| > 2\sqrt{\beta}\cdot n$, there exists a vertex $w_4\in V_{\sigma(4)}'$ such that both $w_3w_4$ and $w_3'w_4$ are normal (and hence in $E(J)$). As remarked before every edge of $J$ is in some copy of $K_4$ in $J$. Since $1-6\beta > \beta$ as $\beta < \frac{1}{12}$, we find that Claim~\ref{cl:NormalK4}(iv) does not hold for a copy of $K_4$ containing $w_3w_4$ and hence (iii) holds yielding that $\dem(vw_4,f)\ge 1-6\beta$. Similarly Claim~\ref{cl:NormalK4}(iii) does not hold for a copy of $K_4$ containing $w_3'w_4$ and hence (iv) holds yielding that $\dem(uw_4,f)\ge 1-6\beta$. But then $1=|\phi(f)|=\dem(uw_4,f)+\dem(vw_4,f) \ge 2-12\beta > 1$, a contradiction. This proves our claim that the same outcome holds for all $K_4$ in $J$.

We assume without loss of generality that outcome (iii) holds for all $K_4$ in $J$. Then it follows from Claim~\ref{cl:NormalK4}(iii) and the remark from before that every edge of $J$ is in a $K_4$ in $J$ that $E(J[V_{\sigma(i)},V_{\sigma(j)}]) \subseteq \mc{E}$ for all $\{i,j\} \in \bigg\{ \{1,2\},\{1,3\},\{2,4\},\{3,4\}\bigg\}$. Then an identical argument as the one given for the case $m_{13} \frac{1}{2}+\beta$ yields that $f$ is $(36\alpha+4900\sqrt{\beta})$-bipartite as desired. 
\end{proof}

\begin{lem}[Extremal Neighborhood Structure - Charge Version]\label{lem:NeighborhoodStructure2}
There exists $\varepsilon_0$ such that the following holds for each $\varepsilon\in (0,\varepsilon_0]$: Let $\alpha \in (0,1)$ and let $\gamma$ such that $2\varepsilon^{1/4} \ge \gamma\ge \sqrt{\varepsilon}$. Let $(G,\phi)$ be a nonzero untied triangle-covered edge-weighted graph on $n$ vertices with $\delta(G)\ge (\frac{3}{4}-\varepsilon)n$ and let $f=uv$ be a negative edge of $G$. If $\textrm{ch}_F(f) \le \gamma \cdot |\phi(f)|$, then 
$$|N(u)\cap N(v)| \le \left(\frac{1}{2}+2\gamma\right)n$$ 
and $f$ is
\begin{enumerate}
    \item[(i)] $\alpha$-locally-small, or
    \item[(ii)] $(36\alpha+5000\cdot \gamma^{1/2000})$-bipartite.
\end{enumerate}
\end{lem}
\begin{proof}
We choose $\varepsilon_0$ small enough so that $\gamma$ is smaller than the $\varepsilon_0$ in Lemma~\ref{lem:NeighborhoodStructure}. Recall that by Lemma~\ref{lem:ChargePrime}, we have that
$$\textrm{ch}_F(f)\ge \textrm{ch}'_F(f) = \phi(f)\left(1-\frac{4}{n}\cdot |N(u)\cap N(v)|\right) - \frac{4}{n^2} \sum_{e\in T^+_{G,\phi}(f)}~~\sum_{f'\in T^{-}_{G,\phi}(e)}~~\demh_f(e,f').$$
Let $d:= \frac{|N(u)\cap N(v)|}{n}$ and $m:= \frac{1}{n}\cdot \max\{|N(u)\setminus N(v)|,~|N(v)\setminus N(u)|\}$. Since $\delta(G)\ge \left(\frac{3}{4}-\varepsilon\right)n$, we find that $1\ge d\ge \frac{1}{2}-2\varepsilon$ and $m\le \frac{1}{4}+\varepsilon$. Also note that $\frac{1-d}{2}\le m\le 1-d$.

For $e\in T^+_{G,\phi}(f)$, we find that
\begin{align*}\sum_{f' \in T^{-}_{G,\phi}(e): f'\setminus e \not\in N(u)\cap N(v)} \demh_f(e,f') &= \sum_{x\not\in N(u)\cap N(v): e\cup\{x\} \textrm{ is a triangle of $G$}}  \demh_f(e,e\cup\{x\}).\\
&\le |\phi(f)|\cdot |\{x\not\in N(u)\cap N(v): e\cup\{x\} \textrm{ is a triangle of $G$}\}|\\
&\le |\phi(f)|\cdot \max\{|N(u)\setminus N(v)|, |N(v)\setminus N(u)|\} \\
&\le |\phi(f)|\cdot m\cdot n,
\end{align*}
where for the first inequality we used that $\demh_f(e,e\cup{x})\le |\phi(f)|$ 
by Proposition~\ref{prop:demandBound}, the second inequality follows as at exactly one of $u$ or $v$ is an end of $e$, and the third inequality follows by definition of $m$. 

Hence we calculate that  
\begin{align*}
\sum_{e\in T^+_{G,\phi}(f)}~~\sum_{f' \in T^{-}_{G,\phi}(e):~f'\setminus e \not\in N(u)\cap N(v)} \demh_f(e,f') &\le m\cdot n \cdot \sum_{e\in T^+_{G,\phi}(f)} \dem(e,f)\\
&\le m\cdot n\cdot |N(u)\cap N(v)| \cdot |\phi(f)|,
\end{align*}
where we used that $|\{e\in T^{+}_{G,\phi}(f)| \le |N(u)\cap N(v)|$. Substituting this bound yields that (recalling as $f$ is negative that $|\phi(f)|=-\phi(f)$)
\begin{align*}
\textrm{ch}_F(f) &\geq \phi(f) \left(1 - \frac{(4-4m)}{n} \cdot |N(u)\cap N(v)|\right) - \frac{4}{n^2} \cdot \textrm{ch}''_F(f),
\end{align*}
Since $\textrm{ch}_F(f) \le  \gamma\cdot |\phi(f)|$, this rearranges to
$$\textrm{ch}''_F(f) \ge \frac{n^2}{4} \cdot |\phi(f)|\cdot \bigg((4-4m)d-1-\gamma\bigg)$$
Recall from our proof of Conjecture~\ref{conj:FracNW6}, we found that $\textrm{ch}''_F(f) \le |\phi(f)|\cdot \frac{|N(u)\cap N(v)|^2}{2}$ and hence combining these inequalities yields
$$2d^2 \ge (4-4m)d-1-\gamma.$$
First suppose that $d\ge \frac{9}{10}$, then $m\le \frac{1}{10}$ and the right-side above is at least $\frac{36}{10}d-1-\gamma \ge \frac{36}{10}\cdot \frac{9}{10}-1 -\gamma \ge \frac{324}{100}-1-\gamma \ge 2.24-\gamma > 2$, a contradiction (since the left-side is at most $2$ as $d\le 1$ and $\gamma \le \varepsilon^{1/4}$ is small enough). 

So we assume $d < \frac{9}{10}$. Now $m\le \frac{1}{4}+\varepsilon$ and hence 
$$\textrm{ch}''_F(f) \ge \frac{n^2}{4} \cdot |\phi(f)|\cdot \bigg((3-4\varepsilon)d-1-\gamma\bigg),$$
and also that
$$2d^2 \ge (3-4\varepsilon)d-1-\gamma.$$
Suppose for a contradiction that $d \ge \frac{1}{2}+2\gamma.$ Then $d\in [\frac{1}{2}+2\gamma,~0.9]$. Consider the function $f(d):= 2d^2-(3-4\varepsilon)d$. Note that $f$ is minimized at $4d=3-4\varepsilon$ or equivalently $d=\frac{3}{4}-\varepsilon$. Since this lies inside our range of $d$ (as $\varepsilon$ is small enough), it follows as $f(d)$ is a quadratic function that the maximum of $f$ in this range is at one of the two endpoints. Yet $f(0.9) = 2(0.9)^2 - (3-4\varepsilon)\cdot 0.9 \le 1.62 - 2.7 + 3.6\cdot \varepsilon = -1.08+3.6\varepsilon < -1-\gamma$ since $\varepsilon$ is small enough. Similarly 
\begin{align*}
f\left(\frac{1}{2}+2\gamma\right) &= 2\left(\frac{1}{2} + 2\gamma\right)^2-(3-4\varepsilon)\left(\frac{1}{2}+2\gamma\right)\\
&= \frac{1}{2} +4\gamma + 8\gamma^2 -\frac{3}{2} -6\gamma+2\varepsilon+8\varepsilon\cdot \gamma\\
&\le -1 -2\gamma+8\gamma(\gamma+\varepsilon)+2\varepsilon\\
&<-1-\gamma,
\end{align*} 
where the last two inequalities follow since $\gamma\ge \sqrt{\varepsilon}$ and $\varepsilon$ is small enough. Thus for all $d\in [\frac{1}{2}+2\gamma,~0.9]$, we find that $f(d) < -1-\gamma$, a contradiction.

So we assume $d\le \frac{1}{2}+2\gamma$. Thus 
$$n\ge \frac{1}{\frac{1}{2}+2\gamma} \cdot |N(u)\cap N(v)|  = \frac{2}{1+4\gamma} \cdot |N(u)\cap N(v)| \ge 2\cdot (1-4\gamma) \cdot |N(u)\cap N(v)|.$$
Since $d\ge \frac{1}{2}-2\varepsilon$, we find that 
\begin{align*}
\textrm{ch}''_F(f) &\ge \frac{n^2}{4}\cdot |\phi(f)| \cdot \left(\frac{3}{2}-8\varepsilon-1-\gamma\right) \ge \frac{n^2}{4}\cdot |\phi(f)|\cdot \left(\frac{1}{2}-9\gamma\right).\\
&\ge (1-4\gamma)^2\cdot |N(u)\cap N(v)|^2\cdot|\phi(f)|\cdot \frac{1-18\gamma}{2} \ge (1-26\gamma)\cdot \frac{|N(u)\cap N(v)|^2}{2}\cdot |\phi(f)|.
\end{align*}
Applying Lemma~\ref{lem:NeighborhoodStructure} with $26\gamma$ and $\alpha$ yields that $f$ is $\alpha$-locally-small of $f$ is $(36\alpha+4900\cdot (26\gamma)^{0.0005})$-bipartite. Since $4900\cdot (26)^{0.0005} \le 4908 \le 5000$, we find that $f$ is $(36\alpha+5000\cdot \gamma^{1/2000})$-bipartite as desired. 
\end{proof}

\section{Forcing the Global Structure}

Finally we proceed with forcing the global extremal structure.

\subsection{Cleaning Lemma}

\begin{definition}
 Let $(G,\phi)$ be a nonzero untied edge-weighted graph on $n$ vertices and let $f=uv$ be a negative edge of $G$. Let $x\in N(u)\cap N(v)$ be $u$-obtuse. We say $z\in N(u)\setminus N(v)$ is \emph{$\alpha$-unclean} for $x$ with respect to $f$ if either  
\begin{itemize}
    \item[(a)] $xz\not\in E(G)$, or
    \item[(b)] $\sum_{e\in E^-(G)\cap \{xz,uz\}} |\phi(e)| < (1-\alpha)\cdot \phi(ux)$.
\end{itemize}
We say $x$ is \emph{$\alpha$-unclean} with respect to $f$ if there exist at least $\alpha\cdot n$ vertices $z$ that are $\alpha$-unclean for $x$ with respect to $f$; we say $x$ is \emph{$\alpha$-clean} with respect to $f$ if it is not $\alpha$-unclean.
\end{definition}

\begin{lem}[Cleaning]\label{lem:FracCleaning}
Let $\varepsilon \in (0,\frac{1}{4}]$. Let $(G,\phi)$ be a nonzero untied triangle-covered edge-weighted graph on $n$ vertices with $\delta(G)\ge \left(\frac{3}{4}-\varepsilon\right)n$ and let $f=uv$ be a negative edge of $G$. If there exist at least $\alpha\cdot n$ vertices in $N(u)\cap N(v)$ that are $\alpha$-unclean with respect to $f$, then $\textrm{ch}'
_F(f)\ge |\phi(f)|\cdot (4\alpha^3-4\varepsilon)$.
\end{lem}

\begin{proof}
Recall that by Lemma~\ref{lem:ChargePrime}, we have that
$$\textrm{ch}'_F(f) = \phi(f)\left(1-\frac{4}{n}\cdot |N(u)\cap N(v)|\right) - \frac{4}{n^2} \sum_{e\in T^+_{G,\phi}(f)}~~\sum_{f'\in T^{-}_{G,\phi}(e)}~~\demh_f(e,f').$$
Yet
$$\sum_{f' \in T^{-}_{G,\phi}(e): f'\setminus e \not\in N(u)\cap N(v)} \demh_f(e,f') = \sum_{z\not\in N(u)\cap N(v): e\cup\{z\} \textrm{ is a triangle of $G$}}  \demh_f(e,e\cup\{z\}).$$
Suppose that $x\in N(u)\cap N(v)$ is $u$-obtuse and $\alpha$-unclean with respect to $f$. Then
\begin{align*}\sum_{f' \in T^{-}_{G,\phi}(ux): f'\setminus ux \not\in N(u)\cap N(v)} \demh_f(ux,f') &= \sum_{z\not\in N(u)\cap N(v): uxz \textrm{ is a triangle of $G$}}  \demh_f(ux,uxz).\\
&\le \dem(ux,f)\cdot (|N(u)\setminus N(v)| - \alpha^2\cdot n)\\
&\le |\phi(f)|\cdot (|N(u)\setminus N(v)| - \alpha^2\cdot n)\}\\
&\le |\phi(f)|\cdot \left(n\left(\frac{1}{4}+\varepsilon\right) - \alpha^2\cdot n\right)\\
&= |\phi(f)|\cdot n\cdot \left(\frac{1}{4} +\varepsilon - \alpha^2\right),
\end{align*}
where we used that $|N(u)\setminus N(v)| \le n\cdot \left(\frac{1}{4}+\varepsilon\right)$ since $\delta(G) \ge (\frac{3}{4}-\varepsilon)n$. 
Since there are at least $\alpha\cdot n$ vertices in $N(u)\cap N(v)$ that are $\alpha$-unclean with respect to $f$, we find that
\begin{align*}
\sum_{e\in T^+_{G,\phi}(f)}~~\sum_{f' \in T^{-}_{G,\phi}(e):~f'\setminus e \not\in N(u)\cap N(v)} \demh_f(e,f') \le |\phi(f)| \cdot \bigg( |N(u)\cap N(v)| \cdot n\left(\frac{1}{4}+\varepsilon\right) - \alpha^3 n^2\bigg).
\end{align*}
Substituting this into the equation from Lemma~\ref{lem:ChargePrime} yields that
\begin{align*}
\textrm{ch}'_F(f) &\ge \phi(f) \left(1 - \frac{3-4\varepsilon}{n} \cdot |N(u)\cap N(v)|\right) + 4\alpha^3 - \frac{4}{n^2} \cdot \textrm{ch}''_F(f).
\end{align*}
We let $H:= \mathrm{Filter}_{G,\phi}(f)$ and define an edge-weighting $\varphi$ on $E(H)$ as $\varphi(g) := \frac{\textrm{Val}_f(g)}{2\cdot |\phi(f)|}$ for each $g\in E(H)$. Then it follows from Proposition~\ref{prop:EdgeWeightUpperBound} that $\varphi$ is a $[0,1]$-edge-weighting of $H$. Similarly, it follows from Lemma~\ref{lem:EdgeWeightTriangleThin} that $\varphi$ is triangle-thin. Hence by Lemma~\ref{lem:Value}, we have that $|\varphi| \le \frac{v(H)^2}{4}$. Since $v(H)=|N(u)\cap N(v)|$, it follows that
$$\textrm{ch}''_F(f) = \sum_{g\in E(G[N(u)\cap N(v)])} \textrm{Val}_f(g) = 2\cdot |\phi(f)| \cdot |\varphi| \le |\phi(f)| \cdot \frac{1}{2} \cdot |N(u)\cap N(v)|^2.$$
Let $d:= \frac{|N(u)\cap N(v)|}{n}$. Substituting this and combining the previous inequalities yields that 
\begin{align*}
\textrm{ch}'_F(f) &\ge -|\phi(f)| \cdot \bigg(1 -(3 -4\varepsilon)d\bigg) + 4\cdot |\phi(f)| \cdot \frac{1}{2} \cdot d^2 + 4\alpha^3 =  |\phi(f)|\cdot \bigg(-1 +(3-4\varepsilon)d -2d^2 + 4\alpha^3\bigg).
\end{align*}
The derivative of $-2d^2+(3-4\varepsilon)d-1$ is $-4d + 3-4\varepsilon$. Since $1 \ge d \ge \frac{1}{2} - 2\varepsilon$ as $\delta(G) \ge (\frac{3}{4}-\varepsilon)n$, the quadratic function is minimized on this range at which ever end point is larger (since the maximum is at $d=\frac{3-4\varepsilon}{4}=\frac{3}{4}-\varepsilon$ which lies inside the range). At $d=1$, the function is $-4\varepsilon$. At $d=\frac{1}{2}-2\varepsilon$, the function evaluates to $-4\varepsilon$. Thus we find that 
\begin{align*}
\textrm{ch}'_F(f) \ge |\phi(f)|\cdot \bigg(4\alpha^3-4\varepsilon\bigg),
\end{align*}
as desired.
\end{proof}

For our main proof, we need the following corollary of Lemma~\ref{lem:FracCleaning}.

\begin{cor}\label{cor:FracCleaning}
Let $\varepsilon \in [0,\frac{1}{4}]$. Let $(G,\phi)$ be a nonzero untied triangle-covered edge-weighted graph on $n$ vertices with $\delta(G)\ge \left(\frac{3}{4}-\varepsilon\right)n$. If $f=uv$ is a negative edge of $G$ with $\textrm{ch}_F(f)\le \sqrt{\varepsilon}\cdot |\phi(f)$, then there exist at most $2\cdot \varepsilon^{1/6}\cdot n$ vertices in $N(u)\cap N(v)$ that are $2\cdot \varepsilon^{1/6}$-unclean with respect to $f$. 
\end{cor}
\begin{proof}
Follows from Lemma~\ref{lem:FracCleaning} with $\alpha=2\cdot\sqrt{\varepsilon}^{1/6}$.         
\end{proof}

For our main proof, we also need a lower bound on the final charges of negative edges as follows.

\begin{lem}\label{lem:MinDegFracLowerBound}
Let $\varepsilon \in [0,\frac{1}{4}]$. Let $(G,\phi)$ be a nonzero untied triangle-covered edge-weighted graph on $n$ vertices with $\delta(G)\ge \left(\frac{3}{4}-\varepsilon\right)n$. If $f=uv$ is a negative edge of $G$, then $\textrm{ch}'
_F(f)\ge -|\phi(f)|\cdot 4\varepsilon$.
\end{lem}
\begin{proof}
Follows from Lemma~\ref{lem:FracCleaning} with $\alpha=0$.     
\end{proof}

\subsection{Fractional Stability Proof}

It is mildly convenient in our next major proofs if the bipartitions from $\alpha$-bipartiteness are balanced. Hence the following definitions. Since technically a balanced bipartition would necessitate the use of floors and ceilings (as in the following definition) but these do not materially affect the proof (which lies in the asymptotic setting), in our proofs in this section we will omit the floors and ceilings for readability.

\begin{definition}
Let $(G,\phi)$ be a nonzero untied triangle-covered edge-weighted graph on $n$ vertices. Let $f=uv$ be a negative edge of $G$. Let $P=(A,B)$ be a bipartition of $N(u)\cap N(v)$. We say a pair $a\in A$, $b\in B$ is \emph{$\alpha$-normal} for $\mc{P}$ if $ab\in E_{\phi}^-(G)$ and $|\phi(ab)| \ge (1-\alpha)\cdot |\phi(uv)|$. We say $\mc{P}$ is an \emph{$\alpha$-bipartition} for $f$ if at least $(1-\alpha)\frac{|N(u)\cap N(v)|^2}{4}$; we say $\mc{P}$ is an \emph{$\alpha$-balanced-bipartition} for $f$ if in addition $|A|=\left\lfloor \frac{|N(u)\cap N(v)|}{2} \right\rfloor$ and $|B|=\left\lceil \frac{|N(u)\cap N(v)|}{2} \right\rceil$. We say a vertex $w\in A\cup B$ is \emph{$\alpha$-normal} for $\mc{P}$ if it is in at least $(1-\sqrt{\alpha})\cdot \frac{\lfloor |N(u)\cap N(v)|\rfloor}{2}$ $\alpha$-normal pairs. We say $f$ is \emph{$\alpha$-balanced-bipartite} if there exists an $\alpha$-balanced-bipartition for $f$.
\end{definition}

The following proposition shows the two concepts are functionally equivalent.

\begin{proposition}\label{prop:balanced}
Let $(G,\phi)$ be a nonzero untied triangle-covered edge-weighted graph on $n$ vertices. Let $f=uv$ be a negative edge of $G$ and $\alpha\in (0,1)$. If $f$ is $\frac{\alpha^2}{2}$-bipartite, then it is $\alpha$-balanced-bipartite.
\end{proposition}
\begin{proof}
Since $f$ is $\frac{\alpha^2}{2}$-bipartite, there exists an $\frac{\alpha^2}{2}$-bipartition $\mc{P}=(A,B)$ of $f$. We assume without loss of generality that $|A|\le |B|$. Let $m:= |N(u)\cap N(v)|$. We claim that $|A|\ge (1-\frac{\alpha}{\sqrt{2}})\cdot  \frac{m}{2}$. Suppose not. Then the number of $\frac{\alpha^2}{2}$-normal pairs for $\mc{P}$ is at most $|A|\cdot |B| < \left(1-\frac{\alpha}{\sqrt{2}}\right)\left(1+\frac{\alpha}{\sqrt{2}}\right)\cdot \frac{m^2}{4}=(1-\frac{\alpha^2}{2})\cdot \frac{m^2}{4}$, a contradiction. 

Let $S\subseteq B$ with $|S|  =\left\lfloor \frac{|N(u)\cap N(v)|}{2} \right\rfloor - |A|$. Hence $|S|\le \frac{\alpha}{\sqrt{2}}\cdot\frac{m}{2} \le \frac{\alpha}{2}\cdot n$. Let $\mc{P}'=(A',B')$ with $A':= A\cup S$ and $B':= B\setminus S$. It follows that $\mc{P}'$ is an $(\frac{\alpha^2}{2} + \frac{\alpha}{2})$-balanced-bipartition of $f$. Since $\alpha\le 1$, $\mc{P}'$ is also then an $\alpha$-balanced-bipartition of $f$ as desired.
\end{proof}

We also require the next proposition asserting there are not too many abnormal vertices for $\alpha$-balanced-bipartition.

\begin{proposition}\label{prop:NotTooManyAbnormal}
Let $(G,\phi)$ be a nonzero untied triangle-covered edge-weighted graph on $n$ vertices. Let $f=uv$ be a negative edge of $G$ and $\alpha\in (0,1)$. If $\mc{P}$ is an $\alpha$-balanced-bipartition for $f$, then there are at most $\sqrt{\alpha}\cdot n$ vertices in $N(u)\cap N(v)$ that are not $\alpha$-normal for $\mc{P}$.
\end{proposition}
\begin{proof}
Let $m:= |N(u)\cap N(v)|$ and let $(A,B):= \mc{P}$. Then there are most $\alpha\cdot \frac{m^2}{4}$ pairs $a\in B, b\in B$ such that $ab$ is not $\alpha$-normal for $\mc{P}$. Hence there are at most $2\alpha\cdot \frac{m^2}{4} \cdot \frac{1}{\sqrt{\alpha}\cdot n}$ vertices in $N(u)\cap N(v)$ that are not $\alpha$-normal for $\mc{P}$. But this is at most $\sqrt{\alpha}\cdot n$ as desired.     
\end{proof}

Here is another useful proposition.

\begin{proposition}\label{prop:BalancedObtuse}
Let $(G,\phi)$ be a nonzero untied triangle-covered edge-weighted graph. If a negative edge $f=uv$ of $G$ is $\alpha$-balanced-bipartite, then there are at most $(\frac{1}{2}+\sqrt{\alpha})\cdot |N(u)\cap N(v)|$ vertices that are $u$-obtuse with respect to $f$. 
\end{proposition}
\begin{proof}
Suppose not. Let $m=|N(u)\cap N(v)|$. Let $S$ be the set of vertices that are $u$-obtuse with respect to $f$. Thus $|S| > (\frac{1}{2}+\sqrt{\alpha})\cdot m$. Let $S_A := S\cap A$ and $S_B:= S\cap B$. Since $|A|=|B|=m$, it follows $|S_A|,|S_B| > \sqrt{\alpha}\cdot m$ and hence $|S_A|\cdot |S_B| > \alpha\cdot m^2$. But for each pair $a\in S_A, b\in S_b$, we find that $ab$ is a positive edge (since $av$ and $bv$ are negative as $a$ and $b$ are $u$-obtuse and $\phi$ is triangle-covered). But then $ab$ is not $\alpha$-normal for $\mc{P}$. So there are at least $\alpha\cdot m^2$ pairs that are not $\alpha$-normal for $\mc{P}$, contradicting that $\mc{P}$ is an $\alpha$-balanced-bipartition.
\end{proof}

We are now prepared to our main fractional stability result which is a slightly stronger version of Theorem~\ref{thm:fracNWstability} (as it asserts that our specific discharging works).

\begin{thm}\label{thm:Charge2Stability}
There exists $\varepsilon_0$ such that the following holds for all $\varepsilon \in (0,\varepsilon_0]$: Let $\sigma:= \varepsilon^{1/20000}$. Let $(G,\phi)$ be a nonzero untied triangle-covered edge-weighted graph on $n$ vertices with $\delta(G) \ge \left(\frac{3}{4}-\varepsilon\right)n$. If there exists some edge $f^*$ of $G$ such that $\textrm{ch}_2(f^*)<0$, then both of the following hold:
\begin{itemize}
    \item[(i)] at least $(1-\sigma)n$ vertices have degree at most $(\frac{3}{4}+\sigma)n$, and
    \item[(ii)] the maximum cut of $G$ is at least $(1-\sigma)\cdot \frac{n^2}{4}$.
\end{itemize}     
\end{thm}
\begin{proof}
Suppose not. Let $\mc{E}_0 := \{f\in E^{-}_{\phi}(G):~\textrm{ch}_F(f) \ge \sqrt{\varepsilon}\cdot |\phi(f)|\}$, that is the set of edges that discharge in Rule 1. Let $\mc{E}_1 := \{f\in E^{-}_{\phi}(G):~\textrm{ch}_1(f) \ge \sqrt{\varepsilon}\cdot |\phi(f)|\}$, that is the set of edges that that discharge in Rule 2. Let $\mc{E}_2:= E(G)\setminus (\mc{E}_0\cup \mc{E}_1)$. Since there exists $f^*$ with $\textrm{ch}_2(f^*)<0$, we find that $\mc{E}_2\ne \emptyset$. Let $f_0=u_0v_0 \in \mc{E}_2$ such that $|\phi(f_0)|$ is maximized. 

Let $\mc{Q}:= \{f\in E^{-}_{\phi}(G): |\phi(f)| > |\phi(f_0)|\}$. We let $\mc{H}(\gamma) := \{E^{-}_{\phi}(G): |\phi(f)| \ge \gamma\}$ and $\mc{H}:= \mc{H}(\frac{3}{4}\cdot |\phi(f_0)|)$. Let $\mc{Q}':= \mc{H} \setminus \mc{E}_2$. Note by our choice of $f_0$, we have that $\mc{Q}\subseteq \mc{E}_0\cup \mc{E}_1$. Thus $\mc{Q}\subseteq \mc{Q}'$. Since $\mc{Q}'\cap \mc{E}_2=\emptyset$ by definition, we find that each $f\in \mc{Q}'$ discharges in either Rule 1 or Rule 2. 

\begin{claim}\label{cl:QSmall}
$|\mc{Q}|\le |\mc{Q}'|\le 6\sqrt{\varepsilon}\cdot n^2$.
\end{claim}
\begin{proofclaim}
Suppose not. Since each $f\in \mc{Q}'$ discharges in either Rule 1 or Rule 2, we have that each such $f\in \mc{Q}'$ sends to $f^*$ an amount of charge at least 
$$\sqrt{\varepsilon}\cdot |\phi(f)|\cdot \frac{1}{2}\cdot \frac{1}{e(G)} \ge \sqrt{\varepsilon}\cdot \frac{1}{n^2}\cdot \frac{3}{4}\cdot  |\phi(f_0)|,$$
where we used that $|\phi(f)|\ge \frac{3}{4}\cdot |\phi(f_0)|$ by definition of $\mc{Q}'$ and that $2\cdot e(G)\le n^2$. Note that $|\phi(f^*)|\le |\phi(f_0)|$ as otherwise $\textrm{ch}_2(f^*)\ge 0$, a contradiction. But then $f^*$ receives in the additional discharging rules an amount of charge at least
$$|Q|\cdot \sqrt{\varepsilon}\cdot \frac{1}{n^2}\cdot \frac{3}{4}\cdot |\phi(f_0)| \ge 4\varepsilon\cdot |\phi(f^*)|.$$
Yet by Lemma~\ref{lem:MinDegFracLowerBound}, $\textrm{ch}_F(f^*) \ge -4\varepsilon\cdot |\phi(f^*)|$ and hence $\textrm{ch}_2(f^*)\ge 0$, a contradiction.
\end{proofclaim}

Let $S$ denote the set of vertices $v\in V(G)$ such that $v$ is incident with at least $3\cdot\varepsilon^{1/4}\cdot n$ edges of $\mc{Q}'$. It follows from Claim~\ref{cl:QSmall} that $|S|\le 4\cdot \varepsilon^{1/4}\cdot n$. Let $G'$ be the graph with $V(G'):= V(G)\setminus S$ and $E(G'):= G[V(G')]\cap \mc{E}_2\cap \mc{H}$. 

We now prove a very useful claim characterizing the structure of edges in $\mc{E}_2\cap H$ as follows. To that end, for the rest of the proof let $\alpha:= \varepsilon^{1/18000}$.

\begin{claim}\label{cl:BipartiteStructure}
If $f=uv \in \mc{E}_2\cap \mc{H}$, then both of the following hold:
\begin{itemize}
\item[(i)] $|N(u)\cap N(v)|\le \left(\frac{1}{2}+2\sqrt{\varepsilon}\right)n$ and hence $d_G(u),d_G(v)\le \left(\frac{3}{4}+3\sqrt{\varepsilon}\right)n$, and 
\item[(ii)] $f$ is $\alpha$-balanced-bipartite.
\end{itemize}
\end{claim}
\begin{proofclaim}
Suppose not. Since $\textrm{ch}'_F(f) \le \textrm{ch}_F(f) \le \textrm{ch}_1(f) \le \sqrt{\varepsilon}\cdot |\phi(f)|\le \sqrt{\varepsilon}\cdot |\phi(f)|$, we have by Lemma~\ref{lem:NeighborhoodStructure2} $\gamma_{\ref{lem:NeighborhoodStructure2}}:=2\cdot \varepsilon^{1/4})$ that $N(u)\cap N(v)| \le \left(\frac{1}{2} + 2\sqrt{\varepsilon}\right)n$ and hence (i) holds. 

Let $\gamma := \varepsilon^{1/4}$. Now we also apply Lemma~\ref{lem:NeighborhoodStructure2} with $\gamma_{\ref{lem:NeighborhoodStructure2}}:=\gamma$ and $\alpha_{\ref{lem:NeighborhoodStructure2}}:=4\gamma$ to find that $f$ is $\gamma$-locally-small or $(144\gamma+5000\cdot \gamma^{1/2000})$-bipartite. First suppose the latter holds. Since $144\gamma + 5000\cdot \gamma^{1/2000}\le \frac{1}{2}\cdot \varepsilon^{1/9000}$ as $\varepsilon$ is small enough, we then find that $f$ is $\frac{1}{2}\cdot \varepsilon^{1/9000}$-bipartite and hence by Proposition~\ref{prop:balanced}, $f$ is $\alpha$-balanced-bipartite, that is (ii) holds and thus the whole claims holds, a contradiction. 

Hence we assume $f$ is $4\gamma$-locally-small, that is $f$ is incident with at least $4\gamma\cdot n$ negative edges $f'$ with $|\phi(f')|\ge 2\cdot |\phi(f)|$. Since $|\phi(f)| > \frac{1}{2}\cdot |\phi(f_0)|$ as $f\in \mc{H}$, it follows that each such $f'\in Q$. Let $R$ be the set of edges incident with $f$ that are in $Q$; hence $|R|\ge 4\gamma\cdot n$. Let $R' := \{f'\in R: \textrm{ch}_F(f') \ge 2\gamma\cdot |\phi(f')|\}$. 

First suppose that $|R'| \ge 2\gamma\cdot n$. But then each $f'\in R'$ sends in Rule $1$ at least $2\gamma\cdot |\phi(f')|\cdot \frac{1}{2}\cdot \frac{1}{2n}$ charge to $f$. Since $|\phi(f')|\ge 2\cdot |\phi(f)|$, this implies that $f$ receives charge under Rule 1 of at least 
$$(2\gamma\cdot n) \cdot \frac{2\gamma\cdot2\cdot |\phi(f)|}{4n}=2\gamma^2\cdot |\phi(f)| =2\sqrt{\varepsilon}\cdot |\phi(f)| \ge (4\varepsilon+\sqrt{\varepsilon})\cdot |\phi(f)|,$$
where the last inequality follows since $\varepsilon$ is small enough. But $\textrm{ch}_F(f) \ge -4\varepsilon\cdot |\phi(f)|$ by Lemma~\ref{lem:MinDegFracLowerBound}. Thus we find that $\textrm{ch}_1(f) \ge \textrm{ch}_F(f) + (4\varepsilon+\sqrt{\varepsilon})\cdot |\phi(f)| \ge \sqrt{\varepsilon}\cdot |\phi(f)|$, contradicting that $f\in \mc{E}_2$.

So we assume $|R'|\le 2\gamma\cdot n$. Thus $|R\setminus R'|\ge 2\gamma\cdot n$. By Lemma~\ref{lem:NeighborhoodStructure2} with $\gamma_{\ref{lem:NeighborhoodStructure2}}:=2\gamma$ and $\alpha_{\ref{lem:NeighborhoodStructure2}}:=\frac{1}{200}$, we have that for each $f'\in R\setminus R'$ that $f'$ is $\frac{1}{200}$-locally-small or $(\frac{36}{200} + 5000\cdot \gamma^{1/2000})$-bipartite. First suppose there exists $f'=u'v'\in R\setminus R'$ that is $(\frac{36}{200} + 5000\cdot \gamma^{1/2000})$-bipartite and hence is $\frac{1}{5}$-bipartite since $\varepsilon$ is small enough. Thus by definition, there exist at least $\frac{4}{5}\cdot \frac{|N(u')\cap N(v')|^2}{4}$ negative edges $ab$ with $|\phi(ab)|\ge \frac{4}{5}\cdot |\phi(f')|$. But $|\phi(f')|\ge 2\cdot |\phi(f)| \ge \frac{3}{2}\cdot |\phi(f_0)|$ since $|\phi(f)|\ge \frac{3}{4}\cdot |\phi(f_0)|$ s $f\in\mc{H}$; hence $\frac{4}{5}\cdot |\phi(f')| \ge \frac{4}{5}\cdot \frac{3}{2} |\phi(f_0)| > |\phi(f_0)|$. Thus we find that $|Q| \ge \frac{4}{5}\cdot \frac{|N(u')\cap N(v')|^2}{4} \ge \frac{4}{5} \cdot \frac{(n/3)^2}{4} \ge \frac{n^2}{45}$, contradicting Claim~\ref{cl:QSmall} since $\varepsilon$ is small enough.

So we assume every edge in $R\setminus R'$ is $\frac{1}{200}$-locally-small. Recalling crucially that locally small means \emph{both} endpoints are incident with many edges that are more negative, we find there are at least $2\gamma n$ vertices each incident with at least $\frac{n}{200}$ edges in $Q$ (noting that for each vertex $y$ only one of $u'y$ or $v'y$ is in $R$ as $\phi$ is triangle-covered). Thus $|Q| \ge \frac{1}{2}\cdot 2\gamma n\cdot \frac{n}{200} \ge \frac{\gamma}{200}\cdot n^2 > 6\sqrt{\varepsilon}\cdot n^2$ contradicting Claim~\ref{cl:QSmall} where we used that $\varepsilon$ is small enough such for this inequality to hold. 
\end{proofclaim}

We also deduce the following claim that more characterizes the common neighborhood of negative edges in $\mc{E}_2$ with slightly larger negative weight than $\frac{3}{4}\cdot |\phi(f_0)|$ as follows.

\begin{claim}\label{cl:BipartiteStructure2}
Suppose that $f=uv \in \mc{E}_2\cap \mc{H}$ with $|\phi(f)|\ge \frac{9}{10}\cdot |\phi(f_0)|$. If $\mc{P}$ is an $\alpha$-balanced-bipartition of $f$ and $\mc{N}$ denotes the set of vertices of $N(u)\cap N(v)\cap V(G')$ that are both $\alpha$-normal for $\mc{P}=(V_1,V_2)$ and are not $2\cdot \varepsilon^{1/6}$-unclean with respect to $f$, then all of the following hold:
\begin{itemize}
    \item[(i)] $|(N(u)\cap N(v))\setminus \mc{N}| \le (\sqrt{\alpha}+6\cdot \varepsilon^{1/6})\cdot n,$
    \item[(ii)] for each $i\in \{1,2\}$ and $w\in V_i\cap \mc{N}$, we have that
    $$\bigg|N_{G'}(w)\cap V_{3-i}\cap \mc{N} \cap \mc{H}\bigg((1-\alpha)\cdot |\phi(f)|\bigg)\bigg|\ge \left(\frac{1}{4}-2\sqrt{\alpha}-10\cdot \varepsilon^{1/6}\right)n,$$
    \item[(iii)] every vertex in $\mc{N}$ has degree at most $(\frac{3}{4}+3\sqrt{\varepsilon})n$, 
    \item[(iv)] if $u\in V(G')$, then for each vertex $w\in \mc{N}$ that is $u$-obtuse such that $|\phi(vw)|\ge \frac{9}{10}\cdot |\phi(f_0)|$, we have that 
    $$\bigg|N_{G'}(w) \cap (N(u)\setminus N(v))\cap \mc{H}\bigg((1-\alpha)\cdot (\phi(f)|+|\phi(vw)|) -|\phi(f_0)|\bigg)\bigg| \ge \left(\frac{1}{4} -18\cdot\varepsilon^{1/6}\right)\cdot n.$$
\end{itemize}
\end{claim}
\begin{proofclaim}
By Proposition~\ref{prop:NotTooManyAbnormal}, at most $\sqrt{\alpha}\cdot n$ vertices of $N(u_0)\cap N(v_0)$ are not $\alpha$-normal for $\mc{P}$. By Corollary~\ref{cor:FracCleaning} since $\textrm{ch}_F(f) \le \sqrt{\varepsilon}\cdot |\phi(f)|$ as $f\in \mc{E}_2$, there are most $2\cdot \varepsilon^{1/6}\cdot n$ vertices in $N(u)\cap N(v)$ that $2\cdot \varepsilon^{1/6}$-unclean with respect to $f$. Since $|V(G)\setminus V(G')|=|S|\le 4\cdot \varepsilon^{1/4}\cdot n$, combining all this we find that 
$$|(N(u_0)\cap N(v_0))\setminus \mc{N}_0| \le \sqrt{\alpha}\cdot n +  4\cdot \varepsilon^{1/4}\cdot n + 2\cdot \varepsilon^{1/6}\cdot n,$$
which proves (i).

Next we proceed towards proving (ii). Without loss of generality let $w\in V_1\cap \mc{N}$. Thus for at least $(1-\alpha)\cdot \frac{|N(u)\cap N(v)|}{2}$ vertices $w'\in V_2$ we have that $ww'$ is $\alpha$-normal; recall this means that $ww'\in E^{-}_{\phi}(G)$ and $|\phi(ww')|\ge (1-\alpha)|\phi(f)| \ge \frac{3}{4}\cdot |\phi(f_0)|$ and hence $ww'\in \mc{H}$. Also since $w\in V(G')$, we have that $w \not\in S$ and hence $w$ is incident with at most $3\cdot\varepsilon^{1/4}\cdot n$ edges of $\mc{Q}'=\mc{H}\setminus \mc{E}_2$. Thus we calculate using (i) that
\begin{align*}
\bigg|\bigg\{w'\in &N_{G'}(w)\cap V_{2}\cap \mc{N}: |\phi(ww')|\ge (1-\alpha)\cdot |\phi(f)|\bigg\}\bigg| \\
&\ge (1-\alpha)\cdot \frac{|N(u)\cap N(v)|}{2} - 3\cdot\varepsilon^{1/4}\cdot n - |(N(u)\cap N(v))\setminus \mc{N}|\\
&\ge (1-\alpha)\cdot \left(\frac{1}{4}-\varepsilon\right)n - 3\cdot\varepsilon^{1/4}\cdot n - (\sqrt{\alpha}+6\cdot \varepsilon^{1/6})\cdot n\\
&\ge \left(\frac{1}{4}-2\sqrt{\alpha}-10\cdot \varepsilon^{1/6}\right)n
\end{align*}
which proves (ii).

Since $\varepsilon$ is small enough, it follows from (ii) that for each $w\in \mc{N}$, we have that $|N_{G'}(w)| > 0$ and hence there exists $w'\in N_{G'}(w)$; but then by Claim~\ref{cl:BipartiteStructure} applied to $ww'$, we find that $w$ has degree at most $(\frac{3}{4}+3\sqrt{\varepsilon})n$. This proves (iii).

Finally we proceed towards proving (iv). Let $w\in \mc{N}$ such that $w$ is $u$-obtuse and $vw\in E(G')$. Since $w\in \mc{N}$, we have by definition that $w$ is not $2\cdot \varepsilon^{1/6}$-unclean for $f$. Let $Z$ denote the set of vertices in $N(u)\setminus N(v)$ that are $2\cdot \varepsilon^{1/6}$-unclean for $w$ with respect to $f$. Recall by definition then $|Z| \le 2\cdot\varepsilon^{1/6}\cdot n$. Recall also that $z\in Z$ means either that $wz\not\in E(G)$ or $\sum_{e\in E^-(G)\cap \{wz,uz\}} |\phi(e)| < (1-\alpha)\cdot \phi(uw)$. 

Since $vw$ is negative as $w$ is $u$-obtuse and by assumption $|\phi(vw)|\ge \frac{9}{10}\cdot |\phi(f_0)|$, we find as $\phi$ is triangle-covered that
$$\phi(uw) \ge |\phi(vw)|+|\phi(f)| \ge \frac{9}{5} \cdot |\phi(f_0)|.$$

Let $Y$ denote the set of vertices $y\in N(u)\setminus N(v)$ such that $uy\in \mc{Q}'$ or $wy\in \mc{Q}'$. Since $u,w\in \mc{N}$, we have that $|Y|\le 8\cdot\varepsilon^{1/6}\cdot n$. Consider $w'\in (N(u)\setminus N(v))\setminus (S\cup Y\cup Z)$. 

We claim that $ww'\in E(G')\cap \mc{H}\bigg((1-\alpha)\cdot (\phi(f)|+|\phi(vw)|) -|\phi(f_0)|\bigg)$. Since $w'\not\in S$, we have $w'\in V(G')$. Since $w'\not\in Y$, we find that $uy,wy\not\in \mc{Q}'$ and hence if $uy\in \mc{H}$, then $uy\in \mc{H}\cap \mc{E}_2$ and so $|\phi(uy)|\le |\phi(f_0)|$ (and similarly for $wy$). On the other hand since $w'\not\in Z$, we find that $ww'\in E(G)$ and furthermore that $$\sum_{e\in E^-(G)\cap \{ww',uw'\}} |\phi(e)| \ge (1-\alpha)\cdot \phi(uw) \ge (1-\alpha)\cdot \frac{9}{5} \cdot |\phi(f_0)| > \frac{7}{4}\cdot |\phi(f_0),$$
where the last inequality follows since $\varepsilon$ (and hence $\alpha$) is small enough. Thus it follows from our comments above that $ww',uw'\in E^-(G)$ and furthermore since $|\phi(uw')|\le |\phi(f_0)|$ as noted above, then we have that 
\begin{align*}
|\phi(ww')| &\ge  (1-\alpha)\cdot|\phi(uw)| - |\phi(uw')| \\
&\ge (1-\alpha)\cdot (\phi(f)|+|\phi(vw)|) -|\phi(f_0)|\\
&\ge (1-\alpha) \frac{4}{5}\cdot |\phi(f_0)|\\
&\ge \frac{3}{4}|\phi(f_0)|,    
\end{align*}
where the last inequality follows since $\alpha$ is small enough and hence $ww'\in E(G')\cap \mc{H}\bigg((1-\alpha)\cdot (\phi(f)|+|\phi(vw)|) -|\phi(f_0)|\bigg)$ as claimed. 

Note that since $|N(u)\cap N(v)|\le (\frac{1}{2}+3\sqrt{\varepsilon})n$ and $|N(u)|\ge (\frac{3}{4}-\varepsilon)n$ , we find that $|N(u)\setminus N(v)| \ge (\frac{1}{4}-4\sqrt{\varepsilon})n$. Thus we calculate that
\begin{align*}
\bigg|N_{G'}(w) &\cap (N(u)\setminus N(v))\cap \mc{H}\bigg((1-\alpha)\cdot (\phi(f)|+|\phi(vw)|) -|\phi(f_0)|\bigg)\bigg| \\
&\ge |N(u)\setminus N(v)| - |S|-|Y|-|Z|\\
&\ge \left(\frac{1}{4}-4\sqrt{\varepsilon}\right)n - 4\cdot\varepsilon^{1/6}n - 8\cdot\varepsilon^{1/6}\cdot n - 2\cdot\varepsilon^{1/6}\cdot n\\
&\ge \left(\frac{1}{4} -18\cdot\varepsilon^{1/6}\right)\cdot n,   
\end{align*}
which proves (iv).
\end{proofclaim}

Since $f_0\in \mc{E}_2$, it follows that Claim~\ref{cl:BipartiteStructure} applies to $f_0$ and hence we find that $|N(u_0)\cap N(v_0)|\le (\frac{3}{4}+3\sqrt{\varepsilon})n$ and $f_0$ is $\alpha$-balanced-bipartite. Thus by definition there exists an $\alpha$-balanced-bipartition $\mc{P}_0=(A_0,B_0)$ of $N(u_0)\cap N(v)$. We let $\mc{N}_0$ denote the set of vertices in $N(u_0)\cap N(v_0) \cap V(G')$ that are $\alpha$-normal for $\mc{P}_0$ and are not $2\cdot \varepsilon^{1/6}$-unclean with respect to $f_0$. 

For the rest of the proof let $\beta := 6\sqrt{\alpha}$.

\begin{claim}\label{cl:NotSmallInterior}
$\max\{e(G[A_0], e(G[B_0])\} \ge \beta\cdot n^2$.    
\end{claim}
\begin{proofclaim}
Suppose not, that is $e(G[A_0]),e(G[B_0]) < \beta\cdot n^2$. By Claim~\ref{cl:BipartiteStructure2}, we find that every vertex of $\mc{N}_0$ has degree at most $(\frac{3}{4}+3\sqrt{\varepsilon})n$. Recall that $|N(u_0)\cap N(v_0)| \le (\frac{1}{2}+2\sqrt{\varepsilon})n$ and hence $|N(u_0)\cap N(v_0)|^2 \le (\frac{1}{4} + 6\sqrt{\varepsilon})n^2$. Thus we calculate that 
$$e(G[A_0\cup B_0]) \le |A_0|\cdot |B_0| + 2\beta\cdot n^2 \le \frac{|N(u_0)\cap N(v_0)|^2}{4} + 2\beta\cdot n^2 \le \left(\frac{1}{16} + 2\sqrt{\varepsilon}+2\beta\right)\cdot n^2.$$

We now consider the bipartition $(A,B)$ of $V(G)$ with $A:= A_0\cup B_0$ and $B:=V(G)\setminus A$. Note that $|A|=|N(u_0)\cap N(v_0)|\ge (\frac{1}{2}-2\varepsilon)n$. Since $\delta(G)\ge (\frac{3}{4}-\varepsilon)n$, it follows that
\begin{align*}
e(G[A,B]) &\ge \delta(G)\cdot |A| - 2\cdot e(G[A])\\
&\ge \left(\frac{3}{4}-\varepsilon\right)n\cdot \left(\frac{1}{2}-2\varepsilon\right)n -2\left(\frac{1}{16} + 2\sqrt{\varepsilon}+2\beta\right)\cdot n^2 \\
&\ge \left(\frac{1}{4} - 6\sqrt{\varepsilon}-4\beta\right)\cdot n^2\\
&\ge \left(\frac{1}{4}-\sigma\right)n^2,
\end{align*}
where the last inequality follows since $\varepsilon$ is small enough.

Let $T_{A_0}:=\{a\in A_0: |N_G(a)\cap A_0| > 16\beta\cdot n\}$ and $T_{B_0}:= \{b\in B_0: |N_G(b)\cap B_0| > 16\beta\cdot n\}$. Since $e(G[A_0])\le \beta\cdot n^2$, we find $|T_{A_0}|\le \frac{n}{8}$; similarly since $e(G[B_0])\le \beta\cdot n^2$, we find $|T_{B_0}|\le \frac{n}{8}$. Note that 
\begin{align*}
|(A_0\cap \mc{N}_0)\setminus T_{A_0}| &\ge |A_0| - |(N(u_0)\cap N(v_0))\setminus \mc{N}_0| - |T_{A_0}| \\
&\ge \frac{1}{2} \cdot |N(u_0)\cap N(v_0)| - (\sqrt{\alpha}+6\cdot\varepsilon^{1/6})\cdot n - 4\cdot\varepsilon^{1/4}\cdot n-\frac{n}{8}    \\
&\ge \frac{1}{2} \cdot \left(\frac{1}{2}-2\varepsilon\right)n - \sqrt{\alpha}\cdot n - 4\cdot\varepsilon^{1/4}\cdot n - \frac{n}{8}\\
&> 0,
\end{align*}
where in the first inequality we used Claim~\ref{cl:BipartiteStructure2}(i) and the last inequality follows since $\varepsilon$ (and hence $\alpha$) are small enough.

Thus there exists $u_0'\in (A_0\cap \mc{N}_0)\setminus T_{A_0}$. Let 
$$U_0 := \bigg|\bigg\{w'\in N_{G'}(u_0')\cap B_0\cap \mc{N}_0: |\phi(u_0w')|\ge (1-\alpha)\cdot |\phi(f_0)|\bigg\}\bigg|.$$ 
Since $u_0'\in A_0\cap \mc{N}_0$, we have by Claim~\ref{cl:BipartiteStructure2}(ii) that
\begin{align*}
|U_0\setminus T_{B_0}| &\ge |U_0| - |T_{B_0}|\\
&\ge\left(\frac{1}{4}-2\sqrt{\alpha}-10\cdot \varepsilon^{1/6}\right)n - \frac{n}{8}\\
&> 0,
\end{align*}
where the last inequality follows since $\varepsilon$ (and hence $\alpha$) are small enough.

Then there exists $v_0'\in U_0\setminus T_{B_0}$. Since $v_0'\not\in T_{B_0}$ we have that $|N_G(v_0')\cap B_0| \le 16\beta\cdot n$. Since $v_0'\in N_{G'}(u_0')$, we find that $u_0'v_0'\in \mc{E}_2\cap \mc{H}$. Let $f_0'=u_0'v_0'$. Since $v_0'\in U_0$, we have that $|\phi(f_0')| \ge (1-\alpha)\cdot |\phi(f_0)| \ge \frac{9}{10}\cdot |\phi(f_0)|$ where the last inequality follows since $\varepsilon$ is small enough. Now we calculate that
$$|N(u_0')\cap N(v_0')\cap B| \ge |N(u_0')\cap N(v_0')| - |N(u_0')\cap A_0| - |N(v_0')\cap B_0| \ge \left(\frac{1}{2}-2\varepsilon-32\beta\right)n.$$ 
Notice that $|B| = n-|A| \le (\frac{1}{2}+2\varepsilon)n$. Hence
$$|B\setminus (N(u_0')\cap N(v_0'))| = |B|- |N(u_0')\cap N(v_0')\cap B|\le (4\varepsilon+32\beta)n.$$
By Claim~\ref{cl:BipartiteStructure} since $f_0'\in \mc{E}_2\cap \mc{H}$, we find that $|N(u_0')\cap N(v_0')| \le (\frac{1}{2}+3\sqrt{\varepsilon})n$ and that $f_0'$ is $\alpha$-balanced-bipartite. Let $\mc{P}_0'=(A_0',B_0')$ be an $\alpha$-balanced-bipartition for $f_0'$. Let $\mc{N}_0'$ be the set of vertices of $N(u_0')\cap N(v_0')\cap V(G')$ that are $\alpha$-normal for $\mc{P}_0'$ and are not $2\cdot \varepsilon^{1/6}$-unclean with respect to $f_0$. By Claim~\ref{cl:BipartiteStructure2}(i) since $f_0'\in \mc{E}_2\cap \mc{H}$ and $|\phi(f_0')|\ge \frac{9}{10}\cdot |\phi(f_0)|$, we have that $|(N(u_0')\cap N(v_0'))\setminus \mc{N}_0'| \le (\sqrt{\alpha}+6\cdot \varepsilon^{1/6})\cdot n$. Similarly by Claim~\ref{cl:BipartiteStructure2}(iii), we have that every vertex in $\mc{N}_0'$ has degree at most $(\frac{3}{4}+3\sqrt{\varepsilon})n$. Hence at most $\sqrt{\alpha}\cdot n+4\cdot\varepsilon^{1/4}\cdot n$ of the vertices of $N(u_0')\cap N(v_0')\cap B$ have degree strictly more than $(\frac{3}{4}+3\sqrt{\varepsilon})n$. Hence at most $2(\sqrt{\alpha}+4\cdot\varepsilon^{1/4})+4\varepsilon+32\beta)\cdot n \le \sigma\cdot n$ vertices of $G$ have degree more than $(\frac{3}{4}+3\sqrt{\varepsilon})n$ (since $\varepsilon$ is small enough for this inequality to hold) and hence at most $\sigma\cdot n$ vertices of $G$ have degree more than $(\frac{3}{4}+\sigma)n$. Thus $G$ is as desired, a contradiction.
\end{proofclaim}

By Claim~\ref{cl:NotSmallInterior}, we have that $\max\{e(G[A_0]),e(G[B_0])\} \ge \beta\cdot n^2$. We assume without loss of generality that $e(G[A_0]) \ge \beta\cdot n^2$.
 
\begin{claim}\label{cl:F1}
There exists a negative edge $f_1=u_1v_1 \in E(G')$ with $|\phi(f_1)|\ge (1-\alpha)|\phi(f_0)|$ such that there exist at least $\frac{\beta}{2}\cdot n$ $u_1$-obtuse vertices $w$ with $|\phi(wv_1)|\ge (1-\alpha)\cdot |\phi(f_0)|$.
\end{claim}
\begin{proofclaim}
Let $U:= \{w\in A_0: |N_G(w)\cap A_0|\ge \beta\cdot n\}$. Since $e(G[A_0])\ge \beta\cdot n^2$, it follows that $|U|\ge \beta\cdot n$ as otherwise $2\cdot e(G[A_0]) = \sum_{v\in A_0} d_{G[A_0]}(v) \le |U|\cdot n + \beta\cdot n\cdot |A_0| < 2\beta\cdot n^2$, a contradiction. By Claim~\ref{cl:BipartiteStructure2}(i) applied to $f_0$, we have that $|(N(u)\cap N(v))\setminus \mc{N}_0|\le (\sqrt{\alpha}+6\cdot \varepsilon^{1/6})n$. 

Thus there exists $u_1\in U\cap \mc{N}_0$ since $\beta \ge 2\sqrt{\alpha}$ and $\varepsilon$ is small enough. Let 
$$U_1:= \bigg|\bigg\{w'\in N_{G'}(u_1)\cap B_0\cap \mc{N}_0: |\phi(u_1w')|\ge (1-\alpha)\cdot |\phi(f_0)|\bigg\}\bigg|$$
By Claim~\ref{cl:BipartiteStructure2}(ii) applied to $f_0$ and $u_1\in A_0\cap \mc{N}_0$, we have that 
$$|U_1|\ge \left(\frac{1}{4}-2\sqrt{\alpha}-10\cdot \varepsilon^{1/6}\right)n > 0,$$
where the last inequality follows since $\varepsilon$ is small enough. Thus there exists a vertex $v_1\in N_{G'}(u_1)\cap B_0\cap \mc{N}_0$ such that $|\phi(u_1v_1)|\ge (1-\alpha)|\phi(f_0)|$. Let $f_1:=u_1v_1$. Let 
$$V_1:= N_{G'}(v_1)\cap A_0\cap \mc{N}_0\cap \mc{H}\bigg((1-\alpha)\cdot |\phi(f_0)|\bigg).$$
By Claim~\ref{cl:BipartiteStructure2}(ii) applied to $f_0$ and $v_1\in B_0\cap \mc{N}_0$, we have that 
$$|V_1|\ge \left(\frac{1}{4}-2\sqrt{\alpha}-10\cdot \varepsilon^{1/6}\right)n > 0,$$
where the last inequality follows since $\varepsilon$ is small enough. Since $V_1\subseteq A_0$ and $|A_0|=\frac{|N(u_0)\cap N(v_0)|}{2} \le (\frac{1}{4}+\sqrt{\varepsilon})n$, we find that
$$|A_0\setminus V_1| \le (2\sqrt{\alpha}+11\cdot \varepsilon^{1/6})n $$
But now we calculate that 
$$|V_1\cap N_{G}(w)| \ge |N_G(w)\cap A_0| - |A_0\setminus V_1| \ge \beta\cdot n - \left(2\sqrt{\alpha}+11\cdot \varepsilon^{1/6}\right)n \ge \frac{\beta}{2}\cdot n,$$
where the last inequality follows since $\beta\ge 6\sqrt{\alpha}$ and $\varepsilon$ is small enough. Since $\phi$ is triangle-covered and $\alpha$ is small enough, it follows that for every $w'\in V_1\cap N_{G}(w)$, we have that $u_1w'$ is a positive edge and hence $w'$ is $u_1$-obtuse with respect to $f_1$. The claim follows.
\end{proofclaim}

Here is a useful claim that will help us to finish the proof.

\begin{claim}\label{cl:NotTooMuchBlueDegree}
If $w\in V(G')$, then $\Delta(G[N_{G'}(w)])\le (\frac{1}{4} + \sqrt{\alpha})n$ and every vertex $w'\in N_{G'}(w)$ has degree in $G$ at most $(\frac{3}{4}+3\sqrt{\varepsilon})n$.   
\end{claim}
\begin{proofclaim}
Let $w'\in N_{G'}(w)$. By definition of $G'$, we have that $ww'\in \mc{E}_2\cap \mc{H}$. Hence by Claim~\ref{cl:BipartiteStructure}, we find that $|N(w)\cap N(w')|\le (\frac{1}{2}+3\sqrt{\varepsilon})n$ and $d_G(w')\le (\frac{3}{4}+3\sqrt{\varepsilon})n$ and that $ww'$ is $\alpha$-balanced bipartite. By Proposition~\ref{prop:BalancedObtuse}, we have that the number of vertices that are $w$-obtuse with respect to $ww'$ is at most 
$$\left(\frac{1}{2}+\sqrt{\alpha}\right)\cdot |N(u)\cap N(v)| \le \left(\frac{1}{2}+\sqrt{\alpha}\right)\cdot\left(\frac{1}{2}+3\sqrt{\varepsilon}\right)n \le \left(\frac{1}{4} + \sqrt{\alpha}\right)n,$$
where the last inequality follows since $\varepsilon$ is small enough. Since $\phi$ is triangle-covered, we find for every $w''\in N_G(w')\cap N_{G'}(w)$ that $w'w''$ is a positive edge and hence $w''$ is $w$-obtuse with respect to $ww'$. The claim follows.
\end{proofclaim}

Let $f_1$ be as guaranteed by Claim~\ref{cl:F1}. Since $\alpha \le \frac{1}{4}$ as $\varepsilon$ is small enough, we find by Claim~\ref{cl:BipartiteStructure} that $|N(u_1)\cap N(v_1)|\le \left(\frac{1}{2}+3\sqrt{\varepsilon}\right)n$ and that $f_1$ is $\alpha$-balanced-bipartite. Let $\mc{P}_1=(A_1,B_1)$ be an $\alpha$-balanced-bipartition of $N(u_1)\cap N(v_1)$. Let $\mc{N}_1$ denote the set of vertices of $N(u_1)\cap N(v_1) \cap V(G')$ that are $\alpha$-normal for $\mc{P}_1$ and are not $2\cdot \varepsilon^{1/6}$-unclean with respect to $f_1$.

\begin{claim}\label{cl:u2}
There exists $u_2 \in \mc{N}_1$ that is $u_1$-obtuse with respect to $f_1$ such that $|\phi(u_2v_1)|\ge (1-\alpha)|\phi(f_0)|$. 
\end{claim}
\begin{proofclaim}
Since $\alpha$ is small enough, we have that $|\phi(f_1)|\ge (1-\alpha)\cdot |\phi(f_0)|\ge \frac{9}{10}\cdot |\phi(f_0)|$. Hence by Claim~\ref{cl:BipartiteStructure2}(i), we find that $|(N(u_1)\cap N(v_1))\setminus \mc{N}_1|\le (\sqrt{\alpha}+6\cdot \varepsilon^{1/6})\cdot n$. By Claim~\ref{cl:F1}, there exists at least $\frac{\beta}{2}\cdot n$ $u_1$-obtuse vertices $w$ with $|\phi(wv_1)|\ge (1-\alpha)\cdot |\phi(f_0)|$. Since $\frac{\beta}{2} = 3\sqrt{\alpha}\cdot n > (\sqrt{\alpha}+6\cdot \varepsilon^{1/6})\cdot n$ as $\varepsilon$ is small enough, we find that there exists $u_2 \in \mc{N}_1$ that is $u_1$-obtuse with respect to $f_1$ such that $|\phi(u_2v_1)|\ge (1-\alpha)\cdot |\phi(f_0)|$ as desired.
\end{proofclaim}

Let $u_2$ be as guaranteed by Claim~\ref{cl:u2}. We assume without loss of generality that $u_2\in A_1$. Let $f_2:= u_2v_1$. Since $\alpha \le \frac{1}{4}$ as $\varepsilon$ is small enough, we find by Claim~\ref{cl:BipartiteStructure} that $|N(u_2)\cap N(v_1)|\le \left(\frac{1}{2}+3\sqrt{\varepsilon}\right)n$ and that $f_2$ is $\alpha$-balanced-bipartite. Let $\mc{P}_2=(A_2,B_2)$ be an $\alpha$-balanced-bipartition of $N(u_2)\cap N(v_1)$. Let $\mc{N}_2$ denote the set of vertices of $N(u_2)\cap N(v_1) \cap V(G')$ that are $\alpha$-normal for $\mc{P}_2$ and are not $2\cdot \varepsilon^{1/6}$-unclean with respect to $f_2$.

Since $\alpha$ is small enough, we have that $|\phi(u_2v_1)|\ge (1-\alpha)\cdot |\phi(f_0)| \ge \frac{9}{10}\cdot |\phi(f_0)|$. Thus by Claim~\ref{cl:BipartiteStructure2})(iii) applied to $f_1$ and $u_2\in \mc{N}_1$ since $u_1\in V(G')$, we find that 
$$\bigg|N_{G'}(u_2) \cap (N(u_1)\setminus N(v_1)) \cap \mc{H}\bigg( (1-4\alpha)\cdot |\phi(f_0)|\bigg)\bigg| \ge \left(\frac{1}{4} -18\cdot\varepsilon^{1/6}\right)\cdot n,$$
where we used that $(1-\alpha)\cdot (|\phi(f_1)| + |\phi(u_2v_1)|) - |\phi(f_0)| \ge (1-\alpha)(2-2\alpha)\cdot |\phi(f_0)|-|\phi(f_0)|\ge (1-4\alpha)\cdot |\phi(f_0)|$.

Let $B:= B_1\cup (N(u_1)\setminus N(v_1))$ and let $A:= V(G)\setminus B = A_1\cup (N(v_1)\setminus N(u_1))$. Since $u_2\in \mc{N}_1$, it follows from Claim~\ref{cl:BipartiteStructure2}(ii) that  
$$\bigg|N_{G'}(u_2)\cap B_1\cap \mc{H}\bigg((1-2\alpha)\cdot |\phi(f_0)|\bigg)\bigg|\ge \left(\frac{1}{4}-2\sqrt{\alpha}-10\cdot \varepsilon^{1/6}\right)n,$$
where we used that $(1-\alpha)\cdot |\phi(f_1)| \ge (1-\alpha)^2\cdot |\phi(f_0)| \ge (1-2\alpha)\cdot |\phi(f_0)|$.
Combining these yields that
\begin{align*}
\bigg|N_{G'}(u_2)&\cap B \cap \mc{H}\bigg((1-4\alpha)\cdot |\phi(f_0)|\bigg)\bigg| \\
&\ge \bigg|N_{G'}(u_2)\cap B_1\cap \mc{H}\bigg((1-2\alpha)\cdot |\phi(f_0)|\bigg)\bigg|+\bigg|N_{G'}(u_2)\cap (N(u_1)\setminus N(v_1))\cap \mc{H}\bigg((1-4\alpha)\cdot |\phi(f_0)|\bigg)\bigg|\\
&\ge \left(\frac{1}{2}-2\sqrt{\alpha}-18\cdot \varepsilon^{1/6}\right)n.   
\end{align*}
Since $u_2\in V(G')$ as $u_2\in \mc{N}_1$, it follows by Claim~\ref{cl:NotTooMuchBlueDegree} that
$$\Delta(G[N_{G'}(u_2)\cap B])\le \left(\frac{1}{4} + \sqrt{\alpha}\right)n$$
and every vertex $w'\in N_{G'}(u_2)\cap B$ has degree in $G$ at most $(\frac{3}{4}+3\sqrt{\varepsilon})n$. 
Note that 
$$|B|\le |N(u_1)\setminus N(v_1)| + |B_1| \le \frac{|N(u_1)\cap N(v_1)|}{2} + n-\delta(G) \le \left(\frac{1}{4}+\sqrt{\varepsilon}\right)n + \left(\frac{1}{4}+\varepsilon\right)n \le \left(\frac{1}{2}+2\sqrt{\varepsilon}\right)n.$$
Moreover
\begin{align*}
e(G[A,B]) &\ge  |N_{G'}(u_2)\cap B|\cdot \left( \left(\frac{3}{4}-\varepsilon\right)n - |B\setminus N_{G'}(u_2)| - \Delta(G[N_{G'}(u_2)\cap B|)\right)\\
&\ge \left(\frac{1}{2}-2\sqrt{\alpha}-18\cdot \varepsilon^{1/6}\right)n \cdot \left(\left(\frac{3}{4}-\varepsilon\right) - \left(\frac{1}{2}+2\sqrt{\varepsilon}\right) + \left(\frac{1}{2}-2\sqrt{\alpha}-18\cdot \varepsilon^{1/6}\right) - \left(\frac{1}{4} + \sqrt{\alpha}\right)\right)n \\
&\ge \left(\frac{1}{2}-2\sqrt{\alpha}-18\cdot \varepsilon^{1/6}\right)\cdot \left(\frac{1}{2}-3\sqrt{\alpha}-21\cdot\varepsilon^{1/6} \right)\cdot n^2\\
&\ge \left(\frac{1}{4}-5\sqrt{\alpha}-39\cdot \varepsilon^{1/6}\right)\cdot n^2\\
&\ge \left(\frac{1}{4}-\sigma\right)\cdot n^2,
\end{align*}
where the last inequality follows since $\varepsilon$ is small enough.

Finally we proceed to show most vertices in $A$ have degree at most $(\frac{3}{4}+3\sqrt{\varepsilon})n$ as follows.

\begin{claim}
There exists $v_2\in \mc{N}_2$ that is $v_1$-obtuse with respect to $f_2$ such that $|\phi(u_2v_2)|\ge (1-4\alpha)\cdot |\phi(f_0)|$.    
\end{claim}
\begin{proofclaim}
Let $V_2:= N_{G'}(u_2)\cap B \cap \mc{H}\bigg((1-4\alpha)\cdot |\phi(f_0)|\bigg)$. Recall from before that 
$|V_2| \ge \left(\frac{1}{2}-2\sqrt{\alpha}-18\cdot \varepsilon^{1/6}\right)n.$   Since $d_G(v_1) \ge (\frac{3}{4}-\varepsilon)n$, we find that
$$|V_2\cap N_G(v_1)| \ge |V_2|+d_G(v_1) - n \ge \left(\frac{1}{4}-2\sqrt{\alpha}-19\cdot \varepsilon^{1/6}\right)n.$$
Note that since $\phi$ is triangle-covered, we have that every vertex of $V_2\cap N_G(v_1)$ is in $N(u_2)\cap N(v_1)$ and is $v_1$-obtuse. Since $\alpha$ is small enough, we have that $|\phi(f_2)|\ge (1-\alpha)\cdot |\phi(f_0)|\ge \frac{9}{10}\cdot |\phi(f_0)|$. Hence by Claim~\ref{cl:BipartiteStructure2}(i), we find that $|(N(u_2)\cap N(v_1))\setminus \mc{N}_2|\le (\sqrt{\alpha}+6\cdot \varepsilon^{1/6})\cdot n$. Thus
$$|V_2\cap N_G(v_1)\cap \mc{N}_2| \ge |V_2\cap N_G(v_1)| - |(N(u_2)\cap N(v_1))\setminus \mc{N}_2| \ge \left(\frac{1}{4}-3\sqrt{\alpha}-25\cdot \varepsilon^{1/6}\right)n > 0,$$
where the last inequality follows since $\varepsilon$ (and hence $\alpha$) is small enough. But then there exists $v_2\in V_2\cap N_G(v_1)\cap \mc{N}_2$ as desired.
\end{proofclaim}

Let $v_2$ be as guaranteed by Claim~\ref{cl:u2}. We assume without loss of generality that $v_2\in B_2$. Since $\alpha$ is small enough, we have that $|\phi(u_2v_2)|\ge \frac{9}{10}\cdot |\phi(f_0)|$. Thus by Claim~\ref{cl:BipartiteStructure2}(iii) applied to $f_2$ and $v_2\in \mc{N}_2$ since $v_1\in V(G')$, we find that 
$$|N_{G'}(v_2) \cap (N(v_1)\setminus N(u_2))| \ge \left(\frac{1}{4} -18\cdot\varepsilon^{1/6}\right)\cdot n.$$
Since $v_2\in \mc{N}_2$, it follows from Claim~\ref{cl:BipartiteStructure2}(ii) that  
$$|N_{G'}(v_2)\cap A_2|\ge \left(\frac{1}{4}-2\sqrt{\alpha}-10\cdot \varepsilon^{1/6}\right)n.$$
Since $A_2\subseteq N(u_2)\cap N(v_1)$, we find that 
$$|N_{G'}(v_2)| \ge |N_{G'}(u_2)\cap A_1|+|N_{G'}(u_2)\cap (N(v_1)\setminus N(u_1))| \ge \left(\frac{1}{2}-2\sqrt{\alpha}-28\cdot \varepsilon^{1/6}\right)n.$$
Recall that $E(G')\subseteq \mc{E}_2\cap \mc{H}\subseteq E^-(G)$. Thus since $\phi$ is triangle-covered and $u_2v_2\in E^-(G)$, we find that $N_{G'}(u_2)\cap N_{G'}(v_2)=\emptyset$. Hence
\begin{align*}
|N_{G'}(v_2)\cap B| &\le |B|-|N_{G'}(u_2)\cap B|   \\
&\le \left(\frac{1}{2}+2\sqrt{\varepsilon}\right)n - \left(\frac{1}{2}-2\sqrt{\alpha}-18\cdot \varepsilon^{1/6}\right)n\\
&\le (2\sqrt{\alpha}+20\cdot \varepsilon^{1/6})n.
\end{align*}
Combining these yields that
\begin{align*}
|N_{G'}(v_2)\cap A| &= |N_{G'}(v_2)|-|N_{G'}(v_2)\cap B|\\
&\ge \left(\frac{1}{2}-2\sqrt{\alpha}-28\cdot \varepsilon^{1/6}\right)n - (2\sqrt{\alpha}+20\cdot \varepsilon^{1/6})n  \\
&\ge \left(\frac{1}{2}-4\sqrt{\alpha}-48\cdot \varepsilon^{1/6}\right)n
\end{align*}
Since $v_2\in V(G')$ as $v_2\in \mc{N}_2$, it follows by Claim~\ref{cl:NotTooMuchBlueDegree} that every vertex $w'\in N_{G'}(v_2)\cap A$ has degree in $G$ at most $(\frac{3}{4}+3\sqrt{\varepsilon})n$. Since $(N_{G'}(u_2)\cap B) \cap (N_{G'}(v_2)\cap A)=\emptyset$ and every vertex of $N_{G'}(u_2)\cap B$ and every vertex of $N_{G'}(v_2)\cap A$ has degree at most $(\frac{3}{4}+3\sqrt{\varepsilon})n$ (which is at most $(\frac{3}{4}+\sigma)n$), we find that the number of vertices of $G$ with degree at most $(\frac{3}{4}+\sigma)n$ is at least  
$$|N_{G'}(u_2)\cap B| + |N_{G'}(v_2)\cap A| \ge \left(\frac{1}{2}-2\sqrt{\alpha}-18\cdot \varepsilon^{1/6}\right)n + \left(\frac{1}{2}-4\sqrt{\alpha}-48\cdot \varepsilon^{1/6}\right)n \ge (1-\sigma)n,$$
where the last inequality follows since $\varepsilon$ (and hence $\alpha$) are small enough. Thus $G$ is as desired, a contradiction.
\end{proof}

We are now prepared to prove Theorem~\ref{thm:fracNWstability}.

\begin{proof}[Proof of Theorem~\ref{thm:fracNWstability}]
Follows from Theorem~\ref{thm:Charge2Stability} and Lemma~\ref{lem:FarkasK3}.
\end{proof}

\section{A Robust Fractional Decomposition via Property Testing and Inheritance}

In this section, we prove Theorem~\ref{thm:RegularFracStability}. First we recall our fractional stability theorem but to make our proof easier we make the following definition (which we also use extensively in Part III).

\begin{definition}
Let $\sigma \ge 0$. A graph $G$ on $n$ vertices is \emph{$\sigma$-extremal} if both of the following hold:    
\begin{itemize}
    \item[(i)] at least $(1-\sigma)n$ vertices of $G$ have degree at most $(\frac{3}{4}+\sigma)n$, and
    \item[(ii)] the maximum cut of $G$ is at least $(1-\sigma)\cdot \frac{n^2}{4}$.
\end{itemize}     
\end{definition}

With this definition, we may restate our fractional stability theorem as follows.

\begin{thm}[Fractional Stability Restated]\label{thm:NWstability2} 
For each sufficiently small $\sigma > 0$, there exists $\varepsilon > 0$ such that the following holds for all sufficiently large $n$: If $G$ is a $K_3$-divisible graph on $n$ vertices with $\delta(G) \ge \left(\frac{3}{4}-\varepsilon\right)n$ and $G$ has no fractional $K_3$-decomposition, then $G$ is $\sigma$-extremal.
\end{thm}

For our proof, we need a result from property testing. To that end, we recall the following definition from the literature.

\begin{definition}
A graph parameter $\gamma$ (i.e. a function $\gamma$ for all graphs $G$ to real numbers) is \emph{estimable} if for all reals $\varepsilon, \rho > 0$ there exists $k_0> 0$ such that for all $n \ge k \ge k_0$, we have that if $G$ is a graph on $n$ vertices and $G[k]$ is an induced subgraph of $G$ on $k$ vertices chosen uniformly at random, then 
$$\Prob{ |\gamma(G[k]) -\gamma(G)| > \varepsilon }\le \rho.$$
\end{definition}

The key fact that we need is the following but first a definition. For a graph $G$, we let $\textrm{MaxCut}(G)$ denote the number of edges in a maximum cut of $G$. We let $t_{\textrm{cut}}(G) := \frac{\textrm{MaxCut}(G)}{v(G)^2}$. 

\begin{thm}\label{thm:MaxCutEstimable}
$t_{\textrm{cut}}$ is estimable.    
\end{thm}

Theorem~\ref{thm:MaxCutEstimable} was originally established in some of the earliest papers on property testing, namely by Arora, Karger, and Karpinski~\cite{AKK95} from 1995 and independently Goldreich, Goldwasser, and Ron~\cite{GGR96} from 1996 (see~\cite{AKK99} and~\cite{GGR98} for the journal versions of both articles respectively).  Since then, general characterizations of estimable parameters have been established from which it is fairly straightforward to see that max cut is estimable (indeed, it is used as one of the canonical parameters having this property). See Example 15.2 in the book by Lov\'asz~\cite{Lbook} on graph limits for one of these more straightforward derivations.

We also need a lemma regarding degrees in a sampled subgraph as follows. 

\begin{lem}\label{lem:Sampled}
For all reals $\varepsilon,\rho  > 0$, there exists an integer $k_0 > 0$ such that the following holds for all integers $n\ge k\ge k_0$: if $G$ is a graph on $n$ vertices and $S$ is a subset of $V(G)$ of size $k$ chosen uniformly at random, then for each $e\in E(G)$ we have that
$$\Prob{\exists~v \in S \textrm{ such that } \frac{d_{G[S]}(v)}{k} < \frac{d_G(v)}{n} - \varepsilon~\bigg|~ e\subseteq S} < \rho.$$ 
\end{lem}
\begin{proof}
Follows from the Chernoff bounds and the union bound.
\end{proof}

Now we use Theorem~\ref{thm:MaxCutEstimable} and Lemma~\ref{lem:Sampled} to prove the following lemma.

\begin{lem}\label{lem:Sampled2}
For all reals $\varepsilon,\rho, \sigma >0$, there exists an integer $k_0 > 0$ such that the following holds for all integers $n\ge k\ge k_0$: Let $G$ be a graph on $n$ vertices with $\delta(G)\ge (\frac{3}{4}-\varepsilon)n$ that is not $\sigma$-extremal. If $S$ is a subset of $V(G)$ of size $k$ chosen uniformly at random, then for each $e\in E(G)$ we have that
$$\Prob{\delta(G[S])\ge \left(\frac{3}{4}-2\varepsilon\right)k \text{ and $G[S]$ is not $\sigma/2$-extremal}~\bigg|~ e\subseteq S} \ge 1-\rho.$$ 
\end{lem}
\begin{proof}
We choose $k_0$ large enough as needed throughout the proof. By Lemma~\ref{lem:Sampled} since $k_0$ is large enough and $\delta(G)\ge (\frac{3}{4}-\varepsilon)n$, we find that 
$$\Prob{\delta(G[S]) < \left(\frac{3}{4}-2\varepsilon\right)k~\bigg|~ e\subseteq S} < \frac{\rho}{2}.$$ 
Since $G$ is not $\sigma$-extremal, $G$ violates condition (i) or condition (ii) of the definition. 

First suppose $G$ violates condition (i), that is, there are fewer than $(1-\sigma)n$ vertices of $G$ with degree at most $(\frac{3}{4} + \sigma)$n or equivalently there are strictly more than $\sigma n$ vertices of $G$ with degree more than $(\frac{3}{4}+\sigma)n$.  Let $B:= \{v\in V(G): d_G(v) > (\frac{3}{4}+\sigma)n\}$. Thus $|B| > \sigma n$. By Lemma~\ref{lem:Sampled} since $k_0$ is large enough, we find that
$$\Prob{G[S]\text{ is $\sigma/2$-extremal}~\bigg|~e\subset S} \le \Prob{\exists v\in B \text{ such that } \frac{d_{G[S]}}{k} < \frac{d_G(v)}{n} - \frac{\sigma}{2}~\bigg|~e\subset S} < \frac{\rho}{2}.$$
Combining the two above inequalities then yields the desired result.

So we assume $G$ violates condition (ii), that is, $\textrm{MaxCut}(G) < (1-\sigma)\cdot \frac{n^2}{4}$ and hence $t_{\textrm{cut}}(G) < \frac{1-\sigma}{4}$. Let $u,v\in V(G)$ such that $e=uv$. Let $G':= G\setminus \{u,v\}$. Let $S'$ be a subset of $V(G')$ of size $k-2$ chosen uniformly at random. Note that $S'\cup \{u,v\}$ has the same probability distribution as $S$ conditioned on the even $e\subseteq S$. Since $n$ is large enough, it follows that $t_{\textrm{cut}}(G') < \frac{1 - 0.9\sigma}{4}$. By Theorem~\ref{thm:MaxCutEstimable}, we find that 
$$\Prob{t_{\textrm{cut}}(G'[S']) < \frac{1-0.7\sigma}{4}} < \frac{\rho}{2}.$$
Note that if $t_{\textrm{cut}}(G'[S']) < \frac{1-0.7\sigma}{4}$, then since $n$ is large enough, we have that $t_{\textrm{cut}}(G[S'\cup \{u,v\}]) < \frac{1-0.5\sigma}{4}$. Thus $\Prob{G[S]\text{ is $\sigma/2$-extremal}~\bigg|~e\subset S} < \frac{\rho}{2}$ and combining this with our first inequality yields the desired result. 
\end{proof}

Lemma~\ref{lem:Sampled2} would be enough when combined with Theorem~\ref{thm:fracNWstability} to show that if $G$ is a graph on $n$ vertices with $\delta(G)\ge \frac{3}{4}n$, then there exists a $C$-low-weight $\rho$-almost fractional triangle decomposition of $G$ (for some fixed integer $C\ge 1$). To show that we can also obtain such a fractional packing that is $\sigma$-seeded (for some $\sigma \in (0,1]$) and furthermore to show that the packing can be further corrected to a fractional triangle decomposition, we require the following lemma which employs ideas from our paper~\cite{DHLP25} which in turn was inspired by ideas of Montgomery in~\cite{montgomery_fractional_2019}.

First let us recall some extra terminology about fractional packings as follows.

\begin{definition}
Let $F$ and $G$ be graphs. If $\phi$ is a fractional $F$-packing of $G$ and $e\in E(G)$, we write $\partial \phi(e) := \sum_{H\subseteq G: e\in E(H)} \phi(H)$.     
\end{definition}

Here is a useful proposition about achieving any target fractional packing for a subset $R$ of $E(G)$ in a graph $G$ provided that $G$ admits a fractional decomposition for $G-R'$ for any $R'\subseteq R$.

\begin{proposition}\label{prop:FracFixing}
Let $F$ and $G$ be graphs. If $R\subseteq E(G)$ such that for each $R'\subseteq R$ we have that $G-R'$ admits a fractional $F$-decomposition, then for any $\varphi:E(R)\rightarrow [0,1]$ there exists a fractional $F$-packing of $G$ with $\partial \phi(e)=1$ for $e\in E(G)\setminus R$ and $\partial \phi(e)=\varphi(e)$ for $e\in R.$ 
\end{proposition}
\begin{proof}
Let $k:= |R|$ and let $e_1,\ldots e_k$ be an enumeration of $R$ such that $\varphi(e_1)\le \varphi(e_2)\le \ldots \le \varphi(e_k)$. Let $R_0:=\emptyset$ and for each $i\in [k]$, let $R_i:=\{e_1,\ldots,e_i\}$. Let $w_0:=\varphi(e_1)$, $w_k:=1-\varphi(e_k)$ and for each $i\in [k-1]$, let $w_i:=\varphi(e_{i+1})-\varphi(e_i)$. By assumption for each $i\in \{0,1,\ldots,k\}$, we have that $G-R_i$ admits a fractional $F$-decomposition $\phi_i$. We define $\phi:= \sum_{i=0}^k w_i\cdot \phi_i$. Thus $\phi$ is a fractional $F$-packing of $G$ such that for each $e\in E(G)\setminus R$, we have that $\partial \phi(e) = \sum_{i=0}^k w_i = 1$. On the other hand, we have for each $i\in [k]$ that $\phi(e_i) = \sum_{j=0}^i w_i = \varphi(e_i)$. Thus $\phi$ is as desired.     
\end{proof}

We now use the above proposition to prove the following lemma which asserts that if when deleting low degree subgraphs of a graph $G$ we still admit fractional decompositions, then $G$ admits any fractional packing of target weights provided they are close enough to $1$.

\begin{lem}\label{lem:FracFixing}
Let $F$ be a graph and let $\varepsilon \in (0,1]$ be real. Let $G$ be a graph on $n$ vertices such that for all $H\subseteq G$ with $\Delta(H)\le \lceil 2\varepsilon n \rceil$, we have that $G-E(H)$ admits a fractional $F$-decomposition. If $\varphi:E(G)\rightarrow [1-\varepsilon,1]$, then there exists a fractional $F$-packing $\phi$ of $G$ such that $\partial \phi(e) = \varphi(e)$ for all $e\in E(G)$. \end{lem}
\begin{proof}
Since $\Delta(G)\le n-1$, we have by Vizing's Theorem that there exist an edge-coloring of $G$ with at most $\Delta(G)+1\le n$ colors. Let $d:= \lceil 2\varepsilon n \rceil$ and $C:= \lfloor \frac{1}{\varepsilon} \rfloor$. Note that $C\le \frac{1}{\varepsilon}$ and hence $\varepsilon C\le 1$. Furthermore, $n \le Cd$ since $2\varepsilon\cdot \lfloor 1/\varepsilon\rfloor \ge 1$ for all $\varepsilon\in (0,1]$. By taking the unions of $d$ color classes, we find that there exists a partition $R_1,\ldots, R_C$ of $E(G)$ such that $\Delta(R_i)\le d$ for each $i\in [C]$. By our assumption on $G$, we have for each $i\in [C]$ and $R'\subseteq R_i$ that $G-R'$ admits a fractional $F$-decomposition. For each $i\in [C]$, let $\varphi_i(e):=1$ for $e\in E(G)\setminus R_i$ and let $\varphi_i(e):= 1 - C\cdot (1-\varphi(e))$. Note that since $\varphi:E(G)\rightarrow [1-\varepsilon,1]$, we find that $\varphi_i:E(G)\rightarrow [1-\varepsilon C,1]$ for each $i\in [C]$. Since $\varepsilon C \le 1$, we further find that $\varphi_i:E(G)\rightarrow [1-\varepsilon C,1]$ for each $i\in [C]$.

Thus by Proposition~\ref{prop:FracFixing}, we have for each $i\in [k]$ that there exists a fractional $F$-packing $\phi_i$ of $G$ such that $\partial \phi_i(e)=1$ for $e\in E(G)\setminus R_i$ and $\partial\phi_i(e) = \varphi_i(e)$ for $e\in R_i$. Now let $\phi:= \frac{1}{C}\cdot \sum_{i\in [C]}\phi_i$. Since $\phi_i$ is a fractional $F$-packing of $G$ for each $i\in [C]$, it follows that $\phi$ is a fractional $F$-packing of $G$. Furthermore, for each $e\in E(G)$, we find that
$$\partial \phi(e) = \frac{1}{C}\cdot \sum_{i\in C} \partial \phi_i(e) = 1 + \frac{1}{C} \cdot (-C(1-\varphi(e)) = \varphi(e)$$
as desired.
\end{proof}

We are now prepared to prove a general lemma which asserts that if a graph $G$ has a family of subgraph $\mc{H}$ where every edge is in roughly the same number of elements of $\mc{H}$ and each graph in $\mc{H}$ admits a fractional $F$-decomposition under deletion of small maximum degree subgraphs, then $G$ admits a $C$-regular fractional $F$-decomposition for some integer $C\ge 1$. 

First we introduce some more notation. If $\mc{H}$ is a family of subgraphs of a graph $G$ and $e\in E(G)$, then we let $\mc{H}(e):= \{H\in \mc{H}: e\in E(H)\}$. 

\begin{lem}\label{lem:RegularFixing}
Let $F$ be a graph and let $C\ge 1$ and let $\varepsilon:= \frac{1}{2C}$. If $G$ is a graph and $\mc{H}$ is a family of subgraphs of $G$ such that all of the following hold: 
\begin{itemize}
    \item[(a)] $v(H)= k$ for all $H\in \mc{H}$,
    \item[(b)] for each $e\in E(G)$, 
    $$\frac{1-\varepsilon}{C}\cdot n^{k-2}\le |\mc{H}(e)| \le \frac{1}{C}\cdot n^{k-2},$$
    \item[(c)] for all $H\in \mc{H}$ and $H_0\subseteq H$ with $\Delta(H_0)\le \lceil4\varepsilon k\rceil$, we have that $H-E(H_0)$ admits a fractional $F$-decomposition,
\end{itemize}
then $G$ has a $2C$-regular fractional $F$-decomposition.
\end{lem}
\begin{proof}
Let $\mc{F}$ be the set of copies of $F$ in $G$. In order to seed the fractional decomposition, we first define a fractional $F_0$-packing as follows. For each $T\in \mc{F}$, we let $\phi_0(T) := w_0:= \frac{\varepsilon}{n^{V(F)-2}}$. Let $h_0:= \frac{1-\varepsilon}{C}\cdot n^{k-2}$. For each $e\in E(G)$, we let
$$\varphi(e) := (1- w_0\cdot |\mc{F}(e)|)\cdot \frac{h_0}{|\mc{H}(e)|}.$$
By condition (b), we find that
$$1\ge \frac{h_0}{|\mc{H}(e)|} \ge 1-\varepsilon.$$
Since $|\mc{F}(e)|\le n^{V(F)-2}$, we also find that $w_0\cdot |\mc{F}(e)|\le \varepsilon$. Since $(1-\varepsilon)^2 \ge 1-2\varepsilon$, we find that
$$\varphi(e) \in [1-2\varepsilon,1].$$
Thus we find that $\varphi:E(G)\rightarrow [1-2\varepsilon,1]$. Note $\varepsilon = \frac{1}{2C} \le \frac{1}{2}$ and hence we also have that $\varphi:E(G)\rightarrow [0,1]$. Now let $\varepsilon':= 2\varepsilon$. By (c), we have for all $H\in \mc{H}$ and $H_0\subseteq H$ with $\Delta(H_0)\le \lceil 2\varepsilon'k\rceil$ that $H-H_0$ admits a fractional $F$-decomposition. Hence by Lemma~\ref{lem:FracFixing}, there exists a fractional $F$-packing $\phi_H$ of $G$ such that $\partial \phi_H(e)= \varphi(e)$ for all $e\in E(G)$. Let 
$$\phi:= \phi_0 +\frac{1}{h_0} \cdot\sum_{H\in \mc{H}} \phi_H.$$ Thus for each $e\in E(G)$ we find that
\begin{align*}
\partial \phi(e) &= \partial \phi_0(e) + \frac{1}{h_0} \cdot\sum_{H\in \mc{H}} \partial \phi_H(e) = w_0\cdot |\mc{F}(e)| +\frac{1}{h_0}\cdot \varphi(e)\cdot |\mc{H}(e)|\\
&= w_0\cdot |\mc{F}(e)| +(1-w_0\cdot |\mc{F}(e)|) = 1
\end{align*}
and hence $\phi$ is a fractional $F$-decomposition of $G$. Note that for each $T\in \mc{F}$, we have that $\phi(T)\ge \phi_0(T) \ge w_0$. Thus $\phi$ is $\frac{1}{2C}$-seeded. On the other hand for each $T\in \mc{F}$, we have that $|\{H\in \mc{H}: T\subseteq H\}| \le n^{k-v(F)}$ and hence we calculate that
\begin{align*}
\phi(T) &= \phi_0(T) + \frac{1}{h_0} \cdot\sum_{H\in \mc{H}} \phi_H(T)\le w_0 + \frac{C}{(1-\varepsilon)\cdot n^{V(F)-2}}\\
&\le \frac{\varepsilon}{n^{v(F)-2}} + \frac{C(1+\varepsilon)}{n^{v(F)-2}}= \frac{1}{n^{v(F)-2}} \cdot \left(\frac{1}{2C} + C + \frac{1}{2}\right)\\
&\le \frac{1}{n^{v(F)-2}}\cdot (C+1) \le \frac{2C}{n^{v(F)-2}},
\end{align*}
and hence $\phi$ is $2C$-low-weight. Altogether then, we find that $\phi$ is a $C$-regular fractional $F$-decomposition of $G$ as desired.
\end{proof}

We are now prepared to prove Theorem~\ref{thm:RegularFracStability}.

\begin{proof}[Proof of Theorem~\ref{thm:RegularFracStability}]
We choose $\varepsilon$ small enough as needed throughout the proof. Let $\varepsilon':=7\varepsilon$ and let $\sigma':=\frac{\sigma}{2}$. We crucially note that since $\varepsilon$ is small enough, we have that $\varepsilon'$ is small enough to satisfy Theorem~\ref{thm:fracNWstability} for $\sigma'$.

Suppose for a contradiction that $G$ is not $\sigma$-extremal. Set $\rho:= \varepsilon$. Let $k_0$ be as in Lemma~\ref{lem:Sampled2} for $\varepsilon, \rho,\sigma$. We choose $k\ge k_0$ such that $\lceil 4\varepsilon k \rceil \le 5\varepsilon k$. Let 
$$\mc{H}:= \bigg\{G[S]: S\subseteq V(G),~|S|=k,~\delta(G[S])\ge \left(\frac{3}{4}-2\varepsilon\right)k \text{ and $G[S]$ is not $\sigma/2$-extremal}\bigg\}.$$
By Lemma~\ref{lem:Sampled2}, we find for each $e\in E(G)$ that
$$(1-\rho)\cdot \binom{n-2}{k-2} \le |\mc{H}(e)| \le \binom{n-2}{k-2}.$$
Let $C:=\frac{1}{2}\cdot n^{k-2} /\binom{n-2}{k-2}$. Let $F:=K_3$. Thus $G$ and $\mc{H}$ satisfy Lemma~\ref{lem:RegularFixing}(a) and (b) for $F$, $\varepsilon$ and $C/2$. 

Further note that for each $H\in\mc{H}$ and $H_0\subseteq H$ with $\Delta(H_0)\le \lceil 4\varepsilon k\rceil$, we have that $\Delta(H-E(H_0)) \ge (\frac{3}{4} - 7\varepsilon)k$ and hence $H-E(H_0)$ admits a triangle decomposition by Theorem~\ref{thm:fracNWstability} as $H$ is not $\sigma'$-extremal. Thus $G$ and $\mc{H}$ satisfy Lemma~\ref{lem:RegularFixing}(c) for $F$, $\varepsilon$ and $C/2$.
By Lemma~\ref{lem:RegularFixing}, it follows that $G$ has a $C$-regular fractional $K_3$-decomposition, a contradiction.
\end{proof}

\part{Proof of Nash-Williams' Conjecture}

\section{Proof Overview}

Before breaking our proof of Nash-Williams' Conjecture into cases, we need the following notation and recall the definition of extremal for Nash-Williams' Conjecture.

\begin{definition}
We write $T_2(n) := K_{\lfloor \frac{n}{2} \rfloor, \lceil \frac{n}{2} \rceil}$ for the Tur\'an graph $T(n,2)$. 
\end{definition}

\begin{definition}
Let $\sigma \ge 0$. A graph $G$ on $n$ vertices is \emph{$\sigma$-extremal} if both of the following hold:    
\begin{itemize}
    \item[(i)] at least $(1-\sigma)n$ vertices of $G$ have degree at most $(\frac{3}{4}+\sigma)n$, and
    \item[(ii)] the maximum cut of $G$ is at least $(1-\sigma)\cdot \frac{n^2}{4}$.
\end{itemize}     
\end{definition}

We recall from the literature that graphs $G$ and $H$ on the same number of vertices are \emph{$\alpha$-close} if one can add and/or delete at most $\alpha\cdot v(G)^2$ edges to transform $G$ into a graph isomorphic to $H$. We also recall that the \emph{join} of graphs $H_1$ and $H_2$, denoted $H_1+H_2$, is the graph obtained by taking the vertex-disjoint union of $H_1$ and $H_2$ and adding a complete bipartite graph between $V(H_1)$ and $V(H_2)$. Also recall that $\overline{G}$ denotes the complement of a graph $G$.
\vskip.1in
\noindent Recall our strategy is to divide the proof of Nash-Williams' Conjecture into five cases as follows: $G$ is
\begin{enumerate}
    \item[(1)] not $\beta$-extremal
    \item[(2)] $\beta$-extremal but not in (3)-(5)
    \item[(3)] $\alpha$-close to $\overline{T(n/2)}+\overline{T(n/2)}$
    \item[(4)] $\alpha$-close to $\overline{T(n/2)}+T(n/2)$
    \item[(5)] $\alpha$-close to $T(n/2)+T(n/2)$.
\end{enumerate}

\begin{remark}
Note in the cases above we should write $\lfloor \frac{n}{2} \rfloor$ for one graph and $\lceil \frac{n}{2} \rceil$ for the other as $n$ might not be even. However, for the sake of the reader we omit the floors and ceilings since in the actual proofs we will fix a partition that is $\alpha$-close to a balanced complete blow-up of the associated family of partitioned triangles (see Definition~\ref{def:Abnormal}) rather than an edit-distance close mapping to a specific graph.    
\end{remark}

\noindent Here are the corresponding theorems for the five cases.

\begin{thm}\label{thm:nonextremal}
There exists $\beta_0 > 0$ such that for every real $\beta \in (0,\beta_0]$, there exists an integer $n_0$ such that for all $n\ge n_0$ the following holds: If $G$ is a $K_3$-divisible graph on $n\ge n_0$ vertices with $\delta(G)\ge \frac{3}{4}n$ and $G$ is not $\beta$-extremal, then $G$ has a $K_3$-decomposition. 
\end{thm}

\begin{thm}\label{thm:neitherside}
For every real $\alpha > 0$, there exists a real $\beta > 0$ and an integer $n_0$ such that for all $n\ge n_0$ the following holds: If $G$ is a $K_3$-divisible graph on $n\ge n_0$ vertices with $\delta(G)\ge \frac{3}{4}n$ and $G$ is $\beta$-extremal but $G$ is not $\alpha$-close to any of $T_2(n/2) + T_2(n/2)$, $\overline{T_2(n/2)} + T_2(n/2)$, or $\overline{T_2(n/2}) + \overline{T_2(n/2)}$, then $G$ has a $K_3$-decomposition. 
\end{thm}

\begin{thm}\label{thm:4cliques}
There exists a real $\alpha > 0$ and an integer $n_0$ such that for all $n\ge n_0$ the following holds: If $G$ is a $K_3$-divisible graph on $n\ge n_0$ vertices with $\delta(G)\ge \frac{3}{4}n$ and $G$ is $\alpha$-close to $\overline{T_2(n/2)} + \overline{T_2(n/2)}$, then $G$ has a $K_3$-decomposition. 
\end{thm}

\begin{thm}\label{thm:2cliqueandbip}
There exists a real $\alpha > 0$ and an integer $n_0$ such that for all $n\ge n_0$ the following holds: If $G$ is a $K_3$-divisible graph on $n\ge n_0$ vertices with $\delta(G)\ge \frac{3}{4}n$ and $G$ is $\alpha$-close to $\overline{T_2(n/2)} + T_2(n/2)$, then $G$ has a $K_3$-decomposition. 
\end{thm}

\begin{thm}\label{thm:K4blowup}
There exists a real $\alpha > 0$ and an integer $n_0$ such that for all $n\ge n_0$ the following holds: If $G$ is a $K_3$-divisible graph on $n\ge n_0$ vertices with $\delta(G)\ge \frac{3}{4}n$ and $G$ is $\alpha$-close to $T_2(n/2) + T_2(n/2)$, then $G$ has a $K_3$-decomposition. 
\end{thm}

\subsection{Outline of Part III}

In Section~\ref{s:Partitioned}, we introduce partitioned graphs. Section~\ref{ss:Partitioned} contains the requisite definitions. Section~\ref{ss:Decomp} states our main decomposition theorem which is a special case of the main `black-box' decomposition theorem from~\cite{DP26}. Section~\ref{ss:Div} characterizes the  fractional and integral divisibility conditions for our four special  cases.

In Section~\ref{s:AbsStabDecomp}, we prove our main absorber stability lemma and idealized decomposition theorems. Namely, in Section~\ref{ss:AbsorberStability}, we prove our absorber stability lemma, Lemma~\ref{lem:ManyHinges}. In Section~\ref{ss:Idealized}, we prove our various idealized decomposition theorems for the different cases.

In Section~\ref{s:Cleaning}, we prove the theorems of the various cases. In Section~\ref{ss:Cleaning}, we prove a few cleaning lemmas that show how to clean a graph to the divisibility conditions for our cases so as to apply the corresponding idealized decomposition theorems. Then we proceed in the remainder of the section to prove our various cases of Nash-Williams' Conjectures, with the cases for $\mc{F}_1$ and $\mc{F}_2$ requiring parity fixers (cases (3) and (4) above) and the case of $\mc{F}_3$ (case (5) above) requiring more careful attention due to its stringent divisibility conditions. 

\subsection{Proof of Nash-Williams' Conjecture}

Assuming the theorems above, we are now prepared to prove Theorem~\ref{thm:NW}, that is Nash-Williams' Conjecture.

\begin{proof}[Proof of Nash-Williams' Conjecture.]
We choose $n$ sufficiently large as needed throughout the proof. Let $\alpha_1$ be as in Theorem~\ref{thm:4cliques}, $\alpha_2$ be as in Theorem~\ref{thm:2cliqueandbip}, and $\alpha_3$ be as in Theorem~\ref{thm:K4blowup}. Let $\alpha:= \min\{\alpha_1,\alpha_2,\alpha_3\}$. Let $\beta$ be as in Theorem~\ref{thm:neitherside} for $\alpha$. Let $\beta_0$ be as in Theorem~\ref{thm:nonextremal} and let $\beta':= \min\{\beta_0,\beta\}$. Recall that $G$ is a $K_3$-divisible graph on $n$ vertices with $\delta(G)\ge \frac{3}{4}n$. 

First suppose that $G$ is $\alpha$-close to $T_2(n/2) + T_2(n/2)$. Then by Theorem~\ref{thm:K4blowup}, $G$ has a $K_3$-decomposition as desired. So we  assume $G$ is not $\alpha$-close to $T_2(n/2) + T_2(n/2)$.

Next suppose that $G$ is $\alpha$-close to $\overline{T_2(n/2)} + T_2(n/2)$. Then by Theorem~\ref{thm:2cliqueandbip}, $G$ has a $K_3$-decomposition as desired. So we assume $G$ is not $\alpha$-close to $\overline{T_2(n/2)} + T_2(n/2)$.

Next suppose that $G$ is $\alpha$-close to $\overline{T_2(n/2)} + \overline{T_2(n/2)}$. Then by Theorem~\ref{thm:4cliques}, $G$ has a $K_3$-decomposition as desired. So we assume $G$ is not $\alpha$-close to $\overline{T_2(n/2)} + \overline{T_2(n/2)}$.

Next suppose that $G$ is $\beta'$-extremal. Since $\beta'\le \beta$, we find that $G$ is also $\beta$-extremal. Since $G$ is not $\alpha$-close to any of $T_2(n/2) + T_2(n/2)$, $\overline{T_2(n/2)} + T_2(n/2)$, or $\overline{T_2(n/2}) + \overline{T_2(n/2)}$, we have by Theorem~\ref{thm:neitherside} that $G$ has a $K_3$-decomposition as desired.

Hence we assume that $G$ is not $\beta'$-extremal. Since $\beta'\le \beta_0$ and $n$ is large enough, we have by Theorem~\ref{thm:nonextremal} that $G$ has a $K_3$-decomposition as desired.
\end{proof}

\section{Partitioned Graphs, Divisibility and a Decomposition Theorem}\label{s:Partitioned}

In this section, we proceed to provide the many definitions we need about partitioned graphs, state our `black-box' partitioned decomposition theorem which is a special case of the main decomposition theorem from~\cite{DP26}, and then characterize integral and fractional divisibility for our various families of partitioned triangles.

\subsection{Partitioned Graphs and Divisibility}\label{ss:Partitioned}

\begin{definition}[Partitioned Graph]
Let $t \ge 1$ be an integer. A \emph{$t$-partitioned} graph is a graph $G$ together with a partition $\mc{P}_G=(V_1,\ldots,V_t)$ of $V(G)$. We use $e_{ij}(G)$ for $i\ne j \in [t]$ to denote $e(G[V_i,V_j])$ and $e_{ii}(G)$ to denote $e(G[V_i])$; we let $\vec{e}(G)$ denote the vector $(e_{ij}(G) : i\le j\in [t])$. For $v\in V(G)$ and $i\in [t]$, we let $d_i(v) := |N_G(v)\cap V_i|$ and let $\vec{d}_G(v)$ denote the vector $(d_i(v): i\in [t])$. 
\end{definition}

\begin{definition}[$\mc{F}$-decomposition]
Let $t \ge 1$ be integers. Let $F$ and $G$ be $t$-partitioned graphs. A \emph{copy} of $F$ in $G$, or \emph{$F$-copy} for short, is a subgraph $H$ of $G$ such that there exists an isomorphism $\phi$ from $H$ to $F$ such that $\mc{P}_F(\phi(v)) = \mc{P}_G(v)$ for all $v\in V(H)$ for all $e\in H$ and $v\in V(e)$. An \emph{$F$-decomposition} of $G$ is a partition of the edges of $G$ into copies of $F$. Let $\mc{F}=\{F_1,\ldots, F_k\}$ be a family of $t$-partitioned graphs. An \emph{$\mc{F}$-copy} of $G$ is an $F_i$-copy of $G$ for some $i\in [k]$. An \emph{$\mc{F}$-decomposition} of $G$ is a partition of the edges of $G$ into $\mc{F}$-copies.
\end{definition}

\begin{definition}[$\mc{F}$-blowups]
Let $t \ge 1$ be integers. Let $F$ and $G$ be $t$-partitioned graphs. We say $G$ is an \emph{$F$-blowup} if for all $i\le j \in [t]$, we have that $e_{ij}(G)\ne 0$ if and only if $e_{ij}(F)\ne 0$. Let $\mc{F}=\{F_1,\ldots,F_k\}$ be a family of $t$-partitioned graphs. We say $G$ is an \emph{$\mc{F}$-blowup} if for all $i\le j \in [t]$, we have that $e_{ij}(G)\ne 0$ if there exists $\ell \in [k]$ such that $e_{ij}(F_{\ell})\ne 0$.
\end{definition}

To characterize the necessary divisibility conditions for an $\mc{F}$-decomposition (what is called $\mc{F}$-divisible) and similarly for integral $\mc{F}$-decompositions (what is called integrally $\mc{F}$-divisible) and fractional $\mc{F}$-decompositions (what is called fractionally $\mc{F}$-divisible), we require the following standard concepts from linear algebra, convex geometry, and optimization.

\begin{definition}
Let $v_1,\ldots, v_k \in \mathbb{Z}_{\ge 0}^d$ and let $A$ be the matrix whose columns are comprised of the $v_i$. The \emph{cone} generated by $v_1,\ldots, v_k$, denoted $\textrm{cone}(A)$, is $\{ Ax: x\in \mathbb{R}^k,~x\ge 0\}$. The \emph{lattice} of $A$, denoted $\textrm{lattice}(A)$, is $\{Ax: x\in \mathbb{Z}^k\}$. The \emph{span} of $A$, denoted $\textrm{span}(A)$, is $\{Ax: x\in \mathbb{R}^k\}.$ The \emph{subgroup} generated by $A$ is $\{Ax: x\in \mathbb{Z}^k,~x\ge 0\}$.       
\end{definition}

\begin{definition}[Divisibility]
Let $t \ge 1$ be integers. Let $\mc{F}=\{F_1,\ldots, F_k\}$ be a family of $t$-partitioned graphs and let $G$ be a $t$-partitioned graph. We say $G$ is
\begin{itemize}
    \item \emph{$\mc{F}$-divisible} if for every $e\in E(G)$, $\vec{e}(G)$ is in the subgroup generated by $(\vec{e}(F_i): i\in [k])$, and for every $v\in V(G)$, $\vec{d}_G(v)$ is in the subgroup generated by $(\vec{d}_{F_i}(u): i\in [k],~u\in V(F_i) \textrm{ with } \mc{P}_{F_i}(u) = \mc{P}_G(v))$,
    \item \emph{integrally $\mc{F}$-divisible} if the same holds with subgroup replaced by lattice,
    \item \emph{fractionally $\mc{F}$-divisible} if the same holds with subgroup replaced by cone,
    \item \emph{rationally $\mc{F}$-divisible} if the same holds with subgroup replaced by span.
\end{itemize}
\end{definition}

Here we make some remarks for the interested reader about how to characterize lattices and cones as follows.

\begin{remark}
The problem of determining the span is simply determining the equations that define it. The problems of characterizing lattices and cones are also classical and well-understood. To characterize a lattice, one seeks to determine the \emph{congruence relations} that define the lattice; these include the equations that define the span with possibly some modular relations in addition. The congruence relations can be found by finding the \emph{Smith normal form} of an integral matrix $A$ which is $UDV$ where $U,V$ are totally unimodular matrices and $D$ is a diagonal matrix. Then the congruence relations that $y\in \textrm{lattice}(A)$ must satisfy are that $U_iy \equiv 0 \mod d_i$ where $d_1,\ldots, d_k$ are the diagonal entries of $D$ and $U_1,\ldots, U_k$ are the rows of $U$. A computer can fairly quickly calculate the Smith normal form for small matrices (and indeed we initially did this for our various cases). That said, $U,V,D$ in fact provide a computer-free certificate that can be checked via matrix multiplication.
\end{remark}

\begin{remark}
To characterize a cone, one seeks to determine the \emph{facet-defining inequalities}. This can be done algorithmically via the \emph{double description method}. A computer can fairly quickly calculate the facet-defining inequalities for small matrices (and indeed we initially did this for our various cases). However, a computer-free certificate (other than verifying that each step of the algorithm was correctly implemented) is more elusive and can be complicated to write. That said, when the columns of $A$ are linearly independent, characterizing the solution set via taking inverses suffices to provide a computer-free proof.
\end{remark}

\begin{remark}
Thankfully for our four cases, three of the cases have linearly independent generators; thus the solution point can easily be characterized and this yields both the facet-defining inequalities directly but this also yields the congruence relations as well (without needing to calculate smith normal form). For the fourth case, the solution set (which is $1$-dimensional in that case) can also be characterized, again avoiding the need to calculate the Smith normal form; as for the facet-defining inequalities, the cone in this case can be recast as a max flow problem and hence the facet-defining inequalities will be equivalent to the corresponding cut constraints (thus providing a computer-free proof via the max-flow/min-cut theorem).    
\end{remark}

\subsection{A Decomposition Theorem}\label{ss:Decomp}

We require a black-box decomposition theorem which is only a very special case of our black-box decomposition theorem from~\cite{DP26} which is our paper on refined absorption for partitioned graphs (and beyond). In particular, we only require our theorem to hold for partitioned families of graphs (instead of hypergraphs possibly with colored and oriented edges); and indeed, we only require our theorem to hold for families of partitioned triangles. So we will restrict to these settings as necessary to reduce the complexity for the reader. Let us start with some definitions for fractional decompositions into families of partitioned graphs as follows. 

\begin{definition}
Let $\mc{F}=\{F_1,\ldots,F_k\}$ be a family of $t$-partitioned graphs. A \emph{fractional $\mc{F}$-packing} of a $t$-partitioned graph $G$ is an assignment $\phi$ of non-negative weights to the $\mc{F}$-copies of $G$ such that for each $e\in E(G)$, $\partial \phi(e):= \sum_{H: e\in H} \phi(H) \le 1$; we say $\phi$ is a \emph{fractional $\mc{F}$-decomposition} if for each $e\in E(G)$, $\phi(e) \le 1$. 
\end{definition}

\begin{definition}
Let $\mc{F}=\{F_1,\ldots,F_k\}$ be a family of $t$-partitioned graphs and let $C \ge 1 \ge \sigma$. Let $\phi$ be a fractional $\mc{F}$-packing of a $t$-partitioned graph $G$ on $n$ vertices. We say $\phi$ is
\begin{itemize}
    \item \emph{$C$-low-weight} if for each $i\in [k]$, every $F_i$-copy in $G$ has weight at most $\frac{C}{n^{v(F_i)-2}}$,
    \item \emph{$\sigma$-seeded} if for each $i\in [k]$, every $F_i$-copy in $G$ has weight at least $\frac{\sigma}{n^{v(F_i)-2}}$,
    \item \emph{$C$-regular} if $\phi$ is $C$-low-weight and $\frac{1}{C}$-seeded,
    \item \emph{$\varepsilon$-almost} if for each $e\in E(G)$, $\partial \phi(e)\ge 1-\varepsilon$.
\end{itemize} 
\end{definition}

We next recall the definitions of anti-edge and hinges along with the corresponding notion for a partitioned graph having many of these gadgets as follows. For anti-edges, we only define this for (families of) partitioned triangles for simplicity. First a notation definition.

\begin{definition}[Partitioned Triangle Notation]
We let $T_{ijk}$ denote the triangle whose vertices have parts $i, j,$
 and $k$.    
\end{definition}

\begin{definition}[Anti-Edge]\label{def:AntiEdge}
Let $F=T_{ijk}$ be a $t$-partitioned triangle. Let $u\in V_i,~v\in V_j$. An \emph{$F$-anti-edge on $uv$}, denoted ${\rm AntiEdge}_{F}(uv)$, is a subgraph $H$ such that $H\cup \{uv\}$ is an $F$-copy of $G\cup \{uv\}$ (or equivalently is a path $uwv$ where $w\in V_k$).
\end{definition}

\begin{definition}[Anti-Edgeful]\label{def:antiedgeful}
Let $\mc{F}$ be a family of $t$-partitioned triangles and let $\alpha \in (0,1)$. We say a $t$-partitioned graph $G$ on $n$ vertices is \emph{$\alpha$-$\mc{F}$-anti-edgeful} if for every $T_{ijk}\in \mc{F}$, $u\in V_i,~v\in V_j$ there exist at least $\alpha\cdot n$ $T_{ijk}$-anti-edges on $uv$ in $G$.
\end{definition}

\begin{definition}[Hinge]\label{def:hinge}
Let $\mc{F}=\{F_1,\ldots,F_k\}$ be a family of $t$-partitioned graphs. Let $S$ and $S'$ be two $\mc{F}$-copies such that $S\cap S'=\{e\}$ for some edge $e$. An \emph{$\mc{F}$-hinge} for $S$ and $S'$ is a graph $H$ such that $H$ is edge-disjoint from $S\cup S'$, $H\cup (S\setminus e)$ has a $\mc{F}$-decomposition and $H\cup (S'\setminus e)$ has a $\mc{F}$-decomposition
\end{definition}

\begin{definition}[Hingeful]\label{def:hingeful}
Let $\mc{F}=\{F_1,\ldots,F_k\}$ be a family of $t$-partitioned graphs and let $\alpha \in (0,1)$ and $C\ge 1$. We say a $t$-partitioned graph $G$ on $n$ vertices is \emph{$(\alpha,C)$-$\mc{F}$-hingeful} if for every two $\mc{F}$-copies $S$ and $S'$ of $G$ such that $S\cap S'$ is an edge, there exists $i\in [C]$ such that there exist at least $\alpha\cdot n^{i-2}$ $\mc{F}$-hinges $H$ for $S$ and $S'$ in $G$ such that $v(H)=i$.
\end{definition}

We are now prepared to state the required decomposition theorem which is a special case of our main decomposition theorem from~\cite{DP26}.

\begin{thm}\label{thm:DecompThm}
For every real $\alpha > 0$, integer $C\ge 1$ and family $\mc{F}$ of $t$-partitioned triangles, there exists a real $\varepsilon>0$ such that the following holds for sufficiently large $n$: If $G$ is a $t$-partitioned graph on $n$ vertices that satisfies both of the following
\begin{itemize}
    \item[(i)] $G$ is fractionally-$\mc{F}$-divisible and there exists a $C$-regular $\varepsilon$-almost fractional $\mc{F}$-packing of $G$, and
    \item[(ii)] $G$ is integrally-$F$-divisible, $\alpha$-$\mc{F}$-anti-edgeful and $(\alpha, C)$-$\mc{F}$-hingeful,
\end{itemize}
then there exists an $\mc{F}$-decomposition of $G$.
\end{thm}

\subsection{Integer and Fractional Divisibility Results}\label{ss:Div}

We now proceed to characterize integral and fractional divisibility for our various partitioned settings. The following proposition is useful for characterizing the cones and lattices of the degree vectors in the various cases.

\begin{proposition}\label{prop:LatticeCone}
Let $A$ be the matrix with columns $(1, 1, 0)$ and $(1,0,1)$ and $(0,1,1)$. Then $\textrm{cone}(A) = \{ (x_1,x_2,x_3)\in \mathbb{R}^3 : 0\le x_i \le \sum_{j\in [3]\setminus \{i\}} x_j \text{ for all } i\in [3]\}$ and $\textrm{lattice}(A) = \{(x_1,x_2,x_3)\in \mathbb{Z}^3: x_1+x_2+x_3 \equiv 0 \pmod 2\}$.

Let $A'$ be the matrix with columns $(1,0,1)$, $(0,1,1)$, $(2,0,0)$ and $(0,2,0)$. Then $\textrm{cone}(A') = \{ (x_1,x_2,x_3)\in \mathbb{R}_{\ge 0}^3 : x_3\le x_1+x_2\}$ and $\textrm{lattice}(A') = \{(x_1,x_2,x_3)\in \mathbb{Z}^3: x_1+x_2+x_3 \equiv 0 \pmod 2\}$.
\end{proposition}

The following lemma characterizes divisibility needed to prove Theorem~\ref{thm:neitherside}.
 
\begin{lem}\label{lem:F0div}
Let $\mc{F}_0 :=\{T_{112}, T_{122}\}$. An $\mc{F}_0$-blowup $G$ is fractionally $\mc{F}_0$-divisible if and only if $G$ satisfies all of the following:
\begin{enumerate}
    \item[(i)] $e_{12}(G) = 2\cdot (e_{11}(G)+e_{22}(G))$,
    \item[(ii)] $d_i(v) \le d_{3-i}(v)$ for each $i\in [2]$ and $v\in V_i$.\end{enumerate}
An $\mc{F}_0$-blowup $G$ is integrally $\mc{F}_0$-divisible if and only if $G$ satisfies the following:
\begin{enumerate}
    \item[(i')] $e_{12}(G) = 2\cdot (e_{11}(G)+e_{22}(G))$,
    \item[(ii')] $G$ is Eulerian.
\end{enumerate}
\end{lem}
\begin{proof}
First we consider the degree constraints. For $v\in V(T_{112})$ with $\mc{P}(v)=1$, we have that $\vec{d}_{T_{112}}(v) = (1,1)$; for $v\in V(T_{122})$ with $\mc{P}(v)=1$, we have that $\vec{d}_{T_{122}}(v) = (0,2)$ and symmetrically $\vec{d}_{T_{112}}(v)(2,0)$ and $\vec{d}_{T_{122}}(v)=(1,1)$ for $v$ with $\mc{P}(v)=2$. Now $\textrm{cone}( (1,1), (2,0)) = \{(x_1,x_2)\in \mathbb{R}_{\ge 0}^2: x_2\le x_1\}$ and $\textrm{lattice}( (1, 1), (2, 0) )=\{(x_1,x_2)\in \mathbb{Z}^2: x_1+x_2\equiv 0\pmod 2\}$. Hence condition (ii) is both necessary and sufficient for degree vectors for fractional $\mc{F}_0$-divisibility while condition (ii') is both necessary and sufficient for degree vectors for integral $\mc{F}_0$-divisibility.

Finally we consider the edge constraints. For this proof when considering a $2$-partitioned graph $H$, we write $\vec{e}(H)$ as $(e_{11}(H),e_{22}(H),e_{12}(H))$. Thus $\vec{e}(T_{112})=(1,0,2)$ and $\vec{e}(T_{122})=(0,1,2)$. We let $A$ be the matrix whose columns are these generators. Then $G$ is fractionally $\mc{F}_0$-divisible if and only if $\vec{e}(G)\in \textrm{cone}(A)$ and (ii) holds. Similarly $G$ is integrally $\mc{F}_0$-divisible if and only if $\vec{e}(G)\in \textrm{lattice}(A)$ and (ii') holds.

Note $\vec{e}(G) \in \textrm{span}(A)$ if and only if $e_{12}(G) = 2\cdot (e_{11}(G)+e_{22}(G))$ (to verify this note that both generators satisfy this relation and they are linearly independent). Quotienting out by this relation and then calculating, we find that $x=(t_{112},t_{122})$ is a solution to $Ax=\vec{e}(G)$ if and only if $t_{112} = e_{11}(G)$ and $t_{122}=e_{22}(G)$ and $G$ satisfies the above equality condition. But then it follows that $x\ge 0$ and $x\in \mathbb{Z}^2$ since $G$ is a graph.  
\end{proof}

Our next lemma characterizes divisibility needed to prove Theorem~\ref{thm:4cliques}.

\begin{lem}\label{lem:F1div}
Let $\mc{F}_1=\{T_{113}, T_{114}, T_{223}, T_{224},T_{133},T_{233},T_{144},T_{244}\}$. An $\mc{F}_1$-blowup $G$ is fractionally $\mc{F}_1$-divisible if and only if $G$ satisfies all of the following:
\begin{enumerate}
    \item[(i)] $e_{13}(G)+e_{14}(G)+e_{23}(G)+e_{34}(G) = 2\cdot (e_{11}(G)+e_{22}(G)+e_{34}(G))$,
    \item[(ii)] for all $i\in [4]$ and $v\in V_i$, $d_i(v) \le d_{\ell}(v)+d_k(v)$ where $\{\ell,k\}=\{3,4\}$ if $i\in\{1,2\}$ and $\{\ell,k\}= \{1,2\}$ if $i\in \{3,4\}$,
    \item[(iii)] all of ``the cut constraints'':
    \begin{itemize}
        \item[(a)] $2\cdot e_{ii}(G) \le e_{ik}(G)+e_{i\ell}(G)$ for $i\in [4]$ where $\{\ell,k\}=\{3,4\}$ if $i\in\{1,2\}$ and $\{\ell,k\}= \{1,2\}$ if $i\in \{3,4\}$,
        \item[(b)] $e_{ik}(G) \le 2\cdot (e_{ii}(G)+e_{kk}(G))$ where $i\in \{1,2\}$ and $k\in \{3,4\}$,
        \item[(c)] $2\cdot (e_{ii}(G)+e_{kk}(G))\le e_{ik}(G)+e_{jk}(G)+e_{i\ell}(G)$ where $i\in \{1,2\}, k\in\{3,4\}$ and $j:=[2]\setminus \{i\}$ and $\ell:=\{3,4\}\setminus \{k\}$.
    \end{itemize}
\end{enumerate}
An $\mc{F}_1$-blowup $G$ is integrally $\mc{F}_1$-divisible if and only if $G$ satisfies all of the following:
\begin{enumerate}
    \item[(i')] $e_{13}(G)+e_{14}(G)+e_{23}(G)+e_{34}(G) = 2\cdot (e_{11}(G)+e_{22}(G)+e_{34}(G))$
    \item[(ii')] $G$ is Eulerian
    \item[(iii')] $e_{13}(G),e_{23}(G),e_{24}(G),e_{34}(G)$ are all even.
\end{enumerate}
\end{lem}
\begin{proof}
First we consider the degree constraints. For $v\in V(T_{113})$ with $\mc{P}(v)=1$, we have that $\vec{d}_{T_{113}}(v) = (1,0,1,0)$; for $v\in V(T_{133})$ with $\mc{P}(v)=1$, we have that $\vec{d}_{T_{133}}(v) = (0,0,2,0)$. Similarly, for $v\in V(T_{114})$ with $\mc{P}(v)=1$, we have that $\vec{d}_{T_{114}}(v) = (1,0,0,1)$; for $v\in V(T_{144})$ with $\mc{P}(v)=1$, we have that $\vec{d}_{T_{144}}(v) = (0,0,0,2)$. By Proposition~\ref{prop:LatticeCone}, we find that $\textrm{cone}( (1,0,1,0), (0,0,2,0), (1,0,0,1), (0,0,0,2) ) = \{(x_1,0,x_3,x_4)\in \mathbb{R}_{\ge 0}^4: x_1\le x_3+x_4\}$ and $\textrm{lattice}( (1,0,1,0), (0,0,2,0), (1,0,0,1), (0,0,0,2) )=\{(x_1,0,x_3,x_4)\in \mathbb{Z}^4: x_1+x_3+x_4\equiv 0\pmod 2\}$. These are then the necessary divisibility conditions for $\vec{d}(v)$ for $v\in V_1$. By symmetry similar symmetric conditions hold for all $v\in V_i$ for $i\in\{2,3,4\}$. Hence condition (ii) is both necessary and sufficient for degree vectors for fractional $\mc{F}_1$-divisibility while condition (ii') is both necessary and sufficient for degree vectors for integral $\mc{F}_1$-divisibility. 

Finally we consider the edge constraints. For this proof when considering a $4$-partitioned graph $H$, we write $\vec{e}(H)$ as $(e_{11}(H),e_{22}(H),e_{33}(H),e_{44}(H),e_{13}(H),e_{14}(H),e_{23}(H),e_{24}(H))$. Thus $\vec{e}(T_{113})=(1,0,0,0,2,0,0,0)$, $\vec{e}(T_{114})=(1,0,0,0,0,2,0,0)$, $\vec{e}(T_{223})=(0,1,0,0,0,0,2,0)$, $\vec{e}(T_{224})=(0,2,0,0,0,0,0,2)$, $\vec{e}(T_{133})=(0,0,1,0,2,0,0,0)$, $\vec{e}(T_{233})=(0,0,1,0,0,0,2,0)$, $\vec{e}(T_{144})=(0,0,0,1,0,2,0,0)$, and $\vec{e}(T_{244})=(0,0,0,1,0,0,0,2)$. We let $A$ be the matrix whose columns are these generators. Then $G$ is fractionally $\mc{F}_1$-divisible if and only if $\vec{e}(G)\in \textrm{cone}(A)$ and (ii) holds. Similarly $G$ is integrally $\mc{F}_1$-divisible if and only if $\vec{e}(G)\in \textrm{lattice}(A)$ and (ii') holds.

We note that (i)/(i') is a necessary condition for $\vec{e}(G)\in \textrm{span}(A)$(which lie in $\mathbb{R}^8$) since it is satisfied by all the generator vectors; since $\textrm{rank}(A)=7$, the condition (i)/(i') is also a sufficient condition for $\vec{e}(G)$ to lie in $\textrm{span}(A)$. Similarly one verifies the remaining conditions are all necessary as they are satisfied by the generator vectors. Quotienting out by this equation, we proceed to describe the solution space namely of the vector $\vec{t}$ whose coordinates $t_{ijk}$ where this represents the number of triangles of type $T_{ijk}$. Note there is one degree of freedom in the solution space, let it be $t_{113}$. One then calculates that 
\begin{align*}
t_{114} &=e_{11}(G)-t_{113}\\
t_{144} &= \frac{1}{2} ( e_{14}(G) - 2\cdot t_{114}) = \frac{1}{2} \cdot e_{14}(G) - t_{114}\\
t_{244} &= e_{44}(G)-t_{144} \\
t_{224} &= \frac{1}{2} \cdot e_{24}(G) - t_{244}\\
t_{223} &= e_{22}(G)-t_{224}\\
t_{233} &= \frac{1}{2}\cdot e_{23}(G)-t_{223} \\
t_{133} &=e_{33}(G)-t_{233}    
\end{align*} For integral divisibility, the only subsequent condition is (iii'). Thus we turn to fractional divisibility and the cone. We characterize this as a max-flow/min-cut problem as follows. We construct an auxiliary digraph $H$ as follows. Let $s$ be a source and $t$ be a sink and let 
$$V(H):=\{s,t,e_{11},e_{22},e_{33},e_{44},e_{13},e_{14},e_{23},e_{24}\}.$$ Then we let
\begin{align*}
E(H):= &\Big\{\overrightarrow{se_{ii}}:i\in [4] \Big\}\cup \Big\{\overrightarrow{e_{ij}t}:i\in \{1,2\},~j\in\{3,4\}\Big\} \\
&\cup \Big\{\overrightarrow{e_{ii}e_{ij}}: i\in \{1,2\}, j\in \{3,4\}\Big\}\cup \Big\{\overrightarrow{e_{jj}e_{ij}}:j\in\{3,4\},~i\in\{1,2\}\Big\}.
\end{align*}
Then we define a capacity function $c$ for $H$ as follows. For $i\in [4]$, we let $c(se_{ii}) :=e_{ii}(G)$. For $i\in \{1,2\},~j\in\{3,4\}$, we let $c(e_{ij}t) := \frac{1}{2} \cdot e_{ij}(G)$ and we let $c(e_{ii}e_{ij})=c(e_{jj}e_{ij}) :=\infty$. Thus there exists a flow $f$ of $H$ from $s$ to $t$ of value $a:=\sum_{i\in [4]}e_{ii}(G)$ if and only if $t_{iij}:=f(e_{ii}e_{ij})$ is in the cone of $A$. 

By the Max-Flow/Min-Cut Theorem such $f$ exists if the minimum capacity over all $s-t$ cuts $(S,T)$ of $H$ is at least $a$, that is $c(S,T)\ge a$ for all $S,T$ with $s\in S$ and $t\in T$. Now we show these cut constraints are equivalent to (iii) as follows. 

Note that it suffices only to consider cuts $(S,T)$ that do not contain an edge from $A$ to $B$ (since those edges have unlimited capacity). Let $A:= \{e_{ii}:i\in [4]\}$ and $B:=\{e_{ij}:i\in\{1,2\},~j\in\{3,4\}\}$. This holds by construction for $S=\{s\}$ and $S=V(H)\setminus \{t\}$. First we note that (iii)(a) is equivalent to the cut condition $c(S,T)\ge a$ when $S:=\{s,e_{ii},e_{ik},e_{i\ell}\}$ where $\{\ell,k\}=\{3,4\}$ if $i\in\{1,2\}$ and $\{\ell,k\}= \{1,2\}$ if $i\in \{3,4\}$. Second we note that (iii)(b) is equivalent to the cut condition $c(S,T)\ge a$ when $T:= \{t,e_{ii},e_{kk},e_{ik}\}$ and $i\in\{1,2\}$, $k\in\{3,4\}$. Finally we note that (iii)(c) is equivalent to the cut condition $c(S,T)\ge a$ when $S:= \{s,e_{ii},e_{kk},e_{ik},e_{jk},e_{i\ell}\}$ where $i\in \{1,2\}, k\in\{3,4\}$ and $j:=[2]\setminus \{i\}$ and $\ell:=\{3,4\}\setminus \{k\}$. Since these are all the cuts not containing an edge from $A$ to $B$, this concludes the proof.
\end{proof}

The following lemma characterizes divisibility needed to prove Theorem~\ref{thm:2cliqueandbip}.

\begin{lem}\label{lem:F2div}
Let $\mc{F}_2=\{T_{113}, T_{114}, T_{223}, T_{224},T_{134},T_{234}\}$. An $\mc{F}_2$-blowup $G$ is fractionally $\mc{F}_2$-divisible if and only if $G$ satisfies all of the following:
\begin{enumerate}
    \item[(i)] $e_{13}(G)+e_{14}(G)+e_{23}(G)+e_{34}(G) = 2\cdot (e_{11}(G)+e_{22}(G)+e_{34}(G))$
    \item[(ii)] For $i\in \{1,2\}$ and $v\in V_i$, 
    \begin{itemize}
        \item[(a)] $d_i(v) \le d_{3}(v)+d_4(v)$,
        \item[(b)] $d_{3}(v)\le d_i(v)+d_4(v)$,
        \item[(c)] $d_4(v)\le d_i(v)+d_3(v)$.
    \end{itemize} 
    For $i\in \{3,4\}$ and $v\in V_i$, $d_{7-i}(v)\le d_1(v)+d_2(v)$.     \item[(iii)]  $2\cdot e_{11}(G) \le e_{13}(G)+e_{14}(G) \le 2\cdot (e_{11}(G)+e_{34}(G)) $
    \item[(iv)] $|e_{13}(G)-e_{14}(G)|\le 2\cdot e_{11}(G)$
    \item[(v)] $|e_{23}(G)-e_{24}(G)|\le 2\cdot e_{22}(G)$.
\end{enumerate}
An $\mc{F}_2$-blowup $G$ is integrally $\mc{F}_2$-divisible if and only if $G$ satisfies all of the following:
\begin{enumerate}
    \item[(i')] $e_{13}(G)+e_{14}(G)+e_{23}(G)+e_{34}(G) = 2\cdot (e_{11}(G)+e_{22}(G)+e_{34}(G))$
    \item[(ii')] $G$ is Eulerian
    \item[(iii')] $2e_{11}(G) + e_{13}(G) \equiv e_{14}(G) \pmod 4$.
\end{enumerate}
\end{lem}

\begin{remark}
We note for the reader that some of the conditions for fractional are redundant. Namely (iii),(iv) and (v) are all redundant by summing degrees as follows: $(ii)(a)$ for $i=1$ implies the lower bound in (iii) and $(ii)(a)$ for $i=2$ and (i) imply the upper bound; similarly (ii)(b)-(c) for $i=1$ implies (iv); and (ii)(b)-(c) for $i=2$ implies (v); nevertheless we keep the redundant conditions for the understanding of the reader.
\end{remark}

\begin{proof}
First we consider the degree constraints. For $v\in V(T_{113})$ with $\mc{P}(v)=1$, we have that $\vec{d}_{T_{113}}(v) = (1,0,1,0)$; for $v\in V(T_{114})$ with $\mc{P}(v)=1$, we have that $\vec{d}_{T_{114}}(v) = (1,0,0,1)$; for $v\in V(T_{134})$ with $\mc{P}(v)=1$, we have that $\vec{d}_{T_{134}}(v) = (0,0,1,1)$. By Proposition~\ref{prop:LatticeCone}, we find that $\textrm{cone}( (1,0,1,0), (1,0,0,1), (0,0,1,1) ) = \{(x_1,0,x_3,x_4)\in \mathbb{R}_{\ge 0}^4: x_1\le x_3+x_4, x_3\le x_1+x_4, x_4\le x_1+x_3\}$ and $\textrm{lattice}( (1,0,1,0), (1,0,0,1), (0,0,1,1) )=\{(x_1,0,x_3,x_4)\in \mathbb{Z}^4: x_1+x_3+x_4\equiv 0\pmod 2\}$. These are then the necessary divisibility conditions for $\vec{d}(v)$ for $v\in V_1$. By symmetry similar symmetric conditions hold for $v\in V_2$. 

So we consider $v\in V_3$ as follows.  For $v\in V(T_{113})$ with $\mc{P}(v)=3$, we have that $\vec{d}_{T_{113}}(v) = (2,0,0,0)$; for $v\in V(T_{223})$ with $\mc{P}(v)=3$, we have that $\vec{d}_{T_{223}}(v) = (0,2,0,0)$; for $v\in V(T_{134})$ with $\mc{P}(v)=3$, we have that $\vec{d}_{T_{134}}(v) = (1,0,0,1)$; for $v\in V(T_{234})$ with $\mc{P}(v)=3$, we have that $\vec{d}_{T_{234}}(v) = (0,1,0,1)$. By Proposition~\ref{prop:LatticeCone}, we find that $\textrm{cone}( (2,0,0,0), (0,2,0,0), (1,0,0,1), (0,1,0,1) ) = \{(x_1,x_2,0,x_4)\in \mathbb{R}_{\ge 0}^4: x_4\le x_1+x_2\}$ and $\textrm{lattice}( (2,0,0,0), (0,2,0,0), (1,0,0,1), (0,1,0,1) )=\{(x_1,x_2,0,x_4)\in \mathbb{Z}^4: x_1+x_2+x_4\equiv 0\pmod 2\}$. These are then the necessary divisibility conditions for $\vec{d}(v)$ for $v\in V_3$. By symmetry similar symmetric conditions hold for $v\in V_4$.

Hence condition (ii) is both necessary and sufficient for degree vectors for fractional $\mc{F}_2$-divisibility while condition (ii') is both necessary and sufficient for degree vectors for integral $\mc{F}_2$-divisibility. 

Finally we consider the edge constraints. this proof when considering a $4$-partitioned graph $H$, we write $\vec{e}(H)$ as $(e_{11}(H),e_{22}(H),e_{34}(H),e_{13}(H),e_{14}(H),e_{23}(H),e_{24}(H))$. Thus $\vec{e}(T_{113})=(1,0,0,2,0,0,0)$, $\vec{e}(T_{114})=(1,0,0,0,2,0,0)$, $\vec{e}(T_{223})=(0,1,0,0,0,2,0)$, $\vec{e}(T_{224})=(0,2,0,0,0,0,2)$, $\vec{e}(T_{134})=(0,0,1,1,1,0,0)$, $\vec{e}(T_{234})=(0,0,1,0,0,1,1)$. We let $A$ be the matrix whose columns are these generators. Then $G$ is fractionally $\mc{F}_2$-divisible if and only if $\vec{e}(G)\in \textrm{cone}(A)$ and (ii) holds. Similarly $G$ is integrally $\mc{F}_2$-divisible if and only if $\vec{e}(G)\in \textrm{lattice}(A)$ and (ii') holds.

We note that (i)/(i') is a necessary condition for $\vec{e}(G)\in \textrm{span}(A)$(which lie in $\mathbb{R}^7$) since it is satisfied by all the generator vectors; since $\textrm{rank}(A)=6$, the condition (i)/(i') is also a sufficient condition for $\vec{e}(G)$ to lie in $\textrm{span}(A)$. Similarly one verifies the remaining conditions are all necessary as they are satisfied by the generator vectors. Quotienting out by this equation, we proceed to describe the solution space namely of the vector $\vec{t}$ whose coordinates $t_{ijk}$ where this represents the number of triangles of type $T_{ijk}$. We the solve the associated matrix equation (which then has full rank), to find that 
\begin{align*}
t_{134} &= \frac{e_{13}(G)+e_{14}(G) -2\cdot e_{11}(G)}{2} \\
t_{234} &= \frac{e_{23}(G)+e_{24}(G) -2\cdot e_{22}(G)}{2} \\
t_{113} &= \frac{1}{2} (e_{13}(G) - t_{134}) = \frac{e_{13}(G) + 2\cdot e_{11}(G) - e_{14}(G)}{4} \\ 
t_{114} &= \frac{1}{2} (e_{14}(G) - t_{134}) = \frac{e_{14}(G) + 2\cdot e_{11}(G) -  e_{13}(G)}{4} \\ 
t_{223} &= \frac{1}{2} (e_{23}(G) - t_{234}) = \frac{e_{23}(G) + 2\cdot e_{22}(G) - e_{24}(G)}{4} \\ 
t_{224} &= \frac{1}{2} (e_{24}(G) - t_{234}) = \frac{e_{24}(G) + 2\cdot e_{22}(G) - e_{23}(G)}{4}
\end{align*}
For the lattice, $t_{134}, t_{234}$ being integral follows from  (iii') while for $t_{113},t_{114},t_{223},t_{224}$, their integrality follows from (ii') and (iii'). For the cone, $t_{134}\ge 0$ by (iii), $t_{234}\ge 0$ by (i) and (iii), $t_{113}, t_{114}\ge 0$ by (iv), and $t_{223}, t_{224}\ge 0$ by (v). 
\end{proof}

Finally our last lemma in this vein characterizes the divisibility conditions which we require in our proof of Theorem~\ref{thm:K4blowup}.

\begin{lem}\label{lem:F3div}
Let $\mc{F}_3=\{T_{123}, T_{124}, T_{134}, T_{234}\}$. An $F$-blowup $G$ is fractionally $\mc{F}_3$-divisible if and only if $G$ satisfies both of the following:
\begin{enumerate}
    \item[(i)] $e_{12}(G)+e_{34}(G) = e_{13}(G) + e_{24}(G) = e_{14}(G)+e_{23}(G)$, and
    \item[(ii)] $d_{j}(v) \le d_k(v)+d_{\ell}(v)$ for $v\in V_i$, $i\in [4]$ and $\{i,j,k,\ell\}=[4]$.
\end{enumerate}
An $\mc{F}_3$-blowup $G$ is integrally $\mc{F}_3$-divisible if and only if $G$ satisfies both of the following:
\begin{enumerate}
    \item[(i')] $e_{12}(G)+e_{34}(G) = e_{13}(G) + e_{24}(G) = e_{14}(G)+e_{23}(G)$, and
    \item[(ii')] $G$ is Eulerian.
\end{enumerate}
\end{lem}
\begin{proof}
First we consider the degree constraints. So we consider $v\in V_1$ as follows.  For $v\in V(T_{123})$ with $\mc{P}(v)=1$, we have that $\vec{d}_{T_{123}}(v) = (0,1,1,0)$; for $v\in V(T_{124})$ with $\mc{P}(v)=1$, we have that $\vec{d}_{T_{124}}(v) = (0,1,0,1)$; for $v\in V(T_{134})$ with $\mc{P}(v)=1$, we have that $\vec{d}_{T_{134}}(v) = (0,0,1,1)$. By Proposition~\ref{prop:LatticeCone}, we find that $\textrm{cone}( (0,1,1,0), (0,1,0,1), (0,0,1,1) ) = \{(0,x_2,x_3,x_4)\in \mathbb{R}_{\ge 0}^4: x_2\le x_3+x_4,~x_3\le x_2+x_4,~x_4\le x_2+x_3 \}$ and $\textrm{lattice}( (0,1,1,0), (0,1,0,1), (0,0,1,1) )=\{(0,x_2,x_3,x_4)\in \mathbb{Z}^4: x_2+x_3+x_4\equiv 0\pmod 2\}$. By symmetry similar symmetric conditions hold for $v\in V_i$ for $i\in \{2,3,4\}$. Hence condition (ii) is both necessary and sufficient for degree vectors for fractional $\mc{F}_3$-divisibility while condition (ii') is both necessary and sufficient for degree vectors for integral $\mc{F}_3$-divisibility.

Finally we consider the edge constraints. this proof when considering a $4$-partitioned graph $H$, we write $\vec{e}(H)$ as $(e_{12}(H),e_{13}(H),e_{14}(H),e_{23}(H),e_{24}(H),e_{34}(H))$. Thus $\vec{e}(T_{123})=(1,1,0,1,0,0)$, $\vec{e}(T_{124})=(1,0,1,0,1,0)$, $\vec{e}(T_{134})=(0,1,1,0,0,1)$, $\vec{e}(T_{234})=(0,0,0,1,1,1)$. We let $A$ be the matrix whose columns are these generators. Then $G$ is fractionally $\mc{F}_3$-divisible if and only if $\vec{e}(G)\in \textrm{cone}(A)$ and (ii) holds. Similarly $G$ is integrally $\mc{F}_3$-divisible if and only if $\vec{e}(G)\in \textrm{lattice}(A)$ and (ii') holds.

We note that (i)/(i') is a necessary condition for $\vec{e}(G)\in \textrm{span}(A)$(which lie in $\mathbb{R}^6$) since it is satisfied by all the generator vectors; since $\textrm{rank}(A)=4$, the condition (i)/(i') is also a sufficient condition for $\vec{e}(G)$ to lie in $\textrm{span}(A)$. Similarly one verifies the remaining conditions are all necessary as they are satisfied by the generator vectors. Quotienting out by this equation, we proceed to describe the solution space namely of the vector $\vec{t}$ whose coordinates $t_{ijk}$ where this represents the number of triangles of type $T_{ijk}$. We the solve the associated matrix equation (which then has full rank), to find that 
\begin{align*}
t_{123} &= \frac{e_{12}(G)+e_{13}(G)-e_{14}(G)}{2} \\
t_{124} &= \frac{e_{12}(G)+e_{14}(G)-e_{13}(G)}{2} \\
t_{134} &= \frac{e_{13}(G)+e_{14}(G)-e_{12}(G)}{2} \\
t_{234} &= \frac{e_{23}(G)+e_{24}(G)-e_{12}(G)}{2}
\end{align*}    
For the lattice, these being integral follows from (ii'). For the cone, these being nonnegative follows from (ii). 
\end{proof}

\section{Absorber Stability and Idealized Decomposition Theorems}\label{s:AbsStabDecomp}

\subsection{Absorber Stability}\label{ss:AbsorberStability}

In this section, we characterize which of our graphs have many absorbers (specifically many hinges, as anti-edges will just follow). First we require the following useful lemma characterizing when there exist many paths of length $3$ between two vertices in a graph of minimum degree roughly half the vertices. We note this lemma is similar in spirit to Lemma 26 of K\"uhn, Lapinskas, and Osthus~\cite{KLO13} from 2013 in their work on optimal packings of Hamilton cycles in graphs of high minimum degree. Said lemma was also used in the work of Csaba, K\"uhn, Lo, Osthus, and Treglown~\cite{CKLOT16} in their proof of the $1$-Factorization and Hamilton Decomposition Conjectures. Nevertheless, we include a proof for completeness. 

\begin{lem}\label{lem:Many3Paths}
For each $\alpha \in (0,1)$, the following holds for sufficiently large $n$: Let $\varepsilon \in [0,\frac{\alpha}{3}]$, let $G$ be a graph on $n$ vertices with $\delta(G) \ge (\frac{1}{2}-\varepsilon)n$, and let $u\ne v\in V(G)$. If there do not exist at least $\frac{a^2}{6}\cdot n^2$ distinct paths of length $3$ from $u$ to $v$ in $G$, then $d_G(u),d_G(v) \le (\frac{1}{2}+2\alpha)n$ and at least one of the following holds:  
\begin{itemize}
    \item[(i)] $|N(u)\cap N(v)|\le \alpha\cdot n$ and $e_G(N(u)\setminus N(v),~ N(v)\setminus N(u)) \le \frac{\alpha^2}{6}\cdot n^2$ and $G$ is $4\alpha$-close to $\overline{T_2(n)}$, or 
    \item[(ii)] $|N(u)\cup N(v)|\le (\frac{1}{2}+\alpha)n$ and $e(G[N(u)\cap N(v)])\le \frac{\alpha^2}{6}\cdot n^2$.
\end{itemize}    
\end{lem}
\begin{proof}
Suppose not. First suppose $|N(u)\cap N(v)|\le \alpha\cdot n$. If $e_G(N(u)\setminus N(v),~ N(v)\setminus N(u)) \ge \frac{\alpha}{6}\cdot n^2$, then there exist at least $\frac{\alpha^2}{6}\cdot n^2$ paths of length $3$ from $u$ to $v$ in $G$, a contradiction. So we assume $e_G(N(u)\setminus N(v),~ N(v)\setminus N(u)) \le \frac{\alpha^2}{6}\cdot n^2$. Note then $|N(u)\setminus N(v)|, |N(v)\setminus N(u)|\ge (\frac{1}{2}-\varepsilon-\alpha)n$. Let $A$ be a subset of $N(u)\setminus N(v)$ of size $(\frac{1}{2}-\varepsilon-\alpha)n$ and let $B$ be a subset of $N(v)\setminus N(u)$ of size $(\frac{1}{2}-\varepsilon-\alpha)n$. Thus we find that $e_G(A,B)\le e_G(N(u)\setminus N(v),~N(v)\setminus N(u)) \le \frac{\alpha^2}{6}\cdot n^2$. This also implies that 
$$e(G[A]) \ge \frac{1}{2}\cdot |A|\bigg( \left(\frac{1}{2}-\varepsilon\right)n - (n-|A|-|B|)\bigg) - e_G(A,B) \ge \frac{1}{2}\cdot |A|\left(\frac{1}{2}-3\varepsilon-2\alpha\right)n - \frac{\alpha^2}{6}n^2.$$ 
Hence $A$ can be made a clique by adding at most $(\alpha+2\varepsilon)n\cdot \frac{|A|}{2} \le \frac{\alpha+2\varepsilon}{4}\cdot n^2$ edges and symmetrically the same for $B$. Let $m:= (\frac{1}{2}-\varepsilon-\alpha)n$. But then $G[A\cup B]$ can be transformed to $K_{m,m}$ by adding and/or deleting at most $(\frac{\alpha+2\varepsilon}{2} + \frac{\alpha^2}{6}) n^2 \le \alpha\cdot n^2$ edges. Hence $G$ can be transformed to $\overline{T_2(n)}$ by adding and/or deleting at most $\alpha\cdot n^2 + n(n-|A|-|B|) \le (3\alpha+2\varepsilon)n^2 \le 4\alpha\cdot n^2$ edges. Thus $G$ is $4\alpha$-close to $\overline{T_2(n)}$. Finally we note that $d_G(u) \le n-|N(v)\setminus N(u)| \le (\frac{1}{2}+\varepsilon+\alpha)n\le (\frac{1}{2}+2\alpha)n$ and similarly for $v$. Thus (i) holds as desired.

So we assume $|N(u)\cap N(v)| > \alpha\cdot n$. Next suppose $|N(u)\cup N(v)|\le (\frac{1}{2}+\alpha)n$. If $e(G[N(u)\cap N(v)])\ge \frac{\alpha^2}{6}\cdot n^2$, then there exist at least $\frac{\alpha^2}{6}\cdot n^2$ paths of length $3$ from $u$ to $v$ in $G$, a contradiction. So we assume $e(G[N(u)\cap N(v)])\le \frac{\alpha^2}{6}\cdot n^2$. But then (ii) holds as desired. 

So we assume $|N(u)\cup N(v)|\ge (\frac{1}{2}+\alpha)n$. Let $R_1:= \{ w\in N(u)\cap N(v): |N(w)\cap N(u)|\ge \frac{\alpha}{3}n\}$ and $R_2:= \{ w\in N(u)\cap N(v): |N(w)\cap N(v)|\ge \frac{\alpha}{3}n \}$. For $w\in N(u)\cap N(v)$, we have that $|N(w)\cap (N(u)\cup N(v))| \ge |N(w)|+|N(u)\cup N(v)|-n \ge (\frac{1}{2}-\varepsilon)n + (\frac{1}{2}+\alpha)n - n \ge \frac{2}{3}\alpha n$; but then $\max\{|N(w)\cap N(u)|~|N(w)\cap N(v)|\} \ge \frac{\alpha}{3}n$ and so $w\in R_1\cup R_2$. It follows that $N(u)\cap N(v) = R_1\cup R_2$. Hence there exists $i\in [2]$ such that $|R_i|\ge \frac{|N(u)\cap N(v)|}{2} \ge \frac{\alpha}{2}\cdot n$. But then there exists at least $|R_i|\cdot \frac{\alpha}{3}\cdot n \ge \frac{\alpha^2}{6} \cdot n^2$ paths of length three from $u$ to $v$ in $G$, a contradiction. 
\end{proof}

We now use the above lemma to characterize when there exist many paths of length $i$ for each $i\ge 4$ between two vertices in a graph of minimum degree roughly half the vertices as follows.

\begin{lem}\label{lem:ManyPaths}
For each $i\ge 4$ and $\alpha \in (0,\frac{1}{8}]$, the following holds for sufficiently large $n$: Let $\varepsilon \in [0,\frac{\alpha}{100}]$ and let $G$ be a graph on $n$ vertices with $\delta(G) \ge (\frac{1}{2}-\varepsilon)n$. Then for all distinct $u,v \in V(G)$, one of the following holds:
\begin{itemize}
    \item[(i)] there exist at least $\frac{\alpha^3}{50} \cdot (\frac{n}{3})^{i-1}$ paths of length $i$ from to $u$ to $v$ in $G$, or 
    \item[(ii)] $G$ is $\alpha$-close to $T_2(n)$ and $d_G(u), d_G(v) \le (\frac{1}{2}+\alpha)n$, or
    \item[(iii)] $G$ is $\alpha$-close to $\overline{T_2(n)}$ and $d_G(u),d_G(v) \le (\frac{1}{2}+\alpha)n$.
\end{itemize} 
\end{lem}
\begin{proof}
We proceed by induction on $i$. First we prove the base case that $i=4$.
\vskip.1in
\noindent \textbf{Base case $i=4$~:}
\vskip.1in
First suppose that at least one of $u$ or $v$ has degree more than $(\frac{1}{2}+\frac{\alpha}{2})n$. Suppose without loss of generality that $v$ does. But then for each $u'\in N(u)$, we have by Lemma~\ref{lem:Many3Paths} with $u'$, $v$ and $\alpha'$ that there exist at least $\frac{(\alpha')^2}{6}\cdot n^2$ paths of length $3$ from $u'$ to $v$ in $G$; note that $u$ is in at most $2n$ of these paths for a fixed $u'$. Hence there exist at least $\left(\frac{(\alpha')^2}{6}\cdot n^2-2n\right)\cdot (\frac{1}{2}-\varepsilon)n \ge \frac{\alpha^2}{400} \cdot n^3 \ge \frac{\alpha^2}{15} \cdot (\frac{n}{3})^3$ paths of length $4$ from $u$ to $v$ in $G$ and (i) holds as desired. 

So we assume $d_G(u), d_G(v) \le (\frac{1}{2}+\frac{\alpha}{2})n$. Let $\alpha':= \frac{\alpha}{6}$. Let $U:= \{u'\in N(u): |N(u')\cup N(v)|\le (\frac{1}{2}+\alpha')n\textrm{ and } e(G[N(u)\cap N(v)])\le \frac{(\alpha')^2}{6}\cdot n^2\}$. First suppose that $|N(u)\setminus U|\ge \frac{\alpha}{6}\cdot n$.  Apply Lemma~\ref{lem:Many3Paths} with $\alpha'$ and $v$ for each $u'\in U$. Since $u'\in U$, we have that Lemma~\ref{lem:Many3Paths}(ii) does not hold. If Lemma~\ref{lem:Many3Paths}(i) holds for any $u'\in U$, then $G$ is $4\alpha'$-close (and hence $\alpha$-close) to $\overline{T_2(n)}$ and (iii) holds as desired. Thus by Lemma~\ref{lem:Many3Paths}, there exist at least $\frac{(\alpha')^2}{6}\cdot n^2$ paths of length $3$ from $u'$ to $v$ in $G$; note that $u$ is in at most $2n$ of these paths for a fixed $u'$. Hence there exist at least $\left(\frac{(\alpha')^2}{6}\cdot n^2-2n\right)\cdot \frac{\alpha}{6} n \ge \frac{\alpha^3}{1300} \cdot n^3 \ge \frac{\alpha^3}{50} \cdot (\frac{n}{3})^3$ paths of length $4$ from $u$ to $v$ in $G$ where the last inequality follows since $n$ is large enough. Thus (i) holds as desired. 

So we assume $|N(u)\setminus U|\le \frac{\alpha}{6}\cdot n$. Hence 
$$|U| \ge |N(u)|- |N(u)\setminus U| \ge \left(\frac{1}{2} - \varepsilon - \frac{\alpha}{6}\right)n.$$ 
Since $\varepsilon \le \frac{\alpha}{12}$, we then find that $U\ne \emptyset$. Note that for each $w\in U$, we have that 
$$|N(w)\cap N(v)|\ge |N(w)|+|N(v)| - |N(w)\cup N(v)| \ge 2\left(\frac{1}{2}-\varepsilon\right)n - \left(\frac{1}{2}+\alpha'\right)n \ge \left(\frac{1}{2} - 2\varepsilon-\alpha'\right)n .$$
Fix $u'\in U$ and let $A$ be a subset of $N(u')\cap N(v)$ of size $(\frac{1}{2} - 2\varepsilon-\alpha')n$ (such exists from above). Since $u'\in U$, we have that 
$$e(G[A])\le e(G[N(u')\cap N(v)]) \le \frac{(\alpha')^2}{6}\cdot n^2 = \frac{\alpha^2}{216}\cdot n^2$$
and hence $A$ can be made independent by deleting at most that many edges. 

Recall that for each $w\in U$, we have that $|N(w)\cap N(v)|\ge (\frac{1}{2} - 2\varepsilon-\alpha')n$. But then 
$$|N(w)\cap A| \ge |A|-(|N(v)|-|N(w)\cap N(v)|)\ge |A| - (2\varepsilon+3\alpha')n\ge \left(\frac{1}{2}-4\varepsilon - 4\alpha'\right)n.$$ Since $\varepsilon\le \alpha' \le \frac{1}{32}$, we also have that $|N(w)\cap A|\ge \frac{n}{4}$. Thus $e(G[A])\ge \frac{1}{2}\cdot |U\cap A|\cdot \frac{n}{4}$. Since $e(G[A])\le \frac{(\alpha')^2}{6}\cdot n^2$, these together imply that 
$$|U\cap A|\le \frac{\alpha^2}{25}\cdot n.$$ 
Thus 
$$|U\setminus A| \ge |U|-|U\cap A| \ge  \left(\frac{1}{2} - \varepsilon - \frac{\alpha}{6} - \frac{\alpha^2}{25}\right)n.$$ 
Let $B$ be a subset of $U\setminus A$ of size $\left(\frac{1}{2} - \varepsilon - \frac{\alpha}{6} - \frac{\alpha^2}{25}\right)n$. Recalling that for each $w\in B$ we have that $|N(w)|\le |N(w)\cup N(v)|\le (\frac{1}{2}+\alpha')n$, we find that
$$e(G[B])\le \frac{1}{2}\cdot \sum_{w\in B} |N(w)|-|N(w)\cap A| \le \frac{1}{2}\cdot |B| \cdot (4\varepsilon+5\alpha')n \le \frac{4\varepsilon+5\alpha'}{4}\cdot n^2,$$
where we used that $|B|\le \frac{n}{2}$. Thus $B$ can be made independent by deleting at most the above number of edges.

Furthermore
$$|A|\cdot |B|-e_G(A,B)= \sum_{w\in B}|A|-|N(w)\cap A| \le |B|\cdot (2\varepsilon+3\alpha')n \le \frac{2\varepsilon+3\alpha'}{2}\cdot n^2.$$
Thus the cut $(A,B)$ can be made complete by adding at most the above number of edges. Combining we find that $G[A\cup B]$ can be transformed to $K_{|A|,|B|}$ by adding and/or deleting at most the following number of edges:
$$e(G[A])+e(G[B]) + (|A|\cdot |B| - e_G(A,B)) \le \left(\frac{\alpha^2}{216} + \left(\varepsilon+\frac{5}{4}\alpha'\right) + \left(\varepsilon+\frac{3}{2}\alpha'\right) \right)\cdot n^2 \le \frac{\alpha}{2}\cdot n^2,$$
where we used that $\varepsilon \le \frac{\alpha}{100}$ and that $\alpha' = \frac{\alpha}{6}$.

On the other hand 
$$n-|A|-|B|\le \left(2\varepsilon+2\cdot \frac{\alpha}{6}+\frac{\alpha^2}{25}\right)n^2 \le \frac{\alpha}{2}\cdot n^2,$$
where we used that $\varepsilon\le \frac{\alpha}{100}$. 
Thus $G$ can be transformed to $T_2(n)$ by adding and/or deleting a number of edges at most 
$$e(G[A])+e(G[B])+(|A|\cdot |B|-e_G(A,B)) + (n-|A|-|B|)n \le \frac{\alpha}{2}\cdot n^2 + \frac{\alpha}{2}\cdot n^2 \le \alpha\cdot n^2.$$ 
Hence $G$ is $\alpha$-close to $T_2(n)$ and (ii) holds as desired.

\vskip.1in
\noindent \textbf{Inductive case $i\ge 5$~:}
\vskip.1in

So we assume $i\ge 5$. First suppose that at least one of $u$ or $v$ has degree more than $(\frac{1}{2}+\frac{\alpha}{2})n$. Suppose without loss of generality that $v$ does. But then for each $u'\in N(u)$, we have by induction with $u'$, $v$ and $\alpha$ that there exist at least $\frac{\alpha^3}{50}\cdot (\frac{n}{3})^{i-2}$ paths of length $i-1$ from $u'$ to $v$ in $G$; note that $u$ is in at most $(i-2)n$ of these paths for a fixed $u'$. Hence there exist at least $\left(\frac{\alpha^3}{50}\cdot n^{i-2}-(i-2)n\right)\cdot (\frac{1}{2}-\varepsilon)n \ge \frac{\alpha^3}{50} \cdot (\frac{n}{3}) \ge \frac{\alpha^3}{50} \cdot (\frac{n}{3})^{i-1}$ paths of length $i$ from $u$ to $v$ in $G$ where the last inequality follows since and (i) holds as desired. 

So we assume $d_G(u), d_G(v) \le (\frac{1}{2}+\frac{\alpha}{2})n$. Thus if $G$ is $\alpha$-close to $T_2(n)$ then (ii) holds as desired. Similarly if $G$ is $\alpha$-close to $\overline{T_2(n)}$, then (iii) holds as desired. So we assume that $G$ is not $\alpha$-close to either $T_2(n)$ or $\overline{T_2(n)}$. But then for each $u'\in N(u)$, we have by induction with $u'$, $v$ and $\alpha$ that there exist at least $\frac{\alpha^3}{50}\cdot (\frac{n}{3})^{i-2}$ paths of length $i-1$ from $u'$ to $v$ in $G$; note that $u$ is in at most $(i-2)n$ of these paths for a fixed $u'$. Hence there exist at least $\left(\frac{\alpha^3}{50}\cdot n^{i-2}-(i-2)n\right)\cdot (\frac{1}{2}-\varepsilon)n \ge \frac{\alpha^3}{50} \cdot (\frac{n}{3}) \ge \frac{\alpha^3}{50} \cdot (\frac{n}{3})^{i-1}$ paths of length $i$ from $u$ to $v$ in $G$ where the last inequality follows since and (i) holds as desired. 
\end{proof}

Now we are prepared to characterize when there exist many $K_3$-hinges in graph $G$ on $n$ vertices with $\delta(G)\ge \frac{3}{4}n$ as follows, but first a definition that will help yield our desired hinges.

\begin{definition}
 Let $g\ge 2$ be an integer. The \emph{$g$-sphere} is the graph $\overline{K_2}+C_{2g}$.    
\end{definition}

We note that the $g$-sphere has two $K_3$-decompositions $\mc{Q}_1,\mc{Q}_2$ such that $\mc{Q}_1\cap \mc{Q}_2=\emptyset$ (i.e.~no triangle is in both decompositions).  Note then that for any two triangles $T_1, T_2$ in the $g$-sphere $S_g$ such that $T_1\cap T_2$ is an edge, we have that $S_g - (T_1\cup T_2)$ is a hinge for $T_1$ and $T_2$.

\begin{lem}\label{lem:ManyHinges}
There exists $\alpha_0 > 0$ such that the following holds for large enough $n$: Let $\alpha \in (0,\alpha_0]$ and let $G$ be a graph on $n$ vertices with $\delta(G) \ge \frac{3}{4}\cdot n$. If $G$ is not $\alpha$-close to $T_2(n/2)+T_2(n/2)$ or $\overline{T_2(n/2)}+T_2(n/2)$ or $\overline{T_2(n/2)}+\overline{T_2(n/2)}$, then $G$ is $(\frac{\alpha^3}{10^{14}},8)$-$K_3$-hingeful. 
\end{lem}
\begin{proof}
Suppose not. Let $\{u_1,u_2\}:= V(T_1)\cap V(T_2)$, $v_1 \in V(T_1)\setminus V(T_2)$ and $v_2\in V(T_2)\setminus V(T_1)$. Let $A$ be a subset of $N(v_1)\cap N(v_2)$ of size $\frac{n}{2}$ such that $u_1,u_2\in A$. Let $H_1:= G[A]$ (note such $A$ exists since $|N(v_1)\cap N(v_2)|\ge \frac{n}{2}$ as $\delta(G) \ge \frac{3}{4}n$). Furthermore, $\delta(H_1)\ge v(H_1) -\frac{n}{4} \ge \frac{v(H_1)}{2}$. 

Let $\alpha':=\frac{\alpha}{800}$. If there exist at least $\frac{(\alpha')^3}{50} \cdot \left(\frac{v(H_1)}{3}\right)^4$ paths of length five in $H_1$ from $u_1$ to $u_2$, then each of these yields a $K_3$-hinge for $T_1$ and $T_2$ (namely the $8$-sphere); hence $G$ is $(\frac{\alpha^3}{50\cdot 800^3\cdot 6^4},8)$-$K_3$-hingeful and hence is $(\frac{\alpha^3}{10^{14}},8)$-$K_3$-hingeful as desired. Thus by Lemma~\ref{lem:ManyPaths} with $i=5$, we find that $H_1$ is $\alpha'$-close to either $T_2(n/2)$ or $\overline{T_2(n/2)}$ and that $d_{H_1}(u_1),d_{H_1}(u_2) \le (\frac{1}{2}+\alpha')\frac{n}{2}$. 

Let $U_1:= \{ u_1':\in N_{H_1}(u_1): d_{H_1}(u_1')\le (\frac{1}{2}+4\alpha')\cdot \frac{n}{2}\}$. Note that $|U_1|\ge \frac{n}{4}$ as otherwise one has to delete strictly more than $\frac{1}{2}\cdot (\frac{n}{2}-|U_1|)\cdot 4\alpha'\cdot v(H_1) \ge \alpha'\cdot v(H_1)^2$ edges to transform $H_1$ to $T_2(n/2)$ or $\overline{T_2(n/2)}$, contradicting that $H_1$ is $\alpha'$-close to at least one of $T_2(n/2)$ or $\overline{T_2(n/2)}$.

Let $\alpha'':= \frac{\alpha}{2}$. First suppose that there exists $u_1'\in U_1$ such that there does not exist at least $\frac{(\alpha'')^3}{50}\cdot (\frac{n}{9})^3$ paths of length four from $v_1$ to $v_2$ in 
$$H_2:=G[N(u_1')\cap N(u_2)\cap (V(G)\setminus V(H_1))].$$
Since $\delta(G)\ge \frac{3}{4}n$ we find that $|N_G(u_2)\setminus V(H_1)| \ge (\frac{1}{2}-\frac{\alpha'}{2})n$ and similarly that $|N_G(u_1')\setminus V(H_1)| \ge (\frac{1}{2}-2\alpha')n$. Hence 
$$v(H_2)\ge \left(\frac{1}{2}-3\alpha'\right)n \ge \frac{n}{3},$$
where we used that $\alpha'\le \frac{1}{18}$. Furthermore, 
$$\delta(H_2)\ge v(H_2) - \frac{n}{4} \ge \left(\frac{1}{2}-4\alpha'\right)\cdot v(H_2),$$
where the last inequality follows since $\frac{n}{4} \le (\frac{1}{2}+4\alpha')v(H_2)$ as $(\frac{1}{2}+4\alpha')(\frac{1}{2}-3\alpha') \ge \frac{1}{4}$ (since $\alpha'\le \frac{1}{24}$). Let $\varepsilon:=4\alpha'$. Note that $\varepsilon \le \frac{\alpha''}{100}$. 
Thus by Lemma~\ref{lem:ManyPaths}, we find that $H_2$ is $\alpha''$-close to either $T_2(v(H_2))$ or $\overline{T_2(v(H_2))}$. But then $G$ is $(2\alpha'+\alpha'')$-close (and hence $\alpha$-close) to $T_2(n/2)+T_2(n/2)$ or $\overline{T_2(n/2)}+T_2(n/2)$ or $\overline{T_2(n/2)}+\overline{T_2(n/2)}$, a contradiction.

So we assume for each $u_1'\in U_1$ there exist at least $\frac{(\alpha'')^3}{50}\cdot (\frac{n}{9})^3$ paths of length four from $v_1$ to $v_2$ in $G[N(u_1')\cap N(u_2)\cap (V(G)\setminus  V(H_1))]$.  then there exist at least $\frac{(\alpha'')^3}{50}\cdot (\frac{n}{9})^3\cdot |U_1| \ge \frac{\alpha^3}{50\cdot 2^3\cdot 9^3\cdot 4}\cdot n^4$ $K_3$-hinges for $T_1$ and $T_2$ (namely $8$-spheres) and hence $G$ is $(\frac{\alpha^3}{10^{14}},8)$-$K_3$-hingeful as desired. 
\end{proof}

\subsection{Proof of Non-Extremal Case}\label{ss:NonExtremal}

We are now prepared to prove Theorem~\ref{thm:nonextremal}.

\begin{proof}[Proof of Theorem~\ref{thm:nonextremal}]
Since $\delta(G)\ge 3n/4$, every edge is in at least $n/2$ triangles, that is $\frac{1}{2}$-$K_3$-anti-edgeful. By Lemma~\ref{lem:ManyHinges} as $G$ is not $\beta$-close to any of $T_2(n/2)+T_2(n/2)$ or $\overline{T_2(n/2)}+T_2(n/2)$ or $\overline{T_2(n/2)}+\overline{T_2(n/2)}$ as $G$ is not $\beta$-extremal, we find as $n$ is sufficiently large that $G$ is $(\frac{\beta^3}{10^{14}},8)$-$K_3$-hingeful. Since $\beta$ is sufficiently small, Theorem~\ref{thm:RegularFracStability} holds for $\beta$; let $C$ be as in Theorem~\ref{thm:RegularFracStability} for $\beta$. Since $G$ is not $\beta$-extremal and $\delta(G)\ge \frac{3}{4}n$, we have by Theorem~\ref{thm:RegularFracStability} that there exists a $C$-regular fractional $K_3$-decomposition of $G$. Note also that $G$ is $K_3$-divisible by assumption (that is integrally $K_3$-divisible and is trivially fractionally $K_3$-divisible). Given all of this and the fact that $n$ is sufficiently large, we find that $G$ satisfies the hypotheses of Theorem~\ref{thm:DecompThm} for $\alpha :=\frac{\beta^3}{10^{14}}$, $C$ and $\mc{F}=\{K_3\}$. Hence by Theorem~\ref{thm:DecompThm}, there exists a $K_3$-decomposition of $G$.
\end{proof}

\subsection{Idealized Decomposition Theorems}\label{ss:Idealized}

\begin{definition}
Let $\mc{F}=\{F_1,\ldots,F_k\}$ be a family of $t$-partitioned graphs for some integer $t\ge 1$. Let $\varepsilon > 0$ be real. Let $G$ be an $\mc{F}$-blowup on $n$ vertices with partition $(V_1,\ldots,V_t)$. We say $G$ is \emph{$\varepsilon$-balanced} if $|V_i| \ge (1-\varepsilon)\cdot \frac{n}{t}$ for each $i\in [t]$. We say $G$ is \emph{$\varepsilon$-nearly-complete} if for each $i\in [t]$, $j\in \{j'\in [t]: \exists \ell\in[k] \textrm{ such that } e_{ij'}(F_k)\ne 0 \}$ and $v\in V_i$, we have that $d_j(v)\ge (1-\varepsilon)\cdot |V_j|$.      
\end{definition}

\subsubsection{Idealized $\mc{F}_0$-decomposition Theorem}

\begin{thm}[Idealized $\mc{F}_0$-decomposition Theorem]\label{thm:F0Idealized}
For every real $\alpha > 0$, there exists $\varepsilon_0 > 0$ such that the following holds for all large enough $n$: Let $\mc{F}_0:= \{T_{112}, T_{122}\}$. Let $G$ be a $2$-partitioned graph on $n$ vertices with partition $(V_1,V_2)$. If $G$ is integrally $\mc{F}_0$-divisible and all of the following hold: 
\begin{enumerate}
    \item[(a)] for each $i\in [2]$, we have $|V_i|\ge (1-\varepsilon_0)n/2$,
    \item[(b)] for each $i\in [2]$ and $v\in V_i$, we have $(\frac{1}{4}-\varepsilon_0)n \le d_i(v) \le (\frac{1}{4}+\varepsilon_0)n$ and $(\frac{1}{2}-\varepsilon_0)n \le d_{3-i}(v)$, and
    \item[(c)] there exists $i\in [2]$ such that $V_i$ is not $\alpha$-close to either $T_2(|V_i|)$ or $\overline{T_2(|V_i|)}$,
\end{enumerate}
then $G$ has an $\mc{F}_0$-decomposition.
\end{thm}
\begin{proof}
\begin{claim}\label{cl:F0fracalmost}
There exists a $14\varepsilon_0$-almost $3$-regular fractional $\mc{F}_0$-packing of $G$.  
\end{claim}
\begin{proofclaim}
Let $w:= \frac{2}{(1+2\varepsilon_0)n}$. For each $\mc{F}_0$-copy $T$ of $G$, let $\phi(T) := w$.  Note that $\frac{3}{n}\ge \phi(T) \ge \frac{1}{3n}\ge 0$ for all $T$ and hence $\phi$ is $3$-regular. Moreover, for each $i\in [2]$ and internal edge $e=uv$ with $u,v\in V_i$, we find that $\phi(e) \le \max\{|V_1|,|V_2|\}\cdot w \le 1$ and similarly that 
$$\phi(e)\ge d_{3-i}(u)+d_{3-i}(v)-|V_{3-i}| \cdot w \ge \left(\frac{1}{2}-3\varepsilon_0\right )n\cdot \frac{2}{(1+2\varepsilon_0)n}\ge 1-8\varepsilon_0.$$
On the other hand, for each cross edge $e=uv$ with $u\in V_1,~v\in V_2$, we find that $\phi(e)\le (d_1(u)+d_2(v) )\cdot w\le 2\cdot (\frac{1}{4}+\varepsilon_0)n \cdot \frac{2}{(1+2\varepsilon_0)n} \le 1$ and similarly 
$$\phi(e) \ge (d_1(u) - |V_1\setminus N_1(v)| +d_2(v)-|V_2\setminus N_2(u)|)\cdot w \ge \left(\frac{1}{2} - 6\varepsilon_0\right)n \cdot \frac{2}{(1+2\varepsilon_0)n}\ge 1-14\varepsilon_0.$$ 
Hence $\phi$ is a $14\varepsilon_0$-almost $3$-regular fractional $\mc{F}_0$-packing of $G$ as desired. 
\end{proofclaim}

\begin{claim}\label{cl:F0fracdiv}
$G$ is fractionally $\mc{F}_0$-divisible.    
\end{claim}
\begin{proofclaim}
For each $i\in[2]$ and $v\in V_i$, we find by our assumptions on $G$  that $$d_i(v) \le \left(\frac{1}{4}+\varepsilon_0\right)n \le \left(\frac{1}{2}-\varepsilon_0\right)n \le d_{3-i}(v),$$
where the middle inequality follows since $\varepsilon_0\le \frac{1}{8}$. Hence $G$ satisfies Lemma~\ref{lem:F0div}(i). Since $G$ is integrally $\mc{F}_0$-divisible, we have by Lemma~\ref{lem:F0div} that $G$ satisfies Lemma~\ref{lem:F0div}(i') and (ii'). Thus $G$ satisfies Lemma~\ref{lem:F0div}(ii). Hence by Lemma~\ref{lem:F0div}, we find that $G$ is $\mc{F}_0$-divisible as desired.
\end{proofclaim}

\begin{claim}\label{cl:F0antiedgeful}
$G$ is $\frac{1}{5}$-$\mc{F}_0$-anti-edgeful.  
\end{claim}
\begin{proofclaim}
For $i\in [2]$, $u,v\in V_i$, the number of $T_{ii(3-i)}$-anti-edges on $uv$ is at least $|N_{3-i}(u)\cap N_{3-i}(v)| \ge d_{3-i}(u)+d_{3-i}(v) - |V_{3-i}| \ge (\frac{1}{2}-3\varepsilon_0)n \ge \frac{n}{5}$ since $\varepsilon_0 \le \frac{1}{10}$. Similarly for $u\in V_1,~v\in V_2$, the number of $T_{112}$-anti-edges on $uv$ is at least $|N_{1}(u)\cap N_{2}(v)| \ge d_{1}(u)+d_{2}(v) - |V_{1}| \ge (\frac{1}{4}-3\varepsilon_0)n \ge \frac{n}{5}$ since $\varepsilon_0 \le \frac{1}{60}$. Symmetrically, we also find that the number of T$_{122}$-anti-edges on $uv$ is at least $\frac{n}{5}$. Thus $G$ is $\frac{1}{5}$-$\mc{F}_0$-anti-edgeful as desired.
\end{proofclaim}

\begin{claim}\label{cl:F0hingeful}
$G$ is $\left(\frac{\alpha^3}{10^8},8\right)$-$\mc{F}_0$-hingeful.  
\end{claim}
\begin{proofclaim}
By (c), there exists $i^*\in [2]$ such that $V_{i^*}$ is not $\alpha$-close to either $T_2(|V_i|)$ or $\overline{T_2(|V_i|)}$. 

Let $T_1=u_1u_2v_1$ and $T_2=u_1u_2v_2$ be two $\mc{F}_0$-copies of $G$. For $i\in [2]$, let $a_i:= \mc{P}_G(u_i)$ and $b_i:= \mc{P}_G(v_i)$.  We assume without loss of generality that $a_1=1$. 
\vskip.1in
\noindent \textbf{Case 1: $a_2=1$}.
\vskip.1in 
Since $a_1=a_2=1$, we find that $b_1=b_2=2$. Additionally suppose that $i^*=1$. Let $H:= G[N(v_1)\cap N(v_2)\cap V_1]$. Note that $v(H)\ge (\frac{1}{2} - 3\varepsilon_0)n \ge \frac{n}{3}$ where the latter follows since $\varepsilon_0$ is small enough. Since $\varepsilon_0$ is small enough, it follows that $H$ is not $\frac{\alpha}{2}$-close to $T_2(v(H))$ or $\overline{T_2(v_H)}$. Moreover, $\delta(H)\ge (\frac{1}{2}-8\varepsilon_0)\cdot (\frac{1}{2} - 9\varepsilon_0)\cdot v(H)$. Since $\varepsilon_0$ is small enough, we have by Lemma~\ref{lem:ManyPaths} for $H$, $\alpha/2$, $9\varepsilon_0$ and $i=5$ that there exist at least $\frac{(\alpha/2)^3}{50}\cdot \left(\frac{v(H)}{3}\right)^4$ paths of length $5$ from $u_1$ to $u_2$ in $H$. Since each of these yields a distinct $\mc{F}_0$-hinge for $T_1$ and $T_2$, it follows that $G$ is $\left(\frac{\alpha^3}{50\cdot 16\cdot 9^4},8\right)$-$\mc{F}_0$ hingeful as desired.  

So we assume $i^*=2$. Let $U_1:= N(u_1)\cap N(v_1)\cap N(v_2)\cap V_1$. Note that $|U_1|\ge (\frac{1}{4} - 3\varepsilon_0)n \ge \frac{n}{5}$ since $\varepsilon_0$ is small enough. Fix $u_1'\in U_1$. Consider $H:= G[N(u_1')\cap N(u_2)\cap V_2]$. Note that $v(H)\ge (\frac{1}{2} - 3\varepsilon_0)n$. Since $\varepsilon_0$ is small enough, it follows that $H$ is not $\frac{\alpha}{2}$-close to $T_2(v(H))$ or $\overline{T_2(v(H))}$. Moreover, $\delta(H)\ge (\frac{1}{2}-8\varepsilon_0)\cdot (\frac{1}{2} - 9\varepsilon_0)\cdot v(H)$. Since $\varepsilon_0$ is small enough, we have by Lemma~\ref{lem:ManyPaths} for $H$, $\alpha/2$, $9\varepsilon_0$ and $i=4$ that there exist at least $\frac{(\alpha/2)^3}{50}\cdot \left(\frac{v(H)}{3}\right)^3$ paths of length $4$ from $u_1$ to $u_2$ in $H$. Since each of these yields a distinct $\mc{F}_0$-hinge for $T_1$ and $T_2$, it follows that there are at least $\frac{(\alpha/2)^3}{50}\cdot \left(\frac{v(H)}{3}\right)^3 \cdot |U_1|$ $\mc{F}_0$ hinges for $T_1$ and $T_2$ in $G$. Hence $G$ is $\left(\frac{\alpha^3}{5\cdot 50\cdot 8\cdot 9^3},8\right)$-$\mc{F}_0$-hingeful as desired.  
\vskip.1in
\noindent \textbf{Case 2: $a_2=2$}.
\vskip.1in 

We assume without loss of generality that $b_1=1$. Further suppose that $b_2 = 1$. Let $U_2:= N(u_2)\cap N(v_1)\cap N(v_2)\cap V_2$. Since $u_2,v_2\in V_1$, we find that $|U_1|\ge (\frac{1}{4}-5\varepsilon_0)n \ge \frac{n}{5}$ since $\varepsilon_0$ is small enough. Fix $u_2'\in U_2$. Let $U_1:=N(u_2')\cap N(u_1)\cap N(v_1)\cap N(v_2)\cap V_2$. Then each $u_2'\in U_2$ yields a $\mc{F}_0$-hinge with $u_1'$. Since $u_1,v_1,v_2\in V_1$, we find that $|U_1| \ge (\frac{1}{4} - 7\varepsilon_0)n \ge \frac{n}{5}$ since $\varepsilon_0$ is small enough. Thus there are at least $|U_1|\cdot |U_2| \ge \frac{n^2}{25}$ $\mc{F}_0$-hinges on $6$ vertices for $T_1$ and $T_2$ in $G$. Thus $G$ is $(\frac{1}{25},8)$-$\mc{F}_0$-hingeful as desired.

So we assume $b_2 = 2$. Let $U_1':= N(u_1)\cap N(v_1)\cap N(v_2)\cap V_2$. Since $u_1,v_1\in V_1$, we find that $|U_1'|\ge (\frac{1}{4}-5\varepsilon_0)n \ge \frac{n}{5}$ since $\varepsilon_0$ is small enough. Fix $u_1'\in U_1'$. Let $U_2':=N(u_1')\cap N(u_2)\cap N(v_1)\cap N(v_2)\cap V_1$. Then each $u_2'\in U_2'$ yields a $\mc{F}_0$-hinge with $u_1'$. Since $u_1',u_2,v_2\in V_2$, we find that $|U_2'| \ge (\frac{1}{4} - 7\varepsilon_0)n \ge \frac{n}{5}$ since $\varepsilon_0$ is small enough. Thus there are at least $|U_1'|\cdot |U_2'| \ge \frac{n^2}{25}$ $\mc{F}_0$-hinges on $6$ vertices for $T_1$ and $T_2$ in $G$. Thus $G$ is $(\frac{1}{25},8)$-$\mc{F}_0$-hingeful as desired.
\end{proofclaim}

By Claims~\ref{cl:F0fracalmost} and~\ref{cl:F0fracdiv}, we have that Theorem~\ref{thm:DecompThm}(i) holds for $G$. Since $G$ is integrally $\mc{F}_0$-divisible and using Claims~\ref{cl:F0antiedgeful} and~\ref{cl:F0hingeful}, we have that Theorem~\ref{thm:DecompThm}(ii) holds for $G$. Thus we find by Theorem~\ref{thm:DecompThm} since $\varepsilon_0$ is small enough that $G$ admits an $\mc{F}_0$-decomposition as desired.
\end{proof}

\subsubsection{Idealized $\mc{F}_i$-decomposition Theorem for each $i\in [3]$}

\begin{thm}[Idealized $\mc{F}_i$-decomposition Theorem for each $i\in \{1,2,3\}$]\label{thm:Idealized}
Let 
\begin{itemize}
    \item $\mc{F}_1:=\{T_{113}, T_{114}, T_{223}, T_{224},T_{133},T_{233},T_{144},T_{244}\}$,
    \item $\mc{F}_2:=\{T_{113}, T_{114}, T_{223}, T_{224},T_{134},T_{234}\}$
    \item $\mc{F}_3:= \{T_{123}, T_{124}, T_{134}, T_{234}\}$.
\end{itemize}
For each $i\in [3]$, there exists $\varepsilon_i > 0$ such that the following holds for all large enough $n$:  If $G$ is an $\varepsilon_i$-balanced $\varepsilon_i$-nearly-complete $\mc{F}_i$-blowup graph on $n$ vertices that is integrally $\mc{F}_i$-divisible, then $G$ has an $\mc{F}_i$-decomposition. 
\end{thm}
\begin{proof}
The key property to note is that for each $i\in [3]$, we have that for each $j,k \in [4]$ such that there exists at least one $\ell\in [4]$ such that $T_{jk\ell}\in \mc{F}_i$, there exists exactly two such indices $\ell$.

\begin{claim}\label{cl:Fifracalmost}
There exists an $8\varepsilon_i$-almost $3$-regular fractional $\mc{F}_i$-packing of $G$.  
\end{claim}
\begin{proofclaim}
Let $w:= \frac{2}{(1+2\varepsilon_i)n}$. For each $\mc{F}_i$-copy $T$ of $G$, let $\phi(T) := w$.  Note that $\frac{3}{n}\ge \phi(T) \ge \frac{1}{3n}\ge 0$ for all $T$ and hence $\phi$ is $3$-regular. 

We claim that for each $e=uv\in E(G)$ we have that $1\ge \phi(e)\ge 1-8\varepsilon_1$ and hence $\phi$ is a $8\varepsilon_i$-almost $3$-regular fractional $\mc{F}_i$-packing of $G$ as desired. Let $j$ be such that $u\in V_j$ and $k$ be such that $v\in V_k$. Let $\ell \ne m\in [4]$ such that $T_{jk\ell}, T_{jkm}\in \mc{F}_i$ (recall there are exactly two such numbers).  Now we calculate that
$$\phi(e) \le (|V_\ell|+|V_m|)\cdot w \le \frac{1}{2}(1+\varepsilon_i)n\cdot \frac{2}{(1+2\varepsilon_i)n} \le 1$$
and similarly that 
$$\phi(e)\ge \bigg((d_{\ell}(u)+d_{\ell}(v)-|V_{\ell}|)+(d_m(u)+d_m(v)-|V_{m}|)\bigg) w \ge \left(\frac{1}{2}-3\varepsilon_i\right )n\cdot \frac{2}{(1+2\varepsilon_i)n}\ge 1-8\varepsilon_i$$
as desired.
\end{proofclaim}

\begin{claim}\label{cl:Fifracdiv}
$G$ is fractionally $\mc{F}_i$-divisible.    
\end{claim}
\begin{proofclaim}
First we note that for all $i\in [3]$ for $j\in [4]$ and $j'\in [4]$ such that there exists $k$ with $T_{jj'k}\in \mc{F}_i$, we have for all $v\in V_j$ that 
$$(1-\varepsilon_i)|V_{j'}|\le d_{j'}(v)\le |V_{j'}|,$$
and hence $d_{j'}(v) \in [(1-2\varepsilon_i)\frac{n}{4},(1+3\varepsilon_i)\frac{n}{4}]$.

Thus for $j \in [4]$ such that there exists $k$ with $T_{jjk}\in \mc{F}_i$, we have that 
$$\frac{(1-\varepsilon_i)}{2}\cdot |V_j|^2 \le e_{jj}(G)\le |V_j|^2$$
and hence $e_{jj}(G) \in [(1-3\varepsilon_i)\cdot \frac{n^2}{32}, (1+6\varepsilon_i)\frac{n^2}{32}]$. Similarly $j\ne j'\in [4]$ such that there exists $k$ with $T_{jj'k}$, we have that
$$\frac{(1-\varepsilon_i)}{2}\cdot |V_j|\cdot |V_k| \le e_{jk}(G)\le |V_j|\cdot |V_k|$$
and hence $e_{jk}(G) \in [(1-3\varepsilon_i)\cdot \frac{n^2}{16}, (1+6\varepsilon_i)\frac{n^2}{16}]$.

First suppose that $i=1$. For all $j\in [4]$ and $v\in V_j$, we find by our assumptions on $G$ that $$d_j(v) \le \left(\frac{1}{4}+3\varepsilon_1\right)n \le \left(\frac{1}{2}-4\varepsilon_1\right)n \le d_{\ell}(v)+d_k(v)$$ 
where $\{\ell,k\}=\{3,4\}$ if $j\in\{1,2\}$ and $\{1,2\}$ if $j\in \{3,4\}$ - where the middle inequality follows since $\varepsilon_1\le \frac{1}{56}$. Hence $G$ satisfies Lemma~\ref{lem:F1div}(ii). Since $G$ is integrally $\mc{F}_1$-divisible, we have by Lemma~\ref{lem:F1div} that $G$ satisfies Lemma~\ref{lem:F1div}(i') and thus $G$ satisfies Lemma~\ref{lem:F1div}(i). Furthermore, it follows from our calculations at the beginning of the proof and since $\varepsilon_i$ is small enough that Lemma~\ref{lem:F1div}(iii) holds. Hence by Lemma~\ref{lem:F1div}, we find that $G$ is fractionally $\mc{F}_1$-divisible as desired.

We omit the proofs of the cases $i=2$ and $i=3$ since they are similar (with (i) following from (i') by integral-divisibility and the remaining conditions following from the calculations at the start of the proof). 
\end{proofclaim}

\begin{claim}\label{cl:Fiantiedgeful}
$G$ is $\frac{1}{5}$-$\mc{F}_i$-anti-edgeful.  
\end{claim}
\begin{proofclaim}
Let $j,k,\ell \in [4]$ such that $T_{jk\ell}\in \mc{F}_i$. Let $u\in V_j$ and $v\in V_k$. Then the number of $T_{jk\ell}$-anti-edges on $uv$ is at least $|N_{\ell}(u)\cap N_{\ell}(v)| \ge d_{\ell}(u)+d_{\ell}(v) - |V_{\ell}| \ge (1-2\varepsilon_i)\cdot |V_\ell| \ge (1-3\varepsilon_i)\cdot \frac{n}{4} \ge \frac{n}{5}$ since $\varepsilon_i \le \frac{1}{15}$. Thus $G$ is $\frac{1}{5}$-$\mc{F}_i$-anti-edgeful as desired.
\end{proofclaim}

\begin{claim}\label{cl:Fihingeful}
$G$ is $(\frac{1}{20},6)$-$\mc{F}_i$-hingeful.  
\end{claim}
\begin{proofclaim}
Let $T_1=u_1u_2v_1$ and $T_2=u_1u_2v_2$ be two $\mc{F}_i$-copies of $G$. Let $j$ be such that $u_1\in V_j$ and $k$ be such that $u_2\in V_k$. Choose $u_1'\in N(u_2)\cap N(v_1)\cap N(v_2) \cap V_j$. Note there are at least $(1-3\varepsilon_i)\cdot |V_j| \ge (1-4\varepsilon_i)\cdot \frac{n}{4}$ choices of $u_1'$. Now choose $u_2'\in N(u_1)\cap N(u_1')\cap N(v_1)\cap N(v_2)\cap V_k$. Note there are at least $(1-4\varepsilon_i)\cdot |V_k| \ge (1-5\varepsilon_i)\cdot \frac{n}{4}$ choices of $u_1'$. But $u_1',u_2'$ together with the edges $u_1'u_2',u_1'u_1,u_1'v_1,u_1'v_2,u_2'v_1,u_2'v_2,u_2'u_1$ form an $\mc{F}_i$-hinge for $T_1\cup T_2$ in $G$. Hence there are at least $(1-4\varepsilon_i)\cdot (1-5\varepsilon_i)\cdot \frac{n^2}{16}$ $\mc{F}_i$-hinges for $T_1\cup T_2$ in $G$. Thus $G$ is $(\frac{1}{20},6)$-$\mc{F}_i$-hingeful as desired since $\varepsilon_i\le \frac{1}{45}$. 
\end{proofclaim}

By Claims~\ref{cl:Fifracalmost} and~\ref{cl:Fifracdiv}, we have that Theorem~\ref{thm:DecompThm}(i) holds for $G$. Since $G$ is integrally $\mc{F}_i$-divisible and using Claims~\ref{cl:Fiantiedgeful} and~\ref{cl:Fihingeful}, we have that Theorem~\ref{thm:DecompThm}(ii) holds for $G$. Thus we find by Theorem~\ref{thm:DecompThm} since $\varepsilon_i$ is small enough that $G$ admits an $\mc{F}_i$-decomposition as desired.
\end{proof}

\section{Cleaning the Extremal Cases}\label{s:Cleaning}

It will be useful to prove a `cleaning' lemma that we can apply inductively to three out of the four cases (i.e.~(2)-(4)) to clean the graph towards its nearly-balanced nearly-complete version so as to apply the idealized decomposition theorems from the previous section. 

First we require some additional terminology about $2$-partitioned graphs as follows.

\begin{definition}
Let $G$ be a $2$-partitioned graph on $n$ vertices with partition $(V_1,V_2)$. An edge $uv$ of $G$ is \emph{cross} if $\mc{P}_G(u)\ne \mc{P}_G(v)$ and \emph{interior} otherwise. For each $v\in V(G)$, we define the \emph{interior degree} of $v$, denoted $\textrm{IntDeg}_G(v)$ as the number of interior edges incident with $v$, and we define the \emph{cross degree} of $v$, denoted $\textrm{CrossDeg}_G(v)$ as the number of cross edges incident with $v$. 

A triangle $T$ of $G$ is \emph{cross} if $V(T)\cap V_i\ne \empty$ for each $i\in [2]$ and \emph{interior} otherwise. We define the \emph{internal surplus} of $G$, denoted $\textrm{IntSurplus}(G)$, as $2(e_{11}(G)+e_{22}(G))-e_{12}(G)$. We say $G$ is \emph{cross-divisible} if $\textrm{IntSurplus}(G)=0$.     
\end{definition}

We note for the sake of the reader that if a $2$-partitioned graph $G$ has triangle decomposition $\mc{T}$, then $\mc{T}$ must include exactly $\frac{1}{6}\cdot \textrm{IntSurplus}(G)$ interior triangles of $G$. In particular, if the $\textrm{IntSurplus}(G)$ is negative or not divisible by $6$, then $G$ will not admit a triangle decomposition. We note that for a graph $G$ on $n$ vertices with $\delta(G)\ge \frac{3}{4}n$, we have that $\textrm{IntSurplus}(G)\ge 0$; for triangle-divisible $G$, we find that $\textrm{IntSurplus}(G)$ is divisible by $6$ as follows.

\begin{proposition}\label{prop:CrossDiffParity}
If $G$ is a $2$-partitioned triangle-divisible graph with partition $(V_1,V_2)$, then 
$$\textrm{IntSurplus}(G)\equiv 0 \pmod 6.$$
\end{proposition}
\begin{proof}
Since $e_{11}(G)+e_{22}(G) = e(G)-e_{12}(G)$, it follows that $\textrm{IntSurplus}(G)= 2\cdot e(G) - 3\cdot e_{12}(G)$. Since $G$ is triangle divisible, we have that $3~|~e(G)$ and hence $6~|~2e(G)$. Since $G$ is triangle divisible. we have that $G$ is Eulerian. But then $\sum_{v\in V_1} d(v) = 2\cdot e_{11}(G)+e_{12}(G)$ is even and hence $e_{12}(G)$ is even. That is $2~|~e_{12}(G)$ and hence $6~|~e_{12}(G)$. Thus $6~|~\textrm{IntSurplus}(G)$. 
\end{proof}

\subsection{Proof Overview for Cleaning}

Now we provide a high-level proof overview for our cleaning steps as the proof is fairly complicated.

In Cases (2)-(4) (the families $\mc{F}_0, \mc{F}_1,\mc{F}_2$ respectively), the integral divisibility conditions by Lemmas~\ref{lem:F0div},~\ref{lem:F1div}, and~\ref{lem:F2div} are essentially that $G$ is Eulerian and $G$ is cross-divisible (where for $\mc{F}_1$ and $\mc{F}_2$ we mean after repartitioning the graph with bipartition $(V_1\cup V_2,V_3\cup V_4)$) where additionally for $\mc{F}_1$ and $\mc{F}_2$ they each have an additional `parity' condition. The three cases mostly admit a common proof then as follows.

First, we prove a `cleaning lemma', Lemma~\ref{lem:OddballCleaning}, that shows how to decompose the triangles around \emph{one oddball vertex} $v$ (meaning it does not conform nicely to the partition of the case). Said lemma says that if the $2$-partitioned graph is extremal then this is possible using only cross triangles provided $\textrm{IntDeg}_G(v)$ is slightly smaller than $\frac{3}{8}n$; on the other hand, when having larger interior degree we can decompose using not too many interior triangles (at most $0.01n$ interior triangles). The latter actually only works provided the interior degree is not too much larger than the cross degree (when all vertices have this property we will call the graph cross-balanced). The latter property will follow by taking the bipartition to be a max cut (instead of a balanced max cut as might have been more natural).

Next we prove a second cleaning lemma, Lemma~\ref{lem:ManyOddballCleaning}, that shows how to decompose \emph{many oddballs} while additionally \emph{leaving a cross-divisible graph} so as to apply an idealized decomposition theorem. To ensure we do not create a negative internal surplus when decomposing the oddballs, we first show that the internal surplus is linear in the number of `big' internal degree vertices where big is as in the first cleaning lemma. Then we apply Lemma~\ref{lem:OddballCleaning} inductively to decomposes all oddball while leaving non-negative internal surplus. However, to ensure cross-divisibility, we then need to remove interior triangles. 

Thus to ensure the final internal surplus is \emph{exactly} zero, we in all three cases we find a large triangle packing among the non-oddballs. Such exists for Case (2) by Tur\'an stability as at least one side is not close to complete bipartite while for Cases (3) and (4) this follows from having a large near-clique subgraph. Then we sparsify this triangle packing to a still triangle packing but of small maximum degree and \emph{reserve} this sparsified triangle packing (before decomposing oddballs). Hence to complete the many oddball cleaning lemma proof, we simply remove interior triangles from this `reserve' until the surplus is exactly zero (which is possible via Proposition~\ref{prop:CrossDiffParity}).

The proof then of Case (2), i.e.~Theorem~\ref{thm:neitherside}, is straightforward from the arguments above, simply checking the various conditions. Cases (3) and (4) require the use of \emph{parity fixers}, gadgets (namely $4$-wheels) which we will use to ensure the extra parity condition is satisfied in each case. For each case then, we first prove a lemma showing these exists a parity fixer in our graph. Then for the proofs of Theorems~\ref{thm:4cliques} and~\ref{thm:2cliqueandbip}, we first set aside this parity fixer, then find a large triangle packing and sparsify to use as a reserve, and then apply the many oddball cleaning lemma. Note here the oddballs will include not only those that are abnormal with respect to the cross partition but also those that are abnormal with respect to the original $4$-partition. After this, we also find a small maximum degree triangle packing to decompose all \emph{oddball edges} among normal vertices (edges not conforming to the partition). Finally, we flip the decomposition of the parity fixer to ensure the parity condition and apply the idealized decomposition theorem.

Lastly, we discuss the final case, Theorem~\ref{thm:K4blowup} and cleaning toward integral $\mc{F}_3$-divisibility. Here our graph is nearly balanced and has few oddballs. Yet here the integral divisibility conditions requires \emph{two equality constraints}, that is for Lemma~\ref{lem:F3div}(i') to hold which says the number of edges in each of the three `matchings' $M_1=\{12,34\}$, $M_2=\{13,24\}$, $M_3=\{14,23\}$ is equal. Thus this case proves a bit trickier and does not unify with the other three. More problematically, we also have \emph{no reservoir of interior triangles} to use to achieve our equality condition. Thus we seek to decompose the interior edges each into their own triangle with $2$ cross edges to one of the other parts and in such a way that leaves the equality condition satisfied. Note we can categorize each such triangle by the matching the cross edges reside in (we call these $i$-matched if they reside in the $i$th matching). 

Thus we first prove in Proposition~\ref{prop:MatchingNumbers} that this is possible in the `abstract', meaning there are numbers for each $i\in [3]$ such that using that number of $i$-matched triangles would yield the equality condition. To decompose, we then color the interior edges with $3$ colors to denote which type of matched triangle they should use and set up an embedding problem to show such matched triangles can be found in a low degree manner. However, care must be taken to color the interior edges around oddball vertices (namely with as balanced a coloring as possible from their neighborhood) since otherwise it may not be possible to decompose the oddball neighborhoods afterwards. 

\subsection{Cleaning Lemmas}\label{ss:Cleaning}

In this subsection, we prove our cleaning lemmas. First we extend the notion of extremal to $2$-partitioned graphs as follows.

\begin{definition}
Let $\sigma > 0$. We say a $2$-partitioned graph $G$ on $n$ vertices with partition $(V_1,V_2)$ is \emph{$\sigma$-extremal} if both of the following hold:    
\begin{itemize}
    \item[(i)] at least $(1-\sigma)n$ vertices of $G$ have degree at most $(\frac{3}{4}+\sigma)n$, and
    \item[(ii)] $e_{G}(V_1,V_2)\ge (1-\sigma)\cdot \frac{n^2}{4}$.
\end{itemize}
We say $G$ is \emph{$\sigma$-cross-balanced} if for each $v\in V(G)$, we have that $\textrm{IntDeg}_G(v) \le \textrm{CrossDeg}_G(v) + \sigma\cdot n$.  
\end{definition}

Next we need a definition of when a triangle packing decomposes and/or saturates a vertex (or more generally a set of vertices) as follows.

\begin{definition}
Let $G$ be a graph. For $v\in V(G)$, we say a triangle packing $\mc{T}$ \emph{saturates} $v$ if every edge incident with $v$ is contained in some triangle of $\mc{T}$; similarly we say $\mc{T}$ is a \emph{triangle decomposition of $v$} in $G$ if $\mc{T}$ saturates $v$ and every triangle of $\mc{T}$ contains $v$. More generally for $S\subseteq V(G)$, a we say $\mc{T}$ \emph{saturates} $S$ if $\mc{T}$ saturates every of $S$ and $\mc{T}$ is a \emph{triangle decomposition of $S$} in $G$ if $\mc{T}$ saturates $S$ and every triangle of $\mc{T}$ contains at least one vertex of $S$.  
\end{definition}

\vskip.1in

For decomposing an oddball neighborhood, we need P\'osa's Theorem which we recall as follows.

\begin{thm}[P\'osa 1962~\cite{P62}]\label{thm:Posa}
Let $G$ be a graph on $n$ vertices and let $v_1,v_2,\ldots,v_n$ be an ordering of $V(G)$ such that $d_G(v_i)\le d_G(v_j)$ for all $i<j\in [n]$. If for each $1\le k \le n/2$ we have $d_G(v_k) > k$, then $G$ has a Hamiltonian cycle. 
\end{thm}

Next we prove a useful lemma about some of the properties of an extremal graph as follows.

\begin{lem}\label{lem:neithersideRobust}
Let $\gamma > 0$. Let $G$ be a $2$-partitioned graph on $n$ vertices with $\delta(G)\ge \frac{3}{4}n$ and partition $(V_1,V_2)$. If $G$ is $\gamma$-extremal and we let $\mc{A} := \bigg\{v\in V(G): \textrm{IntDeg}_G(v) \ge \left(\frac{1}{4}+2\sqrt{\gamma}\right)n\bigg\}$, then 
$$|\mc{A}| \le 2\sqrt{\gamma}\cdot n$$ 
and $G$ satisfies both of the following:
\begin{enumerate}
    \item[(a)] for each $i\in [2]$, we have $|V_i|\ge (1-\sqrt{\gamma})n/2$,
    \item[(b)] for each $i\in [2]$ and $v\in V_i\setminus \mc{A}$, we have $(\frac{1}{4}-\sqrt{\gamma})n \le d_i(v) \le (\frac{1}{4}+2\sqrt{\gamma})n$ and $(\frac{1}{2}-\sqrt{\gamma})n \le d_{3-i}(v)$,
\end{enumerate}
\end{lem}
\begin{proof}

\begin{claim}\label{cl:aRobust}
(a) holds.
\end{claim}
\begin{proofclaim}
Suppose not. But then $e_G(V_1,V_2)\le |V_1|\cdot |V_2| < (\frac{1}{2}-\sqrt{\gamma}) n \cdot (\frac{1}{2}+\sqrt{\gamma})n = (\frac{1}{4}-\gamma)n^2$, contradicting that $G$ is $\gamma$-extremal.
\end{proofclaim}

\begin{claim}\label{cl:SizeOfAbnormal}
$|\mc{A}|\le 2\sqrt{\gamma}\cdot n$.    
\end{claim}
\begin{proofclaim}
Suppose not. Let $B:= \{v\in V(G): d_G(v) > (\frac{3}{4}+\gamma)n)\}$. Since $G$ is $\gamma$-extremal, we have by definition that $|B| \le \gamma\cdot n$. Thus $|\mc{A}\setminus B| > \frac{3}{2} \sqrt{\gamma}n$. For $v\in \mc{A}\setminus B$, we have that 
$$\textrm{CrossDeg}_G(v) = d_G(v) - \textrm{IntDeg}_G(v) \le \left(\frac{3}{4}+\gamma\right)n - \left(\frac{1}{4}+2\sqrt{\gamma}\right)n < \left(\frac{1}{2}-\frac{3}{2}\sqrt{\gamma}\right)n,$$
where the last inequality follows since $\beta$ is small enough. Using Claim~\ref{cl:aRobust}, we find that for each $i\in [2]$ and $v\in V_i\cap  \mc{A}\setminus B$ is non-adjacent to at least $\frac{1}{2}\sqrt{\gamma}n$ vertices in $V_{3-i}$. We assume without loss of generality that $|V_1\cap \mc{A}\setminus B| > \frac{1}{2} \sqrt{\gamma}\cdot n$. But then
$e_G(V_1,V_2) < |V_1|\cdot |V_2| - \frac{1}{2} \cdot \sqrt{\gamma}\cdot n \cdot \frac{1}{2}\sqrt{\gamma} \cdot n \le \frac{n^2}{4} - \frac{\gamma}{4} \cdot n^2$, contradicting that $G$ is $\gamma$-extremal.
\end{proofclaim}

\begin{claim}\label{cl:bRobust}
(b) holds.
\end{claim}
\begin{proofclaim}
Let $i\in [2]$ and $v\in V_i$. Since $|V_{3-i}|\le (\frac{1}{2}+\sqrt{\beta})n$ by Claim~\ref{cl:aRobust} and $\delta(G)\ge \frac{3}{4}n$, we have that $d_i(v) \ge (\frac{1}{4}-\sqrt{\beta})n$. Since $v\not\in A$, we have by definition that $d_i(v) \le (\frac{1}{4}+2\sqrt{\beta})n$. Since $\delta(G) \ge \frac{3}{4}n$, it then follows that $d_{3-i}(v) \ge (\frac{1}{2}-2\sqrt{\beta})n$ as desired.    
\end{proofclaim}

\noindent The lemma now follows from Claims~\ref{cl:aRobust},~\ref{cl:SizeOfAbnormal},~\ref{cl:bRobust}.
\end{proof}

We are finally prepared to state and prove our cleaning lemma for decomposing one oddball as follows.

\begin{lem}[One Oddball Cleaning Lemma]\label{lem:OddballCleaning}
There exists $\varepsilon \in (0,\frac{1}{2000}]$ such that the following holds for all large enough $n$: Let $G$ be a $2$-partitioned graph on $n$ vertices where $\delta(G)\ge (\frac{3}{4}-\varepsilon)n$. If $G$ is $\varepsilon$-extremal and $\varepsilon$-cross-balanced, then the following holds for $v\in V(G)$:

\begin{itemize}
    \item[(i)] if $\textrm{IntDeg}_G(v) \le (\frac{3}{8}-7\varepsilon)n$, then there exists a triangle decomposition of $v$ in $G$ using only cross triangles of $G$, and
    \item[(ii)] if $\textrm{IntDeg}_G(v) > (\frac{3}{8}-7\varepsilon)n$, then there exists a triangle decomposition of $v$ in $G$ using at most $10\varepsilon n$ interior triangles of $G$. 
\end{itemize}
\end{lem}
\begin{proof}
Let $(V_1,V_2)$ be the partition of $G$. Let 
$$\mc{A} := \bigg\{v\in V(G): \textrm{IntDeg}_G(v) \ge \left(\frac{1}{4}+2\sqrt{\varepsilon}\right)n\bigg\}.$$ 
Since $G$ is $\varepsilon$-extremal, we have by Lemma~\ref{lem:neithersideRobust} that $|\mc{A}|\le 2\sqrt{\varepsilon}\cdot n$ and that both of Lemma~\ref{lem:neithersideRobust}(a) and (b) for hold $G$.

We assume without loss of generality that $v\in V_1$. Let $H:= G[N_G(v)] - E(N_G(v)\cap V_1)$. Let $m:=V(H)$. Since $d_G(v)\ge(\frac{3}{4}-\varepsilon)n$, it follows that $m=d_G(v)\ge (\frac{3}{4}-\varepsilon)n$. Since $G$ is $\varepsilon$-cross-balanced, we find for each $w\in V(G)$ that $2\textrm{CrossDeg}_G(w) \ge \textrm{CrossDeg}_G(w) + \textrm{IntDeg}_G(w) - \varepsilon n\ge d_G(w)  - \varepsilon n$ and hence 
$$\textrm{CrossDeg}_G(w) \ge \frac{1}{2} (d_G(w) - \varepsilon) \ge \left(\frac{3}{8} - \varepsilon\right)n,$$ 
since $\delta(G)\ge (\frac{3}{4}-\varepsilon)n$.
By Lemma~\ref{lem:neithersideRobust}(a), we have $|V_i|\ge (1-\sqrt{\varepsilon})n/2$ for each $i\in [2]$ and hence for each $i\in [2]$ and $w\in V_i$, we find that
$$\textrm{IntDeg}_G(w) \ge d_G(w) - |V_{3-i}| \ge \left(\frac{3}{4}-\varepsilon\right)n - \left(\frac{1}{2} - \sqrt{\varepsilon}\right)n \ge \left(\frac{1}{4}-2\sqrt{\varepsilon}\right)n.$$

First suppose that $\textrm{IntDeg}_G(v) \le (\frac{3}{8}-7\varepsilon)n$. Hence $$d_2(v) = d_G(v) - d_1(v) \ge \left(\frac{3}{4}-\varepsilon\right)n - \left(\frac{3}{8}-7\varepsilon\right)n \ge \left(\frac{3}{8}+6\varepsilon\right)n.$$
For each $i\in [2]$ and $w\in V(H)\cap \mc{A}\cap V_1$, we find by Lemma~\ref{lem:neithersideRobust}(b) that  
\begin{align*}
d_H(w) &\ge |N(w)\cap V_{2}\cap V(H)| \ge |N(w)\cap V_{2}| +|N(v)\cap V_{2}| - |V_{2}|\\
&\ge \left(\frac{3}{8}-\frac{\varepsilon}{2}\right)n + \left(\frac{1}{4}-2\varepsilon\right)n - \left(\frac{1}{2}+\varepsilon\right)n \ge \left(\frac{1}{8}-4\varepsilon\right)n \ge  \frac{n}{10},    
\end{align*}
where the last inequality follows since $\varepsilon \le \frac{1}{160}$.  On the other hand, for each $w\in V(H)\cap V_1 \setminus \mc{A}$, we find that
\begin{align*}
d_H(w) &\ge |N(w)\cap V_2\cap V(H)| \ge |N(w)\cap V_2| +|N(v)\cap V_2| - |V_2|\\
&\ge \left(\frac{1}{2}-2\varepsilon\right)n + \left(\frac{3}{8}-2\varepsilon\right)n - \left(\frac{1}{2}+\varepsilon\right)n \ge \left(\frac{3}{8}-5\varepsilon\right)n \ge  d_1(v) = |V(H)\cap V_1|,    
\end{align*}
Finally for each $w\in V(H)\cap V_2$, we find that
$$d_H(v) \ge |N(w)\cap V(H)| = |N(w)\cap N(v)| \ge \left(\frac{1}{2}-2\varepsilon\right)n \ge \frac{17}{32}\cdot v(H),$$
where we used that $v(H) \le \frac{7}{8}n$ and $\varepsilon \le \frac{1}{100}$. 

Note that $|\mc{A}|\le 2\sqrt{\varepsilon}\cdot n \le \frac{n}{10}$. Also note that $d_2(v) \ge (\frac{3}{8}+6\varepsilon)n \ge (\frac{3}{8}-7\varepsilon)n \ge d_1(v)$. Thus by P\'osa's theorem there exists a perfect matching $M$ of $H$ which in turn combined with the edges incident with $v$ yields a triangle decomposition of $v$ in $G$ using only cross triangles of $G$. This proves (a).

So we assume that $\textrm{IntDeg}_G(v) \ge (\frac{3}{8}-7\varepsilon)n$. But then $$\delta(G[N(v)\cap V_1]) \ge \left(\frac{3}{8}-\frac{\varepsilon}{2}\right)n + \left(\frac{1}{4}-2\varepsilon\right)n - \left(\frac{1}{2}+\varepsilon\right)n \ge \left(\frac{1}{8}-4\varepsilon\right)n \ge  \frac{n}{10},$$
where for the last inequality we used that $\varepsilon \le \frac{1}{160}$. Since $\varepsilon \le \frac{1}{200}$, it then follows that there exists a matching $M$ of size at least $10\varepsilon \cdot n$ in $G[N(v)\cap V_1]$ (by greedily removing edges). Let $H:= (G[N_G(v)] - E(N_G(v)\cap V_1))\setminus V(M)$. 

Thus $v(H) = d(v)-20\varepsilon \cdot n$. Note that $d_2(v) \ge d_1(v)-\varepsilon\cdot n$. Hence $d_2(v) \ge \frac{1}{2}(d(v) - \varepsilon n) \ge \frac{1}{2}(v(H)+19\varepsilon n) \ge \frac{v(H)}{2} + 9\varepsilon n$. Thus $|V(H)\cap V_1|\le \frac{v(H)}{2} - 9\varepsilon n$.
For each $i\in [2]$ and $w\in V(H)\cap \mc{A}\cap V_1$, we find by Lemma~\ref{lem:neithersideRobust}(b) that  
\begin{align*}
d_H(w) &\ge |N(w)\cap V_{2}\cap V(H)|-20\varepsilon n \ge |N(w)\cap V_{2}| +|N(v)\cap V_{2}| - |V_{2}|\\
&\ge \left(\frac{3}{8}-\frac{\varepsilon}{2}\right)n + \left(\frac{1}{4}-2\varepsilon\right)n - \left(\frac{1}{2}+\varepsilon\right)n -20\varepsilon n \ge \left(\frac{1}{8}-24\varepsilon\right)n \ge  \frac{n}{10},    
\end{align*}
where the last inequality follows since $\varepsilon \le \frac{1}{960}$.  On the other hand, for each $w\in V(H)\cap V_1 \setminus \mc{A}$, we find that
\begin{align*}
d_H(w) &\ge |N(w)\cap V_2\cap V(H)| \ge |N(w)\cap V_2| +|N(v)\cap V_2| - |V_2|\\
&\ge \left(\frac{1}{2}-2\varepsilon\right)n + \left(\frac{3}{8}-2\varepsilon\right)n - \left(\frac{1}{2}+\varepsilon\right)n \ge \left(\frac{3}{8}-5\varepsilon\right)n \ge |V(H)\cap V_1|,    
\end{align*}
Finally for each $w\in V(H)\cap V_2$, we find that
$$d_H(v) \ge |N(w)\cap V(H)| = |N(w)\cap N(v)| -20\varepsilon n\ge \left(\frac{1}{2}-22\varepsilon\right)n \ge \frac{17}{32}\cdot v(H),$$
where we used that $v(H) \le \frac{7}{8}n$ and $\varepsilon \le \frac{1}{1100}$. 

Note that $|\mc{A}|\le 2\sqrt{\varepsilon} n \le \frac{n}{10}$ since $\varepsilon \le \frac{1}{400}$. Also note that $|V(H)\cap V_2| \ge  (\frac{3}{8}-\varepsilon)n \ge (\frac{3}{8}-19\varepsilon)n \ge |V(H)\cap V_1|$. Thus by P\'osa's theorem there exists a perfect matching $M'$ of $H$; but then $M\cup M'$ combined with the edges incident with $v$ yields a triangle decomposition of $v$ in $G$ using at most $10\varepsilon\cdot n$ interior triangles of $G$. This proves (b). 
\end{proof}

As a point of notation as we delve deeper into the proof, we consider a triangle packing $\mc{T}$ of a graph $G$ as both a set of triangles (so that $|\mc{T}|$ denotes the number of triangles in the packing) and as a subgraph of $G$ (so that $E(\mc{T})$ denotes its set of edges $\bigcup_{T\in \mc{T}} E(T)$ and $\Delta(\mc{T})$ denotes its maximum degree that is $\Delta(G[E(\mc{T})])$).

Here is a useful proposition which shows how to find a still large but small maximum degree triangle packing inside a larger triangle packing.

\begin{proposition}\label{prop:SmallDegSubpacking}
Let $G$ be a graph on $n$ vertices and let $d\in (0,1]$ such that $dn\ge 2$. If $\mc{T}$ is a triangle packing of $G$, then there exists $\mc{T}'\subseteq \mc{T}$ such that $|\mc{T}'|\ge \frac{d}{6}\cdot |\mc{T}|$ and $\Delta(\mc{T}')\le dn+2$.     
\end{proposition}
\begin{proof}
Let $\mc{T}'\subseteq \mc{T}$ be chosen such that $\Delta(\mc{T}')\le dn+2$ and subject to that $|\mc{T}'|$ is maximized. Let $S:= \{v\in V(G): d_{\mc{T}'}(v)\ge dn \}$. Since $\mc{T}'$ is maximized, it follows that every triangle in $\mc{T}\setminus \mc{T}'$ is contains a vertex of $S$. Let $\mc{T}''$ be the set of triangles in $\mc{T}'$ incident that contain a vertex of $S$. Note that 
$$|\mc{T}''| = \frac{1}{3}\cdot E(\mc{T}'') \ge \frac{1}{6} \sum_{v\in S} d_{\mc{T}'}(v) \ge |S|\cdot \frac{dn}{6}.$$
But $|\mc{T}\setminus \mc{T}'| + |\mc{T}''| \le |S|\cdot \frac{n}{2}$. Thus
$$\frac{|\mc{T}'|}{|\mc{T}|} \ge \frac{|\mc{T}''|}{|\mc{T}\setminus \mc{T}'| + |\mc{T}''|} \ge \frac{|S|\cdot (dn/6)}{|S|\cdot n} = \frac{d}{6}$$
as desired.
\end{proof}

We can now prove our second cleaning lemma which decomposes all of the oddballs while ensuring the leftover graph is cross-divisible as follows.

\begin{lem}[Many Oddballs Cleaning Lemma]\label{lem:ManyOddballCleaning}
For each sufficiently small $\gamma$ the following holds for all large enough $n$: Let $G$ be a triangle-divisible $2$-partitioned graph on $n$ vertices where $\delta(G)\ge (\frac{3}{4}-\gamma)n$. Suppose that $G$ is $\gamma$-extremal, $\gamma$-cross-balanced and contains a triangle-packing $\mc{T}_0$ consisting only of interior triangles such that $|\mc{T}_0| \ge 2\sqrt{\gamma}\cdot n^2$. Let 
$$\mc{B} := \bigg\{v\in V(G): \textrm{IntDeg}_G(v) \ge 0.371n\bigg\}.$$
If $\mc{A}\subseteq V(G)$ such that $\mc{B}\subseteq \mc{A}$ and $|\mc{A}|\le 2\sqrt{\gamma} n$ and $\mc{E}\subseteq E(G\setminus \mc{A})$ with $|\mc{E}|\le 2$, then there exists a triangle packing $\mc{T}$ of $G-\mc{E}$ that saturates $\mc{A}$ such that $G':=(G-E(\mc{T}))\setminus \mc{A}$ is cross-divisible and $\delta(G')\ge (\frac{3}{4}-7\sqrt{\gamma})\cdot v(G')$.
\end{lem}
\begin{proof}
We choose $\beta$ sufficiently small to satisfy various inequalities throughout the proof. Let $(V_1,V_2)$ be the partition of $G$. Since $G$ is $\gamma$-extremal, we have by Lemma~\ref{lem:neithersideRobust}(a) that $|V_i| \ge (1-\sqrt{\gamma})n/2$.

Now we turn to cleaning. First, it will behoove us to upper bound $\textrm{IntSurplus}(G)$ robustly so that we know we only have to remove a small triangle packing for our surplus reserve.

\begin{claim}\label{cl:SurplusBounded}
$2\gamma\cdot n^2\ge \textrm{IntSurplus}(G)$.
\end{claim}
\begin{proofclaim}
Since $G$ is $\gamma$-extremal, we have that 
$$e(G) \le \frac{3}{8}n^2 + \frac{1}{2}\left(\gamma n\cdot \frac{n}{4} + n\cdot \gamma n\right) \le \left(\frac{3}{8}-\gamma\right)n^2,$$ 
and 
$$e_{12}(G) \ge (1-\gamma)\cdot \frac{n^2}{4}$$ and hence 
$$\textrm{IntSurplus}(G) = 2\cdot e(G) - 3\cdot e_{12}(G) \le 2\gamma\cdot n^2$$ 
as desired.
\end{proofclaim}

Next we show that $\textrm{IntSurplus}(G)$ is lower bounded linearly by $|\mc{B}|$ as follows.

\begin{claim}
$\textrm{IntSurplus}(G) \ge 0.06n\cdot |\mc{B}|$.
\end{claim}
\begin{proofclaim}
Recall that $\textrm{IntSurplus}(G) = 2\cdot e(G)-3\cdot e_{12}(G)$. Let $m:= |V_1|\cdot |V_2|-e_{12}(G)$ (that is the number of `missing' edges in the cut $(V_1,V_2)$). For each $i\in [2]$ and $v\in V_i$, let $m(v) := |V_{3-i}|-\textrm{CrossDeg}_G(v)$ (that is the number of vertices in $V_{3-i}$ non-adjacent to $v$. By the Handshaking Lemma applied to $\overline{G}$, we have that $2m=\sum_{v\in v(G)} m(v)$. Now let $\textrm{Surplus}(v) := d_G(v)-\frac{3}{4}n + \frac{m(v)}{2}$.

Using that $|V_1|\cdot |V_2|\le \frac{n^2}{4}$, we thus find that
\begin{align*}
\textrm{IntSurplus}(G) &= 2\cdot e(G) - 3 (|V_1|\cdot |V_2|-m) \ge \sum_{v\in v(G)} d_G(v) - \frac{3}{4}n^2 + m \\
&\ge \sum_{v\in V(G)} \bigg(d_G(v) - \frac{3}{4}n + \frac{m(v)}{2}\bigg) = \sum_{v\in V(G)} \textrm{Surplus}(v).
\end{align*}
Since $\delta(G)\ge \frac{3}{4}n$, we have that $d_G(v)\ge \frac{3}{4}n$ for all $v\in V(G)$. Also note that $m(v)\ge 0$ for all $v\in V(G)$. Combining we find that $\textrm{Surplus}(v)\ge 0$ for all $v\in V(G)$.

We claim that for each $v\in \mc{B}$, we have that $\textrm{Surplus}(v) \ge 0.06$. Let $v\in V_i$. But then
\begin{align*}
\textrm{Surplus}(v) &= d_G(v) - \frac{3}{4}n + \frac{|V_{3-i}| - \textrm{CrossDeg}_G(v)}{2} = \frac{\textrm{IntDeg}_G(v)+d_G(v)}{2} -\frac{3}{4}n + \frac{|V_{3-i}|}{2} \\
&\ge \frac{.371n+(\frac{3}{4}-\gamma)n}{2} - .75n + \frac{(\frac{1}{2} - \sqrt{\gamma})n}{2} \ge (0.0605n - 2\sqrt{\gamma})n \ge 0.06n,      
\end{align*}
as claimed where the last inequality follows since $\gamma$ is small enough. 

But then $\textrm{IntSurplus}(G) \ge \sum_{v\in V(G)} \textrm{Surplus}(v) \ge \sum_{v\in \mc{B}'} \textrm{Surplus}(v) \ge 0.06n\cdot |\mc{B}|$ as desired.
\end{proofclaim}

Next we apply Proposition~\ref{prop:SmallDegSubpacking} to find a large enough subset $\mc{T}'$ of $\mc{T}_0$ with small maximum degree (tied to $\gamma$). 

\begin{claim}
There exists a triangle packing $\mc{T}'$ of $(G-\mc{E})\setminus \mc{A}$ consisting only of interior triangles such that $|\mc{T}'| = \lceil \gamma n^2 \rceil$ and $\Delta(\mc{T})\le \sqrt{\gamma} \cdot n$.     
\end{claim}
\begin{proofclaim}
Let $d:= \frac{\sqrt{\gamma}}{2}$. By Proposition~\ref{prop:SmallDegSubpacking}, there exists $\mc{T}\subseteq \mc{T}_0$ such that $|\mc{T}|\ge d\cdot |\mc{T}_0| \ge \frac{\sqrt{\gamma}}{2} \cdot 2\sqrt{\gamma}\cdot n\ge \gamma\cdot n^2$ and $\Delta(\mc{T}) \le dn+2 \le \frac{\sqrt{\gamma}}{2}\cdot n + 2 \le \sqrt{\gamma}\cdot n$ where the last inequality follows since $n$ is large enough.
\end{proofclaim}

Let $G_1:=G-E(\mc{T}')$. Let $\varepsilon\in \left(0,\frac{1}{2000}\right]$ be as in Lemma~\ref{lem:OddballCleaning}. Next via our oddball cleaning lemma, Lemma~\ref{lem:OddballCleaning}, we find a triangle decomposition $\mc{T}_1$ of $\mc{A}$ (the oddballs) using at most $\frac{1}{6} \cdot 0.06n\cdot |\mc{B}|$ (so as to achieve non-negative surplus afterwards).

\begin{claim}
There exists a triangle decomposition $\mc{T}_1$ of $\mc{A}$ in $G_1-\mc{E}$ using at most $0.01n\cdot |\mc{B}|$ interior triangles of $G$.
\end{claim}
\begin{proofclaim}
Let $v_1,\ldots, v_k$ be an enumeration of $\mc{A}$. Let $S_0:=\emptyset$. For each $i\in [k]$ ,let $S_i:=\{v_1,\ldots,v_i\}$. We show by induction for $i\in \{0,\ldots,k\}$ that there exists a triangle decomposition $\mc{T}'_i$ of $S_i$ in $G_1-\mc{E}$ using at most $0.01n\cdot |S_i\cap \mc{B}'|$ interior triangles of $G$. For $i=0$, the empty decomposition suffices. 

So we assume $i\ge 1$. By induction, there exists a triangle decomposition $\mc{T}_{i-1}'$ of $S_{i-1}$ in $G_1-\mc{E}$ using at most $0.01n\cdot |S_{i-1}\cap \mc{B}|$ interior triangles of $G$. Let $G_i':= (G_1-E(\mc{T}_{i-1}'-\mc{E}))\setminus S_{i-1}$. Let $n':=v(G_i')$. Note that $n' = n-(i-1) \ge n-|\mc{A}| \ge (1-4\sqrt{\gamma})n$. Let $(V_1',V_2')$ be the partition of $G_i'$ given by setting $V_j':=V_j\cap V(G_i')$ for each $j\in [2]$. Recall that $\varepsilon\in \left(0,\frac{1}{2000}\right]$ is as in Lemma~\ref{lem:OddballCleaning}. Since $\gamma$ is small enough, it follows that $|V_j'| \ge |V_j| - |\mc{A}| \ge (\frac{1}{2}-4\sqrt{\gamma})n \ge (1-\varepsilon)n/2$.

Let $\mc{A}_i:= \{v\in V(G_i'): \textrm{IntDeg}_{G_i'}(v) \ge (\frac{1}{4}+\varepsilon)n\}$. Since $\gamma$ is small enough, it follows that $\mc{A}_i\subseteq \mc{A}\setminus S_{i-1}$. Hence $|\mc{A}_i|\le |\mc{A}| \le 2\sqrt{\gamma} \cdot n \le \varepsilon\cdot n'$ where the last inequality follows since $\gamma$ is small enough.  

Finally for each $j\in [2]$ and $v\in V_j'$, we have since $G$ is $\gamma$-cross-balanced that 
$$d_{j,G_i'}(v) \le d_{j,G}(v) \le d_{3-j,G}(v) +\gamma n\le d_{3-j,G_i'}(v)+|\mc{A}| +\gamma n\le d_{3-j,G_i'}(v) + 2\sqrt{\gamma} n +\gamma n \le  d_{3-j,G_i'}(v) +\varepsilon \cdot n',$$
where the last inequality follows since $\gamma$ is small enough.

Let $\mc{B}_i:= \{v\in V(G_i'): \textrm{IntDeg}_{G_i'}(v) > (\frac{3}{8}-7\varepsilon)n'\}$. Since $\gamma$ is small enough, it follows that $\mc{B}_i\subseteq \mc{B}$. By Lemma~\ref{lem:OddballCleaning}(i), we have that if $v_i\not\in \mc{B}_i$, then there exists a triangle decomposition $\mc{Q}_i$ of $v_i$ in $G_i'$ using only cross triangles of $G_i$; similarly by Lemma~\ref{lem:OddballCleaning}(i), we have that if $v_i\in \mc{B}_i$, then there exists a triangle decomposition $\mc{Q}_i$ of $v_i$ in $G_i'$ using at most $10\varepsilon n'$ (which is at most $0.01n$ since $\varepsilon\le \frac{1}{1000}$) interior triangles of $G_i'$. Combining, we find in either case that $\mc{T}_i':= \mc{T}_{i-1}'\cup \mc{Q}_i$ is a triangle decomposition of $S_i$ in $G$ using at most $0.01n\cdot |S_i\cap \mc{B}|$ as desired.
\end{proofclaim}

Let $G_2:= (G-E(\mc{T}_1))\setminus \mc{A}$. Note that 
$$\textrm{IntSurplus}(G_2) \ge \textrm{IntSurplus}(G) - 6\cdot 0.01n \cdot |\mc{B}| \ge 0,$$ 
where the last inequality follows since $\textrm{IntSurplus}(G) \ge .06n \cdot |\mc{B}|$.

Next we find a subset $\mc{T}_2$ of $\mc{T}$ whose removal will leave a surplus of exactly zero (i.e.~cross-divisible) as follows.

\begin{claim}
There exists $\mc{T}_2\subseteq \mc{T}$ such that $G_2-E(\mc{T}_2)$ is cross-divisible.    
\end{claim}
\begin{proofclaim}
Choose $\mc{T}_2\subseteq \mc{T}'$ such that $\textrm{IntSurplus}(G_2-E(\mc{T}_2)~) \ge 0$ and subject to that $|\mc{T}_2|$ is maximized. Let $H:= G_2-E(\mc{T}_2)$. If $\textrm{IntSurplus}(H)=0$, then $\mc{T}_2$ is as desired. So we assume that $\textrm{IntSurplus}(H) > 0$.  Note that $H$ is triangle-divisible since $G$ is triangle-divisible and $H$ is obtained from $G$ by the deleting the triangle packing $\mc{T}\cup \mc{T}_1$. Thus by Proposition~\ref{prop:CrossDiffParity} since $G$ is triangle-divisible, we find that $\textrm{IntSurplus}(H) \ge 6$. Since $\mc{T}'$ consists of only internal triangles, it follows that $\textrm{IntSurplus}(H) = \textrm{IntSurplus}(G_2) - 6\cdot |\mc{T}_2|$.  Yet $\textrm{IntSurplus}(G_2) \le 2\gamma\cdot n^2$ by Claim~\ref{cl:SurplusBounded}. Hence 
$$|\mc{T}_2|\le \frac{2\cdot\gamma n^2}{6} < \gamma\cdot n^2 \le |\mc{T}'|.$$ Thus there exists $T\in \mc{T}'\setminus \mc{T}_2$. Let $\mc{T}_3 := \mc{T}_2\cup \{T\}$. But then $\textrm{IntSurplus}(G_2-\mc{T}_3) = \textrm{IntSurplus}(H) - 6 \ge 0$, contradicting the choice of $\mc{T}_2$.
\end{proofclaim}

\noindent Let $\mc{T}:=(\mc{T}_1\cup \mc{T}_2)$. Thus $\mc{T}$ is a triangle packing of $G-\mc{E}$ that saturates $\mc{A}$ such that $G':=(G-E(\mc{T}))\setminus \mc{A}$ is cross-divisible and furthermore we find that 
$$\delta(G') \ge \delta(G) - 2\cdot |\mc{A}| - \Delta(\mc{T}_2) \ge \left(\frac{3}{4}-\gamma\right)n - 4\sqrt{\gamma} n -\sqrt{\gamma}n \ge \left(\frac{3}{4}-7\sqrt{\gamma}\right)\cdot v(G'),$$
where for the last inequality we used that $v(G')\ge (1-\gamma)n$ and that $\gamma$ is small enough. It follows that $\mc{T}$ is as desired.    
\end{proof}

\subsection{Extremal but at least one side is not close to either $T_2(n/2)$ or $\overline{T_2(n/2)}$}

For this case, to create our reserve of triangles, we also require the (weak) stability version of Tur\'an's theorem (for just the triangle case) as follows which was proved by Erd\H{o}s~\cite{E66, E67} and Simonovits~\cite{S68}.

\begin{thm}[Erd\H{o}s-Simonovits 1966~\cite{ES66}]\label{thm:TuranStability}
For every $\alpha > 0$, there exists $\varepsilon > 0$ such that the following holds for all sufficiently large $n$: if $G$ is a graph on $n$ vertices with $e(G) \ge (\frac{1}{4}-\varepsilon)n^2$ and $G$ is not $\alpha$-close to $T_2(n)$, then $G$ contains a triangle.
\end{thm}

We note the following easy corollary of the above theorem.

\begin{cor}\label{cor:TuranStability}
For every $\alpha > 0$, there exists $\varepsilon > 0$ such that the following holds for all sufficiently large $n$: if $G$ is a graph on $n$ vertices with $e(G) \ge (\frac{1}{4}-\varepsilon)n^2$ and $G$ is not $\alpha$-close to $T_2(n)$, then $G$ contains a triangle packing $\mc{T}$ with $|\mc{T}|\ge \frac{\varepsilon}{3}\cdot n^2$.
\end{cor}
\begin{proof}
Let $\varepsilon'$ be as in Theorem~\ref{thm:TuranStability} for $\alpha$ and let $\varepsilon:= \frac{\varepsilon'}{2}$. Now let $\mc{T}$ be a triangle packing of $G$ such that $|\mc{T}|$ is maximized. Let $H:=G-E(\mc{T})$. But then 
$$e(H)\ge e(G) - 3\cdot |\mc{T}| \ge \left(\frac{1}{4}-\varepsilon\right)n^2 - 3\cdot \frac{\varepsilon}{3}\cdot n^2 = \left(\frac{1}{4}-2\varepsilon\right) n^2 = \left(\frac{1}{4}-\varepsilon'\right) n^2.$$
Thus by Theorem~\ref{thm:TuranStability}, $H$ contains a triangle $T$ and hence $\mc{T}\cup \{T\}$ contradicts the choice of $\mc{T}$.
\end{proof}

We are now prepared to prove Theorem~\ref{thm:neitherside}. 

\begin{proof}[Proof of Theorem~\ref{thm:neitherside}]
\renewcommand{\theclaim}{\ref{thm:neitherside}.\arabic{claim}}%
We choose $\beta$ sufficiently small to satisfy various inequalities throughout the proof. Let $(V_1,V_2)$ be a bipartition of $V(G)$ such that $e(G[V_1,V_2])$ is maximized. For the rest of the proof we then consider $G$ as $2$-partitioned graph with partition $(V_1,V_2)$ and hence the $2$-partitioned graph is also $\beta$-extremal by the choice of $(V_1,V_2)$. Given the choice of bipartition, it follows that $\textrm{IntDeg}_G(v) \le \frac{1}{2}\cdot d_G(v)$ for each $v\in V(G)$. Let 
$$\mc{A} := \bigg\{v\in V(G): \textrm{IntDeg}_G(v) \ge \left(\frac{1}{4}+2\sqrt{\beta}\right)n\bigg\}.$$ 
By Lemma~\ref{lem:neithersideRobust}, we have that $|\mc{A}| \le 2\sqrt{\beta}\cdot n$ and $G$ satisfies all of Lemma~\ref{lem:neithersideRobust}(a)-(b). Furthermore, since $\beta$ is small enough and $G$ is not $\alpha$-close to any of $T_2(n/2) + T_2(n/2)$, $\overline{T_2(n/2)} + T_2(n/2)$, or $\overline{T_2(n/2}) + \overline{T_2(n/2)}$, it follows that there exists $i\in [2]$ such that $G[V_i]$ is not $\alpha$-close to either $T_2(|V_i|)$ or $\overline{T_2(|V_i|)}$. 

Now we turn to cleaning. First we find our reserve of interior triangles in $G\setminus \mc{A}$. First, via Corollary~\ref{cor:TuranStability}, we find a large (but with size still tied to $\beta$) triangle packing $\mc{T}_0$ of interior triangles.

\begin{claim}
There exists a triangle packing $\mc{T}_0$ of $G\setminus \mc{A}$ consisting only of interior triangles such that $|\mc{T}_0| \ge 2\sqrt{\beta}\cdot n^2$. 
\end{claim}
\begin{proofclaim}
For $i\in [2]$, let $H_i:= G[V_i]$. By Lemma~\ref{lem:neithersideRobust}(c), there exists $i\in [2]$ such that $H_i$ is not $\frac{\alpha}{3}$-close to $T_2(v(H_i))$. We assume without loss of generality that $H_1$ is not $\frac{\alpha}{3}$-close to $T_2(v(H_1))$. By Lemma~\ref{lem:neithersideRobust}(a), we have that $|V_1|\ge (\frac{1}{2}-\sqrt{\beta})n$ and hence
$$\delta(H_1)\ge v(H_1)-\frac{n}{4} \ge \left(\frac{1}{4} - \sqrt{\beta}\right)n \ge \left(\frac{1}{2} - 4\sqrt{\beta}\right)v(H_1).$$
Let $H_1':= H_1\setminus \mc{A}$. Since $|A|\le 2\sqrt{\beta}\cdot n$, we find that $\delta(H_1') \ge (\frac{1}{2}-6\sqrt{\beta}) v(H_1')$. Thus 
$$e(H_1') \ge \left(\frac{1}{4}-12\sqrt{\beta}\right) \cdot v(H_1)^2.$$
Let $\varepsilon$ be as in Corollary~\ref{cor:TuranStability} for $\frac{\alpha}{3}$. Since $\beta$ is small enough, we have that $12\sqrt{\beta} \le \varepsilon$. Thus by Corollary~\ref{cor:TuranStability}, there exists a triangle packing $\mc{T}_0$ of $H_1$ such that $|\mc{T}_0| \ge \frac{\varepsilon}{3}\cdot n^2$. Since $\frac{\varepsilon}{3}\ge 2\sqrt{\beta}$ as $\beta$ is small enough and every triangle in $\mc{T}_0$ is interior by construction, the claim follows.
\end{proofclaim}

Note that $G$ is triangle-divisible. Thus by Lemma~\ref{lem:ManyOddballCleaning} with $\gamma:=\beta$ and $\mc{E}:=\emptyset$, there exists a triangle packing $\mc{T}$ of $G$ that saturates $\mc{A}$ such that $G':=(G-E(\mc{T}))\setminus \mc{A}$ is cross-divisible and $\delta(G')\ge (\frac{3}{4}-7\sqrt{\beta})\cdot v(G')$.

Finally we apply our idealized decomposition theorem to $G'$, Theorem~\ref{thm:F0Idealized} as follows but first we check the conditions. Let $\alpha':= \frac{\alpha}{4}$. Let $\varepsilon_0$ be as in Theorem~\ref{thm:F0Idealized} for $\alpha'$.

Let $(V_1',V_2')$ be the bipartition of $G'$ given by $V_i':= V(G')\cap V_i$. Note that $v(G')\ge v(G)-|\mc{A}| \ge (1-2\sqrt{\beta})n$. Given that $G$ satisfies Lemma~\ref{lem:neithersideRobust}(a) and that $\beta$ is small enough, it follows that $G'$ satisfies Theorem~\ref{thm:F0Idealized}(a) for $\varepsilon_0$. Similarly given that $G$ Lemma~\ref{lem:neithersideRobust}(b) and that $\beta$ is small enough, it follows that $G'$ satisfies Theorem~\ref{thm:F0Idealized}(b) for $\varepsilon_0$. Finally given that there exists $i\in [2]$ such that $G[V_i]$ is not $\alpha/3$-close to $T_2(|V_i|)$ or $\overline{T_2(|V_i|)}$, it follows since $\beta$ is small enough that $G'$ satisfies Theorem~\ref{thm:F0Idealized}(c) for $\alpha'$.

Hence by Theorem~\ref{thm:F0Idealized}, $G'$ has an $\mc{F}_0$-decomposition $\mc{Q}_0$. But then $\mc{T}\cup \mc{Q}_0$ is a $K_3$-decomposition of $G$ as desired. 
\end{proof}

\subsection{More Terminology for Partitioned Graphs}

For the next three cases, it will be helpful to recast our notion of closeness away from edit distance and instead look to how partitions of a graph $G$ align with a partitioned family $\mc{F}$. To that end, we first define the concept of a complete $\mc{F}$-blowup.

\begin{definition}[Complete $\mc{F}$-blowup]
Let $\mc{F}:=\{F_1,\ldots F_k\}$ be a family of $t$-partitioned graphs for some integer $t\ge 1$. For $i\in [t]$, we let $N_{\mc{F}}(i) := \{ j\in [t] :~\exists \ell\in [k] \text{ such that } e_{ij}(F_{\ell})\ne 0\}$. We say a $t$-partitioned graph $G$ with partition $(V_1,\ldots,V_t)$ is a \emph{complete} $\mc{F}$-blowup if $G$ is an $\emph{F}$-blowup and for all $i\in [t]$ and $j\in N_{\mc{F}}(i)$, we have that $e_{ij}(G) = |V_i|\cdot |V_j|$. We say $G$ is \emph{balanced} if $|V_i| \le |V_j| \le |V_i|+1 $ for all $i < j\in [t]$. We let $\mc{F}(n)$ denote the balanced complete $\mc{F}$-blowup on $n$ vertices.
\end{definition}

In this language, we have that $\mc{F}_1(n) = \overline{T_2(n/2)}+\overline{T_2(n/2)}$, $\mc{F}_2(n)=\overline{T_2(n/2)}+T_2(n/2)$, and $\mc{F}_3(n)=T_2(n/2)+T_2(n/2)$.

Now that we have the concept of a complete $\mc{F}$-blowup, we can also create a notion of how close a partitioned graph $G$ to a complete $\mc{F}$-blowup on the same size of parts as follows.

\begin{definition}[$\mc{F}$-abnormal Graph]\label{def:Abnormal}
Let $\mc{F}:=\{F_1,\ldots F_k\}$ be a family of $t$-partitioned graphs for some integer $t\ge 1$. Let $G$ be a $t$-partitioned graph with partition $(V_1,\ldots,V_t)$. Let $H$ be the complete $\mc{F}$-blowup with partition $(V_1,\ldots,V_t)$ (that is on the same set of vertices of $G$ with the same partition). We say a pair $u\ne v\in V(G)$ is \emph{$\mc{F}$-abnormal} if either $uv \in E(G)$ but $uv\not\in E(H)$, or, $uv\not\in E(G)$ but $uv\in E(H)$. We define the \emph{$\mc{F}$-abnormal graph of $G$}, denoted $\textrm{Abnormal}_{\mc{F}}(G)$, as the graph with vertex set $V(G)$ and edge set the abnormal pairs of $G$. We say a vertex $v\in V(G)$ is \emph{$\alpha$-$\mc{F}$-normal} if $d_{\textrm{Abnormal}_{\mc{F}}}(v) \le \alpha\cdot n$ and \emph{$\alpha$-$\mc{F}$-abnormal} otherwise.
\end{definition}

Finally, we need one last piece of notation.

\begin{definition}{Partitioning an unpartitioned graph}
Let $G$ be a graph. Let $\mc{P}=(V_1,\ldots,V_t)$ be a partition of $V(G)$ for some integer $t\ge 1$. Then we let $G[\mc{P}]$ denote the $t$-partitioned graph $G$ with partition $\mc{P}$. Similarly if $H$ is a subgraph of $G$, we let $H[\mc{P}]$ denote the $t$-partitioned graph $H$ with partition $(V_i\cap V(H): i\in [t])$. 
\end{definition}

\subsection{Close to $\overline{T_2(n/2)}+\overline{T_2(n/2)}$}\label{ss:FourCliques}

We now proceed to the next major case. For this and the next major case, we need an easy proposition to construct our reserve set of triangles as follows.

\begin{proposition}\label{prop:LargeTrianglePacking}
If $G$ is a graph on $n$ vertices with $e(G) \ge \frac{n^2}{4}$, then there exists a triangle packing $\mc{T}_0$ of $G$ with $|\mc{T}_0|\ge \frac{e(G)- (n^2/4)}{3}$.
\end{proposition}
\begin{proof}
Let $\mc{T}_0$ be a triangle packing of $G$ such that $|\mc{T}_0|$ is maximized.  Let $H:=G-E(\mc{T}_0)$. Since $|\mc{T}_0|$ is maximized, it follows that $H$ is triangle-free and hence by Tur\'an's theorem, we have that $e(H)\le \frac{n^2}{4}$. Thus
$$|\mc{T}_0| \ge \frac{e(G)-e(H)}{3}  \ge \frac{e(G) - (n^2/4)}{3}$$
as desired.
\end{proof}

For our next three cases, we also need a lemma about when it is possible to decompose a set of interior edges $\mc{O}$ of a $2$-partitioned graph that is of low maximum degree via disjoint cross triangle also of relatively low maximum degree, but first a definition.

\begin{definition}
Let $G$ be a graph and let $\mc{O}\subseteq E(G)$. A triangle packing $\mc{T}$ of $G$ is \emph{$\mc{O}$-disjoint} if $|\mc{O}\cap E(T)|\le 1$ for all $T\in \mc{T}$. We say $\mc{T}$ is \emph{$\mc{O}$-encompassed} if every $T\in \mc{T}$ contains exactly one edge of $\mc{O}$ and every edge of $\mc{O}$ is in a triangle of $\mc{T}$. 
\end{definition}

\begin{lem}\label{lem:PackingEdges}
There exists $\gamma_0$ such that the following holds for all $\gamma\in (0,\gamma_0]$ and large enough $n$: Let $G$ be a $2$-partitioned graph on $n$ vertices. If $\mc{O}\subseteq E_{11}(G)\cup E_{22}(G)$ with $\Delta(\mc{O})\le \gamma\cdot n$ and for each $i\in [2]$ and $uv\in \mc{O}\cap E_{ii}(G)$ we have that $|N_{3-i}(u)\cap N_{3-i}(v)| \ge \sqrt{\gamma}\cdot n$, then there exists a triangle packing $\mc{T}$ of $G$ such that all of the following hold:
\begin{itemize}
    \item[(i)] $\mc{T}$ consists only of cross triangles of $G$,
    \item[(ii)] $\mc{T}$ is $\mc{O}$-encompassed, and
    \item[(iii)] $\Delta(\mc{T}) \le 4\sqrt{\gamma}\cdot n$.
\end{itemize}
\end{lem}
\begin{proof}
Let $\mc{T}$ be a triangle-packing of $G$ consisting only of cross triangles of $G$ such that every triangle in $\mc{T}$ contains exactly one edge of $\mc{O}$ and that for every $v\in V(G)$, $|\{vwx\in \mc{T}: wx\in \mc{O}\}| \le \sqrt{\gamma}\cdot n$, and subject to those conditions $|\mc{T}|$ is maximized. Let $\mc{O}':= \mc{O}\cap E(\mc{T})$. If $|\mc{O}'|=|\mc{O}|$, then $\mc{T}$ is as desired since $\Delta(\mc{T}) \le 4\sqrt{\gamma} \cdot n$. 

So we assume that $|\mc{O}'| < |\mc{O}|$. Thus there exists $f=xy\in \mc{O}\setminus E(\mc{T}')$. We assume without loss of generality that $x\in V_1$ and $y\in V_1$. By assumption $|N(x)\cap N(y)\cap V_2|\ge \sqrt{\gamma}\cdot n$. Let $W:= \{v\in V(G): |\{vwx\in \mc{T}: wx\in \mc{O}\}| > \sqrt{\gamma}\cdot n-1\}$. It follows that $$|\mc{T}|\ge (\sqrt{\gamma}\cdot n-1)\cdot |W| \ge \frac{1}{2}\sqrt{\gamma}\cdot n \cdot |W|,$$
where the last inequality follows since $\gamma$ is small enough and $n$ is large enough. Yet $|\mc{T}| < |\mc{O}| \le \frac{1}{2}\gamma\cdot n^2$.  Hence 
$$|W| < \frac{\frac{1}{2}\gamma\cdot n^2}{\frac{1}{2}\sqrt{\gamma}\cdot n} = \sqrt{\gamma}\cdot n,$$
and hence $N(x)\cap N(y)\cap V_2 \setminus W\ne \emptyset$. Thus there exists $z\in N(x)\cap N(y)\cap V_2 \setminus W$. But then $\mc{T}\cup \{xyz\}$ contradicts the choice of $\mc{T}$.    
\end{proof}

Next we need the notion of a parity fixer defined as follows. This particular definition is specific to our problem but the same gadget works for both $\mc{F}_1$ and $\mc{F}_2$ (which seems to be mostly happenstance).  

\begin{definition}
Let $G$ be a $4$-partitioned graph. A \emph{parity-fixer for $G$} is subgraph $F$ of $G$ such that $V(F)=\{v, u_1,u_2,u_3,u_4\}$ such that $vu_i\in E(G)$ for all $i\in [4]$ and $u_1u_3u_2u_4$ is a cycle in $G$, and for each $i\in [4]$, $\mc{P}_G(u_i)=i$. We call $v$ the \emph{center} of the fixer and $u_1,u_2,u_3,u_4$ the \emph{spokes}.
\end{definition}

We are now prepared to prove Theorem~\ref{thm:4cliques}.

\begin{proof}[Proof of Theorem~\ref{thm:4cliques}]
\renewcommand{\theclaim}{\ref{thm:4cliques}.\arabic{claim}}%
We choose $\alpha$ sufficiently small to satisfy various inequalities throughout the proof. Since $G$ is $\alpha$-close to $\mc{F}_1(n)$, there exists a partition  $\mc{P}'=(V_1',\ldots,V_4')$ of $V(G)$ such that for each $i\in [4]$, we have $\lfloor n/4 \rfloor \le |V_i'| \le \lceil n/4\rceil$ and $e(\textrm{Abnormal}_{\mc{F}_1}(G[\mc{P'}])) \le \alpha\cdot n^2$. Let 
$$\mc{A} := \bigg\{v\in V(G): d_{\textrm{Abnormal}_{\mc{F}_1}(G[\mc{P}'])}(v)  \ge \sqrt{\alpha}\cdot n\bigg\}.$$ 
Since $e(\textrm{Abnormal}_{\mc{F}}(G[\mc{P}'])) \le \alpha\cdot n^2$, it follows that $|\mc{A}| \le 2\sqrt{\alpha}\cdot n$. 

Let $L:= \{ v\in \mc{A}: |N(v)\cap (V_3'\cup V_4')|\ge |N(v)\cap (V_1'\cup V_2')|\}$ and $R:= \mc{A}\setminus L(\mc{A})$. Let $\mc{A}_1:= \{v\in L: |N(v)\cap V_1'| \ge |N(v)\cap V_2'|\}$ and $\mc{A}_2:= L\setminus \mc{A}_1$; similarly let $\mc{A}_3:= \{v\in R: |N(v)\cap V_3'|\ge |N(v)\cap V_4'|\}$ and $\mc{A}_4:= R\setminus \mc{A}_3$. For each $i\in [4]$, let $V_i:= (V_i'\setminus \mc{A}) \cup \mc{A}_i$. Let $\mc{P}:= (V_i:i\in[4])$. Note $\mc{P}$ is a partition of $V(G)$.

\begin{claim}\label{cl:NearlyBalanced1}
$G[\mc{P}]$ is $2\sqrt{\alpha}$-nearly-balanced.
\end{claim}
\begin{proofclaim}
For each $i\in [4]$, we have that $|V_i|\ge |V_i'|-|\mc{A}|\ge \frac{n}{4} - 2\sqrt{\alpha}\cdot n$ and the claim follows.    
\end{proofclaim}

\begin{claim}\label{cl:NewDegrees1}
Let $i\in [4]$, $v\in V_i$ and $j\in N_{\mc{F}_1}(i)$. If $v\notin \mc{A}$, then $|N_G(v)\cap V_j| \ge (\frac{1}{4}-3\sqrt{\alpha})n$; if $v\in \mc{A}$, then $|N_G(v)\cap V_j| \ge \frac{n}{10}$.    
\end{claim}
\begin{proofclaim}
First suppose $v\not\in \mc{A}$. Hence $|N(v)\cap V_j'| \ge |V_j'| - \sqrt{\alpha}\cdot n = (\frac{1}{4}-\sqrt{\alpha})n$. But then $|N(v)\cap V_j| \ge |N(v)\cap V_j'| - |V_j'\setminus V_j| \ge (\frac{1}{4}-\sqrt{\alpha})n - 2\sqrt{\alpha}\cdot n \ge (\frac{1}{4}-3\sqrt{\alpha})n$ as desired.

So we assume $v\in \mc{A}$. We assume without loss of generality that $i=1$ and hence $v\in \mc{A}_1$. Then by definition of $\mc{A}_1$, we have that  $|N(v)\cap (V_3'\cup V_4')|\ge |N(v)\cap (V_1'\cup V_2')|$ and $|N(v)\cap V_1'| \ge |N(v)\cap V_2'|$. Since $\delta(G)\ge \frac{3}{4}n$, it follows that $|N(v)\cap (V_3'\cup V_4')|\ge \frac{3}{8}n$. First suppose $j\in \{3,4\}$. But then we have $|N(v)\cap V_j'| \ge \frac{n}{8}$ and hence $|N(v)\cap V_j| \ge |N(v)\cap V_j| - |V_j'\setminus V_j| \ge \frac{n}{8} - 2\sqrt{\alpha} \ge \frac{n}{10}$ where the last inequality follows since $\alpha$ is small enough. So we assume $j=1$. Note that $|N(v)\cap (V_1'\cup V_2')| \ge d_G(v) - |V_3'\cup V_4'| \ge \frac{n}{4}$ and hence $|N(v)\cap V_1'| \ge \frac{n}{8}$. Thus $|N(v)\cap V_j| \ge |N(v)\cap V_j'| - |V_j'\setminus V_j| \ge \frac{n}{8} - 2\sqrt{\alpha}\cdot \ge \frac{n}{10}$ where the last inequality follows since $\alpha$ is small enough.
\end{proofclaim}

\begin{claim}\label{cl:ParityFixer4Cliques}
There exists a parity-fixer for $G[\mc{P}]$ whose spokes are in $V(G)\setminus \mc{A}$. 
\end{claim}
\begin{proofclaim}
Recall that $|\mc{A}|\le 2\sqrt{\alpha}\cdot n$. By Claim~\ref{cl:NearlyBalanced1}, we have that $|V_i|\ge (\frac{1}{4}-2\sqrt{\alpha})n$ for each $i\in [2]$. Recall that by Claim~\ref{cl:NewDegrees1}, we have that for all $i\in [4]$, $w\in V_i$ and $j\in N_{\mc{F}_1}(i)$ that $|N_G(w)\cap V_j|\ge (\frac{1}{4}-3\sqrt{\alpha})n$ (a fact we used extensively and implicitly in what follows). 

Since $(V_1,\ldots,V_4)$ is a partition of $V(G)$, there exists $i\in [4]$ such that $|V_i|\ge \frac{n}{4}$. Without loss of generality we assume that $|V_1|\ge \frac{n}{4}$.  Since $\alpha$ is small enough, it follows that $|V_2\setminus \mc{A}| > 0$. Let $u_2\in V_2\setminus \mc{A}$. Since $\delta(G)\ge \frac{3n}{4}$ and $|V_1| \ge \frac{n}{4}$, it follows that $|N(u_2)\cap V_1| > 0$. Let $v\in N(u_2)\cap V_1$. By Claim~\ref{cl:NewDegrees1}, we have that $|N(v)\cap V_i| \ge \frac{n}{10}$ for all $i\in [4]\setminus \{2\}$. Thus there exists $u_3\in V_3\cap N(u_2)\cap N(v) \setminus \mc{A}$. Similarly it follows from Claim~\ref{cl:NewDegrees1} that there exists $u_4\in V_4\cap N(u_2)\cap N(v) \setminus \mc{A}$. Finally by Claim~\ref{cl:NewDegrees1}, there exists $u_1\in V_1\cap N(u_3)\cap N(u_4)\cap N(v) \setminus \mc{A}$. But then $v,u_1,u_2,u_3,u_4$ form a parity-fixer all of whose spokes are in $V(G)\setminus \mc{A}$ as desired. 
\end{proofclaim}

By Claim~\ref{cl:ParityFixer4Cliques}, there exists a parity-fixer $F$ for $G$ with center $v$ and spokes $u_1,u_2,u_3,u_4$ where $u_i\in V_i\setminus \mc{A}$ for each $i\in [4]$.  Let $\mc{Q}_1 := \{vu_1u_3,vu_2u_4\}$ and let $\mc{E}=\{u_1u_4,u_2u_3\}$. Note that $E(F)=\mc{Q}_1\cup \mc{E}$ and also note that $\mc{E}\subseteq E(G\setminus \mc{A})$. Let $G_0:=G-E(\mc{Q}_1)$. 

Now we turn to cleaning.  Let $\mc{P}^*$ be the partition $(V_1\cup V_2, V_3\cup V_4)$ of $V(G)$. 

\begin{claim}
$G_0[\mc{P}^*]$ is $3\sqrt{\alpha}$-cross-balanced.
\end{claim}
\begin{proofclaim}
For each $v\in V(G_0)$, we find that $\textrm{IntDeg}_{G_0[\mc{P}^*]}(v) \le \textrm{CrossDeg}_{G}(v) + |\mc{A}| + 2 \le \textrm{CrossDeg}_{G_0[\mc{P}^*]}(v) + 2\sqrt{\alpha}\cdot n + 2 $ and hence $G^*$ is $\sqrt{\alpha}$-cross-balanced since $\alpha$ is small enough.
\end{proofclaim}

\begin{claim}
$G_0[\mc{P}^*]$ is $15\sqrt{\alpha}$-extremal. 
\end{claim}
\begin{proofclaim}
Note that $|V_1^*|=|V_1|+|V_2|\ge (\frac{1}{2}-4\sqrt{\alpha})n$ and similarly $|V_2^*|=|V_3|+|V_4|\ge  (\frac{1}{2}-4\sqrt{\alpha})n$. For each $v\in V_i\setminus \mc{A}$ and $j\in N_{\mc{F}_3}(i)$, we have by Claim~\ref{cl:NewDegrees1} that $|N_G(v)\cap V_j|\ge (\frac{1}{4}-3\sqrt{\alpha})n$. Thus 
$$e_{G_0[\mc{P^*}]}(V_1^*,V_2^*) \ge \frac{1}{2} \cdot (n-|\mc{A}|) \cdot \left(\frac{1}{2}-6\sqrt{\alpha}\right)n - 2 \ge (1-15\sqrt{\alpha})\cdot \frac{n^2}{4},$$
where for the last inequality we used that $\alpha$ is small enough. On the other hand, for each $i\in[4]$ and $v\in V_i\setminus \mc{A}$, we also find that 
$$d_{G_0}(v) \le d_G(v) \le \sqrt{\alpha}\cdot n +  \sum_{j\in N_{\mc{F}_3}(i)} |V_j'| \le \left(\frac{3}{4} + 2\sqrt{\alpha}\right)n$$ 
and hence $G[\mc{P}^*]$ is $15\sqrt{\alpha}$-extremal.     
\end{proofclaim}

\begin{claim}
There exists a triangle packing $\mc{T}_0$ of $(G_0-\mc{E})\setminus \mc{A}$ consisting only of interior triangles such that $|\mc{T}_0| \ge 30\cdot \alpha^{1/4}\cdot n^2$. 
\end{claim}
\begin{proofclaim}
Let $H:=G[V_1]$. By definition of $\mc{A}$, we find that if $v\in V(H)\cap \mc{A}$, then 
$$d_H(v)\ge v(H) - \sqrt{\alpha}\cdot n \ge \left(1 - \frac{\sqrt{\alpha}\cdot n}{v(H)}\right) \cdot v(H) \ge \left(1 - \frac{4\sqrt{\alpha}}{1-8\sqrt{\alpha}}\right) \ge (1-5\sqrt{\alpha})\cdot v(H),$$
where we used that $v(H)\ge (\frac{1}{4}-2\sqrt{\alpha})$ and $\alpha$ is small enough. Since $|\mc{A}|\le 2\sqrt{\alpha}\cdot n$, it follows that 
$$e(H) \ge \frac{1}{2} (v(H)-2\sqrt{\alpha}n) \cdot (v(H)-\sqrt{\alpha}n) \ge (1-7\sqrt{\alpha})\cdot \frac{v(H)^2}{2}.$$ 
Thus by Proposition~\ref{prop:LargeTrianglePacking}, there exists a triangle packing $\mc{T}_0$ of $G[V_1]$ with 
$$|\mc{T}_0|\ge \frac{e(G[V_1])-\frac{v(H)^2}{4}}{3} \ge \frac{(1-7\sqrt{\alpha}) \cdot \frac{v(H)^2}{2} - \frac{v(H)^2}{4}}{3} \ge \left(\frac{1}{12} - 2\sqrt{\alpha}\right)\cdot v(H)^2 \ge 30\cdot \alpha^{1/4}\cdot n^2,$$
as desired where the last inequality follows since $\alpha$ is small enough.
\end{proofclaim}

Thus by Lemma~\ref{lem:ManyOddballCleaning} applied to $G_0[\mc{P}^*]$ with $\gamma:=15\sqrt{\alpha}$, there exists a triangle packing $\mc{T}_1$ of $G_0-\mc{E}$ that saturates $\mc{A}$ in $G_0$ such that letting $G_1:=(G_0-E(\mc{T}))\setminus \mc{A}$, we have that $G_1[\mc{P}^*]$ is cross-divisible and $\delta(G_1)\ge (\frac{3}{4}-7\sqrt{\alpha})\cdot v(G_1)$.

\begin{claim}
There exists a triangle packing $\mc{T}_2$ of $G_1-\mc{E}$ with $\Delta(\mc{T}_2)\le 8\sqrt{\alpha}n$ such that letting $G_2:= G_1-E(\mc{T}_2)$, we have that $G_2[\mc{P}]$ is an $\mc{F}_1$-blowup and $G_2[\mc{P}^*]$ is cross-divisible. 
\end{claim}
\begin{proofclaim}
Let $\mc{O}:= \{e \in E(G_1): e \in E(\textrm{Abnormal}_{\mc{F}_1}(G_1[\mc{P}]))\}$. Note that $\Delta(\mc{O}) \le \sqrt{\alpha}\cdot n$ by definition of $\mc{A}$. Moreover $|V_i\cap V(G_1)| \ge |V_i'|-|\mc{A}| \ge (1 -2\sqrt{\alpha})n \ge (1-4\sqrt{\alpha}) \cdot v(G_1)$. Thus $G_1$ is $4\sqrt{\alpha}$-nearly-balanced. Let $\gamma:= 4\sqrt{\alpha}$. For each $i\in [2]$ and $uv\in \mc{O}\cap E_{ii}(G_1[\mc{P}^*])$, we find by Claim~\ref{cl:NewDegrees1} that 
$$|N_{G_1}(u)\cap N_{G_1}(v) \cap V_{3-i}^*| \ge \left(\frac{1}{2} -6\sqrt{\alpha}\right)n \ge (1-8\sqrt{\alpha})\cdot \frac{v(G_1)}{2} \ge \sqrt{\gamma}\cdot v(G_1),$$
where the final inequality follows since $\alpha$ is small enough. Hence by Lemma~\ref{lem:PackingEdges}, there exists a triangle-packing $\mc{T}_2$ of $G_1[\mc{P}^*]$ satisfying Lemma~\ref{lem:PackingEdges}(i)-(iii). Condition (i) says that $\mc{T}_2$ consists only of cross triangles of $G_1[\mc{P}^*]$ and hence $G_2[\mc{P}^*]$ is cross-divisible since $G_1[P^*]$ is. Condition (ii) yields that $\mc{O}\subseteq E(\mc{T}_2)$ and hence $G_2[\mc{P}]$ is an $\mc{F}_1$-blowup. Condition (iii) gives that $\Delta(\mc{T}_2) \le 4\sqrt{\gamma}\cdot n \le 8\sqrt{\alpha}\cdot n$. Hence $\mc{T}_2$ is as desired.
\end{proofclaim}

Let $G_2:= G_1-E(\mc{T}_2)$. Next we utilize our parity-fixer. Let $\mc{Q}_2:= \{vu_1u_4,vu_2u_3\}$. Note given our construction, we have $(\mc{T}_1\cup \mc{T}_2)\cap \mc{E}=\emptyset$ and hence $\mc{E}\subseteq E(G_2)$.

\begin{claim}
Either $G_2[\mc{P}]$ or $G_2[\mc{P}]-E(\mc{Q}_1)+E(\mc{Q}_2)$ is integrally $\mc{F}_1$-divisible.  
\end{claim}
\begin{proofclaim}
First suppose $e_{ij}(G_2[\mc{P}])$ is even for all $i\in \{1,2\}$ and $j\in \{3,4\}$. Then $G_2[\mc{P}]$ is integrally $\mc{F}_1$-divisible by Lemma~\ref{lem:F1div} as desired. 

So we assume there exists $i\in \{1,2\}$ and $j\in \{3,4\}$ such that $e_{ij}[G_2[\mc{P}])$ is odd. Since $G_2$ is Eulerian and $G_2[\mc{P}]$ is an $\mc{F}_1$-blowup, it follows that for each $i\in [2]$ that $\sum_{j\in \{3,4\}} e_{ij}(G_2[\mc{P}])$ is even and similarly for each $j\in \{3,4\}$ that $\sum_{i\in [2]} e_{ij}(G_2[\mc{P}])$ is even. But then it follows that for all $i\in \{1,2\}$ and $j\in \{3,4\}$, we have that $e_{ij}[G_2[\mc{P}])$ is odd.

Thus we let $G_2':= G_2-E(\mc{Q}_1)+E(\mc{Q}_2)$. But then that $e_{13}(G_2'[\mc{P}])=e_{13}(G_2[\mc{P}])-1$, $e_{14}(G_2'[\mc{P}])=e_{14}(G_2[\mc{P}])+1,~e_{23}(G_2'[\mc{P}])=e_{23}(G_2[\mc{P}])+1,~e_{24}(G_2'[\mc{P}])=e_{24}(G_2[\mc{P}])-1$ are all even and hence by Lemma~\ref{lem:F1div}, we find that $G_2'[\mc{P}]$ is integrally $\mc{F}_1$-divisible as desired.   
\end{proofclaim}    
 
If $G_2[\mc{P}]$ is integrally $\mc{F}_1$-divisible, let $G_3:= G_2$ and $\mc{T}:=\mc{T}_1\cup \mc{T}_2\cup \mc{Q}_1$; otherwise let $G_3:=G_2-E(\mc{Q}_1)+E(\mc{Q}_2)$ and $\mc{T}:= \mc{T}_1\cup \mc{T}_2\cup \mc{Q}_2$. Finally we apply our idealized decomposition theorem, Theorem~\ref{thm:Idealized}, to $G_3[\mc{P}]$ as follows but first we check the conditions. Let $\varepsilon_1$ be as in Theorem~\ref{thm:Idealized} for $\mc{F}_1$.

Note that $v(G_3)\ge v(G)-|\mc{A}| \ge (1-2\sqrt{\alpha})n$. Given that $G[\mc{P}]$ is $2\sqrt{\alpha}$-nearly-balanced and that $\alpha$ is small enough, it follows that $G_3[\mc{P}]$ is $\varepsilon_1$-nearly-balanced. Similarly given the definition of $\mc{A}$, it follows $G_3[\mc{P}]$ is an $\varepsilon_1$-nearly-complete $\mc{F}_1$-blowup.

Hence by Theorem~\ref{thm:Idealized}, $G_3$ has an $\mc{F}_1$-decomposition $\mc{Q}_0$. But then $\mc{T}\cup \mc{Q}_0$ is a $K_3$-decomposition of $G$ as desired. 
\end{proof}

\subsection{Close to $\overline{T_2(n/2)}+T_2(n/2)$}\label{ss:TwoCliquesOneBip}

We are now prepared to prove Theorem~\ref{thm:2cliqueandbip}.

\begin{proof}[Proof of Theorem~\ref{thm:2cliqueandbip}]
\renewcommand{\theclaim}{\ref{thm:2cliqueandbip}.\arabic{claim}}%
We choose $\alpha$ sufficiently small to satisfy various inequalities throughout the proof. Since $G$ is $\alpha$-close to $\mc{F}_2(n)$, there exists a partition  $\mc{P}'=(V_1',\ldots,V_4')$ of $V(G)$ such that for each $i\in [4]$, we have $\lfloor n/4 \rfloor \le |V_i'| \le \lceil n/4\rceil$ and $e(\textrm{Abnormal}_{\mc{F}_2}(G[\mc{P'}])) \le \alpha\cdot n^2$. Let 
$$\mc{A} := \bigg\{v\in V(G): d_{\textrm{Abnormal}_{\mc{F}_2}(G[\mc{P}'])}(v)  \ge \sqrt{\alpha}\cdot n\bigg\}.$$ 
Since $e(\textrm{Abnormal}_{\mc{F}_2}(G[\mc{P}'])) \le \alpha\cdot n^2$, it follows that $|\mc{A}| \le 2\sqrt{\alpha}\cdot n$. 

Let $L:= \{ v\in \mc{A}: |N(v)\cap (V_3'\cup V_4')|\ge |N(v)\cap (V_1'\cup V_2')|\}$ and $R:= \mc{A}\setminus L(\mc{A})$. Let $\mc{A}_1:= \{v\in L: |N(v)\cap V_1'| \ge |N(v)\cap V_2'|\}$ and $\mc{A}_2:= L\setminus \mc{A}_1$; similarly let $\mc{A}_3:= \{v\in R: |N(v)\cap V_4'|\ge |N(v)\cap V_3'|\}$ and $\mc{A}_4:= R\setminus \mc{A}_3$. For each $i\in [4]$, let $V_i:= (V_i'\setminus \mc{A}) \cup \mc{A}_i$. Let $\mc{P}:= (V_i:i\in[4])$. Note $\mc{P}$ is a partition of $V(G)$.

\begin{claim}\label{cl:NearlyBalanced2}
$G[\mc{P}]$ is $2\sqrt{\alpha}$-nearly-balanced.
\end{claim}
\begin{proofclaim}
For each $i\in [4]$, we have that $|V_i|\ge |V_i'|-|\mc{A}|\ge \frac{n}{4} - 2\sqrt{\alpha}\cdot n$ and the claim follows.    
\end{proofclaim}

\begin{claim}\label{cl:NewDegrees2}
Let $i\in [4]$, $v\in V_i$ and $j\in N_{\mc{F}_2}(i)$. If $v\notin \mc{A}$, then $|N_G(v)\cap V_j| \ge (\frac{1}{4}-3\sqrt{\alpha})n$; if $v\in \mc{A}$, then $|N_G(v)\cap V_j| \ge \frac{n}{10}$.    
\end{claim}
\begin{proofclaim}
First suppose $v\not\in \mc{A}$. Hence $|N(v)\cap V_j'| \ge |V_j'| - \sqrt{\alpha}\cdot n = (\frac{1}{4}-\sqrt{\alpha})n$. But then $|N(v)\cap V_j| \ge |N(v)\cap V_j'| - |V_j'\setminus V_j| \ge (\frac{1}{4}-\sqrt{\alpha})n - 2\sqrt{\alpha}\cdot n \ge (\frac{1}{4}-3\sqrt{\alpha})n$ as desired.

So we assume $v\in \mc{A}$. First suppose that $i\in \{1,2\}$. We assume without loss of generality that $i=1$ and hence $v\in \mc{A}_1$. Then by definition of $\mc{A}_1$, we have that  $|N(v)\cap (V_3'\cup V_4')|\ge |N(v)\cap (V_1'\cup V_2')|$ and $|N(v)\cap V_1'| \ge |N(v)\cap V_2'|$. Since $\delta(G)\ge \frac{3}{4}n$, it follows that $|N(v)\cap (V_3'\cup V_4')|\ge \frac{3}{8}n$. First suppose $j\in \{3,4\}$. But then we have $|N(v)\cap V_j'| \ge \frac{n}{8}$ and hence $|N(v)\cap V_j| \ge |N(v)\cap V_j| - |V_j'\setminus V_j| \ge \frac{n}{8} - 2\sqrt{\alpha} \ge \frac{n}{10}$ where the last inequality follows since $\alpha$ is small enough. So we assume $j=1$. Note that $|N(v)\cap (V_1'\cup V_2')| \ge d_G(v) - |V_3'\cup V_4'| \ge \frac{n}{4}$ and hence $|N(v)\cap V_1'| \ge \frac{n}{8}$. Thus $|N(v)\cap V_j| \ge |N(v)\cap V_j'| - |V_j'\setminus V_j| \ge \frac{n}{8} - 2\sqrt{\alpha}\cdot \ge \frac{n}{10}$ as desired where the last inequality follows since $\alpha$ is small enough.

So we assume that $i\in \{3,4\}$. We assume without loss of generality that $i=3$ and hence $v\in \mc{A}_3$. Then by definition of $\mc{A}_3$, we have that  $|N(v)\cap (V_1'\cup V_2')| \ge |N(v)\cap (V_3'\cup V_4')|$ and $|N(v)\cap V_4'| \ge |N(v)\cap V_3'|$. Since $\delta(G)\ge \frac{3}{4}n$, it follows that $|N(v)\cap (V_1'\cup V_2')|\ge \frac{3}{8}n$. First suppose $j\in \{1,2\}$. But then we have $|N(v)\cap V_j'| \ge \frac{n}{8}$ and hence $|N(v)\cap V_j| \ge |N(v)\cap V_j| - |V_j'\setminus V_j| \ge \frac{n}{8} - 2\sqrt{\alpha} \ge \frac{n}{10}$ where the last inequality follows since $\alpha$ is small enough. So we assume $j=4$. Note that $|N(v)\cap (V_3'\cup V_4')| \ge d_G(v) - |V_1'\cup V_2'| \ge \frac{n}{4}$ and hence $|N(v)\cap V_4'| \ge \frac{n}{8}$. Thus $|N(v)\cap V_j| \ge |N(v)\cap V_j'| - |V_j'\setminus V_j| \ge \frac{n}{8} - 2\sqrt{\alpha}\cdot \ge \frac{n}{10}$ as desired where the last inequality follows since $\alpha$ is small enough.
\end{proofclaim}

\begin{claim}\label{cl:ParityFixer2CliquesBip}
There exists a parity-fixer for $G[\mc{P}]$ whose spokes are in $V(G)\setminus \mc{A}$. 
\end{claim}
\begin{proofclaim}
Recall that $|\mc{A}|\le 2\sqrt{\alpha}\cdot n$. By Claim~\ref{cl:NearlyBalanced2}, we have that $|V_i|\ge (\frac{1}{4}-2\sqrt{\alpha})n$ for each $i\in [2]$. Recall that by Claim~\ref{cl:NewDegrees2}, we have that for all $i\in [4]$, $w\in V_i$ and $j\in N_{\mc{F}_2}(i)$ that $|N_G(w)\cap V_j|\ge (\frac{1}{4}-3\sqrt{\alpha})n$ (a fact we used extensively and implicitly in what follows). 

Since $(V_1,\ldots,V_4)$ is a partition of $V(G)$, there exists $i\in [4]$ such that $|V_i|\ge \frac{n}{4}$. First suppose that there exists $i\in [2]$ such that $|V_i|\ge \frac{n}{4}$. We assume without loss of generality that $|V_1|\ge \frac{n}{4}$.  Since $\alpha_0$ is small enough, it follows that $|V \setminus \mc{A}| > 0$. Let $u_2\in V_2\setminus \mc{A}$. Since $\delta(G)\ge \frac{3n}{4}$ and $|V_1| \ge \frac{n}{4}$ and $\alpha$ is small enough, it follows that $|N(u_2)\cap V_1| > 0$. Let $v\in N(u_2)\cap V_1$. By Claim~\ref{cl:NewDegrees2}, we have that $|N(v)\cap V_i| \ge \frac{n}{10}$ for all $i\in [4]\setminus \{2\}$. Thus there exists $u_3\in V_3\cap N(u_2)\cap N(v) \setminus \mc{A}$. Similarly it follows from Claim~\ref{cl:NewDegrees2} that there exists $u_4\in V_4\cap N(u_2)\cap N(v) \setminus \mc{A}$. Finally by Claim~\ref{cl:NewDegrees2}, there exists $u_1\in V_1\cap N(u_3)\cap N(u_4)\cap N(v) \setminus \mc{A}$. But then $v,u_1,u_2,u_3,u_4$ form a parity-fixer all of whose spokes are in $V(G)\setminus \mc{A}$ as desired. 

So we assume that $|V_1|,|V_2| < \frac{n}{4}$. It follows that there exists $i\in \{3,4\}$ such that $|V_i| > \frac{n}{4}$. Without loss of generality we assume that $|V_3|> \frac{n}{4}$. Since $\alpha_0$ is small enough, it follows that $|V_3\setminus \mc{A}| > 0$. Let $u_3\in V_2\setminus \mc{A}$. Since $\delta(G)\ge \frac{3n}{4}$ and $|V_3| > \frac{n}{4}$, it follows that $|N(u_3)\cap V_3| > 0$. Let $v\in N(u_3)\cap V_3$. By Claim~\ref{cl:NewDegrees2}, we have that $|N(v)\cap V_i| \ge \frac{n}{10}$ for all $i\in [4]\setminus \{3\}$. Thus there exists $u_2\in V_2\cap N(u_3)\cap N(v) \setminus \mc{A}$. Similarly it follows from Claim~\ref{cl:NewDegrees2} that there exists $u_1\in V_1\cap N(u_3)\cap N(v) \setminus \mc{A}$. Finally by Claim~\ref{cl:NewDegrees2}, there exists $u_4\in V_4\cap N(u_1)\cap N(u_2)\cap N(v) \setminus \mc{A}$. But then $v,u_1,u_2,u_3,u_4$ form a parity-fixer all of whose spokes are in $V(G)\setminus \mc{A}$ as desired. 
\end{proofclaim}

By Claim~\ref{cl:ParityFixer2CliquesBip}, there exists a parity-fixer $F$ for $G$ with center $v$ and spokes $u_1,u_2,u_3,u_4$ where $u_i\in V_i\setminus \mc{A}$ for each $i\in [4]$.  Let $\mc{Q}_1 := \{vu_1u_3,vu_2u_4\}$ and let $\mc{E}=\{u_1u_4,u_2u_3\}$. Note that $E(F)=\mc{Q}_1\cup \mc{E}$ and also note that $\mc{E}\subseteq E(G\setminus \mc{A})$. Let $G_0:=G-E(\mc{Q}_1)$. 

Now we turn to cleaning.  Let $\mc{P}^*$ be the partition $(V_1\cup V_2, V_3\cup V_4)$ of $V(G)$. 

\begin{claim}
$G_0[\mc{P}^*]$ is $3\sqrt{\alpha}$-cross-balanced.
\end{claim}
\begin{proofclaim}
For each $v\in V(G_0)$, we find that $\textrm{IntDeg}_{G_0[\mc{P}^*]}(v) \le \textrm{CrossDeg}_{G}(v) + |\mc{A}| + 2 \le \textrm{CrossDeg}_{G_0[\mc{P}^*]}(v) + 2\sqrt{\alpha}\cdot n + 2 $ and hence $G^*$ is $\sqrt{\alpha}$-cross-balanced since $\alpha$ is small enough.
\end{proofclaim}

\begin{claim}
$G_0[\mc{P}^*]$ is $15\sqrt{\alpha}$-extremal. 
\end{claim}
\begin{proofclaim}
Note that $|V_1^*|=|V_1|+|V_2|\ge (\frac{1}{2}-4\sqrt{\alpha})n$ and similarly $|V_2^*|=|V_3|+|V_4|\ge  (\frac{1}{2}-4\sqrt{\alpha})n$. For each $v\in V_i\setminus \mc{A}$ and $j\in N_{\mc{F}_3}(i)$, we have by Claim~\ref{cl:NewDegrees2} that $|N_G(v)\cap V_j|\ge (\frac{1}{4}-3\sqrt{\alpha})n$. Thus 
$$e_{G_0[\mc{P^*}]}(V_1^*,V_2^*) \ge \frac{1}{2} \cdot (n-|\mc{A}|) \cdot \left(\frac{1}{2}-6\sqrt{\alpha}\right)n - 2 \ge (1-15\sqrt{\alpha})\cdot \frac{n^2}{4},$$
where for the last inequality we used that $\alpha$ is small enough. On the other hand, for each $i\in[4]$ and $v\in V_i\setminus \mc{A}$, we also find that 
$$d_{G_0}(v) \le d_G(v) \le \sqrt{\alpha}\cdot n +  \sum_{j\in N_{\mc{F}_3}(i)} |V_j'| \le \left(\frac{3}{4} + 2\sqrt{\alpha}\right)n$$ 
and hence $G[\mc{P}^*]$ is $15\sqrt{\alpha}$-extremal.     
\end{proofclaim}

\begin{claim}
There exists a triangle packing $\mc{T}_0$ of $(G_0-\mc{E})\setminus \mc{A}$ consisting only of interior triangles such that $|\mc{T}_0| \ge 30\cdot \alpha^{1/4}\cdot n^2$. 
\end{claim}
\begin{proofclaim}
Let $H:=G[V_1]$. By definition of $\mc{A}$, we find that if $v\in V(H)\cap \mc{A}$, then 
$$d_H(v)\ge v(H) - \sqrt{\alpha}\cdot n \ge \left(1 - \frac{\sqrt{\alpha}\cdot n}{v(H)}\right) \cdot v(H) \ge \left(1 - \frac{4\sqrt{\alpha}}{1-8\sqrt{\alpha}}\right) \ge (1-5\sqrt{\alpha})\cdot v(H),$$
where we used that $v(H)\ge (\frac{1}{4}-2\sqrt{\alpha})$ and $\alpha$ is small enough. Since $|\mc{A}|\le 2\sqrt{\alpha}\cdot n$, it follows that 
$$e(H) \ge \frac{1}{2} (v(H)-2\sqrt{\alpha}n) \cdot (v(H)-\sqrt{\alpha}n) \ge (1-7\sqrt{\alpha})\cdot \frac{v(H)^2}{2}.$$ 
Thus by Proposition~\ref{prop:LargeTrianglePacking}, there exists a triangle packing $\mc{T}_0$ of $G[V_1]$ with 
$$|\mc{T}_0|\ge \frac{e(G[V_1])-\frac{v(H)^2}{4}}{3} \ge \frac{(1-7\sqrt{\alpha}) \cdot \frac{v(H)^2}{2} - \frac{v(H)^2}{4}}{3} \ge \left(\frac{1}{12} - 2\sqrt{\alpha}\right)\cdot v(H)^2 \ge 30\cdot \alpha^{1/4}\cdot n^2,$$
as desired where the last inequality follows since $\alpha$ is small enough.
\end{proofclaim}

Thus by Lemma~\ref{lem:ManyOddballCleaning} applied to $G_0[\mc{P}^*]$ with $\gamma:=15\sqrt{\alpha}$, there exists a triangle packing $\mc{T}_1$ of $G_0-\mc{E}$ that saturates $\mc{A}$ in $G_0$ such that letting $G_1:=(G_0-E(\mc{T}))\setminus \mc{A}$, we have that $G_1[\mc{P}^*]$ is cross-divisible and $\delta(G_1)\ge (\frac{3}{4}-7\sqrt{\alpha})\cdot v(G_1)$.

\begin{claim}
There exists a triangle packing $\mc{T}_2$ of $G_1-\mc{E}$ with $\Delta(\mc{T}_2)\le 8\sqrt{\alpha}n$ such that letting $G_2:= G_1-E(\mc{T}_2)$, we have that $G_2[\mc{P}]$ is an $\mc{F}_1$-blowup and $G_2[\mc{P}^*]$ is cross-divisible. 
\end{claim}
\begin{proofclaim}
Let $\mc{O}:= \{e \in E(G_1): e \in E(\textrm{Abnormal}_{\mc{F}_1}(G_1[\mc{P}]))\}$. Note that $\Delta(\mc{O}) \le \sqrt{\alpha}\cdot n$ by definition of $\mc{A}$. Moreover $|V_i\cap V(G_1)| \ge |V_i'|-|\mc{A}| \ge (1 -2\sqrt{\alpha})n \ge (1-4\sqrt{\alpha}) \cdot v(G_1)$. Thus $G_1$ is $4\sqrt{\alpha}$-nearly-balanced. Let $\gamma:= 4\sqrt{\alpha}$. For each $i\in [2]$ and $uv\in \mc{O}\cap E_{ii}(G_1[\mc{P}^*])$, we find by Claim~\ref{cl:NewDegrees2} that 
$$|N_{G_1}(u)\cap N_{G_1}(v) \cap V_{3-i}^*| \ge \left(\frac{1}{2} -6\sqrt{\alpha}\right)n \ge (1-8\sqrt{\alpha})\cdot \frac{v(G_1)}{2} \ge \sqrt{\gamma}\cdot v(G_1),$$
where the final inequality follows since $\alpha$ is small enough. Hence by Lemma~\ref{lem:PackingEdges}, there exists a triangle-packing $\mc{T}_2$ of $G_1[\mc{P}^*]$ satisfying Lemma~\ref{lem:PackingEdges}(i)-(iii). Condition (i) says that $\mc{T}_2$ consists only of cross triangles of $G_1[\mc{P}^*]$ and hence $G_2[\mc{P}^*]$ is cross-divisible since $G_1[P^*]$ is. Condition (ii) yields that $\mc{O}\subseteq E(\mc{T}_2)$ and hence $G_2[\mc{P}]$ is an $\mc{F}_1$-blowup. Condition (iii) gives that $\Delta(\mc{T}_2) \le 4\sqrt{\gamma}\cdot n \le 8\sqrt{\alpha}\cdot n$. Hence $\mc{T}_2$ is as desired.
\end{proofclaim}

Let $G_2:= G_1-E(\mc{T}_2)$. Next we utilize our parity-fixer. Let $\mc{Q}_2:= \{vu_1u_4,vu_2u_3\}$. Note given our construction, we have $(\mc{T}_1\cup \mc{T}_2)\cap \mc{E}=\emptyset$ and hence $\mc{E}\subseteq E(G_2)$.

\begin{claim}
Either $G_2[\mc{P}]$ or $G_2[\mc{P}]-E(\mc{Q}_1)+E(\mc{Q}_2)$ is integrally $\mc{F}_2$-divisible.  
\end{claim}
\begin{proofclaim}
First suppose $2\cdot e_{11}(G_2[\mc{P}])+e_{13}(G_2[\mc{P}])-e_{14}(G_2[\mc{P}]) \equiv 0 \pmod 4$. Then $G_2[\mc{P}]$ is integrally $\mc{F}_2$-divisible by Lemma~\ref{lem:F2div} as desired. 

So we assume that $2\cdot e_{11}(G_2[\mc{P}])+e_{13}(G_2[\mc{P}]) - e_{14}(G_2[\mc{P}])\not\equiv 0 \pmod 4$. Since $G_2$ is Eulerian, we have that $e_{13}(G_2[\mc{P}])+e_{14}(G_2[\mc{P}])$ is even and hence so is $e_{13}(G_2[\mc{P}]) - e_{14}(G_2[\mc{P}])$. Since $2\cdot e_{11}(G_2[\mc{P}])$ is also even, we find that $2\cdot e_{11}(G_2[\mc{P}])+e_{13}(G_2[\mc{P}]) - e_{14}(G_2[\mc{P}])\equiv 2 \pmod 4$.

Thus we let $G_2':= G_2-E(\mc{Q}_1)+E(\mc{Q}_2)$. Note that $e_{11}(G_2'[\mc{P}]) = e_{11}(G_2[\mc{P}])$. We also find that $e_{13}(G_2'[\mc{P}])=e_{13}(G_2[\mc{P}])-1$ and  $e_{14}(G_2'[\mc{P}])=e_{14}(G_2[\mc{P}])+1$. But then $2\cdot e_{11}(G_2'[\mc{P}])+e_{13}(G_2'[\mc{P}]) - e_{14}(G_2'[\mc{P}])\equiv 2 -1 - (+1) \equiv 0 \pmod 4$. Hence by Lemma~\ref{lem:F2div}, we find that $G_2'[\mc{P}]$ is integrally $\mc{F}_2$-divisible as desired.   
\end{proofclaim}    
 
If $G_2[\mc{P}]$ is integrally $\mc{F}_2$-divisible, let $G_3:= G_2$ and $\mc{T}:=\mc{T}_1\cup \mc{T}_2\cup \mc{Q}_1$; otherwise let $G_3:=G_2-E(\mc{Q}_1)+E(\mc{Q}_2)$ and $\mc{T}:= \mc{T}_1\cup \mc{T}_2\cup \mc{Q}_2$. Finally we apply our idealized decomposition theorem, Theorem~\ref{thm:Idealized}, to $G_3[\mc{P}]$ as follows but first we check the conditions. Let $\varepsilon_2$ be as in Theorem~\ref{thm:Idealized} for $\mc{F}_2$.

Note that $v(G_3)\ge v(G)-|\mc{A}| \ge (1-2\sqrt{\alpha})n$. Given that $G[\mc{P}]$ is $2\sqrt{\alpha}$-nearly-balanced and that $\alpha$ is small enough, it follows that $G_3[\mc{P}]$ is $\varepsilon_2$-nearly-balanced. Similarly given the definition of $\mc{A}$, it follows $G_3[\mc{P}]$ is an $\varepsilon_2$-nearly-complete $\mc{F}_2$-blowup.

Hence by Theorem~\ref{thm:Idealized}, $G_3$ has an $\mc{F}_2$-decomposition $\mc{Q}_0$. But then $\mc{T}\cup \mc{Q}_0$ is a $K_3$-decomposition of $G$ as desired. 
\end{proof}

\subsection{Close to $T_2(n/2)+T_2(n/2)$}\label{ss:TwoBip}

For this case, we first prove an in-between lemma which yields a decomposition when there are oddballs but they form an independent set and there are no interior edges as follows.

\begin{thm}[Idealized $\mc{F}_3$-decomposition with independent oddballs]\label{thm:IndepOddballCleaning}
There exists $\varepsilon > 0$ such that the following holds for all large enough $n$: Let $\mc{F}_3:= \{T_{123}, T_{124}, T_{134}, T_{234}\}$. Let $G$ be a $4$-partitioned graph on $n$ vertices and suppose that $G$ is integrally $\mc{F}_3$-divisible and $\delta(G) \ge 300\cdot \varepsilon n$. Let $$\mc{A}:= V(G)\setminus \bigcup_{i\in [4]} \Big\{v\in V_i: d_j(v)\ge  (1-2\varepsilon)n/4 \textrm{ for all }j\in [4]\setminus \{i\} \Big\}.$$ 
If all of the following hold:
\begin{enumerate}
    \item[(a)] $|V_i|\ge (1-\varepsilon)n/4$ for all $i\in [4]$,
    \item[(b)] $|\mc{A}| \le \varepsilon\cdot n$ and $\mc{A}$ is an independent set in $G$, and 
    \item[(c)] for each $i\in [4]$, $v\in \mc{A}\cap V_i$, and $j\in [4]\setminus \{i\}$, we have that $d_j(v) \le 0.49\cdot d(v)$,
\end{enumerate}
then $G$ has an $\mc{F}_3$-decomposition.
\end{thm}
\begin{proof}
Let $\mc{P}:=(V_1,\ldots, V_4)$ be the partition of $G$. Let $v_1,\ldots, v_k$ be an enumeration of $\mc{A}$. Let $S_0:=\emptyset$. For each $i\in [k]$, let $S_i:=\{v_1,\ldots,v_i\}$. 

\begin{claim}
For each $i\in \{0,\ldots,k\}$, there exists a triangle decomposition $\mc{T}_i$ of $S_i$ in $G$.    
\end{claim}
\begin{proofclaim}
We proceed by induction on $i$. For $i=0$, the empty decomposition suffices. 

So we assume $i\ge 1$. By induction, there exists a triangle decomposition $\mc{T}_{i-1}$ of $S_{i-1}$ in $G$. Let $G_i:= (G-E(\mc{T}_{i-1}))\setminus S_{i-1}$. Let $n':=v(G_i)$. Note that $n' = n-(i-1) \ge n-|\mc{A}| \ge (1-\varepsilon)n$. We assume without loss of generality that $v_i\in V_1$. Let $H:= G[N_G(v_i)]$ and let $H':=G_{i}[N_{G_i}(v_i)] = (H-E(\mc{T}_{i-1}))\setminus S_{i-1}$. By assumption (b), we have that $\mc{A}$ is an independent set in $G$. Thus $V(H')\subseteq V(G)\setminus \mc{A}$. Also by assumption (c), we have that $|N_G(v)\cap V_j| \le 0.49\cdot d_G(v)$ and hence $|V(H)\cap V_j|\le 0.49\cdot v(H)$ for each $j\in \{2,3,4\}$. By definition of $\mc{A}$, we find that for each $w\in V(H')$, we have
$$d_{H'}(w) \ge  d_H(w) - 2\cdot |S_{i-1}| \ge 0.51\cdot v(H) - 2\varepsilon\cdot n \ge \frac{1}{2}\cdot v(H'),$$
where the last inequality follows since $v(H')\ge v(H)-\varepsilon\cdot n$ and $v(H) = d_G(v_i) \ge \delta(G) \ge 300\varepsilon\cdot n$. Thus by P\'osa's Theorem (in fact Dirac's Theorem suffices here) we find that $H'$ has a Hamilton cycle; as $v(H')$ is even since $G$ is Eulerian (and hence so too is $G_i$), it follows that $H'$ has a perfect matching $M_i$. But then $M_i$ combined with the edges of $G_i$ incident with $v_i$ yields a triangle decomposition $\mc{T}_i'$ of $v_i$ in $G_i$. Hence $\mc{T}_i:=\mc{T}_{i-1}\cup \mc{T}_i'$ is as desired.
\end{proofclaim}

Let $G':= (G-E(\mc{T}_k))\setminus \mc{A}$. Note that $G'$ is integrally $\mc{F}_3$-divisible. Let $\varepsilon_3$ be as in Theorem~\ref{thm:Idealized} for $\mc{F}_3$. Note that $v(G')\ge v(G)-|\mc{A}| \ge (1-\varepsilon)n$. Given that $G$ is $\varepsilon$-nearly-balanced and that $\varepsilon$ is small enough, it follows that $G'$ is $\varepsilon_3$-nearly-balanced. Similarly given the definition of $\mc{A}$, it follows that $G'$ is an $\varepsilon_3$-nearly-complete $\mc{F}_3$-blowup.

Hence by Theorem~\ref{thm:Idealized}, $G'$ has an $\mc{F}_3$-decomposition $\mc{Q}_0$. But then $\mc{T}_k\cup \mc{Q}_0$ is a $K_3$-decomposition of $G$ and in fact an $\mc{F}_3$-decomposition of $G$ as desired (since the only triangles in $G$ are copies of $\mc{F}_3$). 
\end{proof}

Next we prove our key proposition which shows that the numbers of matched triangles required for integral $\mc{F}_3$-divisibility is achievable, but first we define these matched triangles as follows.

\begin{definition}
Let $M_1:= \{ \{1,2\}, \{3,4\}\}$, $M_2:=\{ \{1,3\}, \{2,4\}\}$, $M_3:=\{ \{1,4\}, \{2,3\}\}$. Let $G$ be a $4$-partitioned graph; for each $i\in [3]$, we say a triangle $T$ is \emph{$i$-matched} if $T$ is a copy of $T_{jjk}$ for $\{j,k\}\in M_i$. We consider the $M_i$ also as partitions and hence $G[M_i]$ denotes the $2$-partitioned graph $G[V_{j_1}\cup V_{j_2}, V_{k_1}\cup V_{k_2}]$ where $j_1j_2\ne k_1k_2\in M_i$. For $i\in [3]$, we say an edge $e$ is \emph{$i$-matched} if $e\in E_{jk}(G)$ and $jk\in M_i$; if $e\in E_{jj}(G)$ for some $j\in [4]$, we say $e$ is \emph{$0$-matched}. For $i\in \{0\}\cup [3]$, we let $M_i(G)$ denote the set of $i$-matched edges and let $m_i(G):=|M_i(G)|$.
\end{definition}

\begin{proposition}\label{prop:MatchingNumbers}
Let $G$ be a $4$-partitioned triangle-divisible graph on $n$ vertices with $\delta(G)\ge \frac{3}{4}n$. 
For each $i\in [3]$, let 
$$t_i:= \frac{\textrm{IntSurplus}(G[M_i])}{6}.$$ 
Then $t_i$ is non-negative and integral and $\sum_{i\in [3]}t_i = m_0$. Hence if there exists a triangle packing $\mc{T}$ of $G$ where for each $i\in [3]$ the number of $i$-matched triangles in $\mc{T}$ equals $t_i$, then $G-E(\mc{T})$ is integrally $\mc{F}_3$-divisible.
\end{proposition}
\begin{proof}
For each $i\in \{0\}\cup [3]$, let $m_i:= m_i(G)$. Note that for each $i\in [3]$, we calculate that
$$t_i = \frac{2m_0+2m_i-\sum_{j\in [3]\setminus \{i\}} m_j}{6} = \frac{2\cdot e(G) - 3\cdot \sum_{j\in [3]\setminus \{i\}} m_j}{6}.$$
For each $i\in [3]$, $G[M_i]$ is a $2$-partitioned triangle-divisible graph and hence by Proposition~\ref{prop:CrossDiffParity} we find that $t_i$ is integral; similarly since $\delta(G)\ge \frac{3}{4}n$, it follows that $t_i\ge 0$. One then calculates that $\sum_{i\in [3]} t_i = 3\cdot \frac{m_0}{3} + \sum_{i\in [3]} \frac{2m_i(G) - m_i(G) - m_i(G)}{6} = m_0$ as desired. Finally we prove that $G':=G-E(\mc{T})$ is integrally $\mc{F}_3$-divisible. since $G$ is triangle-divisible and $\mc{T}$ is a triangle packing, it follows that $G'$ is triangle-divisible and hence $G'$ is Eulerian. Moreover, we find that $e_{ii}(G')=0$ for all $i\in [4]$ and hence $G'$ is an $\mc{F}_3$-blowup. Moreover, we find that $e_{12}(G')+e_{34}(G') = m_1 -2\cdot t_1 = \frac{m_1+m_2+m_3-2m_0}{3}$.  Hence by symmetry, we find that  $e_{12}(G')+e_{34}(G') =  e_{13}(G') + e_{24}(G') = e_{14}(G')+e_{23}(G')$. Thus Lemma~\ref{lem:F3div}(i') and (ii') hold for $G'$. So by Lemma~\ref{lem:F3div}, we have that $G'$ is integrally $\mc{F}_3$-divisible as desired.
\end{proof}

\begin{remark}
As per our discussion before about internal surplus, in order for $G[M_i]$ to admit a triangle-decomposition (equivalently for $G$ to admit one), such a decomposition requires the use of precisely $\frac{\textrm{IntSurplus}(G[M_i])}{6} = t_i$ interior triangles of $G[M_i]$. Now an interior triangle of $G[M_i]$ is either an $i$-matched triangle or a triangle inside $G[V_j]$ for some $j\in [4]$. Since in general $G[V_j]$ may not contain any triangles (as in our current major case of being edit distance close to $T_2(n/2)+T_2(n/2)$, we would have that $V_j$ is edit distance close to an independent set), it is natural then to seek a decomposition using $t_i$ $i$-matched triangles.    
\end{remark}

To create our coloring of the internal edges with the desired properties, we use the following proposition. As a point of notation, for a $t$-partitioned graph $G$ and $v\in V(G)$ and $i\in [t]$, we let $\overline{d_i}(v) := |V_i|-d_i(v)$. First we need the following definition.

\begin{definition}
Let $G$ be a $4$-partitioned graph on $n$ vertices with partition $(V_1,\ldots, V_4)$. Let $i\in [3]$, $v\in V_j$ for some $j\in [4]$, and let $k$ be such that $jk\in M_i$; we define $m_i(v):=d_k(v)$ and $\overline{m_i}(v) := |V_k|-d_k(v)$. We define
$$t_i(v):= \left\lfloor 0.99 \cdot \bigg( \frac{d_G(v) - \frac{3}{4}n}{6} + \frac{1}{2}\cdot \sum_{j \in [3]\setminus \{i\}} \overline{m_j}(v)\bigg) \right\rfloor.$$    
\end{definition}

Here then is our coloring proposition. It might be helpful to the reader if we first describe in words what each outcome guarantees and what its use will be as in the following remark. 

\begin{remark}
Recall that we need to color all interior edges while being careful how we color edges incident with abnormal vertices. For interior edges incident with two abnormal vertices, this needed particular care and was the subject of our earlier claims. After coloring those, we will color the edges incident with an abnormal vertex $v$ with $t_i(v)$ edges of color $i$ so as to achieve condition Theorem~\ref{thm:IndepOddballCleaning}(c) and then color the remaining interior edges arbitrarily. This means the number of edges incident with $v$ of color $i$ will be at least $t_i(v)$ and at most $m_0(v)-\sum_{j\in [3]\setminus \{i\}} t_j(v)$. 

With this in mind here is a description of each outcome below. Outcome (i) says that $t_i(v)$ is non-negative which is clearly necessary to show. Outcome (ii) guarantees that our upper bound above on the number of edges colored $i$ is at most $m_i(v)$, a necessary condition if we plan to use that number of $i$-matched triangles at $v$ (where we have an additional buffer of $0.01n$ to account for small degrees losses). Outcome (iii) guarantees that new $i$-degree after matching (which is at most $m_i(v)-t_i(v)$) is at most $0.485$ times the degree after matching (which is $d_G(v)-2\cdot m_0(v)$) so that Theorem~\ref{thm:IndepOddballCleaning}(c) will hold (which was crucial to applying P\'osa's Theorem); note we included a buffer here as well ($0.005n$ and $0.485$ as opposed to $0.49$). Outcome (iv) guarantees the sum of the $t_i(v)$ is at most the internal degree $m_0(v)$, an obviously necessary condition, where again we have a buffer ($0.995$ vs $1$). Finally outcome (v) guarantees that two of the $t_i(v)$ are large which in turn will show that the two of the global parameter $t_i$'s are linear in the number of abnormals which is key to finding a good color for interior abnormal-abnormal edges. This will follow since we will show the sum of $t_i(v)$ over the abnormal vertices is a lower bound for $t_i$ (indeed it would be almost equal if we dropped the $0.99$ in the definition of $t_i(v)$ and changed the $6$ in the denominator to $3$).
\end{remark}

\begin{proposition}\label{prop:ColorProperties}
There exists $\gamma >0$ such that the following holds for all large enough $n$. Let $G$ be a $4$-partitioned graph on $n$ vertices with partition $(V_1,\ldots, V_4)$ and $\delta(G)\ge \frac{3}{4}n$ such that $G$ is $\gamma$-nearly-balanced. Let $G$ be a $4$-partitioned graph on $n$ vertices. Let $v\in V(G)$. Let $i\in [3]$. If $m_0(v) \ge \gamma^{2/3}\cdot n$ and $m_0(v) \le m_j(v) + 2\gamma\cdot n$ for all $j\in [3]$, then all of the following hold:
\begin{itemize}
    \item[(i)] $t_i(v) \ge 0$,
    \item[(ii)]$m_i(v) - 0.01n \ge m_0(v) - \sum_{j\in [3]\setminus \{i\}} t_j(v) $,
    \item[(iii)] $m_i(v)-t_i(v) \le 0.485 (d_G(v) - 2\cdot m_0(v)) -0.005n$,
    \item[(iv)] $\sum_{j\in [3]} t_j(v) \le 0.995 \cdot m_0(v)$,
    \item[(v)]  there exist distinct $j,k \in [3]$ such that $t_{j}(v), t_{k}(v) \ge \frac{1}{12} \cdot m_0(v)$.
\end{itemize}
\end{proposition}
\begin{proof}
We assume without loss of generality that $v\in V_4$. 

First we prove (i). Since $d_G(v)\ge \frac{3}{4}n$ as $\delta(G)\ge \frac{3}{4}n$, we find that $d_G(v)-\frac{3}{4}n\ge 0$. We also note that $\overline{m_j}(v)\ge 0$ for any $j\in [3]$. Hence $t_i(v)\ge 0$ and (i) holds as desired.

Next we prove (ii). Now we calculate that
\begin{align*}
m_i(v) + \sum_{j\in [3]\setminus \{i\}} t_j(v) &\ge 0.99\left(2\cdot \frac{d_G(v)  - \frac{3}{4}n}{6} + m_i(v) + \overline{m_i}(v) + \sum_{j\in [3]\setminus \{i\}} \frac{\overline{m_j}(v)}{2} \right)\\  
&\ge 0.99\left(\frac{d_G(v)}{3}  - \frac{n}{4} + |V_i| + \sum_{j\in [3]\setminus \{i\}} \frac{\overline{m_j}(v)}{2} \right)\\
&\ge 0.99 \cdot \left(\frac{d_G(v)}{3} -2\gamma \cdot n. \right)
\end{align*}
First suppose that $m_0(v) \le 0.2 n$. Since $\delta(G)\ge \frac{3n}{4}$, the above is at least $0.99(0.25n - 2\gamma n)\ge 0.22n \ge m_0(v)+0.01n$ as desired. So we assume $m_0(v) \ge 0.2n$. But then by assumption, we have  for each $j\in [3]$ that $m_j(v) \ge 0.2n-2\gamma \cdot n$. Hence $d_G(v) \ge 0.8n-6\gamma\cdot n$. But then the above is at least $0.99 (0.266n - 4\gamma\cdot n) > |V_4| + 0.01n \ge m_0(v) + 0.01n$ where we used that $\gamma$ is small enough and also that $|V_4|\le \frac{n}{4}+2\gamma\cdot n$ as $G$ is $\gamma$-nearly-balanced. This proves (ii).

Next we prove (iii). Now we calculate that
\begin{align*}
t_i(v) + 0.485(d_G(v)-2\cdot m_0(v)) &\ge 0.99 \cdot \left( \frac{d_G(v) - \frac{3}{4}n}{6} +  \frac{1}{2}\cdot \sum_{j\in [3]\setminus \{i\}} \overline{m_j}(v)\right) + 0.485\left(\sum_{i\in [3]} m_i(v) - m_0(v)\right) \\
&\ge 0.485 \cdot \sum_{j\in [3]\setminus \{i\}}(m_j(v)+\overline{m_j}(v)) + 0.99\cdot \frac{d_G(v) - \frac{3}{4}n}{6} + 0.485\cdot (m_i(v)-m_0(v)) \\
&\ge 0.485\cdot 2\cdot \left(\frac{n}{4}-\gamma \cdot n\right)+0.99\cdot \frac{d_G(v) - \frac{3}{4}n}{6} + 0.485\cdot (m_i(v)-m_0(v))\\
&\ge (0.2425-\gamma)n+0.99\cdot \frac{d_G(v) - \frac{3}{4}n}{6} + 0.485\cdot (m_i(v)-m_0(v)).
\end{align*}
First suppose that $m_i(v)-m_0(v)\ge 0.026n$. Then the above is at $(0.2425-\gamma)n+0.49(0.026n) \ge (0.25524-\gamma)n \ge m_i(v) + 0.005n$ as desired where the last inequality follows since $G$ is $\gamma$-nearly balanced and $G$ is small enough. So we assume $m_i(v) \le m_0(v) + 0.026n$. Next suppose that $m_i(v) < 0.237n$, but then the above is at least $(0.2425-\gamma)n + 0.485(-2\gamma n) \ge (0.2425-2\gamma)n \ge m_i(v) +0.005n$ as desired since $\gamma$ is small enough. So we assume $m_i(v) \ge 0.237n$ and hence $m_0(v) \ge 0.211n$ from above. But then we also have that for $j\in [3]\setminus \{i\}$, $m_j(v) \ge m_0(v)-2\gamma\cdot n$. Thus $d_G(v) \ge 0.237n + 0.211n + 2\cdot (0.211n -2\gamma)n \ge (0.87-4\gamma)n$. But now it follows that the last term in the large equation above is at least $(0.2425-\gamma)n + 0.99 \frac{(0.12-4\gamma)n}{6} + 0.485(-2\gamma n) \ge (0.26-3\gamma)n \ge m_i(v)+0.005n$ as desired where for the last inequality we used that $G$ is $\gamma$-nearly-balanced and $\gamma$ is small enough. This proves (iii).

Next we prove (iv). Thus we calculate that 
\begin{align*}
\sum_{j\in [3]} t_j(v) &\le 0.99 \cdot \left( 3\cdot  \frac{d_G(v) - \frac{3}{4}n}{6} + \sum_{j\in [3]} (|V_j|-m_j(v))\right) \\
&\le 0.99\cdot  \left( d_G(v) - \frac{3}{4}n + \sum_{j\in [3]} (|V_j|-m_j(v))\right) \\
&\le 0.99 \cdot (m_0(v) + 3\gamma \cdot n) \le 0.995\cdot m_0(v),
\end{align*}
where the second to last inequality follows since $G$ is $\gamma$-nearly-balanced and the last inequality follows since $m_0(v)\ge \gamma^{2/3}\cdot n$ and $\gamma$ is small enough. This proves (iv).

Finally we prove (v). First suppose $d_G(v) \ge \frac{3}{4}n + \frac{2}{3} m_0(v)$. Then it follows that for each $j\in [3]$, we have that $t_j(v) \ge \lfloor 0.99 \frac{2m_0(v)/3}{6} \rfloor \ge \frac{1}{12}m_0(v)$ as desired. So we assume $d_G(v) \le \frac{3}{4}n+\frac{2}{3}m_0(v)$. But then there exists $\ell\in [3]$ such that $m_\ell(v) \le \frac{1}{3}(d_G(v)-m_0(v)) \le  \frac{1}{3} (\frac{3}{4}n - \frac{m_0(v)}{3}) = \frac{n}{4} - \frac{m_0(v)}{9}$. Since $G$ is $\gamma$-nearly-balanced, it then follows that $\overline{m_{\ell}}(v) \ge \frac{m_0(v)}{9} - \gamma\cdot n \ge \frac{m_0(v)}{10}$ where the last inequality follows since $m_0(v)\ge \gamma^{2/3}\cdot n$ and $\gamma$ is small enough. But then each $j\in [3]\setminus \{\ell\}$ satisfies $t_j(v) \ge \lfloor 0.99\cdot \overline{m_\ell}(v) \rfloor \ge \frac{1}{12}m_0(v)$ as desired. This proves (v).
\end{proof}

We are now prepared to prove Theorem~\ref{thm:K4blowup}.

\begin{proof}[Proof of Theorem~\ref{thm:K4blowup}]
\renewcommand{\theclaim}{\ref{thm:K4blowup}.\arabic{claim}}%
We choose $\alpha$ sufficiently small to satisfy various inequalities throughout the proof. Since $G$ is $\alpha$-close to $\mc{F}_3(n)$, there exists a partition  $\mc{P}'=(V_1',\ldots,V_4')$ of $V(G)$ with $\lfloor n/4 \rfloor \le |V_i'| \le \lceil n/4\rceil$ for $i\in [4]$ such that $e(\textrm{Abnormal}_{\mc{F}_3}(G[\mc{P'}])) \le \alpha\cdot n^2$. Hence $\sum_{i<j\in [4]} e_{ij}(G[\mc{P}']) \ge (\frac{3}{8}-\alpha)n^2$. Also note that the abnormality inequality also implies that $e(G) \le (\frac{3}{8}+\alpha)n^2$. Let 
$$\mc{A}' := \bigg\{v\in V(G): d_{\textrm{Abnormal}_{\mc{F}_2}(G[\mc{P}'])}(v)  \ge \sqrt{\alpha}\cdot n\bigg\}.$$ 
Since $e(\textrm{Abnormal}_{\mc{F}_2}(G[\mc{P}'])) \le \alpha\cdot n^2$, it follows that $|\mc{A}'| \le 2\sqrt{\alpha}\cdot n$.  For each $i\in [4]$, let 
$$\mc{A}_i' := \{ v\in \mc{A}' \setminus \bigcup_{k<i} \mc{A}_k: |N(v)\cap V_i'|\le |N(v)\cap V_j'| \text{ for all } j\in [4]\setminus \{i\}\}.$$ 
For each $i\in [4]$, let 
$$V_i:= (V_i'\setminus \mc{A}') \cup \mc{A}_i'.$$ 
Let $\mc{P}:= (V_i:i\in[4])$. Note $\mc{P}$ is a partition of $V(G)$. For the rest of the proof we consider $G$ to be the $4$-partitioned graph $G[\mc{P}]$.

\begin{claim}\label{cl:NearlyBalanced3}
For each $i\in [4]$, we have that $(\frac{1}{4}+2\sqrt{\alpha})n\ge |V_i| \ge (\frac{1}{4}-2\sqrt{\alpha})n$ and hence $G$ is $2\sqrt{\alpha}$-nearly-balanced.
\end{claim}
\begin{proofclaim}
For each $i\in [4]$, we have that $|V_i|\ge |V_i'|-|\mc{A}'|\ge \frac{n}{4} - 2\sqrt{\alpha}\cdot n$ and similarly that $|V_i|\le |V_i'|+|\mc{A}'|\le \frac{n}{4} + 2\sqrt{\alpha}\cdot n$ and the claim follows.    
\end{proofclaim}

Now for each $i\in [4]$, we let
$$\mc{A}_i:= \{v\in V_i: |N_{G}(v)\cap V_i| \ge \alpha^{1/4}\cdot n\},$$
and
$$\mc{A}:= \bigcup_{i\in [4]} \mc{A}_i.$$

First we prove a series of claim concerning the degrees of vertices not in $\mc{A}$ and in $\mc{A}$ as well as upper bound the size of $\mc{A}$.

\begin{claim}\label{cl:NewDegrees3Normal}
Let $i\in [4]$. If $v\in V_i\setminus \mc{A}$ and $j\in N_{\mc{F}_3}(i)$, then $|N_G(v)\cap V_j| \ge (\frac{1}{4}-2\cdot \alpha^{1/4})n$.
\end{claim}
\begin{proofclaim}
Since $v\not\in \mc{A}$, we find that $|N(v)\cap V_i|\le \alpha^{1/4}\cdot n$. Note that $d_G(v) \ge \frac{3}{4}n$. By Claim~\ref{cl:NearlyBalanced3}, we have that $|N(v)\cap V_k|\le (\frac{1}{4}+2\sqrt{\alpha})n$ for each $k\in [4]\setminus \{i,j\}$. Hence 
$$|N(v)\cap V_j| \ge d_G(v) - |N_G(v)\cap V_i| - \sum_{k\in [4]\setminus\{i,j\}} |N(v)\cap V_k|\ge \frac{3}{4} n - \alpha^{1/4}\cdot n - 2\cdot \left(\frac{1}{4}+2\sqrt{\alpha}\right)n \ge \left(\frac{1}{4}-2\cdot \alpha^{1/4}\right)n,$$
where the last inequality follows since $\alpha$ is small enough. 
\end{proofclaim}

\begin{claim}
$|\mc{A}|\le 2\alpha^{3/4}\cdot n$.    
\end{claim}
\begin{proofclaim}
Note that for each $v\in \mc{A}$, we find that $|N_G(v)\cap V_i'|\ge |N_G(v)\cap V_i|-|\mc{A}'| \ge (\alpha^{1/4}-2\sqrt{\alpha})n\ge \frac{\alpha^{1/4}}{2} n$ since $\alpha$ is small enough. Hence $d_{\textrm{Abnormal}_{\mc{F}_3}(G[\mc{P}'])} (v) \ge \frac{\alpha^{1/4}}{2}n$. Since $e(\textrm{Abnormal}_{\mc{F}_3}(G[\mc{P'}])) \le \alpha\cdot n^2$, it follows that $|\mc{A}| \le \frac{\alpha\cdot n^2}{\alpha^{1/4}\cdot n /2} = 2\alpha^{3/4}\cdot n$ as desired.
\end{proofclaim}

\begin{claim}\label{cl:NewDegrees3}
Let $i\in [4]$, $v\in V_i$ and $j\in N_{\mc{F}_3}(i)$. If $v\in \mc{A}$, then $|N_G(v)\cap V_i| \le \frac{d_G(v)}{4} + 2\sqrt{\alpha}\cdot n$ and $|N_G(v)\cap V_i| \le |N_G(v)\cap V_j|+4\sqrt{\alpha}\cdot n$.
\end{claim}
\begin{proofclaim}
By definition of $\mc{A}_i'$ and since $\mc{A}\subseteq \mc{A}'$, it follows that $|N(v)\cap V_i'|\le |N(v)\cap V_k'|$ for all $k\in [4]\setminus \{i\}$. Hence $|N(v)\cap V_i'|\le \frac{d_G(v)}{4}$. But then $|N(v)\cap V_i| \le \frac{d_G(v)}{4} + |\mc{A}'| \le \frac{d_G(v)}{4} + 2\sqrt{\alpha}\cdot n$. Similarly $|N_G(v)\cap V_i|- |N_G(v)\cap V_j|\le |N_G(v)\cap V_i'|-|N_G(v)\cap V_j'| + 2|\mc{A}'| \le 2|\mc {A}'| \le 4\sqrt{\alpha}\cdot n$ as desired.
\end{proofclaim}

For cleaning, it will behoove us to categorize the types of edges to be cleaned as follows. Let
\begin{align*}
\mc{O}_1 &:= \bigcup_{i<j\in [4]} \{e=uv \in E_{ij}(G): u,v\in \mc{A}\}, \\
\mc{O}_2 &:=  \bigcup_{i\in [4]} \{e=uv\in E_{ii}(G): u,v\in \mc{A}\},\\
\mc{O}_3 &:=  \bigcup_{i\in [4]} \{e=uv\in E_{ii}(G): u,v\not\in \mc{A}\},\\
\mc{O}_4 &:= \bigcup_{i\in [4]} \{e=uv\in E_{ii}(G): u\in \mc{A} \text{ and } v\not\in \mc{A}\}.
\end{align*}
[That is $\mc{O}_1$ consists of `cross' edges whose both ends are in $\mc{A}$, $\mc{O}_2$ consists of `interior' edges whose ends are both in $\mc{A}$, $\mc{O}_3$ consists of interior edges whose ends are both not in $\mc{A}$, and $\mc{O}_4$ consists of interior edges with one end in $\mc{A}$ and the other end not in $\mc{A}$.]

We next prove a series of claims to establish properties we will need to create our coloring of the internal edges. First we prove that any edge in $\mc{O}_2$ has at least two color choices for which there are many options of matched triangle. To that end, we make the following definition. For each $e=uv\in \mc{O}_2$ such that $u,v\in V_i$ for $i\in [3]$, we define 
$$L(e):= \{ j\in [4]\setminus \{i\}: |N(u)\cap N(v) \cap V_j|,|N(u)\cap N(v)\cap V_k|\ge 0.05n\}.$$

Now we prove this has size at least $2$ as follows.

\begin{claim}\label{cl:ColorListTwo}
If $e \in \mc{O}_2$, then $|L(e)|\ge 2$.
\end{claim}
\begin{proofclaim}
Suppose not. Let $e=uv \in \mc{O}_2$. We assume without loss of generality that $u,v\in V_4$ and that $1,2\not\in L(e)$. Hence for each $i\in \{1,2\}$, we find that $|N(u)\cap N(v)\cap V_i| \le 0.05n$ and hence $d_i(u)+d_i(v) \le |V_i|+|N(u)\cap N(v)\cap V_i|\le (\frac{1}{4}+2\sqrt{\alpha})n + 0.05n$. By Claim~\ref{cl:NewDegrees3}, we have that $d_4(u)\le \frac{d_G(u)}{4}+2\sqrt{\alpha}\cdot n$ and similarly $d_4(v) \le \frac{d_G(v)}{4} + 2\sqrt{\alpha}\cdot n$. 

But then we find that
$$d_G(u)+d_G(v) \le 2\cdot \left(\left(\frac{1}{4}+2\sqrt{\alpha}\right)n + .05n\right)+ 2\cdot \left(\frac{1}{4}+2\sqrt{\alpha}\right)n + \frac{d_G(u)+d_G(v)}{4} + 4\sqrt{\alpha}\cdot n.$$
Rearranging yields that 
$$d_G(u)+d_G(v) \le \frac{4}{3} \cdot (1.1 + 12\sqrt{\alpha})n \le (1.47+16\sqrt{\alpha})n < \frac{3}{2}n,$$
where the last inequality follows since $\alpha$ is small enough. But this is a contradiction since $d_G(u)+d_G(v) \ge \frac{3}{2}n$ as $\delta(G)\ge \frac{3}{4}n$.
\end{proofclaim}

For each $i\in [3]$, let $t_i$ be as in Proposition~\ref{prop:MatchingNumbers} for $G$. Next we prove that if we color edges around abnormal vertices with $t_i(v)$ edges of color $i$ and account for coloring $\mc{O}_2$, we do not exceed $t_i$ (and hence have not `overcolored' and so will be able to complete the coloring to having exactly $t_i$ edges of color $i$). 

\begin{claim}\label{cl:NoColorOverused}
For each $i\in [3]$, we have that $\sum_{v\in \mc{A}} t_i(v) + |\mc{A}|^2 \le t_i$.
\end{claim}
\begin{proofclaim}
We assume without loss of generality that $i=1$. Now we calculate as follows that
\begin{align*}
t_1 &= \frac{\textrm{IntSurplus}(G[M_1])}{6} = \frac{2\cdot e(G) - 3\cdot (m_2(G)+m_3(G))}{6} \\
&= \frac{2\cdot e(G) - \frac{3}{4}\cdot n^2 + \frac{3}{4}\cdot n^2  - 3(|V_1|+|V_2|)(|V_3|+|V_4|) + 3 \cdot \bigg((|V_1|+|V_2|)(|V_3|+|V_4|) - m_2(G)-m_3(G)\bigg)}{6} \\
& = {\sum_{v\in V(G)}\frac{(d_G(v)-\frac{3}{4}n)}{6}} +\frac{\frac{n^2}{4} - (|V_1|+|V_2|)(|V_3|+|V_4|)}{2} + \frac{(|V_1|+|V_2|)(|V_3|+|V_4|) - m_2(G)-m_3(G)}{2}\\
& = {\sum_{v\in V(G)}\frac{(d_G(v)-\frac{3}{4}n)}{6}} +\frac{\frac{n^2}{4} - (|V_1|+|V_2|)(|V_3|+|V_4|)}{2} + \frac{1}{2}\cdot \sum_{v\in V(G)} \frac{\overline{m_2}(v) + \overline{m_3}(v)}{2}\\
&\ge \sum_{v\in \mc{A}} \frac{(d_G(v)-\frac{3}{4}n)}{6} + \sum_{v\in \mc{A}} \frac{\overline{m_2}(v) + \overline{m_3}(v)}{2} -|\mc{A}|^2\\
&\ge 1.005 \sum_{v\in \mc{A}} t_1(v) -|\mc{A}|^2\\
&\ge \sum_{v\in \mc{A}} t_1(v) + 0.005\cdot \frac{\alpha^{1/4}}{12}\cdot n \cdot |\mc{A}| -|\mc{A}|^2 \\ 
&\ge \sum_{v\in \mc{A}} t_1(v) + |\mc{A}|^2 ,
\end{align*}
where the first inequality follows since $\frac{n^2}{4} \ge (|V_1|+|V_2|)(|V_3|+|V_4|)$ and the fact that $\frac{1}{2}\cdot \sum_{v\in V(G)} \frac{\overline{m_2}(v)+\overline{m_3}(v)}{2} \ge \sum_{v\in \mc{A}} \frac{\overline{m_2}(v) + \overline{m_3}(v)}{2} -|\mc{A}|^2$, the second inequality follows from the definition of $t_1(v)$, the third inequality follows from Proposition~\ref{prop:ColorProperties}(v) since $\alpha$ is small enough, and the final inequality follows since $2|\mc{A}|^2 \le 0.005\cdot \frac{\alpha^{1/4}}{12}\cdot n \cdot |\mc{A}|$ as $|\mc{A}|\le 2\sqrt{\alpha}\cdot n \le \frac{1}{2}\cdot 0.005\cdot \frac{\alpha^{1/4}}{12}\cdot n$ since $\alpha$ is small enough. 
\end{proofclaim}

Next we prove a claim to show that at least two of the colors (read as two $t_i$'s) are large (and hence we will be able to color an edge $e\in \mc{O}_2$ with at least one color from $L(e)$ by pigeonhole). 

\begin{claim}\label{cl:TwoGoodColors}
There exist distinct $i,j\in [3]$ such that $t_i, t_j \ge \frac{1}{36}\cdot \alpha^{1/4}\cdot n \cdot |\mc{A}|$.
\end{claim}
\begin{proofclaim}
For each $v\in \mc{A}$, let $B(v):= \{i\in [3]: t_i(v) \ge \frac{1}{36}\cdot \alpha^{1/4}\cdot n\}$. By Proposition~\ref{prop:ColorProperties}(v), we have that $|B(v)|\ge 2$ for all $v\in \mc{A}$. By pigeonhole, that there exists $S\subseteq [3]$ with $|S|=2$ such that for at least $\frac{|\mc{A}|}{3}$ vertices $v$ of $\mc{A}$ we have $S\subseteq B(v)$. But then for each $i\in S$, we find that $t_i \ge \frac{|\mc{A}|}{3}\cdot \frac{1}{12}\cdot \alpha^{1/4}\cdot n$ as desired.\end{proofclaim}

Now we are prepared to prove our key coloring claim which colors the interior edges with the desired properties to achieve integral $\mc{F}_3$-divisibility so as to be able to apply Theorem~\ref{thm:IndepOddballCleaning}.

\begin{claim}\label{cl:3Coloring}
There exists a  (not necessarily proper) $3$-coloring $\varphi$ of $\mc{O}_2\cup \mc{O}_3\cup \mc{O}_4$ such that all of the following hold:
\begin{itemize}
    \item[(i)] for each $i\in [3]$, we have that $|\{e\in \mc{O}_2\cup \mc{O}_3\cup \mc{O}_4: \varphi(e) = i\}| = t_i$,
    \item[(ii)] for each $uv\in \mc{O}_2$, we have $\varphi(e)\in L(e)$,
    \item[(iii)] for each $i\in [3]$ and $v\in \mc{A}$, we have 
    $$m_i(v) - 0.01n \ge |\{e=vw\in \mc{O}_2\cup \mc{O}_4: \varphi(e)=i\}| \ge m_i(v) - 0.485 (d_G(v) - 2 \cdot m_0(v)) +0.005n.$$ 
\end{itemize}
\end{claim}
\begin{proofclaim}
Let $\mc{C}:= \{ k\in [3]: t_k \ge \frac{1}{36}\cdot \alpha^{1/4}\cdot n\cdot |\mc{A}|\}$. By Claim~\ref{cl:TwoGoodColors}, we find that $|\mc{C}|\ge 2$. For each $e\in \mc{O}_2$, we choose $\varphi(e)\in L(e) \cap \mc{C}$; note that such a color exists by the pigeonhole principle since $|\mc{C}|\ge 2$ and by Claim~\ref{cl:ColorListTwo} we have that $|L(e)|\ge 2$. Hence (ii) holds for $\varphi$.

Next for each $v\in \mc{A}$ and $i\in [3]$, we color $t_i(v)$ edges of $\mc{O}_3$ incident with $v$ with color $i$ (in a disjoint manner so that no edge receives two colors). Let $\gamma:= 2\sqrt{\alpha}$ and hence $G$ is $\gamma$-nearly-balanced by Claim~\ref{cl:NearlyBalanced3}. Note that as $v\in\mc{A}$, we have by definition that $m_0(v) \ge \alpha^{1/4}\cdot n\ge \gamma^{2/3}\cdot n$ since $\alpha$ is small enough. Thus by Proposition~\ref{prop:ColorProperties}(i), $t_i(v)\ge 0$ and by definition $t_i(v)$ is integral. Similarly by Proposition~\ref{prop:ColorProperties}(iv), $\sum_{j\in [3]} t_j(v) \le 0.995\cdot m_0(v)$. So this is indeed possible since $\alpha$ is small enough so that $|\mc{A}| \le 0.005\cdot \alpha^{1/4}\cdot n \le 0.005\cdot m_0(v) $.

Now we color the remaining uncolored edges of $\mc{O}_3\cup \mc{O}_4$ arbitrarily such that the total number of edges of color $i$ is exactly $t_i$. Such is possible by Claim~\ref{cl:NoColorOverused}. Then the upper bound in (iii) follows from Proposition~\ref{prop:ColorProperties}(ii) and the lower bound in (iii) follows from Proposition~\ref{prop:ColorProperties}(iii). Hence all of (i), (ii) and (iii) hold and $\varphi$ is as desired.
\end{proofclaim}

\begin{claim}
There exists an $\mc{O}_1$-encompassed triangle-packing $\mc{T}_1$ of $G$ using only $\mc{F}_3$-copies such that $\Delta(\mc{T}_1)\le 24\cdot \alpha^{3/8}\cdot n$.    
\end{claim}
\begin{proofclaim}
First note that for $i\in [3]$ and $v\in \mc{A}\cap V_i$ and $k\ne \ell\in [4]\setminus \{i\}$, we have by Claim~\ref{cl:NewDegrees3} letting $j = [4]\setminus \{i,k,\ell\}$ that
$$|N_G(v)\cap (V_k\cup V_{\ell})| \ge d_G(v) - |N_G(v)\cap V_i| - |V_j| \ge \frac{3}{4}\cdot d_G(v) - 2\sqrt{\alpha}\cdot n - \left(\frac{1}{4}+2\sqrt{\alpha}\right)n \ge \left(\frac{5}{16}-4\sqrt{\alpha}\right)n.$$
Thus for each $u\in V_i, v\in V_j$ such that $uv\in \mc{O}_1$, we find that letting $k\ne \ell \in [4]\setminus \{i,j\}$,
$$|N_G(u)\cap N_G(v)\cap (V_k\cup V_{\ell})| \ge 2\cdot \left(\frac{5}{16}-4\sqrt{\alpha}\right)n - 2\cdot \left(\frac{1}{4}+2\sqrt{\alpha}\right)n \ge \left(\frac{1}{8} -12\sqrt{\alpha}\right)n.$$

For each $i\in \{3\}$, we let $\mc{O}_{1,i}:= \mc{O}_1\cap M_i(G)$.  Note that $\Delta(\mc{O}_{1,i})\le \Delta(\mc{O}_1)\le |\mc{A}|\le 2\cdot \alpha^{3/4}\cdot n$.

We claim that for each $i\in [3]$, there exists an $\bigcup_{j\in [i]} \mc{O}_{1,j}$-encompassed triangle-packing $\mc{T}_{1,i}$ of $G$ using only $\mc{F}_3$-copies such that $\Delta(\mc{T}_{1,i})\le i\cdot 8\cdot \alpha^{3/8}\cdot n$. We proceed by induction. Let $G_1:= G-\mc{O}_{1,2}-\mc{O}_{1,3}$.  For $i=1$, such a $\mc{T}_{1,1}$ exists by Lemma~\ref{lem:PackingEdges} applied to $G_1[M_1]$ with $\gamma := 2\cdot \alpha^{3/4}$ since $\alpha$ is small enough. For $i\in \{2,3\}$, by Lemma~\ref{lem:PackingEdges} applied to $G_i:=G[M_i]-E(\mc{T}_{i-1}) - \bigcup_{j\in [3]\setminus \{i\}} \mc{O}_{1,j}$ and $\gamma$ since $\alpha$ is small enough, there exists an $\mc{O}_{1,i}$-encompassed triangle-packing $\mc{T}_{1,i}'$ of $G_i$ using only $\mc{F}_3$-copies in $G$ such that $\Delta(\mc{T}_{1,i}')\le 8\cdot \alpha^{3/8}\cdot n$ and hence $\mc{T}_{1,i}:= \mc{T}_{1,i-1}\cup \mc{T}_{1,i}'$ is as claimed.

But then $\mc{T}_1:=\mc{T}_{1,3}$ is as desired.
\end{proofclaim}

Let $G_1:= G-E(\mc{T}_1)$.

\begin{claim}\label{cl:T2}
There exists an $(\mc{O}_2\cup \mc{O}_3)$-encompassed triangle-packing $\mc{T}_2$ of $G_1$ such that $\Delta(\mc{T}_2)\le 4\cdot \alpha^{1/8}\cdot n$  and for each $e\in \mc{O}_2\cup \mc{O}_3$, the triangle in $\mc{T}_2$ containing $e$ is $\varphi(e)$-matched.
\end{claim}
\begin{proofclaim}
For each $i\in [3]$, let $\mc{O}_{2,i}:= \{e\in \mc{O}_2\cup \mc{O}_3: \varphi(e)=i\}$. For each $i\in [3]$, let $G_{1,i}$ be the $2$-partitioned graph $V(G_{1,i}):=V(G)$, $E(G_{1,i}):= M_i(G_1)$ and partition $(V_{j_1}\cup V_{k_1}, V_{j_1}\cup V_{k_2})$ where $M_i=\{j_1j_2,k_1k_2\}$. 

We claim that for each $i\in [3]$, $e=uv\in \mc{O}_2\cup \mc{O}_3$ with $u,v\in V_i$ and letting $j$ be such that $ij\in M_{\varphi(e)}$, we have that $|N_{G_1}(u)\cap N_{G_1}(v)\cap V_j| \ge 0.04n$. To see this, first suppose $u,v\in \mc{O}_3$, that is , $u,v\not\in \mc{A}$. Then by Claim~\ref{cl:NewDegrees3} we find that $|N_G(u)\cap V_j|, |N_G(v)\cap V_j|\ge (\frac{1}{4}-2\cdot \alpha^{1/4})n$ and hence $|N(u)\cap N(v)\cap V_j|\ge (\frac{1}{4}-6\cdot \alpha^{1/4})n$. Since $|N_{G}(w)\setminus N_{G_1}(w)|\le \Delta(\mc{T}_1)\le 24\cdot \alpha^{3/8}\cdot n$ for all $w\in V(G)$, it follows that $|N_{G_1}(u)\cap N_{G_1}(v)\cap V_j| \ge 0.04n$ since $\alpha$ is small enough. 

So we assume $uv\in \mc{O}_2$, that is $u,v\in \mc{A}$. Then by Claim~\ref{cl:3Coloring}(ii) we find that $|N(u)\cap N(v)\cap V_j|\ge .05n$. Since $|N_{G}(w)\setminus N_{G_1}(w)|\le \Delta(\mc{T}_1)\le 24\cdot \alpha^{3/8}\cdot n$ for all $w\in V(G)$, it follows that $|N_{G_1}(u)\cap N_{G_1}(v)\cap V_j| \ge 0.04n$ since $\alpha$ is small enough. 

Thus by Lemma~\ref{lem:PackingEdges} applied to $G_{1,i}$ with $\gamma:= \alpha^{1/4}$ since $\alpha$ is small enough, there exists an $\mc{O}_{2,i}$-encompassed triangle-packing $\mc{T}_{2,i}$ of $G_{1,i}$ consisting only of cross triangles of $G_i$ and such that $\Delta(T_{2,i})\le 4\cdot \alpha^{1/8}\cdot n$. By construction of $G_i$, it follows that every $\mc{T}_{2,i}$ is $i$-matched. Note also by construction that the $T_{2,i}$ are pairwise disjoint. Hence $\mc{T}_2:= \bigcup_{i\in [3]} T_{2,i}$ is as desired.
\end{proofclaim}

Let $G_2:= G_1-E(\mc{T}_2)$.

\begin{claim}\label{cl:T3}
There exists an $\mc{O}_4$-encompassed triangle-packing $\mc{T}_3$ of $G_2$ such that for each $e\in \mc{O}_4$, the triangle in $\mc{T}_3$ containing $e$ is $\varphi(e)$-matched.
\end{claim}
\begin{proofclaim}
Let $\mc{T}_3$ be a triangle-packing of $G_2$ such that every triangle in $\mc{T}_3$ contains exactly one edge in $\mc{O}_4$ and for each $e\in \mc{O}_3$ the triangle in $\mc{T}_3$ containing $e$ is $\varphi(e)$-matched, and subject to that $|\mc{T}_3|$ is maximized. If $|\mc{T}_3|=|\mc{O}_4|$, then $\mc{T}_3$ is as desired.

So we assume there exists $e=uv\in \mc{O}_4\setminus E(\mc{T}_3)$. We assume without loss of generality that $u\in \mc{A}\cap V_1$ and $v\in V_1\setminus \mc{A}$ and that $\varphi(e)=1$. By Claim~\ref{cl:NewDegrees3}, we find that $|N_G(v)\cap V_2| \ge (\frac{1}{4}-2\cdot \alpha^{1/4})n$ and hence $|V_2\setminus N_G(v)| \le 3\cdot \alpha^{1/4}\cdot n$. 

Note that for each $w\in V(G)$, we have that $|N_G(w)\setminus N_{G_2}(w)|\le \Delta(\mc{T}_1\cup\mc{T}_2) \le 28\cdot \alpha^{1/8}\cdot n$. Let $G_3:=G_2-E(\mc{T}_3)$. On the other hand by construction, we have that 
\begin{align*}
|N_{G_3}(u)\cap V_2| &\ge |N_{G_2}(u)\cap V_2| - |\{e=uv\in \mc{O}_2\cup \mc{O}_3\cup \mc{O}_4:\varphi(e)= 1\}| \\
&\ge |N_{G}(u)\cap V_2| - 28\cdot \alpha^{1/8}\cdot n - |\{e=uv\in \mc{O}_2\cup \mc{O}_3\cup \mc{O}_4:\varphi(e)= 1\}| \\
&\ge 0.01n - 28\cdot \alpha^{1/8}\cdot n,
\end{align*}
where the last inequality follows from the upper bound in Claim~\ref{cl:3Coloring}(iii). Since $|N_{G_2}(v)\setminus N_{G_3}(v)|\le 2\cdot |\mc{A}|\le 4\cdot \alpha^{3/4}\cdot n$ as $v\not\in \mc{A}$, we find that
\begin{align*}
|N_{G_3}(u)\cap N_{G_3}(v)\cap V_2| &\ge |N_{G_3}(u)\cap V_2| - |N_{G_2}(v)\setminus N_{G_3}(v)| - |N_{G}(v)\setminus N_{G_2}(v)| - |V_2\setminus N_G(v)| \\
&\ge 0.01n - 28\cdot \alpha^{1/8}\cdot n - 4\cdot \alpha^{3/4}\cdot n - 28\cdot \alpha^{1/4}\cdot n - 3\cdot \alpha^{1/4}\cdot n \\
&\ge 0.01n - 63\cdot \alpha^{1/8} \cdot n > 0,
\end{align*}
where the last inequality follows since $\alpha$ is small enough. Hence there exists $z\in N_{G_3}(u)\cap N_{G_3}(v)\cap V_2$. Let $\mc{T}_3':=\mc{T}_3\cup \{uvz\}$. But then $\mc{T}_3'$ is a triangle-packing of $G_2$ such that every triangle in $\mc{T}_3'$ contains exactly one edge in $\mc{O}_4$ and for each $e\in \mc{O}_4$ the triangle in $\mc{T}_3$ containing $e$ is $\varphi(e)$-matched (since $T$ is $1$-matched and $\varphi(uv)=1$). Since $|\mc{T}'| > |\mc{T}|$, this contradicts the choice of $\mc{T}$.
\end{proofclaim}

Let $G_3:= G_2-E(\mc{T}_3)$. Let $\mc{T}:=\bigcup_{i\in [3]} \mc{T}_i$ and hence $G_3=G-E(\mc{T})$. By Claim~\ref{cl:T2}, $\mc{T}_2$ is an $(\mc{O}_2\cup \mc{O}_3)$-encompassed triangle packing of $G$ such that for each $e\in \mc{O}_2\cup \mc{O}_3$, the triangle in $\mc{T}_2$ containing $e$ is $\varphi(e)$-matched. Similarly by Claim~\ref{cl:T3}, $\mc{T}_3$ is an $\mc{O}_4$-encompassed triangle packing of $G$ such that for each $e\in \mc{O}_4$, the triangle in $\mc{T}_3$ containing $e$ is $\varphi(e)$-matched. Since $\mc{T}_2\cap \mc{T}_3=\emptyset$ by construction, it follows from Claim~\ref{cl:3Coloring}(i) that for each $i\in [3]$, the number of $i$-matched triangles in $\mc{T}$ equals $t_i$. Hence by Proposition~\ref{prop:MatchingNumbers}, we find that $G_3$ is integrally $\mc{F}_3$-divisible. 

Let $\varepsilon$ be as in Theorem~\ref{thm:IndepOddballCleaning}. Since $G$ is $2\sqrt{\alpha}$-balanced by Claim~\ref{cl:NearlyBalanced3}, we find that Theorem~\ref{thm:IndepOddballCleaning}(a) holds for $G_3$. Let
$$\mc{A}_3:= V(G)\setminus \bigcup_{i\in [4]} \Big\{v\in V_i: d_{G_3,j}(v)\ge  (1-2\varepsilon)n/4 \textrm{ for all }j\in [4]\setminus \{i\} \Big\}.$$ 

Note that for each $w\in V(G)\setminus \mc{A}$, we have that $|N_{G}(w)\setminus N_{G_3}(w)| \le \Delta(\mc{T}) \le \Delta(\mc{T}_1\cup \mc{T}_2) + |\mc{A}| \le 28\cdot \alpha^{1/8}\cdot n + 2\alpha^{3/4}\cdot n$. Since $\alpha$ is small enough, it then follows that $\mc{A}_3\subseteq \mc{A}$. Thus $|\mc{A}_3|\le |\mc{A}|\le 2\alpha^{3/4}\cdot n\le \varepsilon n$ since $\alpha$ is small enough. Furthermore, $\mc{A}$ is an independent set in $G_3$ as $(\mc{O}_1\cup \mc{O}_2\cup \mc{O}_3\cup\mc{O}_4)\cap E(G_3)=\emptyset$; hence it follows that $\mc{A}_4$ is an independent set in $G$. Thus Theorem~\ref{thm:IndepOddballCleaning}(b) holds for $G_3$ and $\mc{A}_4$.  

Finally Theorem~\ref{thm:IndepOddballCleaning}(c) holds for $G_3$ and $\mc{A}_3$ (in fact for $\mc{A}$) by the lower bound in Claim~\ref{cl:3Coloring}(c) and the fact that $\Delta(T_1\cup T_2)\le 28\cdot \alpha^{1/8}\cdot n$ and $\alpha$ is small enough. Thus by Theorem~\ref{thm:IndepOddballCleaning}, there exists an $\mc{F}_3$-decomposition $\mc{Q}_0$ of $G_3$. But then $\mc{T}\cup \mc{Q}_0$ is a $K_3$-decomposition of $G$ as desired.
\end{proof}

\bibliographystyle{plain}
\bibliography{bibliography}

\end{document}